\documentclass[a4paper,10pt]{amsart}
\usepackage{epsf}  
\usepackage{amsfonts, amssymb, amsthm, amsmath}  
\usepackage{hyperref}
\usepackage{graphicx}
\usepackage{psfrag}
\newtheorem{thm}{Theorem}

\newtheorem{lemma}[thm]{Lemma}  
\newtheorem{remark}[thm]{Remark}  
\newtheorem{defn}[thm]{Definition}  
\newtheorem{prop}[thm]{Proposition}  
  
\newtheorem{example}[thm]{Example}  
\numberwithin{thm}{section}  
\def\pf{\noindent\emph{Proof: }}  
\def\stop{\hfill$\square$}  

\providecommand{\tdbar}[1]{\pi_{Y}(T_{#1}\dbar)}
\providecommand{\Y}[1]{Y(#1)}
\providecommand{\wb}[1]{wb(#1)}

\providecommand{\totl}[1]{\ensuremath{\lceil #1\rceil }}
\providecommand{\totb}[1]{\ensuremath{\underline{ #1}}}
\DeclareMathOperator{\End}{End}
\DeclareMathOperator{\Aut}{Aut}
\newcommand{\ro}{{}^{r}\Omega}
\providecommand{\rof}{{}^{\phantom{f}r}_{fg}\Omega}
\newcommand{\rh}{{}^{r}H}
\newcommand{\ex}{\bold}
\providecommand {\e}[1]{\mathfrak t^{#1}}
\providecommand{\sfp}[3]{\ensuremath \lrb{{#1}}_{#2}^{#3}}
\providecommand{\C}[2]{\ensuremath {C^{#1,\underline{#2}}}}

\newcommand{\Mod}{\mathcal M}
\newcommand{\Ms}{\mathcal M^{\infty,\underline 1}}
\newcommand{\Msw}{\mathcal M^{\omega}}
\newcommand{\bMs}{\bar{\mathcal M}^{\omega}}
\newcommand{\M}{\expl \bar{\mathcal M}}

\providecommand{\fp}[2]{{}_{\hspace{3pt}#1\hspace{-2pt}}\times_{#2}}
  
\DeclareMathOperator{\EV}{EV}
  
\DeclareMathOperator{\id}{id}
\DeclareMathOperator{\expl}{Expl}

\newcommand{\dbar}{\bar{\partial}}  

\newcommand{\dvert}{d_{\text vert}}

\providecommand{\et}[2]{\ensuremath{\bold T^{#1}_{#2}}}
\providecommand{\lrb}[1]{\ensuremath{\left(#1\right)}}
\providecommand{\abs}[1]{\left\lvert #1\right\rvert}

\author{Brett Parker   }
\email{brettdparker@gmail.com}  
\thanks{Partially supported by SNF, No 200020-119437/1. Part of this work was also completed during the authorÕs stay at UC Berkeley and the Mathematical Science Research Institute in Berkeley.}
  
\title{Gromov Witten invariants of exploded manifolds}

\begin{document}
\maketitle

\begin{abstract}This paper describes the structure of the moduli space of holomorphic curves and constructs Gromov Witten invariants in the category of exploded manifolds. This includes  defining Gromov Witten invariants relative to normal crossing divisors and proving the associated gluing theorem which involves summing relative invariants over a count of tropical curves. 

\end{abstract}

\tableofcontents

\newpage

\section{Introduction}

\

In this paper, we shall define Gromov Witten invariants of a basic, complete exploded manifold $\ex B$ with an almost complex structure $J$ and a taming form $\omega$. Definitions for exploded manifolds can be found in \cite{iec}.

\

 Section \ref{structure section} describes the structure of the $\dbar$ equation on a moduli stack $\Msw$ of  $\C\infty1$ curves in $\ex B$, then uses multiperturbations on a neighborhood of the substack of holomorphic curves to define a virtual moduli space $\mathcal M$. Similar results are obtained for a family $\hat{\ex  B}\longrightarrow \ex G$ of exploded manifolds. As this section is quite technical, the results are described in separate sections before they are proved.

\
 
 Section \ref{numerical GW} is concerned with integrating differential forms over the virtual moduli space $\mathcal M$ to obtain Gromov Witten invariants.  To do this we first restrict to a compact component of the virtual moduli space. Convenient notation is as follows:
 
\begin{defn}Given a family $\hat {\ex B}$ of exploded manifolds, almost complex structures and taming forms, choose a genus $g\in\mathbb N$, a tropical curve $\gamma$ in $\totb{\hat{\ex B}}$ a nonnegative  real number $E$ and a linear map $\beta:H^{2}(\hat{\ex B})\longrightarrow \mathbb R$ where $H^{2}(\hat {\ex B})$ is the DeRham cohomology of $\hat{\ex B}$ as defined in \cite{dre}.
\begin{itemize}
\item $\Msw_{g,[\gamma],\beta}(\hat {\ex B})$  is the moduli stack of curves in $\Msw(\hat{\ex B})$ with genus $g$, and with tropical part isotopic to $[\gamma]$, so that integration over the curve gives the map $\beta:H^{2}(\hat {\ex B})\longrightarrow \mathbb R$.
\item \[\Msw_{g,[\gamma],E}=\coprod_{\beta(\omega)=E}\Msw_{g,[\gamma],\beta}\]
\item $\Mod_{g,[\gamma],\beta}(\hat {\ex B})$ and $\Mod_{g,[\gamma],E}(\hat {\ex B})$ are the virtual moduli spaces of holomorphic curves in $\Msw_{g,[\gamma],\beta}(\hat {\ex B})$ and $\Msw_{g,[\gamma],E}(\hat{\ex B})$ respectively.
\item $\overline{\mathcal M}_{g,[\gamma]}$ is the Deligne Mumford space of genus $g$ curves with punctures labeled by the infinite ends of $\gamma$.
\end{itemize}

We shall use $\Msw_{g,[\gamma],\beta}$ in place of $\Msw_{g,[\gamma],\beta}(\hat{\ex B})$ where no ambiguity is present.

\end{defn}

\begin{defn}\label{gcompactness}
Say that Gromov compactness holds for $(\ex B,J)$ if   the substack of   holomorphic curves  in $\Msw_{g,[\gamma],\beta}(\ex B)$ is compact (in the topology on $\Msw_{g,[\gamma],\beta}$ described in \cite{cem}), and  there are only a finite number of $\beta$ with $\beta(\omega)\leq E$  so that there is a holomorphic curve in $\Msw_{g,[\gamma],\beta}$. 

Say that Gromov compactness holds for a family $\hat {\ex B}\longrightarrow\ex G$ if the following holds: given any compact exploded manifold $\hat {\ex G}'$ with a map $\hat {\ex G}'\longrightarrow \ex G$, let $\hat {\ex B}'\longrightarrow  \ex G'$ be the pullback of our original family $\hat{\ex B}\longrightarrow \ex G$. Then the substack of holomorphic curves in $\Msw_{g,[\gamma],\beta}(\hat {\ex B}')$ is compact and there are only a finite number of $\beta$ with $\beta(\omega)\leq E$  so that there is a holomorphic curve in $\Msw_{g,[\gamma],\beta}(\hat {\ex B}')$.

\end{defn}
Some cases in which the results of \cite{cem} imply  Gromov compactness are discussed in appendix \ref{compactness}.

\

  Suppose that Gromov Compactness holds for $(\ex B,J)$.  Given any $\C\infty1$ map \[\psi:\Msw_{g,[\gamma],\beta}\longrightarrow \ex X\] where $\ex X$ is an exploded manifold (or orbifold), and a closed differential form $\theta\in \Omega^{*}(\ex X)$,  the integral \[\int_{\mathcal M_{g,[\gamma],\beta}}\psi^{*}\theta\] is well defined independent of the choices used in defining $\mathcal M_{g,[\gamma],\beta}$. The same result holds for refined differential forms  $\theta\in\ro^{*}(\ex X)$. (See \cite{dre} for a thorough discussion of the different types of differential forms on exploded manifolds.) In fact, Theorem \ref{pd class} on page \pageref{pd class} implies that there is a closed differential form \[\eta_{g,[\gamma],\beta}\in \ro^{*}_{c}(\ex X)\] Poincare dual to $\psi$  in the sense that  given any closed differential form $\theta\in \ro^{*}(\ex X)$, \[\int_{\mathcal M_{g,[\gamma],\beta}}\psi^{*}\theta=\int_{\ex X}\theta\wedge \eta_{g,[\gamma],\beta}\] Theorem \ref{pd class} also shows that the class of $\eta_{g,[\gamma],\beta}$ in  the homology $\rh^*(\ex X)$ of $(\ro^{*}(\ex X),d)$  is an independent of choices made in its construction.
   The same results hold for $\mathcal M_{g,[\gamma],E}$ because for any fixed $g$ and $\gamma$,  there are only a finite number of $\beta$ with $\mathcal M_{g,[\gamma],\beta}$ nonempty and $\beta(\omega)\leq E$. 
  
  \
  
For example, section \ref{construction of EV} describes a complete exploded manifold $\End_{[\gamma]}\ex B$ which parametrizes the location of the ends of curves in $\Msw_{g,[\gamma],\beta}$. There is a natural evaluation map 
 \[\psi:\Msw_{g,[\gamma],\beta}\longrightarrow \overline{\mathcal M}_{g,[\gamma]}\times \End_{[\gamma]} \ex B\] 

In the case that $\ex B$ is a compact symplectic manifold, $\rh^{*}(\ex B)=H^{*}(\ex B)$ is just the usual DeRham cohomology of $\ex B$, and $\End_{[\gamma]}\ex B$ is just $\ex B^{n}$ where $n$ is the number of ends of $\gamma$. In this case  $[\eta_{g,[\gamma],\beta}]\in H^{*}(\overline{\mathcal M}_{g,n}\times \ex B^{n})$ is a version of the familiar Gromov Witten invariants of $(\ex B,\omega)$. These Gromov Witten invariants satisfy all of  Kontsevich and Mannin's  axioms of Gromov Witten invariants stated in \cite{KM} apart from the `Motivic' axiom which does not make sense outside of the algebraic context. 

In some cases, the Gromov Witten invariants defined in this paper can be verified to coincide with previously defined Gromov Witten invariants of symplectic manifolds defined by   Fukaya and Ono in \cite{FO}, McDuff in \cite{McDuff}, Ruan in \cite{Ruanvirtual}, Liu and Tian  in \cite{Liu-Tian}, Siebert in \cite{Siebert}, and Li and Tian in \cite{Tian-Li}. In particular,  Theorem \ref{obstruction bundle} on page \pageref{obstruction bundle} states roughly that if the moduli space of holomorphic curves in $\Msw_{g,[\gamma],\beta}$ is an orbifold with a nice obstruction bundle, then integration of the pullback of a closed form  $\theta$ on $\Mod_{g,[\gamma],\beta}$ is equal to integration of the wedge product of the Euler class of the obstruction bundle with  the corresponding pullback of $\theta$ to the moduli space of holomorphic curves.   Our invariants do not depend on a choice  of tamed almost complex structure $J$, so if $J$ may be chosen so that Theorem \ref{obstruction bundle} holds and a similar theorem holds for the above constructions, then the different definitions of Gromov Witten invariants will agree. A special case is the case of `transversality' when there is no obstruction bundle required, and all definitions of Gromov Witten invariants agree. 
I expect that in the algebraic case, the definition of Gromov Witten invariants given here will agree with the algebraic definition given by Behrend and Fantechi in \cite{BehrendFantechi}, however this is not proved in this paper.

\

The fact that Gromov Witten invariants do not change in families for which Gromov compactness holds is proved in Theorem \ref{family invariance} on page \pageref{family invariance}. (In fact, just as cohomology does not change locally in a family, but there may be monodromy, Gromov Witten invariants do not change locally, but there may be monodromy in a family with smooth part that is not simply connected.) This invariance allows Gromov Witten invariants of a compact symplectic manifold to be calculated by deforming the symplectic manifold in a connected family of exploded manifolds to an exploded manifold in which the calculation is easier. 

Gromov Witten invariants of an exploded manifold $\ex B$ with nontrivial tropical part $\totb{\ex B}$ are often easier to calculate because of  Theorem \ref{gluing formula} on page \pageref{gluing formula}, which is a generalization of the symplectic sum formula for Gromov Witten invariants. In particular, we may break up Gromov Witten invariants into contributions coming from tropical curves. 

\begin{defn} Given some moduli stack or virtual moduli space $\mathcal M$ of curves, define $\mathcal M\rvert_{\gamma}$ to be the subset of $\mathcal M$   consisting of curves with tropical part isomorphic to $\gamma$. 
\end{defn}

Then the integral over the virtual moduli space breaks up into a sum of contributions from each tropical curve,

\[\int_{\Mod_{g,[\gamma],\beta}}\psi^{*}\theta=\sum_{\gamma_{i}}\int_{\Mod_{g,[\gamma],\beta}\rvert_{\gamma_{i}}}\psi^{*}\theta\]
where only a finite number of tropical curves $\gamma_{i}$ have a nonzero contribution to the sum.
If $\theta\in\rof(\ex X)$ is a closed differential form generated by functions in the sense of definition \ref{rof}, then this integral $\int_{\Mod_{g,[\gamma],\beta}\rvert_{\gamma}}\psi^{*}\theta$ is well defined independent of the choices in the construction of $\mathcal M_{g,[\gamma],\beta}$. (All differential forms on a smooth manifold are generated by functions. Examples of  interesting differential forms generated by functions on an exploded manifold are the Poincare duals to points, symplectic taming forms and  Chern classes constructed using the Chern-Weyl construction.) 

Theorem \ref{gluing cobordism} describes the cobordism class of $\Mod_{g,[\gamma],E}\rvert_{\gamma}$ in terms of a fiber product of moduli spaces corresponding to each vertex of $\gamma$. In particular, in the case that $\ex X=\End_{[\gamma]}\ex B$, Theorem \ref{gluing cobordism} implies a generalization of the symplectic sum formula for Gromov Witten invariants, 
  Theorem \ref{gluing formula} on page \pageref{gluing formula}. In that case, 
  \[\int_{\mathcal M_{g,[\gamma],E}\rvert_{\gamma}}\psi^{*}\theta=\frac {\prod m_{e}}{\abs{\Aut \gamma}}\sum_{ g_{v}, E_{v}}\int_{\prod_{e\subset\gamma}\Mod_{\gamma_{e}} }\theta\bigwedge_{v\in\gamma}\eta_{g_{v},[\gamma_{v}],E_{v}}\]
where:
\begin{itemize}
\item The sum is over choices of $g_{v}\in \mathbb N$ and $E_{v}\in[0,\infty)$ for each vertex $v$ of $\gamma$ so that $\sum E_{v}=E$ and $\sum g_{v}$ is equal to $g$ minus the genus of $\gamma$. Gromov compactness assumptions give that only a finite number of terms in this sum will be nonzero.
\item $m_{e}$ indicates a multiplicity of each internal edge of $\gamma$ (and is $0$ if an edge is sent to a point).
\item For each edge $e$ of $\gamma$, $\mathcal M_{\gamma_{e}}$ is a manifold parametrizing the possible maps of $\et 1{(0,l)}$ into $\ex B$ which have tropical part equal to $e$.
\item $\theta$ on the right hand side of the equation indicates the pullback of $\theta\in \rof^{*}(\End_{[\gamma]} \ex B)$ under a map $\prod_{e\subset\gamma}\Mod_{\gamma_{e}}\longrightarrow \End_{[\gamma]}\ex B$.

\item To the strata of $\ex B$ containing each vertex $v$ of $\gamma$, an exploded manifold $\check{\ex B}_{v}$ called a tropical completion of the strata is associated. For this gluing formula to hold, Gromov compactness is assumed to hold for $\check{\ex B}_{v}$.  There is a tropical curve $\gamma_{v}$ in $\totb{\check{\ex B}_{v}}$ with one vertex (corresponding to $v$) and one edge corresponding to each end of an edge $e$ of $\gamma$ attached to $v$. 
\item $\eta_{g,[\gamma_{v}],E_{v}}$ is the Poincare dual to the map 

\[\Mod_{g,[\gamma_{v}],E_{v}}(\check{\ex B}_{v})\longrightarrow \End_{[\gamma_{v}]}\check{\ex B}_{v}\]
pulled back under a map 
\[\prod_{e\subset\gamma}\Mod_{\gamma_{e}}\longrightarrow \End_{[\gamma_{v}]}\check{\ex B}_{v}\]

\end{itemize} 

These $\eta_{g,[\gamma_{v}],E_{v}}$ can be regarded as relative Gromov Witten invariants.  The exploded manifold $\check{\ex B}_{v}$ has tropical part equal to a cone. In the case that the corners of the polytopes in $\totb{\ex B}$ are standard, the smooth part $\totl{\check{\ex B}_{v}}$ of $\check{\ex B}_{v}$ is a symplectic manifold with a tamed almost complex structure $J$, and each strata of $\totl{\check{\ex B}_{v}}$ is a (pseudo)holomorphic submanifold. If $J$ is integrable, $\totl{\check{\ex B}_{v}}$ is a complex manifold with normal crossing divisors, and $\check{\ex B}_{v}$ is some $\ex T^{n}$ bundle over the explosion of $\totl{\check{\ex B}_{v}}$.

Given a compact complex manifold $M$ with normal crossing divisors given by embedded complex codimension $1$ submanifolds which intersect each other transversely and a symplectic taming form $\omega$, \cite{cem} implies that Gromov compactness holds for  the explosion $\expl M$ of $M$ described in \cite{iec}. The Gromov Witten invariants of $M$ relative to its normal crossing divisors can be defined as the Gromov Witten  invariants of $\expl M$. Given $n$  $\mathbb C^{*}$ bundles over $M$, there is a natural corresponding $\ex T^{n}$ bundle $\ex B$ over $\expl M$. Appendix \ref{compactness} implies that given any lift of the  complex structure on $\expl M$ to an almost complex structure on $\ex B$, Gromov Compactness will hold on $\ex B$. The Gromov Witten invariants of $\ex B$ do not depend on this lift of almost complex structure, so we may define the relative Gromov Witten invariants of the $n$ $\mathbb C^{*}$ bundles over  $M$ with its normal crossing divisors to be the Gromov Witten invariants of $\ex B$. 

The symplectic analogue of explosion described in \cite{cem} allows Gromov Witten invariants of (some number of $\mathbb C^{*}$ bundles over) symplectic manifolds relative to orthogonally intersecting codimension $2$ symplectic submanifolds to be defined analogously.  The gluing formula from Theorem \ref{gluing formula} gives a lot of structure to these relative Gromov Witten invariants, which makes their computation easier. Other relationships between relative Gromov Witten invariants follow from Theorem \ref{refinement theorem} on page \pageref{refinement theorem} which implies that Gromov Witten invariants do not change under refinement of exploded manifolds. For example, this means that the Gromov Witten invariants of all toric manifolds of a given dimension relative to their toric boundary divisors can be viewed as being the same.

\section{Structure of the $\dbar$ equation on $\mathcal M^{\omega}$ }\label{structure section}

\

\

$\Msw$, (defined more precisely shortly) is the moduli stack of $\C\infty1$ curves for which the integral of some two form  $\omega$ is positive on each smooth component equal to a sphere with at most one puncture and nonnegative on all other smooth components. We shall describe the structure of $\Msw$ using the notion of a core family below, and describe the $\dbar$ equation on $\Msw$ using the notion of an obstruction model.   

\subsection{the functors $\ex F$ and $\ex C$}

\

\

This paper studies  families of  holomorphic curves in a smooth family of targets in the exploded category, \[\pi_{\ex G}:(\hat {\ex B},J,\omega)\longrightarrow\ex G\] where each fiber $\ex B$ is a complete, basic exploded manifold with a civilized almost complex structure $J$ tamed by a symplectic form $\omega$ (using terminology from \cite{iec} and \cite{cem}). We will often talk about $\C\infty1$ families of curves $\hat f$ in $\hat {\ex B}\longrightarrow \ex G$ which will correspond to commutative diagrams
\[\begin{array}{ccc}\ex C(\hat f)&\xrightarrow {\hat f}&\hat{\ex B}
\\\downarrow\pi_{\ex F(\hat f)}& &\downarrow\pi_{\ex G} 
\\\ex F(\hat f)&\longrightarrow &\ex G\end{array}\]

where $\pi_{\ex F(\hat f)}:\ex C(\hat f)\longrightarrow \ex F(\hat f)$ is a family of curves (as defined in \cite{iec}). Where no ambiguity is present, this family is just referred to as $\pi_{\ex F}:\hat{\ex C}\longrightarrow\ex F$, however, we shall also think of $\ex F$ and $\ex C$ as functors in the following way:
As noted in \cite{iec}, families of $\C\infty1$ curves in $\hat{\ex B}\longrightarrow \ex G$ form a category with morphisms $\hat f\longrightarrow \hat g$ given by commutative diagrams
\[\begin{array}{ccccc}\ex F(\hat f)&\longleftarrow &\ex C(\hat f)&\longrightarrow &\hat {\ex B}
\\\downarrow && \downarrow&&\downarrow
\\ \ex F(\hat g)&\longleftarrow&\ex C(\hat g)&\longrightarrow &\hat {\ex B}\end{array}\]
so that restricted to each fiber of $\pi_{\ex F(\hat f)}$ and $\pi_{\ex F(\hat g)}$, $\ex C(\hat f)\longrightarrow\ex C(\hat g)$ is a holomorphic isomorphism. In this way, $\ex C$ and $\ex F$ can be regarded as functors from the category of $\C\infty1$ families of curves to the category of $\C\infty1$ exploded manifolds. Use the notation $\Ms(\hat {\ex B})$ for the moduli stack of $\C\infty1$ curves, which  is the category of $\C\infty1$ families of curves together with the functor $\ex F$. When no ambiguity is present, simply use the notation $\Ms$.

\subsection{$\omega$-positive curves}

\

\

\begin{defn} Call a $\C\infty1$ curve $f:\ex C(f)\longrightarrow \hat{\ex B}$ $\omega$-positive if for each smooth strata of $\ex C(f)$ which is a sphere with at most one puncture, the integral of $\omega$ over that strata is positive, and the integral of $\omega$ over every smooth strata of $\ex C(f)$ is non-negative.

Let $\Msw$ be the substack of $\Ms$ consisting of families of  $\omega$-positive curves.
\end{defn}

Note that if $\omega$ is a symplectic form that tames the almost complex structure on $\ex B$, every stable holomorphic curve is $\omega$-positive. Restricting to $\omega$ positive curves is a technical assumption chosen for the property that given any $\omega$-positive curve $f$ which is transverse to some other complete  map to $\ex B$, then there exists a neighborhood of $f$ in $\Msw$ so that all curves in  that neighborhood are also transverse to this map.

\subsection{evaluation maps and adding extra marked points to families}

\

\

The following theorem will be proved on page \pageref{moduli}. The notation $\bar{\mathcal M}$ refers to the usual Deligne Mumford space considered as a complex orbifold with normal crossing divisors given by its boundary components, and $\M$ is the explosion of this discussed in \cite{iec}. The map $\M^{+1}\longrightarrow\M$ is the explosion of the usual forgetful map from the moduli space of curves with one extra marked point to the moduli space of curves forgetting that extra marked point. The following theorem implies that $\M$ represents the moduli stack of stable $\C\infty1$ curves, and that the map  $\M^{+1}\longrightarrow \M$  is a universal family of stable $\C\infty1$ curves.

\begin{thm}\label{moduli1} Consider a $\C\infty1$ family of exploded curves $(\hat{\ex C},j)\longrightarrow \ex F$ so that each curve has $2g+n\geq 3$ where $g$ is the genus and $n$ is the number of punctures. Then there exists a unique fiber-wise holomorphic map 
\[\begin{array}{ccc}
(\hat {\ex C},j)&\longrightarrow & (\M^{+1},j)
\\ \downarrow & &\downarrow
\\ \ex F&\xrightarrow{ev^{0}} & \M
\end{array}\]
so that the map on each fiber $\ex C$ factors into a degree one holomorphic map to a stable exploded curve $\ex C'$ and a holomorphic map from $\ex C'$ to a fiber of $\M^{+1}$ given by quotienting $\ex C'$ by its automorphism group. 

The above maps all have regularity $\C\infty 1$.
\end{thm}

\
 
 In what follows, we define an `evaluation map' for a family of curves using a functorial construction of a family of curves $\hat f^{+n}$ with $n$ extra punctures from a given family of curves $\hat f$. 
 
\begin{defn}
Given a submersion $f:\ex D\longrightarrow \ex E$, use the following notation for the fiber product of $\ex D$ over $\ex E$ with itself $n$ times:
\[\sfp {\ex D}{\ex E}n:=\ex D\fp ff \ex D\fp ff \dotsb \fp ff \ex D\]

\end{defn}

\begin{defn}\label{f^{+n}}
Given a family of curves $\hat f$ in $\hat{\ex B}\longrightarrow \ex G$ and $n\in \mathbb N$, the family $\hat f^{+n}$ is a family of curves with $n$ extra punctures
\[\begin{array}{ccc}\ex C(\hat f^{+n})&\xrightarrow{\hat f^{+n}}&\sfp{\hat{\ex B}}{\ex G}{n+1}
\\ \downarrow\pi_{\ex F(\hat f^{+n})} & &\downarrow 
\\ \ex F(\hat f^{+n}) &\xrightarrow{\hat f^{+(n-1)}}&\sfp{\hat {\ex B}}{\ex G}{n}\end{array}\]
 satisfying the following conditions
\begin{enumerate}
\item The family of curves $\hat f^{+0}$ is $\hat f$.
\item The base of the family $\hat f^{+n}$ is the total space of the family $\hat f^{+(n-1)}$.
\[\ex F(\hat f^{+n})=\ex C(\hat f^{+(n-1)})\]
\item The fiber of $\pi_{\ex F(\hat f^{+1})}:\ex C(\hat f^{+1})\longrightarrow \ex F(\hat f^{+1})$ over a point $p\in\ex F(\hat f^{+1})=\ex C(\hat f)$ is equal to the fiber of $\pi_{\ex F(\hat f)}:\ex C(\hat f)\longrightarrow \ex F$ containing $p$ with an extra puncture at the point $p$.

\item The family of curves $\hat f^{+n}$ is $\lrb{\hat f^{+(n-1)}}^{+1}$
\item There exists a fiberwise holomorphic degree $1$ map 
\[\begin{array}{ccc}\ex C(\hat f^{+1})&\longrightarrow &{\ex C}(\hat f)\fp{\pi_{\ex F}}{\pi_{\ex F}}{\ex C}(\hat f)
\\\downarrow &&\downarrow
\\ \ex C(\hat f)&\xrightarrow{\id}&\ex C(\hat f) \end{array}\] so that the following diagram commutes.
\[\begin{array}{ccc}
{\ex C}(\hat f^{+1})&&
\\ \downarrow & \hat f^{+1}\searrow &
\\ {\ex C}(\hat f)\fp{\pi_{\ex F}}{\pi_{\ex F}}{\ex C}(\hat f)  &\xrightarrow{\hat f\times \hat f} & \hat{\ex B}\fp{\pi_{\ex G}}{\pi_{\ex G}}{\hat {\ex B}}
\\ \downarrow & & \downarrow
\\ {\ex C}(\hat f)&\xrightarrow {\hat f}&\hat{\ex B}
\\ \downarrow\pi_{\ex F} & & \downarrow\pi_{\ex G}
\\ \ex F & \longrightarrow & \ex G
\end{array}\]

\end{enumerate}
 
\end{defn}
 
 It is shown in appendix \ref{construct f^{+n}} on page \pageref{construct f^{+n}} that $\hat f^{+n}$ exists and is smooth or $\C\infty 1$ if $\hat f$ is. The above conditions imply that the  map $\hat f^{+n}:\ex C(\hat f^{+n})\longrightarrow \sfp{\hat {\ex B}}{\ex G}{n+1}$ factors as
 \[\ex C(\hat f^{+n})\longrightarrow\sfp{\ex C(\hat f)}{\ex F}{n+1}\longrightarrow \sfp{\hat {\ex B}}{\ex G}{n+1}\]
 where the second map $\sfp{\ex C(\hat f)}{\ex F}{n+1}\longrightarrow \sfp{\hat {\ex B}}{\ex G}{n+1}$ is given by the $(n+1)$-fold product of $\hat f$, and the first  map $\ex C(\hat f^{+n})\longrightarrow \sfp{\ex C(\hat f)}{\ex F}{n+1}$ is a degree one map.

   The construction of $\hat f^{+n}$ is functorial, so given a map of families of curves $\hat f\longrightarrow \hat g$, there is an induced map $\hat f^{+n}\longrightarrow \hat g^{+n}$. This map  $\hat f^{+n}\longrightarrow \hat g^{+n}$ is compatible with the naturally induced map $\sfp{\ex C(\hat f)}{\ex F(\hat f)}{n+1}\longrightarrow \sfp{\ex C(\hat g)}{\ex F(\hat g)}{n+1}$.
 
 \

 Combining $\hat f^{+(n-1)}$ with the map $ev^{0}:\ex F(\hat f^{+n})\longrightarrow \M$ given by Theorem \ref{moduli1} when $n$ is large enough, we get the evaluation map
\[ ev^{+n}(\hat f):=(ev^{0},\hat f^{+n-1}):\ex F(\hat f^{+n})\longrightarrow \M\times\sfp{ \hat{\ex B}}{\ex G}n\]

\subsection{core families}

\

\

 The following notion of a core family gives a way of  locally describing  the moduli stack $\Msw$ of $\omega$-positive $\C\infty1$ curves. A notion such as this is necessary, as the `space' of $\omega$-positive curves with in $\hat{\ex B}\longrightarrow \ex G$ of a given regularity can not in general be locally modeled on even an orbifold version of a Banach space - this is because the domain of curves that we study are not fixed, and because of  phenomena which would be called bubble and node formation in the setting of smooth manifolds.  (The moduli stack of curves satisfying a slightly stronger $\omega$-positivity condition could be described as a `space' by using an adaption to the exploded setting  of the theory of polyfolds being developed by Hofer, Wysocki and Zehnder in a series of papers including \cite{polyfold1}. An adaption of the theory of polyfolds to the exploded setting is a worthwhile direction for further research which is not explored in this paper.)

\begin{defn}\label{core family}
A core family of curves,  $(\hat f/G,\{s_{i}\},F)$ for an open substack $\mathcal O$ of $\Msw\lrb{\hat {\ex B}}$ is: 
\begin{itemize}
 \item a basic $\C\infty 1$ family $\hat f$ of $\omega$-positive, stable curves with a group $G$ of automorphisms,
 \[\begin{array}{ccc}\ex C(\hat f)&\xrightarrow{\hat f}&\hat{\ex B}
 \\ \downarrow & &\downarrow
 \\ \ex F(\hat f)&\longrightarrow &\ex G
 \end{array}\]
 \item a nonempty finite collection of $\C\infty 1$ `marked point' sections $s_{i}:\ex F(\hat f)\longrightarrow{\ex C}(\hat f)$ which do not intersect each other, and which do not intersect the edges of the curves in $ \ex C(\hat f)$.  
 \item a $G$ invariant $\C\infty1$ map,
 \[\begin{array}{ccc} \hat f^{*}{T_{vert}\hat{\ex B}}&\xrightarrow{F} &\hat{\ex B} 
\\ \downarrow &&\downarrow
\\ \ex F(\hat f)&\longrightarrow &\ex G
 \end{array}\]
where $T_{vert}\hat {\ex B}$ indicates the vertical tangent space of the family $\hat{\ex B}\longrightarrow \ex G$.

  \end{itemize}
   so that
  \begin{enumerate}
 
 \item  For all curves $f$ in $\hat f$, there are exactly $\abs G$ maps $f\longrightarrow \hat f$, and the action of $G$ on the set of maps of $f$ into $\hat f$ is free.
 \item For all curves $f$ in $\hat f$,  the smooth part of the domain $\ex C(f)$ with the extra marked points from $\{s_{i}\}$ has no automorphisms.
 \item The action of $G$ preserves the set of sections $\{s_{i}\}$, so there is some  action of $G$ as a permutation group on the set of indices $\{i\}$ so that for all $g\in G$ and $s_{i}$,
 \[s_{i}\circ g=g\circ s_{g(i)}\]
 where  the action of $g$ is on $\ex F(\hat f)$, $\ex C(\hat f)$ or the set of indices $\{i\}$ as appropriate. 
 \item \label{local core crit} \begin{enumerate}\item\label{lcca} There exists a neighborhood $U$ of the image of the section
\[s:\ex F(\hat f)\longrightarrow \ex F(\hat f^{+n})\]
defined by the $n$ sections $\{s_{i}\}$ so that  
\[ev^{+n}(\hat f):\ex F(\hat f^{+n})\longrightarrow \M\times\sfp{\hat {\ex B}}{\ex G}n\]
restricted to $U$ is an equi-dimensional embedding

\item\label{lccb} The tropical part of $ev^{+n}\circ s$ is a complete map, and restricted to any polytope in $\totb{\ex F(\hat f)}$ is an isomorphism onto a strata of the image of $\mathcal O$ in  $\totb{\M\times\sfp{\hat {\ex B}}{\ex G}n}$ under $\totb{ev^{+n}}$.
\end{enumerate}
\item\label{F crit} 
 \begin{enumerate}\item $F$ restricted to the zero section is equal to $\hat f$, 
 \item $TF$ restricted to the natural inclusion of $\hat f^{*}{T_{vert}\hat{\ex B}}$ over the zero section is equal to the identity,
 \item  $TF$ restricted to the vertical tangent space at any point of $\hat f^{*}T_{vert}\hat{\ex B}$ is injective.
 \end{enumerate}

  \item \label{core family main}
   given any $\C\infty1$ family  $\hat f'$ in $\mathcal O$, there exists a unique $\C\infty1$ fiber-wise holomorphic map
  \[\begin{array}{ccc}(\ex C(\hat f'),j)&\xrightarrow{\Phi_{\hat f'}}&({\ex C}(\hat f),j)/G
  \\ \downarrow & & \downarrow
  \\ \ex F(\hat f')&\longrightarrow &\ex F(\hat f)/G\end{array}\]
  and unique $\C\infty 1$ section 
  \[\psi_{\hat f'}:\ex C(\hat f')\longrightarrow\Phi_{\hat f'}^{*}\lrb{\hat f^{*}T_{vert}\hat{\ex B}}\]
    which vanishes on the pullback of marked points, so that 
   
\[\hat f'=F\circ\psi_{\hat f'}\]
   \end{enumerate}
\end{defn}

 The last  condition can be roughly summarized by saying that a $\C\infty1$ family in $\mathcal O$ is equivalent to a $\C\infty1$ map to $\ex F(\hat f)/G$ and a sufficiently small $\C\infty1$ section of some vector bundle associated to this map.
Theorem \ref{core criteria} on page \pageref{core criteria} states that if the first five conditions hold, there exists some $\mathcal O$ so that the  condition holds.
 Proposition \ref{smooth model family} stated on page \pageref{smooth model family} constructs a core family containing any given stable holomorphic curve with at least one smooth component.

\begin{thm} \label{local description}
Given any stable holomorphic $f$ with at least one smooth component in a basic family of exploded manifolds, $\hat {\ex B}\longrightarrow  \ex G$,
 there exists an open neighborhood $\mathcal O$ of $f$ in $\Msw$ and  a core family $(\hat f/G,\{s_{i}\},F)$ for $\mathcal O$ containing $f$.
%
%


%
\end{thm}

Given any core family $(\hat f/G,\{s_{i}\},F)$ there is a canonical orientation of $\hat f$ given as follows. The sections $\{s_{i}\}$ define a section \[s:\ex F(\hat f)\longrightarrow \ex F(\hat f^{+n})\] so that in a neighborhood of the image of this section $s$,  the map 
\[ev^{+n}(\hat f):\ex F(\hat f^{+n})\longrightarrow\M\times \sfp{\hat{\ex B}}{\ex G}n\] is an equi-dimensional embedding. There is a canonical orientation of $\M\times \sfp{\hat{\ex B}}{\ex G}n$ relative to $\ex G$ given by the orientation from the complex structure on $\M$ and the almost complex structure $J$ on the fibers of $\hat {\ex B}\longrightarrow \ex G$. Therefore, there is a canonical  orientation of $\ex F(\hat f^{+n})$ relative to $\ex G$ so that $ev^{+n}(\hat f)$ preserves orientations in a neighborhood of the above section $s$. This in turn gives a natural orientation of $\ex F(\hat f)$ relative to $\ex G$, because the fibers of each of the maps $\ex F(\hat f^{+k})\longrightarrow \ex F(\hat f^{+(k-1)})$ are complex.

%



\subsection{trivializations and (pre)obstruction models}

\

\

Now we shall start describing the $\dbar$ equation on the moduli stack of $\C\infty1$ curves.

\begin{defn}
Given a smooth (or $\C\infty 1$) family,  

  \[\begin{array}{ccc} 
  (\ex {\hat C},j) & \xrightarrow{\hat f} &(\hat{\ex B},J)
  \\ \ \ \ \  \downarrow\pi_{\ex F} & &\ \ \ \ \downarrow\pi_{\ex G}
  \\  \ex F &\longrightarrow &\ex G\end{array}\]
  
\begin{enumerate}
\item Use the notation 
\[T_{vert}\hat {\ex C}:=\ker d\pi_{\ex F}\subset T\hat{\ex C}\ \ \ \ \ \ \ \ \text{ and }T_{vert}\hat {\ex B}:=\ker d\pi_{\ex G}\subset T\hat{\ex B}\]
  Let $T^{*}_{vert}\hat{\ex C}$ be the dual of $T_{vert}\hat{\ex C}$. Of course, $T^{*}_{vert}\hat {\ex C}$ may also be described as $T^{*}\hat{\ex C}/\pi_{\ex F}^{*}T^{*}\ex F$
\item Define
\[\dvert \hat f: T_{vert}\hat {\ex C}\longrightarrow T_{vert}\hat{\ex B}\]
to be $d\hat f$ restricted to the vertical tangent space, $T_{vert}\hat {\ex C}\subset T\hat{\ex C}$.
\item  Define
\[\dbar \hat f: T_{vert}\hat {\ex C}\longrightarrow T_{vert}\hat{\ex B}\]
\[\dbar \hat f:=\frac 12\lrb{ \dvert f +J\circ \dvert f\circ j}\]
 Consider 
 \[\dbar \hat f \in \Gamma \lrb{T^{*}_{vert}\hat{\ex C} \otimes \hat f^{*}T_{vert}\hat{\ex B}}^{0,1}\]

 \end{enumerate}
\end{defn}

 As the above bundle is cumbersome to write out in full, and can be considered as the pull back of a vector bundle $Y$ over the moduli stack of $\C\infty 1$ curves,  use the following notation: 
 \begin{defn}\label{Ydef} Use the notation
 $\Y{\hat f}$ to denote $\lrb{T_{vert}^{*}\hat{\ex C}\otimes {\hat f^{*}T_{vert}\hat{\ex B}}}^{0,1}$, which is the sub vector bundle of $T_{vert}^{*}\hat{\ex C}\otimes_{\mathbb R} {\hat f^{*}T_{vert}\hat{\ex B}}$ consisting of vectors so that the action of $J$ on the second factor is equal to $-1$ times the action of $j$ on the first factor.
 \end{defn}

Note that given any map  of families of curves $\hat f\longrightarrow \hat g$, there is a corresponding map of vector bundles $\Y{\hat f}\longrightarrow \Y{\hat g}$.

\

\begin{defn}\label{trivialization def} Given a $\C\infty1$ family, 
\[\begin{array} {ccc}\hat {\ex C}&\xrightarrow{\hat f}& \hat{\ex B}
\\ \downarrow & &\downarrow\pi_{\ex G}
\\ \ex F&\longrightarrow &\ex G\end{array}\]
 a choice of trivialization $(F,\Phi)$ is 
 \begin{enumerate}
 \item
 a $\C\infty1$ map 
 \[\begin{array}{ccc} \hat f^{*}{T_{vert}\hat{\ex B}}&\xrightarrow{F} &\hat{\ex B} 
\\ \downarrow &&\downarrow
\\ \ex F&\longrightarrow &\ex G
 \end{array}\]
 so that \begin{enumerate}\item $F$ restricted to the zero section is equal to $\hat f$, 
 \item $TF$ restricted to the natural inclusion of $\hat f^{*}{T_{vert}\hat{\ex B}}$ over the zero section is equal to the identity,
 \item  $TF$ restricted to the vertical tangent space at any point of $\hat f^{*}T_{vert}\hat{\ex B}$ is injective.
 \end{enumerate} 
 \item A $\C\infty1$ isomorphism  from the bundle $F^{*}{T_{vert}\hat{\ex B}}$  to the vertical tangent bundle of $\hat f^{*}T_{vert}\hat{\ex B}$ which preserves $J$, and which restricted to the zero section of $\hat f^{*}T_{vert}\hat {\ex B}$ is the identity.
 
  In other words, if $\pi:\hat f^{*}T_{vert}\hat{\ex B}\longrightarrow \hat {\ex C} $ denotes the vector bundle projection, a $\C\infty1$ isomorphism between $F^{*}{T_{vert}\hat{\ex B}}$ and $\pi^{*}\hat f^{*}{T_{vert}\hat{\ex B}}$ which preserves the almost complex structure $J$ on ${T_{vert}\hat{\ex B}}$. This can be written as a $\C\infty1$ vector bundle map 
 \[\begin{array}{ccc}F^{*}{T_{vert}\hat{\ex B}}&\xrightarrow {\Phi}&\hat f^{*}{T_{vert}\hat{\ex B}}
 \\ \downarrow & & \downarrow
 \\ \hat f^{*}{T_{vert}\hat{\ex B}}&\longrightarrow &\hat{\ex C}\end{array}\]
 which is the identity when the vector bundle  $F^{*}{T_{vert}\hat{\ex B}}\longrightarrow \hat f^{*}{T_{vert}\hat{\ex B}}$ is restricted to the zero section of $\hat f^{*}T_{vert}\hat{\ex B}$.
\end{enumerate}
 
 A trivialization allows us to define $\dbar$ of a section $\nu :\hat{\ex C}\longrightarrow\hat f^{*}{T_{vert}\hat{\ex B}}$ as follows: $F\circ \nu$ is a family of maps $\hat {\ex C}\longrightarrow \hat {\ex B}$, so $\dbar (F\circ \nu)$ is a section of $\Y{F\circ\nu}=\lrb{{T_{vert}^{*}\hat{\ex C} }\otimes {(F\circ \nu)^{*}T_{vert}\hat{\ex B}}}^{(0,1)}$. Applying the map $\Phi$ to the second component of this tensor product gives an identification of $\Y{F\circ \nu}$ with $\Y{\hat f}$, so we may consider $\dbar(F\circ \nu)$ to be a section of $\Y{\hat f}$. Define $\dbar \nu$ to be this section of $\Y{\hat f}$.

\end{defn}

For example, we may construct a trivialization by extending $\hat f$ to a map $F$ satisfying the above conditions (for instance by choosing a smooth connection on $T_{vert}\hat {\ex B}$ and reparametrising the exponential map on a neighborhood of the zero section in  $f^{*}T_{vert}\hat {\ex B}$), and letting $\Phi$ be given by parallel transport along a linear path to the zero section using a smooth $J$ preserving connection on $T_{vert}\hat{\ex B}$. 

Given a choice of  trivialization for $\hat f$ and a $\C\infty1$ section $\nu$ of $\hat f$, there is an induced choice of trivialization for the family $F(\nu)$, described in \cite{reg}.

\begin{defn}\label{pre obstruction model}
A  $\C\infty1$  pre obstruction model $(\hat f,V,F,\Phi,\{s_{i}\})$, is given by
\begin{enumerate}
\item
 a $\C\infty1$ family

\[\begin{array}{ccc} 
  (\ex {\hat C},j) & \xrightarrow{\hat f} &(\hat{\ex B},J)
  \\ \ \ \ \  \downarrow\pi_{\ex F} & &\ \ \ \ \downarrow\pi_{\ex G}
  \\  \ex F &\longrightarrow &\ex G\end{array}\]

\item a choice of trivialization $(F,\Phi)$ for $\hat f$ in the sense of definition \ref{trivialization def}.
 
 \item a finite collection $\{s_{i}\}$ of extra marked points on $\hat{\ex C}$ corresponding to $\C\infty1$ sections $s_{i}:\ex F\longrightarrow \hat{\ex C}$, so that restricted to any curve $\ex C$ in $\hat {\ex C}$, these marked points are distinct and contained inside the smooth components of $\ex C$.
 
 \item a vector bundle $V$ over $\ex F$
 \[ \begin{split}  & V \\ & \downarrow \\ &\ex F\end{split}\]
\item  a smooth or $\C\infty1$  map of vector bundles over $\hat{\ex C}$
\[\pi^{*}_{\ex F}(V)\longrightarrow \Y{\hat f}:=\lrb{{T_{vert}^{*}\hat{\ex C}}\otimes {\hat f^{*}T_{vert}\hat{\ex B}}}^{0,1}\]
which vanishes on the edges of curves in  $\hat{\ex C}\longrightarrow \ex F$.
The above map must be non trivial in the sense that for any nonzero vector in $V$, there exists a choice of lift to a vector in $\pi^{*}_{\ex F}(V)$ which is not sent to $0$. (Said differently, a point $(f,v)\in V$ corresponds to a section of $\pi^{*}_{\ex F}(V)$ over the curve $\ex C(f)$ over $f$. This section is  sent by the above map to a section of  the  bundle $\Y{ f}$. This section is the zero section if and only if $v$ is $0$.)
\end{enumerate}
\end{defn}

Note that a $\C\infty1$ section $\ex F\longrightarrow V$ corresponds to a $\C\infty1$ section  $\hat{\ex C}\longrightarrow \Y{\hat f}$. Call such a section  $\hat{\ex C}\longrightarrow\Y{\hat f}$ `a section of $V$'. We shall usually use the notation $(\hat f,V)$ for a pre obstruction bundle. 

\begin{defn}Given any family $\hat f$, with a collection $\{s_{i}\}$ of extra marked points on $\ex C(\hat f)$,  let  $X^{\infty,\underline1}(\hat f)$ indicate the space of $\C\infty1$ sections of $ f^{*}T_{vert}\hat{\ex B}$ which vanish on the extra marked points $\{s_{i}\}$ on ${\ex C}(\hat f)$, and let $Y^{\infty,\underline 1}(\hat f)$ indicate the space of $\C\infty1$ sections of $\Y {\hat f}$ which vanish on edges of curves in  ${\ex C}(\hat f)$.
\end{defn}

Note that both  $X^{\infty,\underline 1}(\hat f)$ and $Y^{\infty,\underline1}(\hat f)$ are complex vector spaces because they consist of sections of complex vector bundles.

We may restrict any pre obstruction bundle $(\hat f,V)$ to a single curve $f$ in $\hat f$. The restriction of $V$ to this curve $f$ can be regarded as a finite dimensional linear subspace  $V(f)\subset Y^{\infty,\underline 1}(f)$. 

 Let
 $D\dbar(f):X^{\infty,\underline 1}(f)\longrightarrow Y^{\infty,\underline 1}(f)$ indicate the (directional) derivative of $\dbar$ at $0\in X^{\infty,\underline 1}(f)$. We are most interested in pre obstruction bundles $(\hat f,V)$ containing curves $f$ that $D\dbar(f)$ is injective and has image complementary to $V(f)$.

\

In what follows, we define obstruction models $(\hat f/G,V)$ which can be regarded as giving a kind of $\C\infty1$ Kuranishi structure to the moduli stack of holomorphic curves. Roughly speaking, an obstruction model $(\hat f/G,V)$  is a $G$-invariant pre obstruction model $(\hat f,V)$ so that $\hat f/G$ is a core family, $\dbar \hat f$ can be regarded as a $\C\infty1$  section $\dbar :\ex F(\hat f)\longrightarrow V$, and so that perturbing the $\dbar$ equation can locally be modeled on perturbing this section of $V$.

\begin{defn}\label{obstruction model}

An obstruction model $(\hat f/G,V)$ is core family $(\hat f/G,\{s_{i}\},F)$ together with a compatible $G$-invariant trivialization $(F,\Phi)$ and obstruction bundle $V$ making a pre obstruction model $(\hat f,V,F,\Phi,\{s_{i}\})$ so that:
\begin{itemize}
\item $\dbar \hat f$ is a section of $V$.
\item \[D\dbar(f):X^{\infty,\underline 1}(f)\longrightarrow Y^{\infty,\underline 1}(f)\]
is injective and has image complementary to $V(f)$ for all $f$ in $\hat f$.
\end{itemize} 

Say that an obstruction model $(\hat f/G,V)$ is extendible if it is the restriction of some larger obstruction model $(\hat f'/G,V')$ to a compactly contained sub family $\hat f$ of $\hat f'$.

\end{defn}

In Definition \ref{obstruction substack} below, we shall say what it means for $(\hat f/G,V)$ to be an obstruction model for a substack $\mathcal O\subset\Msw$.
The existence of obstruction models is proved on page \pageref{construct obstruction model}. It follows from the existence of core families and the results of \cite{reg}.

To describe the importance of obstruction models, we shall need the notion of a simple perturbation below.

\begin{defn}\label{simple def} A simple perturbation parametrized by a family 

  \[\begin{array}{ccc} 
  \ex {\hat C} & \xrightarrow{\hat f} &\hat{\ex B}
  \\ \ \ \ \  \downarrow\pi_{\ex F} & &\ \ \ \ \downarrow\pi_{\ex G}
  \\  \ex F &\longrightarrow &\ex G\end{array}\]
is a section $\mathfrak P$ of the bundle 
 $\lrb{{T_{vert}^{*}\hat{\ex C}}\otimes {T_{vert}\hat{\ex B}}}^{(0,1)}$ over $\hat{\ex C}\times \hat {\ex B}$ with the same regularity as $\hat f$ which vanishes on all edges of the curves which are the fibers of $\hat {\ex C}$. 
 
The topology on the space of simple perturbations parametrized by $\hat f$ is the corresponding topology on the space of sections. Say that a simple perturbation has compact support if the corresponding section has compact support on $\hat{\ex C}\times \hat{\ex B}$.
 
\end{defn}

Let $\mathfrak P$ be a simple perturbation parametrized by  a family of curves $\hat f$ in $\hat{\ex B}$ with a trivialization $(F,\Phi)$. So $\mathfrak P$ is a section of $\lrb{{T_{vert}^{*}{\ex C}(\hat f) }\otimes {T_{vert}\hat{\ex B}}}^{(0,1)}$ over ${\ex C}(\hat f)\times \hat {\ex B}$. A section $\nu$ of $\hat f^{*}T_{vert}\hat{\ex B}$ defines a map $(\id,F(\nu)):\ex C(\hat f)\longrightarrow \ex C(\hat f)\times \hat {\ex B}$. Pulling back the section $\mathfrak P$ over this map gives a section of  $Y(F(\nu))$, which we can identify as a section of $Y(\hat f)$ using the map $\Phi$ from our trivialization. Therefore, we get a modification $\dbar'$ of the usual $\dbar$ equation on sections $\hat f^{*}T_{vert}\hat{\ex B}$ given by the trivialization 

\[\dbar'\nu:=\dbar\nu-\Phi\lrb{(\id,F(\nu))^{*}\mathfrak P}\]

This modification $\dbar'$ of $\dbar$ was what was referred to as a simple perturbation of $\dbar$ in \cite{reg}.

The following theorem is the main theorem of \cite{reg}.

\begin{thm}\label{regularity theorem}
Suppose that $(\hat f,V)$ is a $\C\infty1$ pre obstruction model containing the curve $f$, so that $\dbar f\in V(f)$, and 
\[D\dbar(f):X^{\infty,\underline 1}(f)\longrightarrow Y^{\infty,\underline 1}(f)\]
is injective and has image complementary to $V(f)$.

 Then the restriction $(\hat f',V)$ of $(\hat f,V)$ to some open neighborhood of $f$ satisfies the following:

	There exists a neighborhood $U$ of $\dbar$ in the space of $\C\infty1$ perturbations of $\dbar$  and a neighborhood $O$ of $0$ in  $X^{\infty,\underline 1}(\hat f')$ so that  
	\begin{enumerate}
	\item \label{rt1}Given any curve $f'$ in $\hat f'$, section $\nu\in O$, and simple perturbation $\dbar'$ of $\dbar$ in $U$, 
	\[D\dbar'(\nu(f')):X^{\infty,\underline 1}(\nu(f'))
	\longrightarrow Y^{\infty,\underline 1}(\nu(f'))\]
	is injective and has image complementary to $V(f')$.
		
\item	\label{rt2}For  any $\dbar'\in U$, there exists some $\nu\in O$ and a section $v$ of $V$ so that 
\[\dbar'\nu=v\]
The sections $\nu$ and $v$ are unique in the following sense:  Given any curve $g$ in $\hat f'$, let $\nu(g)$ and $ O(g)$ be the relevant restrictions of $\nu$ and $O$ to $g$. Then $\nu(g)$ is the unique element of $O(g)$ so that $\dbar'\nu(g)\in V(g)$. 

The map $U\longrightarrow O$ which sends $\delta'$ to the corresponding solution $\nu$ is continuous in the $\C\infty1$ topologies on $U$ and $O$.
\end{enumerate}

\end{thm}

In particular, Theorem \ref{regularity theorem} above tells us how the solutions of $\dbar$ equation behave in an open neighborhood of  an obstruction model when the $\dbar$ equation is perturbed by a simple perturbation. In light of this, we make the following definition:

\begin{defn}\label{obstruction substack} An obstruction model $(\hat f/G,V,F,\Phi,\{s_{i}\})$ for an open substack $\mathcal O\subset\Msw$ is an obstruction model  so that $(\hat f/G,F,\{s_{i}\})$ is a core family for $\mathcal O$, and 
there exists a neighborhood $U$ of $0$ in the space of $\C\infty1$ simple perturbations $\dbar'$ of $\dbar$ so that items \ref{rt1} and \ref{rt2} of Theorem \ref{regularity theorem} hold for $U$ and the open neighborhood $O$ of the zero section in $X^{\infty,\underline 1}(\hat f)$ defined so that $\nu\in O$ if and only if $F(\nu)$ is in $\mathcal O$.

\end{defn}

Theorem \ref{regularity theorem} implies that every obstruction model $(\hat f/G,V)$ is an obstruction model for some open  neighborhood  $\mathcal O$ of $\hat f$ in $\Msw$. The following theorem is proved on page \pageref{construct obstruction model}. 

\begin{thm}Given any stable holomorphic curve $f$ with at least one smooth component in a basic family of targets $\hat {\ex B}\longrightarrow \ex G$, there exists an obstruction model $(\hat f/G,V)$ for an  open substack $\mathcal O\subset \Msw$ so that $f$ is isomorphic to a member of the family $\hat f$.
\end{thm}

 Theorem \ref{regularity theorem} is also useful for proving Theorem \ref{Multi solution} on page \pageref{Multi solution} which describes the solutions of the $\dbar$ equation when perturbed  using multiple simple perturbations parametrized by different obstruction models.

\

The following theorem, proved in  \cite{reg} implies that for any simple perturbation $\dbar'$ of $\dbar$, we may treat $D\dbar'(f)$ like it is a Fredholm operator, and that we may define orientations on the kernel and cokernel of $D\dbar'(f)$ by choosing a homotopy of $D\dbar'(f)$ to a complex map.  

\begin{thm}\label{f replacement}Given any $\C\infty1$ family of curves $\hat f$  with a trivialization, a set $\{s_{i}\}$ of extra marked points and a simple perturbation $\dbar'$ of $\dbar$, the following is true:
\begin{enumerate}\item\label{fr1} for every curve $f$ in $\hat f$, 
\[D\dbar'(f):X^{\infty,\underline 1}(f)\longrightarrow Y^{\infty,\underline 1}(f)\]
is a linear map which has a closed image and finite dimensional kernel and cokernel.
\item \label{fr2}
The dimension of the kernel minus the dimension of the cokernel of $D\dbar'(f)$ is a topological invariant \[2c_{1}-2n(g+s-1)\] where $c_{1}$ is the integral of the first Chern class of $J$ over the curve $f$, $2n$ is the relative dimension of $\hat{\ex B}\longrightarrow \ex G$, $g$ is the genus of the domain of $f$, and $s$ is the number of extra marked points on which sections in $X^{\infty,\underline 1}$ must vanish.

\item \label{fr3}If $D\dbar'(f)$ is injective, then for all $f'$ in an open neighborhood of $f$ in $\hat f$, $D\dbar'(f')$ is injective.

\item \label{fr4}If $D\dbar'(f)$ is injective for all $f$ in $\hat f$, then there is a $\C\infty1$ vector bundle $K$ over $\ex F(\hat f)$ with an identification of  the fiber over $f$ with the dual of the cokernel of $D\dbar(f)$, \[K(f)=Y^{\infty,\underline 1}(f)/D\dbar'(f)(X^{\infty,\underline 1}(f))\] so that given any section $\theta$  in $Y^{\infty,\underline 1}(\hat f)$, the corresponding section of $K$ is $\C\infty1$.
\end{enumerate}

The set of maps 
\[X^{\infty,\underline 1}(\hat f)\longrightarrow Y^{\infty,\underline 1}(\hat f)\]
equal to some $D\dbar'$ for some simple perturbation $\dbar'$ of $\dbar$ is convex and contains the complex map \[\frac12(\nabla\cdot+J\circ\nabla \cdot \circ j):X^{\infty,\underline 1}(\hat f)\longrightarrow Y^{\infty,\underline 1}(\hat f)\]
for any $\C\infty1$ connection $\nabla$ on $T_{vert}\hat{\ex B  }$ which preserves $J$.

The set of all such $D\dbar':X^{\infty,\underline 1}(\hat f)\longrightarrow Y^{\infty,\underline 1}(\hat f)$ is independent of choice of trivialization for $\hat f$. 

\end{thm}

\begin{remark}\label{orientation} The kernel of $D\dbar'(f)$ is oriented relative to the cokernel by choosing a homotopy of $D\dbar'(f)$ to a complex operator as follows:
\end{remark}
 Let $\hat f$ be the trivial family consisting of $\mathbb R$ times $f$, and let $\dbar'$ be chosen so that $D\dbar'(f,t)$ is equal to a complex map when $t=1$ and $D\dbar'(f)$ when $t=0$. The properties of the kernel of $D\dbar'(f,t)$ given in Theorem \ref{f replacement} imply that we may choose a finite set of marked points so that the kernel of $D\dbar'(f,t)$ restricted to sections which vanish at those marked points is zero for $t$ in some neighborhood of $[0,1]$. Let $X'\subset X^{\infty,\underline 1}(f)$ denote this set of sections vanishing at these marked points.
Theorem \ref{f replacement} tells us that cokernel of $D\dbar'(f,t)$ restricted to $X'$ gives a finite dimensional smooth vector bundle over $[0,1]$. Give this vector bundle the relative orientation defined by its complex structure when $t=1$. This defines an orientation for the cokernel of $D\dbar'(f)$ restricted $X'$. The quotient $X^{\infty,\underline 1}(f)/X'$ is finite dimensional and complex, and comes with a canonical map to $Y^{\infty,\underline 1}(f)/D\dbar(f)(X')$, which is the dual to the cokernel of $D\dbar'(f)$ restricted to $X'$, and hence oriented by the homotopy above. We may identify the kernel and cokernel of $D\dbar'(f)$ with the kernel and  cokernel of this map \[X^{\infty,\underline 1}(f)/X'\longrightarrow Y^{\infty,\underline 1}(f)/D\dbar(f)(X')\]
As both the domain and target of this map are finite dimensional and oriented, this defines an orientation of the kernel relative to the cokernel, which we may take as the orientation of the kernel of $D\dbar'(f)$ relative to its cokernel.

To see that this construction does not depend on the choice of extra marked points, let $X''$ denote set of sections in $X'$ which vanish on some more chosen marked points. Then the following is a short exact sequence for all $t\in[0,1]$:
\[ X'/X''\longrightarrow Y^{\infty,\underline 1}(f)/D\dbar'(f,t)(X'')\longrightarrow Y^{\infty,\underline 1}(f)/D\dbar'(f,t)(X')\]

 At $t=1$, the inclusion of $X'/X''$ into $Y/D\dbar'(f,1)(X'')$ is complex, so the orientation on $Y^{\infty,\underline 1}(f)/D\dbar'(f,t)(X'')$
given by our construction using this larger set of marked points is the same as the orientation determined by the above short exact sequence
giving $X'/X''$ its complex orientation and $Y^{\infty,\underline 1}(f)/D\dbar'(f,t)(X')$ the orientation defined above. The orientation of $X/X''$ given by its complex structure is also compatible with the orientation on $X'/X''$ and $X/X'$ given by their complex structures in the sense that the following short exact sequence is oriented 
\[X'/X''\longrightarrow X/X''\longrightarrow X/X'\]
Therefore the relative orientation of the kernel and cokernel of $D\dbar'(f)$ using this larger set of marked points is the same as the relative orientation using the original set of marked points. The orientation determined using any other set of marked points is equivalent, because it is equivalent to the orientation obtained using the union of our two sets of marked points. 

Note that the orientation of the kernel of $D\dbar' (f)$  relative to its cokernel thus determined does not depend on the trivialization used to define $D\dbar'(f)$, as with any choice of trivialization, we may still use the linear homotopy between $D\dbar'(f)$ and a complex operator. Note also that this orientation does not depend on the choice of complex operator, because the set of complex operators we may use is convex.

  For any obstruction model  $(\hat f/G,V)$, $D\dbar(f)$ is injective and $V(f)$ is equal to the dual of the cokernel of $D\dbar(f)$, and is therefore oriented as above. In fact, the properties of the kernel and cokernel of $D\dbar'$ in families given by items \ref{fr3} and \ref{fr4} of Theorem \ref{f replacement}  imply that this gives an orientation of the vector bundle $V$.

Therefore, any obstruction model $(\hat f/G,V)$ has a canonical orientation which is an orientation of $\hat f$ relative to $\ex G$ and an orientation of $V$ relative to $\hat f$. This gives an orientation relative to $\ex G$ for the transverse intersection of any two sections of $V$.

\subsection{weighted branched sections of sheaves}

\

\

Obstruction models give a local model for the behavior of $\dbar$ on the moduli stack of curves. For the construction of the virtual moduli space of holomorphic curves, we need some way of dealing with the usual orbifold issues that arise when dealing with moduli spaces of holomorphic curves. I think that the best way of defining the virtual moduli space probably involves the use of Kuranishi structures, first defined by Fukaya and Ono in \cite{FO}. A generalization of the  Kuranishi homology developed by Joyce in \cite{kuranishihomology} should extend to the exploded setting, however this is not done in this paper. Instead, we shall work with weighted branched objects.  There are a few approaches to weighed branched manifolds - our definition below only allows the definition of weighted branched sub objects, and is subtly different from the definition given by Cieliebak, Mundet i Rivera and Salamon in \cite{wbsubmanifolds} or the  intrinsic definition given by Mcduff in \cite{orbifolds}, because our definition has the notion of a `total weight' and  allows for the possibility of an empty submanifold being given a positive weight. I do not know which approach is better.

\begin{defn}\label{multi} The following is a construction of a `weighted branched' version of any sheaf of sets or vector spaces.

Given a vector space $V$, consider the group ring of $V$ over $\mathbb R$. This is  the free commutative ring generated as a $\mathbb R$ module by elements of the form $t^{v}$ where $v\in V$ and $t$ is a dummy variable used to write addition on $V$ multiplicatively. Multiplication on this group ring is given by
\[w_{1}t^{v_{1}}\times w_{2}t^{v_{2}}=(w_{1}w_{2})t^{v_{1}+v_{2}}\]
where $w_{i}\in\mathbb R$, $t$ is a dummy variable, and $v_{i}\in V$. Denote by $\wb V$ the sub semiring of the group ring of $V$ over $\mathbb R$ consisting of elements of the form $\sum_{i=1}^{n}w_{i}t^{v_{i}}$ where $w_{i}\geq 0$. 

There is a homomorphism 
\[\text{Weight}:\wb V\longrightarrow \mathbb R\]
\[\text{Weight}\lrb{\sum_{i=1}^{n}w_{i}t^{v_{i}}}:=\sum_{i=1}^{n}w_{i}\]

Similarly, if  $X$ is a set, consider the free $\mathbb R$ module generated by elements of the form $t^{x}$ for $x\in X$. Define $\wb{X}$ to be the `$\mathbb R^{+}$ submodule'  consisting of elements in  the form  of finite sums $\sum_{i=1}^{n} w_{i}t^{x_{i}}$ where $w_{i}\geq0$. 
The homomorphism $\text{Weight}:\wb X\longrightarrow \mathbb R$ is defined similarly to the case of vector spaces: 
$\text{Weight}(\sum w_{i}t^{x_{i}})=\sum w_{i}$.

Given a sheaf $S$ with stalks $S_{x}$, define the corresponding weighted branched sheaf $\wb S$ to be the sheaf with stalks $\wb {S_{x}}$. Call a section of $\wb S$ a weighted branched section of $S$. The Weight homomorphism gives a sheaf homomorphism of $\wb S$ onto the locally constant sheaf with stalks equal to $[0,\infty)$. (The weight of a section of $\wb S$ is a locally constant, $[0,\infty)$ valued section.) We shall usually just be interested in weighted branched sections of $S$ with weight $1$.

\end{defn}
This construction allows us to talk of the following weighted branched objects:

\begin{enumerate}
\item A smooth weighted branched section of a vector bundle $X$ over  a manifold $M$ is a global section of $\wb{C^{\infty}(X)}$ where $C^{\infty}(X)$ indicates the sheaf on $M$ of smooth sections of $X$. In particular, such a weighted branched section is locally of the form 
\[\sum_{i=1}^{n}w_{i}t^{\nu_{i}}\]
where $\nu_{i}$ is  a smooth section. This section has weight $1$ if $\sum w_{i}=1$.

\item Given a vector bundle $X$ over the total space of a family of curves $\hat {\ex C}\longrightarrow \ex F$, a $C^{k,\delta}$ weighted branched section of $X$  is defined as follows: consider the sheaf $C^{k,\delta}(X)$ over the topological space $\ex F$ which assigns to each open set $U\subset \ex F$, the vector space of $C^{k,\delta}$ sections of the vector bundle $X$ restricted to the inverse image of $U$ in $\hat {\ex C}$. Then a $C^{k,\delta}$ weighted branched section of $X$ is a global section of $\wb{C^{k,\delta}(X)}$. Such a weighted branched  section is equal to $\sum w_{i}t^{\nu_{i}}$ restricted to sufficiently small open subsets of $\ex F$, where $\nu_{i}$ indicates a $C^{k,\delta}$ section of $X$. Note that considering $C^{k,\delta}$ sections of $X$ as a sheaf over different topological spaces allows different branching behavior.

\item \label{wbv} Suppose that we have an obstruction model $(\hat f/G,V)$. Define a weighted branched section of $\hat f^{*}T\ex B\oplus V$ as follows: Consider the sheaf of $\C\infty1$  sections of $\hat f^{*}T\ex B$ and the sheaf of $\C\infty1$ sections of $V\longrightarrow \ex F(\hat f)$ as sheaves of vector spaces over $\ex F(\hat f)$. Let $X$ be the product of these two sheaves. A weighted branched section of $\hat f^{*}T\ex B\oplus V$ is a global section of $\wb{X}$.

\item Given exploded manifolds $\ex A$ and $\ex B$, consider the sheaf of sets  $S(\ex A)$  so that $S(U)$ is the set of maps of $U$ to $\ex B$. A weighted branched map of $\ex A$ into $\ex B$ is a global section of $\wb{S(\ex A)}$.

\item For any topological space $X$, consider the sheaf of sets  $S(X)$, where $S(U)$ is the set of subsets of $U$. A weighted branched subset of $X$ is a global section of $\wb{S(X)}$. (As an example of what is meant by `branching' in this context if $X=\mathbb R$, the global section $t^{(-1,2)}+t^{(0,1)}= t^{(-1,1)}+t^{(0,2)}$, $t^{(-2,-1)}+t^{(0,1)}=t^{(-2,-1)\cup(0,1)}+t^{\emptyset}$, but $t^{(-1,0)}+t^{(0,1)}\neq t^{(-1,0)\cup(0,1)}+t^{\emptyset}$.)
\item For any smooth manifold $M$, consider the sheaf of sets $S(M)$ where $S(U)$ is the set of smooth submanifolds of $U$. Then a smooth weighted branched submanifold of $M$ is a global section of $\wb{S(M)}$. 
 Locally such a weighted branched submanifold is equal to 
\[\sum_{i=1}^{n}w_{i}t^{N_{i}}\]
where each $N_{i}$ is a submanifold. Note that $N_{i}$ might intersect $N_{j}$, and $N_{i}$ might be equal to an empty submanifold. 

Similarly, one can talk of weighted branched submanifolds which are proper, oriented, or have a particular dimension.

\item An $n$-dimensional substack of the moduli stack of $\omega$-positive $\C\infty1$ curves  $\Msw$  is a substack  equal to $\hat f$ where $\hat f$ is an $n$-dimensional $\C\infty1$  family of curves. (In other words, a $\C\infty 1$ family $\hat f'$ in  this substack  is equivalent to a $\C\infty1$ map $\hat f'\longrightarrow\hat f$.) Define a sheaf of sets $\mathcal X$ over $\Msw$ by setting $\mathcal X(U)$ to be the set of complete oriented $n$-dimensional  $\C\infty1$ substacks of $U\subset \Msw$. A complete oriented weighted branched  $n$-dimensional $\C\infty1$ substack of $\Msw$ is a global section of $\wb{\mathcal X}$ with weight 1.
The support of such a weighted branched substack is the set of curves $f$ so that there is no neighborhood of $f$ in which our weighted branched substack is equal to the empty substack with weight 1.

\end{enumerate}

\subsection{multiperturbations}

\

\

\begin{defn}
A multi-perturbation on a substack   $\mathcal O\subset \Msw$ is an assignment to each $\C\infty 1$ family $\hat f$ of curves  in $\mathcal O$,  a choice of $\C\infty1$ weighted branched section of $\Y{\hat f}$ with weight $1$ so that given any map of families of curves in $\mathcal O$, $\hat f\longrightarrow \hat g$,  the weighted branched section of $\Y{\hat f}$ is the pull back of the weighted branched section of $\Y{\hat g}$.
\end{defn}

%

\begin{example}[The multiperturbation defined by a simple perturbation]\label{pullback perturbation}
\end{example}
A simple perturbation $\mathfrak P$ parametrized by $\hat f$ where $\hat f/G$ is a core family  for $\mathcal O$ defines a multi-perturbation on $\mathcal O$ which is a weighted branched section $(\hat f')^{*}\mathfrak P$ of $\Y{\hat f'}$ for all families of curves $\hat f'$ in $\mathcal O$ as follows:

Recall  from the definition of the core family $\hat f/G$ that given any $\C\infty 1$ family $\hat f'$ in $\mathcal O$, there exists a unique map satisfying certain conditions 
\[\begin{array}{ccc}\ex  C(\hat f')&\xrightarrow{\Phi_{\hat f'}}&\ex C(\hat f)/G
\\ \downarrow & & \downarrow
\\ \ex F(\hat f')&\longrightarrow &\ex F(\hat f)/G\end{array}\] 
Given such a map, around any point $p\in \ex F(\hat f')$, there exists a  neighborhood $U$ of $p$ so that $\Phi_{\hat f'}$ restricted to the lift $\tilde U$ of $U$ to $\ex C(\hat f')$ lifts to exactly $\abs{G}$ maps $\Phi_{ U,g}:\tilde U\longrightarrow \ex C(\hat f)$ 

\[\begin{array}{ccc}
 &&\ex C(\hat f)
\\ &\nearrow \Phi_{U,g}&\downarrow
\\\ex  C(\hat f')\supset\tilde U&\xrightarrow{\Phi_{\hat f'}}&\ex C(\hat f)/G
\end{array}
\]

Recall that the simple perturbation $\mathfrak P$ is some $\C\infty 1$ section of the bundle $\lrb{{T_{vert}^{*}\ex C(\hat f)  }\otimes {T_{vert}\hat{\ex B}}}^{(0,1)}$ over $\ex C(\hat f)\times \hat {\ex B}$. We can pull this section $\mathfrak P$ back over the maps $(\Phi_{U,g},\hat f')$ to obtain $\abs G$  sections $(\Phi_{U,g},\hat f')^{*}\mathfrak P$ of the bundle $\Y{\hat f'}$ restricted to $U$. Define our weighted branched section $(\hat f')^{*}\mathfrak P$  of $\Y{\hat f'}$ over $U$ to be 
\[(\hat f')^{*}(\mathfrak P)\rvert_{U}:=\sum_{g\in G}\frac 1{\abs G}t^{(\Phi_{U,g},\hat f')^{*}\mathfrak P}\]

This construction is clearly compatible with any map of families $\hat f''\longrightarrow \hat f'$ in $\mathcal O$, so the above local construction glues together to a global weighted branched section $(\hat f')^{*}(\mathfrak P)$ of $\Y{\hat f'}$, and we get a multi-perturbation defined over all families  in $\mathcal O$. Note that $\hat f^{*}(\mathfrak P)$ is not just $(id,\hat f)^{*}\mathfrak P$, but the weighted branched section which is $\sum_{g\in G}\frac 1{\abs G}t^{(g,\hat f)^{*}\mathfrak P}$ where $g$ indicates the map  $g:\ex C(\hat f)\longrightarrow\ex C(\hat f)$ given by the group action.

\

To extend a multiperturbation on $\mathcal O$ defined as above  to some other  substack $\mathcal U$ of $\Msw$, we need  $\mathcal U$ to meet $\mathcal O$ properly as defined below.

\begin{defn}\label{proper meeting} Say that a family $\hat g$ of curves meets a substack $\mathcal O\subset\Msw$ with a core family $\hat f/G$ properly if the following holds. Let $\hat g\rvert_{\mathcal O}$ indicate the sub family of $\hat g$ consisting of all curves contained in $\mathcal O$, let $\iota:\ex C(\hat g\rvert_{\mathcal O})\longrightarrow\ex C(g) $ indicate the natural inclusion, and let $\Phi_{\hat g\rvert_{\mathcal O}}:\ex C(\hat g\rvert_{\mathcal O})\longrightarrow \ex C(\hat f)/G$ indicate the projection from the definition of a core family. Then the map 
\[(\iota, \Phi_{\hat g\rvert_{\mathcal O}}):\ex C(\hat g\rvert_{\mathcal O})\longrightarrow \ex C(\hat g)\times \ex C(\hat f)/G\]
	is  proper.

Say that a substack $\mathcal U\subset\Msw$ meets $\mathcal O$ properly if every family $\hat g$ in $\mathcal U$ meets $\mathcal O$ properly.
\end{defn}

\begin{example}\label{proper pullback}\end{example} Let $\mathfrak P$ be a compactly supported simple perturbation parametrized by a core family $\hat f/G$ for $\mathcal O$, and let $\mathcal U$ meet $\mathcal O$ properly. Define a multi-perturbation on $\mathcal U$ as follows: If $\hat g$ is a family in $\mathcal U$, and $\hat g\rvert_{\mathcal O}$ indicates the sub family of all curves contained in $\mathcal O$, let $\hat g^{*}\mathfrak P$ indicate the weighted branched section of $\Y{\hat g}$ which when restricted to $\hat g_{\rvert_{\mathcal O}}$ equals the multi perturbation $(\hat g\rvert_{\mathcal O})^{*}\mathfrak P$ 
 from example \ref{pullback perturbation}, and which is equal to the zero section with weight $1$ everywhere else. As $\hat g$ meets $\mathcal O$ properly and $\mathfrak P$ is compactly supported, $\hat g^{*}\mathfrak P$ is a $\C\infty 1$ weighted branched section if $\hat g$ is $\C\infty1$.

\

More generally, given a finite collection of compactly supported simple perturbations $\mathfrak P_{i}$  parameterized by core families $\hat f_{i}/G_{i}$ for substacks  $\mathcal O_{i}$ which meet $\mathcal U$ properly, we may define a multiperturbation on $\mathcal U$  by multiplying together all the multiperturbations from example \ref{proper pullback}, so given a family $\hat g$ in $\mathcal U$, the multiperturbation on $\hat g$ is $\prod_{i}\hat g^{*}\mathfrak P_{i}$. (Note that a product of weighted branched sections as defined in definition \ref{multi} involves adding sections and multiplying weights.)
Theorem \ref{Multi solution} on page \pageref{Multi solution} describes the solution to the $\dbar$  equation  within an obstruction model when perturbed by a multiperturbation constructed in this way.

\subsection{Summary of construction of virtual moduli space}

\

\

The following theorem proved beginning on page \pageref{final theorem proof} should be thought of as outlining  the construction of a `virtual class' for a component of the moduli stack of holomorphic curves. This virtual class is a cobordism class of finite dimensional $\C\infty1$ weighted branched substacks of $\Msw$, oriented relative to $\ex G$. Other approaches such as those  used in \cite{FO}, \cite{kuranishihomology}, \cite{Polyfoldint}, \cite{Tian-Li}, \cite{Liu-Tian}, \cite{Ruanvirtual}, \cite{McDuff} and \cite{Siebert} should also generalize to the exploded setting. The existence of obstruction models imply that given any compact basic family of targets $\hat {\ex B}\longrightarrow \ex G$ on which Gromov compactness holds, the stack of holomorphic curves in  $\Msw_{g,[\gamma],\beta}(\hat {\ex B})$   may be covered by a finite number of extendible obstruction models. The following theorem constructs a virtual moduli space $\Mod_{g,[\gamma],\beta}(\hat{\ex B})$ using these obstruction models.

\begin{thm}\label{final theorem}
Given 
\begin{itemize}
\item a compact, basic family of targets $\hat{\ex  B}\longrightarrow \ex G$ in which Gromov compactness holds in the sense of definition \ref{gcompactness}. 
\item a genus $g$, a tropical curve $\gamma$ in $\totb{\hat{\ex B}}$ and a linear map $\beta :H^{2}(\hat{\ex B})\longrightarrow \ex B$ ,
\item and any finite collection of extendible obstruction models covering  the moduli stack of holomorphic curves  in  $\Msw_{g,[\gamma],\beta}(\hat{\ex B})$, 
\end{itemize}
each obstruction model may be modified by restricting to an open subset covering the same set of holomorphic curves, and satisfying the following: There exists an open $\C\infty1$ neighborhood $U$ of $0$ in the space of collections of compactly supported simple perturbations parametrized by these obstruction models so that for any collection $\{\mathfrak P_{i}\}\in U$ of such perturbations,  the following is true
\begin{enumerate}
\item \label{ft1}There is some open neighborhood $\mathcal O\subset \Msw$ of the set of holomorphic curves in $\Msw_{g,[\gamma],\beta}(\hat{\ex B})$ which meets each of the our obstruction models properly. On $\mathcal O$ there is a $\C\infty1$ multi-perturbation $\theta$ defined by \[\theta(\hat f):=\prod_{i} \hat f^{*}\mathfrak P_{i}\] where each multi-perturbation $\hat f^{*}\mathfrak P_{i}$ is as defined in example \ref{proper pullback}. (Note that the notation of a product of the weighted branched sections involves adding sections and multiplying weights. See definition \ref{multi} on page \pageref{multi}.)

\item \label{ft2}For each of our obstruction models $(\hat f/G,\{s_{i}\},F,V)$ there exists a unique $G$-invariant $\C\infty1$ weighted branched section $(\nu,\dbar'\nu)$ of $\hat f^{*}T_{vert}\hat{\ex B}\oplus V$  with weight $1$  (see example \ref{wbv} below definition \ref{multi}) so that

\begin{enumerate}
\item locally on $\ex F(\hat f)$, 
\[(\nu,\dbar'\nu)=\sum_{k=1}^{n}\frac 1n t^{(\nu_{k},\dbar'\nu_{k})}\]
where  the sections $\nu_{k}$ vanishes on the image of the marked point sections $s_{i}:\ex F(\hat f)\longrightarrow\ex C(\hat f)$ and 
correspond to families $F(\nu_{k})$ in $\mathcal O$

\item
\[\theta(F(\nu_{k}))=\sum_{j=1}^{n}\frac 1nt^{\mathfrak P_{k,j}}\]
 and \[\dbar F(\nu_{k})-\mathfrak P_{k,k}=\dbar'\nu_{k}\]
\item Given any curve $f$ in the above subset of $\hat f$ where $\nu_{k}$ is defined, if $\nu'$ is a section of $f^{*}{T_{vert}\hat{\ex B}}$ vanishing on $\{s_{i}\}$ so that $F(\nu')$ is in $\mathcal O$ and  the multi-perturbation $\theta(F(\nu'))=\sum_{j}w_{j}t^{\mathfrak Q_{j}}$, then $\frac 1n$ times the number of the above $\nu_{k}$ so that $\nu'$ is the restriction of $\nu_{k}$ is equal to the sum of $w_{j}$ so that $\dbar F(\nu')- \mathfrak Q_{j}$ is in V.
\end{enumerate}
\item\label{ft3}
Say that  the multi perturbation $\theta$ is transverse to $\dbar$ on a sub family $C\subset \hat f$ if the sections $\dbar'\nu_{k}$ of $V$ are transverse to the zero section on $C$.

Given any compact subfamily of one of our obstruction models $C\subset\hat f$, the subset of the space $U$ of collections of simple perturbations $\{\mathfrak P_{i}\}$ discussed above  so that  the corresponding multi-perturbation $\theta$ is transverse to $\dbar$ is open and dense in the $\C\infty1$ topology.

\item\label{ft4}
Say that the multi-perturbation $\theta$ is fixed point free on a sub family $C\subset \hat f$ if none of the curves in $F(\nu_{k})$ restricted to $C$ have smooth part with nontrivial automorphism group.

If the relative dimension of $\hat {\ex B}\longrightarrow \ex G$ is greater than $0$, then given any compact subfamily of one of our obstruction models $C\subset\hat f$, the subset of the space $U$ of collections of simple perturbations $\{\mathfrak P_{i}\}$ discussed above  so that $\theta$ is fixed point free is open and dense in the $\C\infty1$ topology.

\item \label{ft5}There exists an open substack $\mathcal O^{\circ}\subset \mathcal O$ which contains  the  holomorphic curves in $\Msw_{g,[\gamma],\beta}(\hat{\ex B})$ and a collection of  compact sub families of our obstruction models $C_{i}\subset \hat f_{i}$ so that
\begin{enumerate}
\item if $f$ is a stable curve in $\mathcal O^\circ$, then there exists some curve $f'$ in one of these sub families $C_{i}$ and section $\nu$ of $(f')^{*}T_{vert}\hat{\ex B}$ vanishing on marked points so that $f=F(\nu)$,
\item if
 $f$ is any curve in $\mathcal O$ so that for any collection of perturbations in $U$, $\theta(f)=wt^{\dbar f}+\dotsc$ where $w>0$, then  $f$ is in $\mathcal O^\circ$.
\end{enumerate}

\item \label{ft6}Say that  $\theta$ is transverse to $\dbar$ and fixed point free if it is transverse to $\dbar$ and fixed point free on each $C_{i}$ from item \ref{ft5} above. Given any such $\theta$, there is a unique complete weighted branched finite dimensional $\C\infty1$ substack $\mathcal M_{g,[\gamma],\beta}\subset\mathcal O^{\circ}$, oriented relative to $\ex G$ which is the solution to $\dbar=\theta$ in the following sense:

\begin{enumerate}
\item
Given a curve $f$ in $\hat f$ in the region where the equations from item \ref{ft2} hold, and given $\psi\in X^{\infty,\underline 1}$ so that $F(\psi)$ is in $\mathcal O^{\circ}$,  restricted to some neighborhood of $F(\psi)$ in $\mathcal O^{\circ}$, $\Mod_{g,[\gamma],\beta}$ is locally equal to
 \[\sum _{k=1}^{n}\frac 1n t^{\hat g_{k}}\]
 where $\hat g_{k}$ is the empty substack if $\nu_{k}(f)\neq \psi$ or $\dbar'\nu_{k}\neq 0$, and  otherwise, $\hat g_{k}$ is given by the restriction of $F(\nu_{k})$ to some neighborhood of $F(\psi)$ in the intersection of $\dbar'\nu_{k}$ with the zero section, oriented using the orientation of $\ex F(\hat f)$ relative to $\ex G$ and the orientation of $V$ relative to $\ex F(\hat f)$.
\item If $\mathcal M_{g,[\gamma],\beta}$ is locally equal to $\sum_{k}w_{k}t^{\hat g_{k}}$ then $\theta(\hat g_{k})=w_{k}t^{\dbar\hat g_{k}}+\dotsc$. 
 \item If $\theta (f)=\sum_{k}w'_{k}t^{\mathfrak Q_{k}}$, then the sum of weights $w_{k}$ so that $f$ is in $\hat g_{k}$
is equal to the sum of weights $w'_{k}$ so that $\dbar f=\mathfrak Q_{k}$.
\item If a family $\hat g$ in $\mathcal O^{\circ}$ containing a curve $f$ satisfies 
\[\theta(\hat g)=w t^{\dbar \hat g}\] 
then if $\mathcal M_{g,[\gamma],\beta}$ is equal to $\sum_{k}w_{k}t^{\hat g_{k}}$ on a neighborhood of $f$, then the sum of the $w_{k}$ so that there is a map from some neighborhood of $f$ in $\hat g$ to $\hat g_{k}$ is at least $w$.  
\end{enumerate}
\item\label{ft7} The support of $\mathcal M_{g,[\gamma],\beta}$ is  a compact subset of $\Msw_{g,[\gamma],\beta}$. 
\item \label{ft8}
Given any construction of virtual moduli space $\mathcal M'_{g,[\gamma],\beta}$ defined using another small enough multi-perturbation $\theta'$ which is fixed point free and  transverse to $\dbar$, defined on some other open neighborhood of the holomorphic curves in $\Msw_{g,[\gamma],\beta}(\hat{\ex B})$ using different choices of obstruction models, $\mathcal M_{g,[\gamma],\beta}$ is cobordant to $\mathcal M'_{g,[\gamma],\beta}$ in the following sense: 
 
 Let $\hat {\ex B}\times S^{1}\longrightarrow \ex G\times S^{1}$ be the product of our original family $\hat {\ex B}\longrightarrow \ex G$ with a circle. We can regard $\gamma$ as a tropical curve in $\totb{\hat{\ex B}\times S^{1}}$ and $\beta$ as a map $H^{2}(\hat{\ex B}\times S^{1})\longrightarrow \mathbb R$. Then we may construct $\Mod_{g,[\gamma],\beta}(\hat{\ex B}\times S^{1})$ so that the map 
 \[\Mod_{g,[\gamma],\beta}(\hat{\ex B}\times S^{1})\longrightarrow S^{1}\]
 is transverse to two points $p_{1}$ and $p_{2}$ in $S^{1}$, and the restriction of $\Mod_{g,[\gamma],\beta}(\hat{\ex B}\times S^{1})$ to $\Msw_{g,[\gamma],\beta}(\hat{\ex B}\times \{p_{1}\})$ is equal to $\Mod_{g,[\gamma],\beta}(\hat{\ex B})$ and the restriction of $\Mod_{g,[\gamma],\beta}(\hat{\ex B}\times S^{1})$ to $\Msw_{g,[\gamma],\beta}(\hat{\ex B}\times \{p_{2}\})$ is equal to $\Mod'_{g,[\gamma],\beta}(\hat{\ex B})$.

\item \label{ft 9}Given any compact exploded manifold $\ex G'$ and $\C\infty1$ map $\ex G'\longrightarrow \ex G$, the virtual moduli space
$\Mod_{g,[\gamma],\beta}(\hat {\ex B})$ can be constructed so that the map \[\Mod_{g,[\gamma],\beta}(\hat{\ex B})\longrightarrow \ex G\] is transverse to the map $\ex G'\longrightarrow \ex G$.

Suppose that this is the case and let $\hat{\ex B}'\longrightarrow \ex G'$ be the pullback of $\hat{\ex B}\longrightarrow \ex G$. Suppose also that there are tropical curves  $\gamma_{i}$ in $\totb{\hat{\ex B}'}$ and  maps  $\beta_{i}:\mathbb H^{2}(\hat{\ex B}')\longrightarrow \mathbb R$ so that a curve is in  $\coprod_{i}\Msw_{g,[\gamma_{i}],\beta_{i}}(\hat{\ex B'})$ if and only if its composition with $\hat {\ex B}'\longrightarrow \hat{\ex B}$ is a curve in  $\Msw_{g,[\gamma],\beta}(\hat {\ex B})$. 

Then  so long as the perturbations used to define $\Mod_{g,[\gamma],\beta}$ are small enough, $\coprod_{i}\Mod_{g,[\gamma_{i}],\beta_{i}}(\hat{\ex B}')$ may be constructed as the inverse image of $\Mod_{g,[\gamma],\beta}(\hat{\ex B})$, and may be considered as the fiber product of  $\ex G'\longrightarrow \ex G$  with $\Msw_{g,[\gamma],\beta}(\hat{\ex B})\longrightarrow \ex G$.


\end{enumerate}
\end{thm}

Item \ref{ft8} above gives that the cobordism class of $\Mod_{g,[\beta],\gamma}$ is independent of the choices made in its construction.  Item \ref{ft 9} can be used to show that Gromov Witten invariants do not change in families.

\subsection{proofs}

\

\

The following theorem states roughly that the explosion of Deligne Mumford space, $\M$ (discussed in \cite{iec}) represents the moduli stack of $\C\infty 1$ families of stable exploded curves. A similar theorem probably holds over the complex version of the exploded category with `smooth and holomorphic' replacing `$\C\infty 1$'.

\begin{thm}\label{moduli} Consider  any $\C\infty 1$ family of exploded curves $(\hat{\ex C},j)\longrightarrow \ex F$ so that each exploded curve is connected and has $2g+n\geq 3$ where $g$ is the genus and $n$ is the number of punctures. Then  there exists a unique fiber wise holomorphic map 
\[\begin{array}{ccc}
(\hat {\ex C},j)&\longrightarrow & (\M^{+1},j)
\\ \downarrow & &\downarrow
\\ \ex F&\longrightarrow & \M
\end{array}\]
so that the map on each fiber $\ex C$ factors into a degree one holomorphic map to a stable exploded curve $\ex C'$ and a map from $\ex C'$ to a fiber of $\M^{+1}$ given by quotienting $\ex C'$ by its automorphism group. 

The above maps all have regularity $\C\infty 1$.
\end{thm}

\pf

We first construct this map for the fiber $\ex C$ over a single point of $\ex F$. The first stage of this is to construct a stable curve $\ex C'$  with a holomorphic degree one map $\ex C\longrightarrow \ex C'$. The idea is to `remove' all unstable components using a series of maps of the following two types: 
\begin{enumerate}
\item
If a smooth component of $\ex C$ is a sphere attached to only one edge, put holomorphic coordinates on a neighborhood of the edge modeled on  an open subset of $\et 1{[0,l]}$ with coordinate $\tilde z$ so that $\totl{\tilde z}$ gives coordinates on the smooth component of $\ex C$ attached to the other end of the edge. Replace this coordinate chart with the corresponding open subset  of $\mathbb C$ with coordinate $z=\totl{\tilde z}$. There is an obvious degree one holomorphic map from our old curve to this new one that is given in this coordinate chart by $\tilde z\mapsto\totl{\tilde z}$, and sends our unstable sphere and the edge attached to it to the point $p$ where $z(p)=0$.  (This map  is the identity everywhere else.)

\item
If a smooth component of $\ex C$ is a sphere attached to exactly two edges, there exists a holomorphic identification of a neighborhood of this smooth component with a refinement of an open subset of $\et 1{[0,l]}$. Replace this open set with the corresponding open subset of $\et 1{[0,l]}$. The degree one holomorphic map from the old exploded curve to the new one is this refinement map. (Refer to \cite{iec} for the definition of refinements.)
\end{enumerate}

Each of the above types of maps removes one smooth component, so after applying maps of the above type a finite number of times, we obtain a connected exploded curve with no smooth components which are spheres with one or two punctures. Our theorem's hypotheses then imply that the resulting exploded curve $\ex C'$ must be stable. It is not difficult to see that the stable curve obtained this way is unique.

  The smooth part of this stable exploded curve $\ex C'$ is a stable nodal Riemann surface with punctures, $\totl{\ex C'}$. This nodal Riemann surface determines a point in Deligne Mumford space $p_{\totl{\ex C'} }\in \bar{\mathcal M}$, and a corresponding map of $\totl{\ex C'}$ to the fiber of $\bar{\mathcal M}^{+1}$ over this point $p_{\totl{\ex C'}}$. Note that $\bar{\mathcal M}$ is the smooth part of $\M$ and $\bar{\mathcal M}^{+1}$ is the smooth part of $\M^{+1}$.  The smooth part of our map $\ex C'\longrightarrow \M^{+1}$ must correspond with this map $\totl{\ex C'}\longrightarrow \bar{\mathcal M}^{+1}$. We must now determine the remaining information.
  
  If we give Deligne Mumford space its usual holomorphic structure, there is a holomorphic uniformising chart $(U,G)$ containing this point $p_{\totl{\ex C'}}$, where $U$ is some subset of $\mathbb C^{n}$ so that the boundary strata of $\bar{\mathcal M}$ correspond to where coordinates $z_{i}=0$, and $G$ is a finite group with a holomorphic action on $U$ which 
preserves the boundary strata. $\M$ is constructed so that it has a corresponding uniformising coordinate chart $(\tilde U,G)$ where $\tilde U$ is an open subset of $\et n{[0,\infty)^{n}}$ which corresponds to the set where $\totl{\tilde z}\in U$, and the action of $G$ on $\tilde U$ induces the action of $G$ on the smooth part $\totl{\tilde U}=U$. The inverse image of $\tilde U$ in $\M^{+1}$ is some exploded manifold $\tilde U^{+1}$ quotiented by $G$, and the smooth part of this is the inverse image of $ U$ in $\bar{\mathcal M}^{+1}$, which is equal to some smooth complex manifold $U^{+1}$ quotiented by $G$. There are $\abs G$ identifications  of $\totl{\ex C'}$ with a fiber of  $U^{+1}\longrightarrow U$, which are permuted by the action of $G$, (so together they correspond to a unique map to $\bar{\mathcal M}^{+1}$). Choose one of these maps. 

Each of the nodes of $\totl{\ex C}$ now correspond to some coordinate $z_{i}$ on $U$ which is equal to $0$. We must determine the value of the corresponding $\tilde z_{i}$. (All other coordinates  are nonzero so $\tilde z_{k}$ is given by $\tilde z_{k}=z_{k}$.) There is a chart $U_{i}^{+1}$ on $U^{+1}$ containing this node which is equal to a convex open subset of $\mathbb C^{n+1}$ with coordinates $ z_{j}$, $j\neq i$ and $ z_{i}^{+}, z_{i}^{-}$, so that the map $U_{i}^{+1}\longrightarrow U$ is given by $z_{i}=z_{i}^{+}z_{i}^{-}$ and $z_{j}=z_{j}$.  The identification of a neighborhood of this node with a fiber of $U_{i}^{+1}$ means that we can use $ z_{i}^{+}$ and $ z_{i}^{-}$ respectively to parametrize the two discs that make up the neighborhood of the node. The open subset of $\ex C'$ with smooth part equal to this neighborhood can then be covered by an open subset of $\et 1{[0,l]}$ with coordinates $\tilde z_{i}^{+}$ and  $\tilde z_{i}^{-}$ so that $\totl{\tilde z_{i}^{\pm}}=z_{i}^{\pm}$. These coordinate are  related  by the equation 
\[\tilde z_{i}^{+}\tilde z_{i}^{-}=c\e l\] In the above, $l$ is the length of our edge, and $c$ is  canonically determined by our choice of holomorphic  coordinate chart on $\bar{\mathcal M}$. Our coordinate $\tilde z_{i}$ must be equal to $ c\e l$.   
To see this consider the corresponding coordinate chart $\tilde U_{i}^{+1}$ with coordinates $\tilde z_{j}$ and $\tilde z_{i}^{\pm}$ so that $\totl{\tilde z_{j}}=z_{j}$ and $\totl{\tilde z_{i}^{\pm}}=z_{i}^{\pm}$. The map $\tilde U_{i}^{+1}\longrightarrow \tilde U$ is given by $\tilde z_{i}=\tilde z_{i}^{+}\tilde z_{i}^{-}$ and $\tilde z_{j}= \tilde z_{j}$. The smooth part of the intersection of our curve with $\tilde U_{i}^{+}$ must be as described above, and the parametrization of the smooth part by $\totl{\tilde z_{i}^{\pm}}$ must also be as above. The fiber is over a point where $\tilde z_{i}=c\e l$ has two coordinates related by the equation $\tilde z_{i}^{+} \tilde z_{i}^{-}=c\e l$. This fiber is equal to the corresponding open subset of our curve $\ex C'$ and parametrized correctly if and only if $\tilde z_{i}=c\e l$.

We have shown that after choosing any one of the $\abs G$ holomorphic maps $\totl {\ex C'}\longrightarrow U^{+1}$ there is a unique holomorphic map $\ex C'\longrightarrow \tilde U^{+1}$ onto a fiber of $U^{+1}$ with smooth part equal to this. Therefore, there is a unique holomorphic map $\ex C'\longrightarrow \M^{+1}$ which factors as an inclusion as a fiber of $\tilde U^{+1}$ followed by  quotiententing by the action of the group $G$. In particular, there is a unique holomorphic map $\ex C'\longrightarrow \M^{+1}$ satisfying the required conditions of our theorem. This completes the construction of our map for each individual fiber. We must now verify that the resulting map on the total space has regularity $\C\infty 1$.  

To verify the regularity of the map we've constructed, we need only to work locally around a fiber. As this is local, we may assume that the base of our family $\ex F$ is covered by a single standard coordinate chart. Start with the map on a single fiber $\ex C\longrightarrow \tilde U^{+1}$ constructed above. We shall prove that this extends to a $\C\infty 1$ fiberwise holomorphic map from a neighborhood of the fiber. The uniqueness of our map on fibers shall then imply that this map must agree with the map constructed above, proving the required regularity. We shall consider $\tilde U^{+1}\longrightarrow \tilde U$ to give a family of targets,  to which we shall first construct a smooth map from a neighborhood of the fiber, and then apply Theorem \ref{regularity theorem} to correct this to a fiber wise holomorphic map. 

Construct a smooth extension of  $\ex C\longrightarrow\ex C'\longrightarrow \tilde U^{+1}$  using local coordinate charts as follows. Cover $\tilde U^{+1}$ with a finite number of charts of the following three types: the charts $\tilde U_{i}^{+1}$ mentioned earlier which cover edges of the image of $\ex C'$; charts covering punctures of the image of $\ex C'$, which are all of the form of some open subset of  $\et 11\times \tilde U$; and charts containing only smooth parts of the image of $\ex C'$, which can all be identified with some open subset of $\mathbb C$ times $\tilde U$. On $\hat{ \ex C}\longrightarrow \ex F$, consider a single coordinate chart  on $\ex F$ which we may assume (without losing generality in this part of the construction) is equal to some subset of $\et mP$ containing all strata of $P$. ( If this is not the case, and our coordinate chart on $\ex F$ is a subset of $\mathbb R^{n}\times \et mP$, we may construct our smooth map to be independent of the $\mathbb R^{n}$ coordinates.) Cover the inverse image of this coordinate chart in $\hat {\ex C}$ with a finite number of coordinate charts $\ex V$ which project to $\et mP$ in one of the standard forms for coordinate charts on families discussed in \cite{iec}. Construct these charts $\ex V$ small enough so that the portion of $\ex C$ contained in any one of these coordinate charts is contained well inside one of the coordinate charts on $\tilde U^{+1}$. In this case say that $\ex V$ is `sent to' the corresponding chart on $\tilde U^{+1}$.

Consider a chart $\ex V$ on $\hat{\ex C}$ which is sent to a chart $\tilde U_{i}^{+1}$ corresponding to the $i$th node of $\ex C'$. We may assume that if any two of these coordinate charts intersect, then the tropical part of the intersection is equal to $P$. 

 Define $\e {\mathbb R}$ valued integral affine functions $h_{i,\ex V}$ and $h_{i,\ex V}^{\pm}$ on $\totb{\ex V}$ as follows: If the tropical part of $\ex C$ intersecting this chart is sent to a single point in $\tilde U_{i}^{+1}$, then set all three of these functions equal to  $\e 0$. If not, define $h_{i,\ex V}$ at $p\in P$ to be the length of the inverse image of $p$ in $\totb{\ex V}$, and define $h_{i,\ex V}^{\pm}$ on $\totb{\ex V}$ to be the distance to either end of the fibers of $\totb{\ex V}\longrightarrow  P$, choosing the relevant `ends' so that on the intersection with $\ex C$ these $h_{i,\ex V}^{\pm}$ are equal to the pull back of $\totb{\tilde z_{i}^{\pm}}$ times some constant. Now define the function $h_{i}$ on $ P$ by multiplying together all the functions $h_{i,\ex V}$ from each of the coordinate charts above.  Note that on the intersection with $\ex C$, this integral affine function is equal to $\totb{\tilde z_{i}}$. We now define integral affine functions $h_{i}^{\pm}$ which correspond to $\totb{\tilde z_{i}^{\pm}}$. Define a partial order on these charts as follows: if $\totb{\tilde z^{+}}$ is greater on the part of $\ex C$ in chart 1 than on the part of $\ex C$ in chart 2, and at some point strictly greater, then say that chart 1 is greater than chart 2. Define $h^{+}_{i}$ on $\ex V$ to be equal to the product of  $h^{+}_{i}\ex V$ with $h_{i,\ex V'}$ for all $\ex V'$ greater that $\ex V$. Similarly define $h^{-}_{i}$. Define the function $\tilde z_{i}$ on $\et mP$ to be equal to the unique monomial so that $\totb{\tilde z_{i}}=h_{i}$ and $\tilde z_{i}$ restricted to $\ex C$ is equal to the pull back of $\tilde z_{i}$ from $\tilde U$. Doing the same for all other nodes and setting the other coordinates constant gives our smooth map from our subset of $\et mP$ to $\tilde U$. 

 Now choose  functions $\tilde z^{\pm}_{i}$ on each coordinate chart $\ex V$ so that 
\begin{enumerate}
\item \[\tilde z_{i}^{+}\tilde z_{i}^{-}=\tilde z_{i}\]
\item \[\totb{\tilde z_{i}^{\pm}}=h_{i}^{\pm}\]
\item Restricted to $\ex C$, $\tilde z_{i}^{\pm}$ is equal to the pull back of $\tilde z_{i}^{\pm}$ from $\tilde U_{i}$.
\end{enumerate} 
Because the tropical part $\totb{\tilde z_{i}^{\pm}}$ is compatible with coordinate changes, and because $\tilde z_{i}^{\pm}$ is compatible with coordinate changes, these functions are compatible with coordinate changes on any fiber with smooth part equal to the smooth part of $\ex C$, and are almost compatible with coordinate changes in a small  neighborhood of $\ex C$. We can therefore modify them to obey all the above conditions, and be compatible with coordinate changes, defining smooth exploded functions $\tilde z_{i}^{\pm}$ on the union of all coordinate charts $V$ which are sent to $\tilde  U^{+1}_{i}$.  These together with the map from our subset of $\et mP$ to $\tilde U$ defined above define smooth maps from these coordinate charts to $\tilde U_{i}^{+1}$ which are compatible with coordinate changes. 

We must also define our map on coordinate charts which are sent to coordinates charts on $\tilde U^{+1}$ which are a product  of $\tilde U$ with some  open subset of $\et 11$. As we already have our map to $\tilde U$, this amounts to constructing a map into $\et 11$, which we shall give a coordinate $\tilde w$. The construction of this map is entirely analogous to the construction of the function $\tilde z_{i}^{+}$ above. Once we have done this, we have smooth maps from each of our coordinate charts into $\tilde U^{+1}$ which are compatible with coordinate changes on any fiber with smooth part  equal to the smooth part of $\ex C$, and which agree with the map constructed earlier on $\ex C$. There is no obstruction to modifying these maps to give a smooth map which is compatible with all coordinate changes and satisfies the required conditions.

We have now shown that there exists a smooth map 
\[\begin{array}{ccc}
\hat{\ex C}&\longrightarrow & \tilde U^{+1}
\\ \downarrow& &\downarrow
\\ \ex F &\longrightarrow &\tilde U
\end{array}\]
so that the restriction to the fiber $\ex C$ is holomorphic. (We proved this under the assumption that $\ex F$ is covered by a single coordinate chart.) We now wish to show that this can be modified to a fiber wise holomorphic $\C\infty1$ map on a neighborhood of $\ex C$.  We shall show below that if the above is considered as a map into a family of targets $\tilde U^{+1}\longrightarrow \tilde U$, the cokernel of the relevant linearized $\dbar$ operator is naturally identified with the cotangent space of $\tilde U$. To deal with this cokernel within the framework of this paper,  we shall  extend our map to a smooth map 
\[\begin{array}{ccc}
\hat{\ex C}\times \mathbb R^{n} &\xrightarrow{\hat g} & \tilde U^{+1}
\\ \downarrow & &\downarrow
\\ \ex F\times \mathbb R^{n} &\longrightarrow &\tilde U
\end{array}\]
where the tangent space of $\tilde U$ is identified with $\mathbb R^{n}\times \tilde U$, and the derivative of this map on the $\mathbb R^{n}$ factor at $0$ is the identity.

Consider the corresponding linearized operator $D\dbar$ at $\hat g$ restricted to our curve $\ex C$. This is just the standard $\dbar$ operator on sections of the pullback of  $T\ex C'$ to $\ex C$. Standard complex analysis tells us that  as $\ex C'$ is stable, this operator is injective, and has a cokernel which we may identify with `quadratic differentials', which are  holomorphic sections $\theta$ of the pull back to $\ex C$ of the symmetric square of the holomorphic cotangent bundle of  $\ex C'$ which vanish at punctures. (This actually corresponds to allowing a simple pole at punctures viewed from the smooth perspective as quadratic differentials on $\et 11$ look like holomorphic functions times $(\tilde z^{-1}d\tilde z)^{2}$). This is proved by showing that the quadratic differentials are the kernel of the adjoint of $\dbar$. All we shall use is that the relationship is as follows: the wedge product of $\dbar v$ with $\theta$ gives a two form which is equal to $d(\theta(v))$. This vanishes at all edges and punctures because $\dbar v$ does, so the integral is well defined. Therefore, as it equals $d(\theta(v))$, and $\theta(v)$ is a one form which vanishes on punctures,  the integral of the wedge product of $\dbar v$ with $\theta$ over $\ex C$ must vanish. As holomorphic sections of the pullback to $\ex C$ of any bundle on $\ex C'$ can be identified with holomorphic sections of the bundle on $\ex C'$, we may identify the cokernel of our $D\dbar$ with the quadratic differentials on $\ex C'$. 

We can also identify the holomorphic cotangent space to $\tilde U$ at the image of the curve $\ex C'$ with the space of quadratic differentials as follows: Refine $\tilde U$ so that $\ex C'$ is the fiber over a smooth point, and trivialize a small neighborhood of $\ex C'$ in the corresponding refinement of $\tilde U^{+1}$.  Given a tangent vector $v$ to $\tilde U$, by differentiating the almost complex structures on fibers using our trivialization, we then obtain a tensor $v(j)$ which is a section of $T^{*}\ex C'\otimes T\ex C'$. The derivative of $\dbar\hat g$ using this trivialization  at the curve  $\ex C$ in the direction corresponding to $v$ is $\frac 12v(j)\circ j$. Then taking the wedge product of $\frac 12v(j)\circ j$ with a quadratic differential gives a two form on $\ex C'$ which we can then integrate over $\ex C'$. The result of this integral does not depend on the choice of trivialization because of the above discussion identifying the cokernel of the restriction of $D\dbar$ to vertical vector fields   with the quadratic differentials. It is a standard fact from Teichmuller theory that this  will give an identification of quadratic differentials with the holomorphic cotangent space to 
$U$ when $\ex C'$ has no internal edges. It follows from this fact that restricting to quadratic differentials that vanish on edges, we get the holomorphic cotangent space to the appropriate strata of the smooth part of $\tilde U$. Using this fact, it is not difficult to prove directly that the above gives an identification of the space of quadratic differentials with the holomorphic cotangent space to $\tilde U$ in general.

Add a bundle $V$ to our smooth map to get a smooth pre obstruction model $(\hat g,V)$ so that the fibers of $V$ are dual to the space of quadratic differentials, and $D\dbar$ is injective and has image complementary to $V$.  Theorem \ref{regularity theorem} implies that we can modify $(\hat g,V)$ on a  neighborhood of $\ex C$ to a $\C\infty 1$ pre obstruction model $(\hat f,V)$ with $\dbar f$ a section of $V$. Referring to this neighborhood as 
$ \mathbb R^{n}\times \ex F$, we have 
\[\begin{split} &V
\\ &\downarrow\uparrow\dbar
\\ \mathbb R^{n} &\times \ex F
\end{split}\]
The differential of $\dbar$ restricted to the $\mathbb R^{n}$ factor at $0$ is surjective due to the identification of the cotangent space of $U$ with the space of quadratic differentials. Therefore, there is a $\C\infty 1$ map from a  neighborhood of $\ex C$ in $\ex F$ to $\mathbb R^{n}\times \ex F$ so that the composition with $\dbar$ is $0$. This constructs a $\C\infty 1$ map from a  neighborhood of $\ex C$ to $\tilde U^{+1}$ which is fiber wise holomorphic, and which is equal to our chosen map on $\ex C$. The uniqueness proved above gives that this must agree with our map to $\M^{+1}$, therefore this map must therefore actually be $\C\infty 1$.

\stop

Recall the definition of a core family given on page \pageref{core family}. The following theorem gives criteria for when a given family with a collection of marked point sections is a core family:

\begin{thm}\label{core criteria}
A $\C\infty1$ family of stable $\omega$-positive curves $\hat f$  with group $G$ of automorphisms, finite, nonempty  set of disjoint sections $s_{i}:\ex F(\hat f)\longrightarrow \ex C(\hat f)$ which do not intersect the edges of the curves in $\ex C(\hat f)$ and a $\C\infty1$ map 
\[\begin{array}{ccc}\hat f^{*}T_{vert}\hat{\ex B}&\xrightarrow{F} &\hat {\ex B}
\\ \downarrow& & \downarrow
\\ \ex F(\hat f)&\longrightarrow &\ex G
\end{array}\]
 is a core family $(\hat f/G,\{s_{i}\},F)$ for some open neighborhood $\mathcal O\subset\Msw$ of $\hat f$  if and only if the following criteria are satisfied:
 \begin{enumerate}
  \item \label{crit1} For all curves $f$ in $\hat f$, there are exactly $\abs G$ maps of $f$ into $\hat f$ and the action of $G$ on the set of maps of $f$ into $\hat f$ is free and transitive.
 \item\label{crit2} For all curves $f$ in $\hat f$, the smooth part of the domain $\ex C(f)$ with the extra marked points from $\{s_{i}\}$ has no automorphisms.
 \item \label{crit3}The action of $G$ preserves the set of sections $\{s_{i}\}$, so there is some  action of $G$ as a permutation group on the set of indices $\{i\}$ so that for all $g\in G$ and $s_{i}$,
 \[s_{i}\circ g=g\circ s_{g(i)}\]
 where the action of $g$ is on $\ex F(\hat f)$, $\ex C(\hat f)$ or the set of indices $\{i\}$ as appropriate. 
 \item\label{crit4} \begin{enumerate}\item There exists a neighborhood $U$ of the image of the section
\[s:\ex F(\hat f)\longrightarrow \ex F(\hat f^{+n})\]
defined by the $n$ sections $\{s_{i}\}$ so that  
\[ev^{+n}(\hat f):\ex F(\hat f^{+n})\longrightarrow \M\times\sfp{\hat {\ex B}}{\ex G}n\]
restricted to $U$ is an equi-dimensional embedding

\item The tropical part of $ev^{+n}\circ s$ is a complete map, and restricted to any polytope in $\totb{\ex F(\hat f)}$ is an isomorphism onto a strata of the image  in  $\totb{\M\times\sfp{\hat {\ex B}}{\ex G}n}$ under $\totb{ev^{+n}}$ of some open neighborhood of $\hat f$ in $\Msw$.
\end{enumerate}

\item \label{crit5} \begin{enumerate}\item $F$ restricted to the zero section is equal to $\hat f$, 
 \item $TF$ restricted to the natural inclusion of $\hat f^{*}{T_{vert}\hat{\ex B}}$ over the zero section is equal to the identity,
 \item  $TF$ restricted to the vertical tangent space at any point of $\hat f^{*}T_{vert}\hat{\ex B}$ is injective.
 \end{enumerate} 
\end{enumerate}
\end{thm}

\pf Throughout this proof, use  $\ex F$ to refer to $\ex F(\hat f)$.
Consider the pullback of the family of curves $\hat f^{+n}$ under the map
 \[s:\ex F\longrightarrow \ex F(\hat f^{+n})\]
 This gives a family of curves $s^{*}(\hat f^{+n})$  over $\ex F$ with $n$ extra punctures.
 For any individual curve, $f\in \hat f$, criterion \ref{crit1}  implies that the family $f^{+n}$ contains exactly $\abs G$ curves which are contained in $s^{*}(\hat f^{+n})$, and criterion \ref{crit4} implies that  $ev^{+n}(f)(\ex F(f^{+n}))$ intersects the image of the section $s:\ex F\longrightarrow \ex F(\hat f^{+n})$ under $ev^{+n}(\hat f)$ transversely at each of the $\abs G$ points in $\ex F(f^{+n})$ corresponding to these curves.  Therefore, for any curve $f'$ in $\Msw$ sufficiently  close to $f$  in $C^{1,\delta}$, $ev^{+n}(f')(\ex F(f'^{+n}))$ intersects the image of $ev^{+n}(\hat f)(s(\ex F))$ transversely $\abs G$ times so that the corresponding $\abs G$ curves in $f'^{+n}$ are  close in $C^{1,\delta}$ to curves in $s^{*}(\hat f^{+n})$. In other words, there exists a $C^{1,\delta}$ neighborhood $\mathcal O$ of $\hat f$ in $\Msw$ and a $C^{1,\delta}$ neighborhood $\mathcal O_{s}$ of $s^{*}(\hat f^{+n})$ so that for any curve $f'$ in $\mathcal O$, $ev^{+n}(f')(\ex F(f'^{+n}\rvert_{\mathcal O_{s}}))$ intersects the image of $ev^{+n}(\hat f)(s(\ex F))$ transversely exactly $\abs G$ times, where $f'^{+n}\rvert_{\mathcal O_{s}}$ indicates the restriction of $f'^{+n}$ to curves in $\mathcal O_{s}$.
 
  It follows that for any family $\hat f'$ of curves in $\mathcal O$, the following fiber product comes with  an equidimensional submersion  of degree $\abs G$ to $\ex F(\hat f')$.

 \begin{equation}\label{fp1}
 \begin{split}&\ex F(\hat f'^{+n})\longleftarrow\ex F(\hat f'^{+n}\rvert_{\mathcal O_{s}})\fp{ev^{+n}(\hat f')}{ev^{+n}(\hat f)\circ s} \ex F(\hat f)
 \longrightarrow \ex F(\hat f)
  \\&\downarrow
 \\ &\ex F(\hat f')
 \end{split}\end{equation}
 
 Actually, criterion \ref{crit4} implies that the above fiber product is locally equal to a subset of $\ex F(\hat f'^{+n}\rvert _{\mathcal O_{s}})$ which is defined by the inverse image of some regular value of a $\C\infty1$, $\mathbb R^{2n}$ valued function, so the above submersion is actually a $\abs G$-fold cover of $\ex F(\hat f')$. 
   We therefore get a map from this $\abs G$-fold cover of $\ex F(\hat f')$ to  $\ex F(\hat f)$. Criterion \ref{crit3} implies that the action of $G$ on $\ex F(\hat f)$ gives an action of $G$ on $ev^{+n}(\hat f)(s(\ex F))$ which does nothing apart from permuting the marked points. As the image of $ev^{+n}(\hat f')$ automatically contains all the results of a permutation of marked points,  this $G$ action gives an action of $G$ on the above fiber product in (\ref{fp1}). This makes the above $\abs G$-fold cover of  $\ex F(\hat f')$ into a $G$-bundle because the action on $ev^{+n}(\hat f)(s(\ex F))$ simply permutes the marked points, so each $G$-orbit is contained within the same fiber of $\ex F(\hat f'^{+n})\longrightarrow \ex F(\hat f')$.  Therefore, the above map from our $G$-bundle to $\ex F$ is equivalent to a map from $\ex F(\hat f') $ to $\ex F/G$.  It follows from Theorem \ref{moduli} that if $\hat f'$ is $\C\infty 1$, this map is actually $\C\infty1$ map.

  \label{lift argument}
 
 There is a unique lift of this map  $\ex F(\hat f')\longrightarrow \ex F(\hat f)/G$ to a fiberwise holomorphic map 
 \[\begin{array}{ccc}\ex C(\hat f')&\xrightarrow{\Phi}& \ex C(\hat f)/G
 \\ \downarrow & & \downarrow 
 \\ \ex F(\hat f')&\longrightarrow &\ex F(\hat f)/G 
 \end{array}\]
 so that $\hat f'$ is equal to $\hat f\circ \Phi$ when restricted to the pullback under $\Phi$ of each of the sections $s_{i}$. 
 This map $\Phi$ is
constructed as follows: Consider the map 
\begin{equation}\label{tilde ev}\tilde{ev}^{+n}(\hat f):\ex C(\hat f^{+n})\longrightarrow \M\times\sfp{\hat{\ex B}}{\ex G}n\end{equation} which is equal to $ev^{+(n+1)}$ composed with a projection $\M\times\sfp{\hat {\ex B}}{\ex G}{n+1}\longrightarrow\M\times\sfp{\hat {\ex B}}{\ex G}n$ forgetting the image of the $(n+1)$st marked point and also equal on the second component to the composition of the projection $\ex C(\hat f^{+n})\longrightarrow \ex F(\hat f^{+n})$ with the map $\hat f^{+(n-1)}$. Criteria \ref{crit2} and \ref{crit4} imply that this evaluation map $\tilde {ev}^{+n}(\hat f)$ is an equidimensional embedding in a neighborhood of $\ex C(s^{*}\hat f^{+n})\subset \ex C(\hat f^{+n})$, and the following is a pullback diagram of families of curves

\[\begin{array}{ccc}\ex C(\hat f^{+n})&\xrightarrow{\tilde{ev}^{+n}(\hat f)}&\M^{+1}\times \sfp{\hat {\ex B}}{\ex G}{n}
\\\downarrow & &\downarrow
\\\ex F(\hat f^{+n})&\xrightarrow{ev^{+n}(\hat f)}&\M\times \sfp{\hat {\ex B}}{\ex G}{n}
\end{array}\]
 
Use the notation $\ex F'$ for the $G$-bundle over  $\ex F(\hat f')$ featured above in (\ref{fp1}) 
\[\ex F':=\ex F(\hat f'^{+n}\rvert_{\mathcal O_{s}})\fp{ev^{+n}(\hat f')}{ev^{+n}(\hat f)\circ s} \ex F(\hat f)\subset \ex F(\hat f'^{+n})\] 
The action of $G$ on $\ex F'$ is some permutation of marked points. This $G$ action extends to a $G$ action  on $\ex F(\hat f'^{+n})$ permuting these marked points, and  lifts to a $G$ action  on $\ex C(\hat f'^{+n})$ which just permutes the same marked points. Let  $\ex C'$ be the subset of $\ex C(\hat f'^{+n})$ over $\ex F'$.  The above mentioned $G$ action almost  makes $\ex C'$ a $G$-bundle over $\ex C(\hat f')$, except instead of locally being composed of $\abs G$ copies of $\ex C(\hat f')$, the `bundle' $\ex C'$ is locally equal to $\abs G$ copies of $\ex C(\hat f')$ with $n$ extra marked points.  

There is a unique  fiberwise holomorphic $\C\infty1$ map $\tilde \Phi$ from $\ex C'$ to $\ex C(s^{*}\hat f^{+n})$ so that $ \tilde{ev}^{+n}(\hat f)\circ\tilde \Phi=\tilde{ev}^{+n}(\hat f')$ on $\ex C'$. The extra marked points on $\ex C'$ are just the pullback of the extra marked points on $\ex C(s^{*}\hat f)$, and as all curves in $\hat f$ are stable, and the extra marked points are distinct points in the smooth components of $\ex C(\hat f)$, for $\hat f'$ close enough to $\hat f$,  the extra marked points on $\ex C'$ shall also correspond to  distinct points in the smooth components of $\ex C'$. It follows that we may forget the extra marked points in the domain and target of $\tilde \Phi$ to obtain a fiberwise holomorphic map $\C\infty1$ map from a $G$ bundle over $\ex C$ to $\ex C(\hat f)$  which corresponds to a $\C\infty 1$ map $\Phi:\ex C(\hat f')\longrightarrow \ex C(\hat f)/G$. The uniqueness of $\tilde\Phi$ implies the uniqueness of such a holomorphic map $\Phi$ so that restricted to the inverse image of the extra marked points, $\hat f'$ is equal to $\hat f$.

 The fact that this map $\Phi:\ex C(\hat f')\longrightarrow \hat{\ex C}/G $ is $\C\infty1$ means the composition of it with $\hat f$ is $\C\infty1$. This is by construction close to our other family of curves $\hat f'$, and is equal to $\hat f'$ on all the marked points coming from $\{s_{i}\}$. Criterion \ref{crit5} states  that $F$ restricted to fibers of $\hat f^{*}T_{vert}\hat{\ex B}$  is an equidimensional embedding into fibers of $\hat {\ex B}\longrightarrow \ex G$. Therefore, there is a unique  $\C\infty 1$  section $v$ of $\Phi^{*}\hat f^{*}T_{vert}\hat{\ex B}$ which vanishes at all marked points so that $f'=F(\Phi_{*}(v))$.  Therefore $(\hat f/G,\{s_{i}\},F)$ is a core family for $\mathcal O$.
 
 \stop

\

The following proposition constructs a core family containing a given stable holomorphic curve which has at least one smooth component (so it isn't $\ex T$).
 
\begin{prop}
\label{smooth model family}

Given a stable,  connected, holomorphic curve $f$ with at least one smooth component in a basic family of exploded manifolds  $\hat{\ex B} $ and a collection of marked points $\{p_{j}\}$ in the interior of  the smooth components of $\ex C(f)$, 
  there exists a $\C\infty1$ core family $(\hat f/G,\{s_{i}\},F)$ with $\hat f$ a family containing $f$ so that the restriction of $\{s_{i}\}$ to $f$ contains the given marked points $\{p_{j}\}$.

\end{prop}

\pf 
The automorphism group $G$ of $\hat f$ shall be equal to the group of automorphisms of $\totl f$, the smooth part of $f$. 

 By restricting to an open subset of our family $\hat {\ex B}\longrightarrow \ex G $ containing the image of our curve, we may assume that $\totb{\ex G}$ is an integral affine polytope.

We shall enumerate the steps of this construction so that we can refer back to them 
\begin{enumerate}

\item Choose extra marked points on the smooth components of $\ex C$ so each smooth component of $\ex C$ contains at least one marked point, the smooth part of $\ex C$ has no automorphisms with these extra marked points, and so that that we can divide the marked points on $\ex C$ into the following types:
\begin{enumerate}
\item \label{mp1} On any smooth component of $\ex C$ which is unstable,  choose enough extra marked points at which $d\totl{f}$ is injective to stabilize the component. Note that $G$ has a well defined action on the smooth part of $\ex C$. Choose the set of marked points of this type to be  preserved by the action of $G$. (Note that the fact that  $f$ is a stable holomorphic curve  implies that each unstable smooth component of $\ex C$ must contain a nonempty open set where $d\totl{f}$ is injective.)

\item\label{mp2} Choose the set of remaining marked points to be preserved by the action of $G$. 
\end{enumerate}

\item The tropical part of the family $\ex F$ will be an integral affine polytope $P=\totb{\ex F}$, or a disjoint union of some number of copies of this polytope $P$.  Construct this polytope $P$ as follows:  

\begin{enumerate}
\item Construct a polytope $\check P$ as follows: The image of each marked point is contained in the interior strata of a coordinate chart with tropical part $\mathcal P(f(p_{i}))$. Construct $\check P$ by taking the fiber product over $\totb{\ex G}$ of the polytopes $\mathcal P(f(p_{i}))$, then take the product of this with  a copy of $\e{[0,\infty)}$ for every internal edge of $\ex C$.
The coordinates on $\check P$ record the tropical position of each marked point, and the length of each internal edge of the tropical curves in our family.  Note that $\check P$ is a convex integral affine polytope. 
 
\item \label{tconstraints}Consider the tropical structure $f_{T}$ of the map $f$ discussed in \cite{iec}. Corresponding to each marked point $p$, there is a point  $f_{T}(\mathcal P(p))\in\mathcal P(f(p))$. Corresponding to each homotopy class $\gamma$ of path between marked points $p_{i}$ and $p_{j}$ on a smooth component, there is a map $\mathcal P(f_{T}(\gamma)): \mathcal P(p_{i})\longrightarrow \mathcal P(p_{j})$. The requirement that $\mathcal P(f_{T}(\gamma))(f_{T}(\mathcal P(p_{i})))=f_{T}(\mathcal P(p_{j}))$
is an integral affine condition. Similarly, there is an integral affine condition corresponding to an internal edge $e$ of $\ex C$ as follows: The image of the edge $e$ is contained in a coordinate chart with tropical part $\mathcal P(f(e))$. Given a path $\gamma_{i}$  joining the marked points $p_{i}$ to the edge $e$, there is an inclusion $\mathcal P(f_{T}(\gamma_{i})):\mathcal P(f(p_{i}))\longrightarrow \mathcal P(f_{T}(e))$. If $p_{i}$ and $p_{j}$ are at opposite ends of the edge $e$, the requirement  that $\mathcal P(\gamma_{i})(\mathcal P(f(p_{i})))$ and $\mathcal P(\gamma_{j})(\mathcal P(f(p_{j})))$ are joined by an edge of the specified length and with the same velocity as the original edge of $f$ is an integral affine condition on $\check P$.
The polytope $P\subset \check P$ is the solution to these integral affine conditions. As $\check P$ was a  convex integral affine polygon, $P\subset\check P$ is too. $P$ must be nonempty because it contains a point corresponding to $\totb f$.

\end{enumerate}

\item \label{strata charts}
Use equivariant coordinate charts on $\hat{\ex B}\longrightarrow \ex G$ as constructed in \cite{cem}. Recall that a coordinate chart $\mathbb R^{n}\times\et mP$ has a $\ex T^{m}$ action which is given by multiplying the coordinates of $\et mP$ with constants in $\ex T^{m}$. So
\[(c_{1},\dotsc,c_{m})*(x,\tilde z_{1},\dotsc,\tilde z_{m})=(x,c_{1}\tilde z_{1},\dotsc,c_{m}\tilde z_{m})
\] (This action is only sometimes defined.) A map $\phi:\mathbb R^{n'}\times \et {m'}{P'}\longrightarrow \mathbb R^{n}\times\et mP$ is equivariant if there exists a homomorphism  $\psi:\ex T^{m'}\longrightarrow \ex T^{m}$ so that 
$\phi(c*p)=\psi(c)*\phi(p)$. Note that every $\C\infty1$ map $\phi$ is equivariant restricted to the interior of $P'$. Coordinate charts on $\hat{\ex B}\longrightarrow \ex G$ are equivariant if each transition map is either equivariant or its inverse is equivariant, and if the projection to $\ex G$ in coordinates is equivariant.

\item Construct $\ex F$ as follows: Coordinates on $\ex F$ shall be given by the position of marked points and the complex structure of our curves. $\ex F$ will be equal to some open subset of $\mathbb R^{n}\times \et mP$. First construct $\check {\ex F}$ in analogy to $\check P$ which will have tropical part $\check P$. $\ex F$ will be an open subset of a refinement of $\check{\ex  F}$ corresponding to $P\subset\check P$.

\begin{enumerate}
\item Construct $\check {\ex F}$ using the following coordinates:
\begin{enumerate}
\item\label{F1}  For each marked point $p$ from part \ref{mp1}, as  $d\totl{f}(p)$ is injective, we can choose a  coordinate chart on $\hat {\ex B}$ which identifies  a neighborhood of $f(p)$ with $\mathbb R^{2}\times \mathbb R^{k}\times \et m{\mathcal P(f(p))}$,  so that the restriction of $f$ to a neighborhood of $p$ is equal to an inclusion $x\in\mathbb R^{2}\mapsto (x,c)$ where $c\in\mathbb R^{k}\times\et m{\mathcal P(f(p))}$ and $p$ is given by $x=0$.  We can also construct our coordinate chart above so that the slices $\mathbb R^{2}\times c'$ are all contained in a fiber of $\hat{\ex B}\longrightarrow\ex G$, so that there is a well defined submersion  $\mathbb R^{k}\times \et m{\mathcal P(f(p))}\longrightarrow \ex G$. We can construct this coordinate chart by taking one of the equivariant coordinate charts from item \ref{strata charts} and reparametrising the $\mathbb R^{k+2}$ factor, which will not affect the equivariance property. Use the same coordinate chart for each marked point in an orbit of $G$. Include in our coordinates for $\check{\ex F}$ the fiber product of  $\mathbb R^{k}\times \et m{\mathcal P(f(p))}\longrightarrow \ex G$ for every marked point $p$ from part \ref{mp1}.  
\item\label{F3}  The image of each marked point from part \ref{mp2} is contained in the interior strata of one of our equivariant coordinate charts $\ex U_{i}$. Take the fiber product over $\ex G$ of the coordinates from part \ref{F1} with a copy of $\ex U_{i}\longrightarrow \ex G$ for each marked point from part $\ref{mp2}$.
\item\label{F4} Each smooth component of  $\ex C$ can be regarded as a stable punctured Riemann surface with labeled punctures determined by the exploded structure of $\ex C$ plus the extra marked points from part \ref{mp1}. Take the product of the above coordinates from part \ref{F1} and \ref{F3}  with a  uniformizing neighborhood of the corresponding point in Deligne Mumford space for each smooth component of $\ex C$.  Do this so that the obvious $G$ action is well defined. 
\item\label{F5} Take the product of the above coordinates with a copy of $\et 11$ for each internal edge of $\ex C$. 
\end{enumerate}
\end{enumerate}

Observe that the tropical part of $\check{\ex F}$, $\totb{\check{\ex F}}$ is equal to $\check P$. If $\check P$ is $m'$ dimensional, there is a (sometimes defined) free $\ex T^{m'}$ action on $\check{\ex F}$ given by multiplication on the correct coordinates from item \ref{F1}, \ref{F3} and \ref{F5}. The (sometimes defined) action of a subgroup $\ex T^{m}\subset \ex T^{m'}$ preserves $P\subset\check P$ where $P$ is $m$ dimensional.  There is a corresponding action of $\ex T^ m$ on each of the coordinate charts $\ex U_{i}$ referred to in item \ref{F3} (which of course is not necessarily free).

 There is a distinguished point  $p_{0}\in \check{\ex F}$ corresponding to our curve $f$, which is the point $f(p)$ in item \ref{F1} and \ref{F3}, the point corresponding to the complex structure on the smooth components of $\ex C$ in item \ref{F4},  and for item \ref{F5},   $1\e l$ where the strata of $\ex C$ corresponding to the internal edge in question is equal to $\et 1{(0,l)}$.
  Roughly speaking, our family $\ex F$ will be some neighborhood of orbit of this point $p_{0}$ under an action of $G$ and the above mentioned $\ex T^{m}$ action.

\item Construct $\ex F$ as follows:
\begin{enumerate}
\item Take any refinement $\check{\ex F}'$ of $\check{\ex F}$ so that the $\totb{\check{\ex F}'}$ includes a strata with closure equal to $P\subset\check P$.  
\item \label{gluing limit} There is an action of $G$ the subset of $\check{\ex F}'$ with tropical part $P$  which shall be defined in item \ref{hat Ai} below.  $\ex F$ is given by the orbit under $G$ of  a small  open neighborhood of  the point $p_{0}\in \check{\ex F}'$ corresponding to $f$ so that the coordinates from item \ref{F5} have absolute value strictly less that some $\epsilon \e 0$. 
\end{enumerate}

\item We shall now construct $(\hat{\ex C},j)\xrightarrow{\pi_{\ex F}}\ex F$. Roughly speaking, the coordinates \ref{F4} and \ref{F5} give a map from $\ex F$ to $\M$ which at $f$ corresponds to the complex structure on $\ex C$ with the extra punctures mentioned in \ref{mp1}. Pulling back $\M^{+1}\longrightarrow \M$ gives $(\hat{\ex C},j)$ with the sections corresponding to marked points from \ref{mp1}, and we just need to extend the other marked points to appropriate sections to define $(\hat {\ex C},j)\longrightarrow \ex F$ with all its sections $s_{i}$.  We shall do this below in an explicit way to enable us to describe more explicitly the extension of  the map $f$ to $\hat f$. 
  
\begin{enumerate}
\item Choose holomorphic identifications of a neighborhood of each internal edge of $\ex C$ with 
\[A_{i}:=\{ c_{i}\e 0>\abs{\tilde z}> 1\e l\}\subset \et 1{[0,l]}\]
so that these neighborhoods $A_{i}$ are disjoint, and all marked points are in the complement of these annuli $A_{i}$. Do this so that the set of images of $A_{i}$ in the smooth part of $\ex C$ are preserved by the action of $G$. Also choose  $ c_{i}>8\epsilon$ where $\epsilon$ is the constant mentioned in part \ref{gluing limit} above. (Of course, to achieve this, we need to choose $\epsilon$ small enough.)
\item \label{hat Ai}
Use the notation $\hat A_{i}$  to refer to the part of $\hat{\ex  C}$ corresponding to $A_{i}$. This is given as follows:

In the construction of $\check{ \ex F}$, replace the factor of $\et 11$ from item \ref{F5} corresponding to the edge $A_{i}$ with $\et 22$. If this $\et 22$ has coordinates $\tilde w_{1},\tilde w_{2}$, then let $\check A_{i}$ be the subset of this $\et 22$ so that 
$\abs{\tilde w_{1}}<c_{i} \e 0$ and $\abs{\tilde w_{2}}<1\e 0$. The map $\check A_{i}\longrightarrow \check{\ex F}$ is given by the map $\tilde w_{1}\tilde w_{2}:\et 22\longrightarrow \et 11$ and is the identity on all other coordinates.

There is a natural action of $G$ on the union of these $\check A_{i}$ given as follows: On the coordinates corresponding to  all coordinates on $\ex F$ apart from part \ref{F5}, there is an obvious action of $G$. If $g\in G$ sends the smooth part of $A_{i}$ to $ A_{j}$, the pull back of smooth part of the coordinates $\totl{\tilde w_{1}}, \totl{\tilde w_{2}}$ on $ A_{j}$ is equal to some constant times the smooth part of the corresponding coordinates on $ A_{i}$. Define the map from $\check A_{i}$ to $\check A_{j}$ by defining the pull back of $\tilde w_{1} $ and $\tilde w_{2}$ simply to be the corresponding coordinate multiplied by the above constant, and the pull back of other coordinate functions as given by the obvious $G$ action on coordinates from parts \ref{F1}, \ref{F3} and \ref{F4}. This induces an action of $G$ on $\check {\ex F}$ so that the map $\bigcup _{i}\check A_{i}\longrightarrow \check{\ex F}$ is $G$ equivariant. 

We have chosen $\ex F$ to be a $G$ equivariant subset of $\check{\ex F}$. Let  $\hat{A_{i}}$ be the restriction of $\check{ A_{i}}$ to the inverse image of $\ex F\subset\check{\ex F}$.

\item Label by $C_{k}$ the connected components of the complement of the sets  $\{\frac {c_{i}}2\e 0>\abs{\tilde z}>2\e l\}\subset A_{i}$ . 
Again, use the notation $\hat C_{k}$ to refer to the corresponding part of $\hat {\ex C}$. This is simply given by 
the product 
\[\hat C_{k}:=C_{k}\times {\ex F} \]
The map $C_{k}\longrightarrow{\ex F}$ is simply projection onto the second component. Note that there is an action of $G$ on the union of these $C_{k}$ given by the action on $\ex F$ defined at the end of item \ref{hat Ai} above, and the action of $G$ on the union of $C_{k}$ as a subset of the smooth part of $\ex C$.
\item The transition maps between $A_{i}$ and $C_{k}$ induce in an obvious way transition maps between $\hat A_{i}$ and $\hat C_{k}$, which defines the family $\hat{\ex C}\longrightarrow\ex F$. Note that the inverse image of our special point $p_{0}\in \ex F$ is equal to $\ex C$. Note also that these transition maps are compatible with the action of $G$ on  the union of the  $\hat A_{i}$ and the union of the $\hat C_{k}$, so there is an action of $G$ on $\hat{\ex C}\longrightarrow \ex F$. 
  
  Remembering the positions of our marked points in $C_{k}$ gives the sections $s_{i}:\ex F\longrightarrow \hat{\ex C}$ referred to in the statement of this proposition.
\item It remains to construct the  complex structure $j$ on the fibers of $\hat{\ex C}$. Recall that the coordinates on $\ex F$ from item \ref{F4} are intended to give the almost complex structure on smooth components of $\ex C$. Choose a smooth family of complex structures $j$ on the smooth components of $\ex C$   parameterised by these coordinates with the correct isomorphism class of complex stucture, so that $j$ at our special point is the original complex structure on $\ex C$, and $j$ restricted to the the subsets $A_{i}$ is also the original complex structure. Do this equivariantly with respect to the action of $G$ on the smooth part of $\ex C$ and the action of $G$ on the coordinates from part \ref{F4}. This gives a family of  complex structures  on the fibers of  $\hat C_{k}$. This is compatible with the standard holomorphic structure on $\hat A_{i}$, so using this gives us our globally defined $(\hat {\ex C}, j)$. Note that the restriction of this to the curve corresponding to our special point $p\longrightarrow \ex F$ will give $\ex C$ with the original complex structure.

\end{enumerate}

\item Construct the family of maps $\hat f:\hat{\ex  C}\longrightarrow \hat {\ex B}$. This will involve translating around pieces of the original map $f$,  modifying this map near marked points as directed by the coordinates of $\ex F$, and gluing together the result. The last `gluing' step  only affects the map $\hat f$ on $\hat A_{i}$, so we shall now perform the first two steps to construct $\hat f$ on $\hat C_{k}$.

Construct $\hat f$ on $\hat C_{k}$ as follows:
\begin{enumerate}
\item   Recall that there is a (sometimes defined) action of $\ex T^{m}$ on $\ex F$ and a corresponding 
action of $\ex T^{m}$ on the coordinate charts on $\hat{\ex B}$ containing the marked points. As we chose our coordinate charts on $\hat{\ex B}$ equivariantly and $\hat {\ex B}$ is basic, this action of $\ex T^{m}$ can be extended to a collection of coordinate charts which contain the image of a smooth part of $\ex C$.
If $p'=\tilde c *p_{0}$, where $p_{0}\in\ex F$ is the special point corresponding to $f$, then define
\[\hat f(z,p'):= c*f(z)\text{ when }z\in C_{k}\]
This defines $\hat f$ on the part of $\hat C_{k}$ which projects to the orbit of $p_{0}\in \ex F$ under the action of $\ex T^{m}$. Note that this map is preserved by the action of  the subgroup of $G$ which sends $p_{0}$ to some $\tilde c*p_{0}$. We may extend the definition of $\hat f$ to be $G$ equivariant on the orbit of $p_{0}$ under the $\ex T^{m}$ action and the action of $G$.
 
\item We must make sure that each of the individual smooth curves in $\hat f$ are contained in the correct fiber of $\hat{\ex  B}\longrightarrow \ex G$. Note that this is automatically true so far, because of the compatibility of our $G$ and $\ex T^{m}$ action with the map $\hat{\ex B}\longrightarrow \ex G$. (In fact, there is a (sometimes defined) action of $\ex T^{m}$ on $\ex G$ so that this map is equivariant.) We shall now extend the definition of $\hat f$ to a subset of $\hat{\ex C}$ which is equal to $\mathbb R^{n'}$ times where it is already defined by `translating in directions coming from $\ex G$'. 

As constructed, the obvious map (trivial on all coordinates apart from those from items \ref{F1} and \ref{F3}), $\ex F \longrightarrow\ex G$ is a submersion which is preserved by the action of $G$. The image of the tropical part $\totb{\ex F}$ under this map is some polytope $Q\subset \totb{\ex G}$ and the image of $\ex F$ under this map is an open subset of some refinement of $\ex G$ that has the interior of  $Q$ as a strata. If $\ex F$ is chosen small enough, this open subset of the refinement of $\ex G$ is isomorphic to $\mathbb R^{n'}\times \et {m'}Q$. We can pull our family $\hat{\ex B}\longrightarrow \ex G$ back to to be a family over $\mathbb R^{n'}\times \et {m'}Q$. If $\ex F$ was chosen small enough, this family will split into a product $\mathbb R^{n'}\times\hat{\ex B}'\longrightarrow\mathbb R^{n'}\times \et {m'}Q$ which is the identity on the $\mathbb R^{n'}$ component, and some family $\hat{\ex B}'\longrightarrow \et {m'}Q$ on the second component. We can choose this splitting so that it is compatible with our local actions of $\ex T^{m}$ on coordinate charts $\ex U_{i}$. This also gives a splitting of $\ex F$ into $\mathbb R^{n'}\times \ex F'$. We can choose this splitting so that the subset of $\ex F$ where we have already defined $\hat f$ is contained inside $0\times\ex F'$. Now we can define $\hat f$ as a map to $\mathbb R^{n'}\times\hat{\ex B'}$ follows:
\[\hat f(z,x,p')=(x, y)\text{ if }f(z,0,p')\text{ is already defined, and }f(z,0,p')=(0,y)\]
Here $(z,x,p')$ denotes coordinates on $\hat C_{k}=C_{k}\times\mathbb R^{n'}\times \ex F'$. This map is defined on a $G$ invariant subset, and is preserved by the action of $G$.
\item Split $\ex F$ further into an open subset of $\mathbb R^{n''}\times\ex F''$, so that our map $\hat f$ is defined so far on the subset of $\hat C_{k}$ which is over $0\times \ex F''$, and the splitting is preserved by the action of $G$. Extend the map defined so far to a smooth map $\hat f$ defined on all of $C_{k}$ so that
\begin{enumerate}
\item $\hat f$ fits into the commutative diagram
\[\begin{array}{ccc} 
  \hat C_{k} & \xrightarrow{\hat f} &\hat{\ex B}
  \\ \ \ \ \  \downarrow\pi_{\ex F} & &\ \ \ \ \downarrow\pi_{\ex G}
  \\  \ex F &\longrightarrow &\ex G\end{array}\]

\item $\hat f$ is preserved by the action of $G$ on $\hat C_{k}$.

\item On the intersection of $\hat A_{i}$ with $\hat C_{k}$ and outside a small neighborhood of all marked points, $\hat f(z,x,y)=\hat f(z,0,y)$. (This uses coordinates $\hat C_{k}=C_{k}\times\mathbb R^{n''}\times \ex F''$.)
\item For each marked point $q$ from part \ref{mp1}, the value of $\hat f$ at the point  $(q,x,y)$ is equal to the corresponding coordinate of $\ex F$ from part \ref{F1}.
\item For each marked point $q$ from part \ref{mp2},  $(x,y)\in \mathbb R^{n''}\times \ex F''$ determines a value for the coordinate on $\ex F$ from part $\ref{F3}$ which is a point some coordinate chart. For such a marked point, $\hat f(q,x,y)$ is equal to this point.
 
\end{enumerate}

\end{enumerate}

\item Now to define $\hat f$ on $\hat A_{i}$. Consider the subset of $\ex F$ obtained by taking the orbit of the point $p$ corresponding to $f$ under the previously mentioned action of $\ex T^{m}$ and the action of $G$.

  Cut $\hat A_{i}$ into two pieces. Translate each piece the same way as the $C_{k}$ it is attached to, and then use a smooth gluing procedure to glue together the result which only modifies $\hat f$ on the region  where $\hat A_{i}$ does not intersect $C_{k}$. ( Examples of such a smooth gluing procedure are given in \cite{cem} and \cite{reg}.) Do this so that $\hat f$ is compatible with $\hat {\ex B}\longrightarrow \ex G$, and $\hat f$ is preserved by the action of $G$. Note that modification is not necessary over the point corresponding to $f$.


\

\item
We have now constructed the required family of stable $\omega$-positive curves.

\[\begin{array}{ccc} 
  (\ex {\hat C},j) & \xrightarrow{\hat f} &(\hat{\ex B},J)
  \\ \ \ \ \  \downarrow\pi_{\ex F} & &\ \ \ \ \downarrow\pi_{\ex G}
  \\  \ex F &\longrightarrow &\ex G\end{array}\]
This map $\hat f$ is smooth or $\C\infty1$ if $f$ is.  The family $\hat f$ with the sections $\{s_{i}\}$ satisfies criteria \ref{crit1}, \ref{crit2} and \ref{crit3} from Theorem \ref{core criteria} on page \pageref{core criteria}, and if $F$ is  any map 
\[\begin{array}{ccc}\hat f^{*}T_{vert}\hat{\ex B}&\xrightarrow{F} &\hat {\ex B}
\\ \downarrow& & \downarrow
\\ \ex F(\hat f)&\longrightarrow &\ex G
\end{array}\]
given by exponentiation using some smooth family of connections on $\hat {\ex B}\longrightarrow \ex G$, then $F$ is $\C\infty1$ if $\hat f$ is, and $F$ satisfies the criterion \ref{crit5} from Theorem \ref{core criteria}. Therefore, it remains to check criterion \ref{crit4} from Theorem \ref{core criteria}.

In the remainder of this proof, let $n$ denote the number of our sections, so $n=\abs{\{s_{i}\}}$.
We need to show that the evaluation map from $\ex F(\hat f^{+n})$ to $\M$ times the fiber product of $\hat{\ex B}\longrightarrow \ex G$ with itself $n$ times 
\[ev^{+n}(\hat f):\ex F(\hat f^{+n})\longrightarrow \M\times \sfp {\hat {\ex B}}{\ex G}{n}\]
 is an equidimensional embedding when restricted to some neighborhood of the section $s:\ex F\longrightarrow\ex F(\hat f^{+n})$ given by our sections $\{s_{i}\}$, and to check a condition on the tropical part of $ev^{+n}(\hat f)\circ s$.  Coordinates on a neighborhood of the image of $\ex F\longrightarrow\ex F(\hat f^{+n})$ are given by coordinates on $\ex F$ and coordinates on a neighborhood of each marked point. The evaluation map from this neighborhood splits into two equidimensional embeddings as follows: By construction, the coordinates from \ref{F4} and \ref{F5} together with the coordinates around all marked point not of type \ref{mp1} describe an equidimensional embedding into $\M$. The  coordinates from \ref{F3}, and the coordinates from \ref{F1} plus the coordinates around marked points of type \ref{mp1} describe an equidimensional embedding to $\sfp {\hat {\ex B}}{\ex G}{n}$ restricted to a small enough neighborhood. (We should restrict to a suitably small $G$ equivariant subset so that this holds.) 
 
 Recall that the tropical part of $\ex F$ is some number of copies of a polytope $P$, which we constructed by subjecting a 
 polytope $\check P$ in the tropical part of $\M\times \sfp{\hat {\ex B}}{\ex G}{n}$ to the conditions that a tropical curve must satisfy, so $P$ is a polytope in the image of $\Msw$ in the tropical part of $\M\times \sfp{\hat {\ex B}}{\ex G}{n}$ under $\totb{ev^{+n}}$. By construction $\totb{ev^{+n}\circ s}$ is an isomorphism from $P$ considered as a polytope in the tropical part of $\ex F$ to $P$ considered as a a polytope in the image of $\Msw$ in the tropical part of $\M\times \sfp{\hat {\ex B}}{\ex G}{n}$ under $\totb{ev^{+n}}$. Therefore, criterion \ref{crit4} holds.
 
 We have now checked that $(\hat f/G,\{s_{i}\})$ satisfies the requirements of Theorem \ref{core criteria}, and may construct the additional map $F$ required for criterion {crit5} using a smooth connection and parallel transport then reparametrizing so that $TF$ remains injective on fibers. Therefore $(\hat f/G,\{s_{i}\}, F)$ is a core family for some open substack of $\Msw$. 
 
\end{enumerate}
\stop

We now prove the existence of obstruction models, defined on page \pageref{obstruction model}.

\begin{thm}\label{construct obstruction model} Any stable  holomorphic curve $f$ with at least one smooth component in a basic exploded manifold $\ex B$ is contained inside some $\C\infty 1$ obstruction model $(\hat f/G,V)$. 
\end{thm}

\pf

 \
 
Proposition \ref{smooth model family} tells us that the curve $f$ must be contained in a $\C\infty 1$ core family $\hat h/G$. We may choose this core family to include any collection of marked points on $f$, so after choosing a ($G$ invariant) trivialization, Theorem \ref{f replacement} implies that we may arrange that $D\dbar:X^{\infty,\underline 1}(f)\longrightarrow Y^{\infty,\underline 1}(f)$ is injective. Theorem \ref{f replacement} implies that we may then choose a finite dimensional  complement $V_{0}(f)$ to $D\dbar(X^{\infty,\underline 1}(f))$ consisting of $\C\infty1$ sections of $Y(f)$. Below, we check that we can make $V_{0}(f)$ compatible with the action of $G$ and extend it to a $G$ invariant pre obstruction model $(\hat h,V)$.

Extend this $V_{0}$ to all curves in $\hat h$ with the same smooth part as $\hat f$. 
 Consider $V_{0}$ as giving a projection onto the ($G$-invariant) image of $D\dbar$. We may average this projection under the action of $G$ to obtain a $G$-invariant projection corresponding to a different, $G$ invariant complement $V'$ to $D\dbar(f)$, so far defined over the curves in $\hat h$ with the same smooth part as $\hat f$. This $V'$ is canonically a trivial $\mathbb R^{n}$ bundle. The action of $G$ gives a linear action $\rho$ of $G$ on $\mathbb R^{n}$ so that $g*( f,x)=(g*f,\rho(g)x)$. In a neighborhood of $f$, choose $\C\infty1$ sections $v'_{1},\dotsc,v'_{n}$ of $\Y{\hat h}$ so that restricted to $f$, $v_{i}$ is close to the section of $\Y f$  corresponding to the $i$th standard basis vector of $\mathbb R^{n}$. Then consider the $n$ sections of $\Y{\hat h}$ defined by
 \[[v_{1},\dotsc,v_{n}]:=\frac 1{\abs{G}}\sum_{g\in G}[g*v'_{1},\dotsc,g^{*}v'_{n}]\rho(g)^{-1}\]
  Note that the action of $G$ will preserve the span of these sections: 
  \[\begin{split}[g'*v_{1},\dotsc,g'*v_{n}]&=\frac 1{\abs{G}}\sum_{g\in G}[(g'g)*v'_{1},\dotsc,(g'g)^{*}v'_{n}]\rho(g)^{-1}
  \\ &=\frac 1{\abs{G}}\sum_{g\in G}[g*v'_{1},\dotsc,g^{*}v'_{n}]\rho((g')^{-1}g)^{-1}
  \\ &=[v_{1},\dotsc,v_{n}]\rho(g')\end{split}\]
 Note also that if restricted to $f$, each $v'_{i}$ was chosen to be exactly equal to the section of $Y(f)$ corresponding to the $i$th standard basis vector of $\mathbb R^{n}$, then $v_{i}'=v_{i}$. Therefore, close to $f$, the span of the sections $v_{1},\dotsc,v_{n}$ gives a $G$-invariant pre obstruction model $(\hat h,V)$ so that restricted to $f$, $V$ is complementary to the image of $D\dbar$. Then Theorem \ref{regularity theorem} gives that restricted to a small enough  open neighborhood of $f$,  we may modify $\hat g$ to a $G$ equivariant $\C\infty1$ family  $\hat f$ so that $\dbar f$ is a section of $V$, and  $(\hat f/G,V)$ is a $\C\infty 1$ obstruction model. 

\stop

\

\

The following theorem describes the `solution' to the  $\dbar$ equation  perturbed by multiple simple perturbations parametrized by different obstruction models. 

\

\begin{thm} \label{Multi solution}
Given
\begin{itemize}
\item a finite collection of  core families $\hat f'_{i}/G_{i}$ for the substacks $\mathcal O'_{i}\subset\Msw$,
\item  an open substack $\mathcal O\subset\Msw$ which  meets $\mathcal O'_{i}$ properly for all $i$ (definition \ref{proper meeting}), 
\item an obstruction model $(\hat f_{0}/G_{0},V)$ for the substack $\mathcal O$
 \item compactly contained $G_{i}$ invariant sub families $\hat f_{i}\subset \hat f'_{i}$,
 \end{itemize}
 then given any collection of  $\C\infty1$ simple perturbations $\mathfrak P_{i}$ parametrized by  $\hat f'_{i}$ which are compactly supported in $\hat f_{i}$ and are small enough in $\C\infty1$, and a sufficiently small simple perturbation $\mathfrak P_{0}$ parametrized by $\hat f_{0}$ there exists a solution mod $V$ on $\hat f_{0}$  which is a $G_{0}$ equivariant $\C\infty1$ weighted branched section $(\nu,\dbar'\nu)$ of $\hat f_{0}^{*}T_{vert}\hat{\ex B}\oplus V$  with weight $1$ (see example \ref{wbv} below definition \ref{multi}) so that the following holds:

Locally on $\ex F(\hat f_{0})/G_{0}$ (so restricted to a small enough $G_{0}$ equivariant neighborhood of any curve in $\ex F(\hat f_{0})$), 
\[(\nu,\dbar'\nu)=\sum_{l=1}^{n}\frac 1n t^{(\nu_{l},\dbar'\nu_{l})}\]
where $\dbar'\nu_{l}$ is a section of $V$, and  $\nu_{l}$ is in $X^{\infty,\underline 1}(\hat f_{0})$ so that $F(\nu_{l})$ is in $\mathcal O$ and in the notation of example \ref{proper pullback} on page \pageref{proper pullback}
\[\prod_{i}F(\nu_{l})^{*}\mathfrak P_{i}=\sum_{j=1}^{n} \frac 1nt^{\mathfrak P_{j,l}}\]
 
 so that 
 \[\dbar'\nu_{l}=\dbar F(\nu_{l})-\mathfrak P_{l,l}\]
 
  The weighted branched section $(\nu,\dbar'\nu)$ is the unique weighted branched section of $\hat f_{0}^{*}T_{vert}\hat {\ex B}\oplus  V$ with weight $1$ satisfying the following two conditions:
 \begin{enumerate}
 \item
   Given any curve $f$ in $\mathcal O$ and choice of $f'\in\hat f$ and $\psi\in X^{\infty,\underline 1}(f')$, so that 
   \[f=F(\psi)\]
    if 
    \[\prod_{i} f^{*}\mathfrak P_{i}=\sum w_{k}t^{Q_{k}}\] and near $f'$,  \[(\nu,\dbar'\nu)=\sum w'_{l}t^{(\nu_{l},\dbar'\nu_{l})}\] then the sum of the weights $w_{k}$ so that $\dbar f-Q_{k}$ is in $V$ is equal to  the sum of the weights $w'_{l}$  so that $\psi$ is equal to $\nu_{l}(f')$.
 \item For any  locally defined section $\psi$ in $X^{\infty,\underline 1}(\hat f_{0})$, if the multi perturbation 
 $\prod_{i} F(\psi)^{*}\mathfrak P_{i}=wt^{Q}+\dotsc$, and $\dbar F(\psi)-Q$ is a section of $V$, then locally, $(\nu,\dbar'\nu)=wt^{(\psi,\dbar F(\psi)-Q)}+\dotsc$.
 \end{enumerate}
 
This weighted branched section determines the solutions to the perturbed $\dbar$ equation on $\mathcal O$ in the following sense:  
Given any family $\hat g$ in $\mathcal O$, if $\prod_{i}\hat g^{*}\mathfrak P_{i}=wt^{\dbar\hat g}+\dotsc$,  then around each curve in $\hat g$ which projects to the region where the above $\nu_{l}$ are defined, there is a connected open neighborhood in $\hat g$ with at least $nw$  different  maps to $\lrb{\coprod_{l} F(\nu_{l})}/G_{0}$.

\ 
 
 If $\{\mathfrak P_{i}'\}$ is another collection of simple perturbations satisfying the same assumptions as $\{\mathfrak P_{i}\}$ then the sections $(\nu_{l}',\dbar'\nu_{l}')$ corresponding to $(\nu_{l},\dbar'\nu_{l})$, with the correct choice of indexing can be forced to be  as close to $(\nu_{l},\dbar'\nu_{l})$ as we like in $\C\infty1$ by choosing $\{\mathfrak P_{i}'\}$ close to $\{\mathfrak P_{i}\}$  in $\C\infty1$. If $\{\mathfrak P_{i,t}\}$ is a $\C\infty1$ family of simple perturbations satisfying the same assumptions as $\{\mathfrak P_{i}\}$, then the corresponding family of solutions mod $V$, $(\nu_{t},\dbar'\nu_{t})$ form a $\C\infty1$ family of weighted branched sections.

\end{thm}

\pf

Use $\mathcal O_{i}$ to denote the restriction of $\mathcal O_{i}'$ to the subset with core $\hat f_{i}/G_{i}$. As $\hat f_{0}'$ meets $\mathcal O'_{i}$ properly for all $i$, and $\hat f_{i}$ is  compactly contained in $\hat f_{i}'$, 
there is some $C^{1,\delta}$ neighborhood $U$ of $0$ in $X^{\infty,\underline 1}(\hat f_{0}')$ and some finite covering of $(\hat f_{0}/G_{0},V)$ by  extendible obstruction models  $(\hat f/G_{0},V)$ so that either
\begin{itemize}
  \item for all $\nu$ which are the restriction to $\hat f$ of sections in $U$, $F(\nu)$ is  contained inside $\mathcal O'_{i}$

or    
  \item  $F(\nu)$
  does not intersect $\mathcal O_{i}$ for any $\nu$ which is the restriction to $\hat f$ of some section in $U$.
  \end{itemize}
  
   Let $I$ indicate the set of indices $i$ so that the first option holds, so $F(\nu)$ is  contained inside $\mathcal O'_{i}$ for $\nu$ small enough.
    
 The main problem that must be overcome in the rest of this proof is that the simple perturbations $\mathfrak P_{i}$ are not parametrized by $\hat f$. We will extend $\hat f$ to a  family $\hat h$ which can be regarded as parametrizing the simple perturbations $\mathfrak P_{i}$ for all $i\in I$ and use the resulting unique solution $\tilde \nu$ to the corresponding perturbed $\dbar$ equation over $\hat h$ to construct the weighted branched section of $\hat f^{*}T_{vert}\hat{\ex B}\oplus V$ which is our `solution' with the required properties. This will involve reexamination of ideas that came up in the proof of Theorem \ref{core criteria}.

 Use the notation 
 \[s^{i}:\ex F(\hat f'_{i})\longrightarrow  \ex F(\hat f_{i}'^{+n_{i}})\]
for the map coming from the extra marked points on the core family $\hat f'_{i}$.

 The map $ev^{+n_{i}}(\hat f_{i}'):\ex F(\hat f_{i}'^{+n_{i}})\longrightarrow \M\times \sfp{\hat{\ex B}}{\ex G}{n_{i}}$ has the property that it is an equidimensional embedding in a neighborhood of the section $s^{i}$. There exists an open neighborhood $\mathcal O_{s^{i}}$ of the family of curves $(s^{i})^{*}\hat f_{i}'^{+n}$ so that  given any curve $f$ in $\mathcal O'_{i}$, if $ f^{+n_{i}}\rvert_{\mathcal O_{s^{i}}}$ indicates the restriction of the family $f^{+n_{i}}$ to $\mathcal O_{s^{i}}$, then $ev^{+n_{i}}(f)(\ex F(f^{+n_{i}}\rvert_{\mathcal O_{s^{i}}}))$ intersects $ev^{+n_{i}}(\hat f')(s^{i}(\ex F(\hat f'_{i})))$ transversely exactly $\abs {G_{i}}$ times, corresponding to the $\abs {G_{i}}$ maps from $\ex C(f)$ into $\ex C(\hat f'_{i})$.

Consider the family $\hat f^{+(n-1)}:\ex F(\hat f^{+n})\longrightarrow \sfp{\hat{\ex B}}{\ex G}{n}$. Use the notation 
$X^{+(n-1)}$ to denote the vector bundle over $\ex F(\hat f^{+n})$ which is the pullback under $\hat f^{+n-1}$ of the vertical tangent space of the family $\sfp{\hat{\ex B}}{\ex G}{n}\longrightarrow \ex G$. The $G_{0}$ action on $\hat f$ gives a $G_{0}$ action on $X^{+(n-1)}$.
Any section $\nu$ of $\hat f^{*}T_{vert}\hat{\ex B}$ corresponds in an obvious way to a section $\nu^{+(n-1)}$ of $X^{+(n-1)}$, and the map $F:\hat f^{*}T_{vert}\hat{\ex B}\longrightarrow \hat{\ex B} $ corresponds to a $G_{0}$ invariant $\C\infty 1$ map  
\[F^{+(n-1)}:X^{+(n-1)}\longrightarrow  \sfp{\hat{\ex B}}{\ex G}n\]
so that
\[F^{+(n-1)}(\nu^{+(n-1)})=\lrb{F(\nu)}^{+(n-1)}\]
Use the notation $\nu^{+(n_{i}-1)}\rvert_{\mathcal O_{s^{i}}}$ to denote the restriction of $\nu^{+(n_{i}-1)}$ to the subset $\ex F(F(\nu)^{+n_{i}}\rvert_{\mathcal O_{s^{i}}})\subset \ex F(\hat f^{+n_{i}})$ 

Define a map \[EV^{+n}:X^{+(n-1)}\longrightarrow  \M\times\sfp{\hat{\ex B}}{\ex G}n\] so 
that $EV^{+n}$ is equal to the natural map coming from the complex structure of curves in $\ex C(\hat f^{+n})\longrightarrow \ex F(\hat f^{+n})$ on the first component, and $F^{+(n-1)}$ on the second component. So
\[EV^{+n}( \nu^{+(n-1)}(\cdot))=ev^{+n}(F(\nu))(\cdot)\] 
The map $EV^{+n}$ is $\C\infty1$ and $G_{0}$ invariant. 

 Use the notation $\nu(g)^{+(n_{i}-1)}\rvert_{\mathcal O_{s^{i}}}$ for the restriction of  $\nu^{+(n_{i}-1)}\rvert_{\mathcal O_{s^{i}}}$ to the inverse image of a curve $g\in \hat f$. For any section $\nu$ small enough in $C^{1,\delta}$,
  the map $EV^{+n}$ restricted to $\nu(g)^{+(n_{i}-1)}\rvert_{\mathcal O_{s^{i}}}$  intersects $ev^{+n_{i}}(\hat f'_{i})(s^{i}(\ex F(\hat f'_{i})))$ transversely in exactly $\abs {G_{i}}$ points.  Denote by $S_{i}$ the subset of $X^{+(n_{i}-1)}$ which is the pullback of the image of the section $s^{i}$:
\[S_{i}:=(EV^{+n_{i}})^{-1}\lrb{ev^{+n_{i}}(\hat f'_{i})(s^{i}(\ex F(\hat f'_{i})))}\subset X^{+(n_{i}-1)}\]

$S_{i}$ is $G_{0}$ invariant, and  close to the zero section in $X^{+(n_{i}-1)}$, $S_{i}$ has regularity $\C\infty1$, and for sections $\nu$ small enough in $C^{1,\delta}$, $S_{i}$ is transverse to $\nu^{+(n_{i}-1)}\rvert_{\mathcal O_{s^{i}}}$ and  $\nu^{+(n_{i}-1)}\rvert_{\mathcal O_{s^{i}}}\cap S_{i}$ is a  $\abs {G_{i}}$-fold  multisection $\ex F(\hat f)\longrightarrow S_{i}\subset X^{+(n_{i}-1)}$.

Use the notation $\tilde X^{+(n-1)}$ for the pullback along the map $\ex C(\hat f^{+n})\longrightarrow \ex F(\hat f^{+n})$ of the vector bundle $X^{+(n-1)}$, $\tilde S_{i}$ for the inverse image of $S_{i}$ in $\tilde X^{+(n_{i}-1)}$, and 
$\tilde \nu^{+(n_{i}-1)}_{\rvert_{\mathcal O_{s^{i}}}}$ for the pullback of $\nu^{+(n_{i}-1)}\rvert_{\mathcal O_{s^{i}}}\subset X^{+(n_{i}-1)}$ to a section of $\tilde X^{+(n_{i}-1)}$. Define a map 
\[\tilde {EV}^{+n}:\tilde X^{+{n-1}}\longrightarrow \M\times \sfp{\hat {\ex B}}{\ex G}n\]
so that recalling the notation $\tilde{ev}^{+n}$ from (\ref{tilde ev}) on page \pageref{tilde ev} in the proof of Theorem \ref{core criteria},
\[\tilde {EV}^{+n}(\tilde\nu^{+(n-1)})=\tilde{ev}^{+n}(F(\nu))\]
 Note that
 \[\tilde S_{i}=\lrb{\tilde{EV}^{+n_{i}}}^{-1}(\tilde{ev}^{+n_{i}}(\hat f'_{i})(\ex C((s^{i})^{*}\hat f_{i}'^{+n_{i}})))\]
Of course, everything in the above construction is $G_{0}$ invariant. The fact that $ev^{+n_{i}}(\hat f'_{i})$ and $\tilde {ev}^{+n_{i}}(\hat f'_{i})$ are embeddings in a neighborhood of $s^{i}$ imply that there are natural maps  
 \begin{equation}\label{S map}\begin{array}{ccccc}\tilde S_{i}&\longrightarrow &\ex C\lrb{(s^{i})^{*}\hat f_{i}'^{+n}}&\longrightarrow &\ex C(\hat f_{i}')
 \\\downarrow &&\downarrow &&\downarrow
 \\ S_{i}&\longrightarrow &\ex F\lrb{(s^{i})^{*}\hat f_{i}'^{+n}}&\longrightarrow &\ex F(\hat f_{i}')
 \end{array}\end{equation}
 so that the map on the left is an isomorphism on each fiber, and the map on the right is the map that forgets the extra marked points. We can define a family of curves $\check S_{i}\longrightarrow S_{i}$ by forgetting the extra marked points in the family  $\tilde S_{i}\longrightarrow S_{i}$ - so $\check S_{i}$ is equal to the domain of the pullback of the family $\hat f_{i}'$ over the map $S_{i}\longrightarrow \ex F(\hat f_{i}')$. The above map of families of curves (\ref{S map}) factors differently into the following diagram 

 \begin{equation}\label{S map2}\begin{array}{ccccc}\tilde S_{i}&\longrightarrow &\check S_{i}&\xrightarrow{\Phi_{i}} &\ex C(\hat f_{i}')
 \\\downarrow &&\downarrow &&\downarrow
 \\ S_{i}&\longrightarrow &S_{i}&\longrightarrow &\ex F(\hat f_{i}')
 \end{array}\end{equation} 
 
Note that $\tilde S_{i}$ and $\Phi_{i}$ are $G_{0}$ invariant.  As discussed in the proof of Theorem \ref{core criteria}, the right hand map $\Phi_{i}$ above determines the maps $\Phi_{F(\nu)}$ from the definition of the core family $\hat f_{i}'/G_{i}$ in the following sense: For $\nu$ small enough, the intersection of $ \nu^{+(n_{i}-1)}\rvert_{\mathcal O_{s^{i}}}$
 with $ S_{i}$ is transverse, and is a $\abs {G_{i}}$-fold cover  of $\ex F(\hat f)$, which lifts to a $\abs {G_{i}}$-fold cover of $\ex C(\hat f)$ which is a subset of $\check S$. Then $\Phi_{i}$ gives a map of our $\abs {G_{i}}$-fold cover of $\ex C(\hat f)$ into $\ex C(\hat f'_{i})$, which corresponds to the map $\Phi_{F(\nu)}:\ex C(\hat f)\longrightarrow \ex C(\hat f'_{i})/G_{i}$.
 
%
%
%

Denote by $X^{+}$ the fiber product of $X^{+(n_{i}-1)}$ over $\ex F(\hat f)$ for all $i\in I$,
and  denote by $\ex F^{+}$ the fiber product of $\ex F(\hat f^{+n_{i}})$ over $\ex F(\hat f)$ for all $i\in I$;
    so $X^{+}$ is a vector bundle over $\ex F^{+}$ with a natural $G_{0}$ action. 
        A $\C\infty1$ section $\nu$ of $\hat f^{*}T_{vert}\hat {\ex B}$ corresponds in the obvious way to a $\C\infty1$ section $\nu^{+}$ of $X^{+}$ which is equal to $\nu^{+(n_{i}-1)}$ on each $X^{+(n_{i}-1)}$ factor. Similarly, denote by $\nu^{+}\rvert_{\mathcal O_{s^{+}}}$ the open subset of $\nu^{+}$ inside $\nu^{+(n_{i}-1)}\rvert_{\mathcal O_{s^{i}}}$ on each $X^{+(n_{i}-1)}$ factor.
          Denote by $ S^{+}\subset X^{+}$ the  subset corresponding to all $S_{i}\subset X^{+(n_{i}-1)}$  restricted to a neighborhood of the zero section small enough that $S^{+}$ is $\C\infty1$.
We can choose $S^{+}$ small enough so that pulling back $(\hat f,V)$ over the map $S^{+}\longrightarrow \ex F(\hat f)$ gives an extendible $G_{0}$ invariant pre obstruction model $(\hat h,V)$. Note that $\ex C(\hat h)$ is some open subset of the fiber product of $\check S_{i}$ over $\ex C(\hat f)$ for all $i\in I$, so the maps $\Phi_{i}$ from (\ref{S map2}) induce $G_{0}$ invariant maps
\[\begin{array}{ccc}\ex C(\hat h)&\xrightarrow{\Phi_{i}}& \ex C(\hat f'_{i})
\\ \downarrow & & \downarrow
\\\ex F(\hat h)&\longrightarrow &\ex F(\hat f'_{i})\end{array}\]
Pulling a simple perturbation $\mathfrak P_{i}$ parametrized by $\hat f_{i}'$ back over the map $\Phi_{i} $ gives a $G_{0}$ invariant simple perturbation $\Phi_{i}^{*}\mathfrak P_{i}$ parametrized by $\hat h$. Use the notation
\[\mathfrak P:=\sum_{i\in I}\Phi_{i}^{*}\mathfrak P_{i}\]
   If $\nu$ is any small enough section of $\hat f^{*}T_{vert}\hat{\ex B}$, then the multi perturbation $\prod_{i\in I}F(\nu)^{*}\mathfrak P_{i}$ defined as in example \ref{proper pullback} on page \pageref{proper pullback} can be constructed as follows: If $\nu$ is small enough, then $\nu^{+}\rvert_{\mathcal O_{s^{+}}}$ is transverse to $S^{+}$, and the intersection of $\nu^{+}\rvert_{\mathcal O_{s^{+}}}$ with $S^{+}$ is a $n$-fold cover of $\ex F(\hat f)$ in $\ex F(\hat h)$ which lifts to an $n$-fold cover  of $\ex C(\hat f)$ inside $\ex C(\hat h)$ (where $n=\prod_{i\in I}\abs{ G_{i}}$). Together these give the domain for a family of curves $F(\nu)'$ which is a  
 $n$-fold multiple cover of $F(\nu)$.  Restricting $\mathfrak P$ to $F(\nu)'$ then gives a section of $\Y{F(\nu)'}$, which corresponds to a $n$-fold multi section of $\Y{F(\nu)}$. Locally, giving each of these $n$ sections a weight $1/n$ gives a weighted branched section of $\Y{F(\nu)}$ with total weight $1$ which is equal to the multi perturbation   $\prod_{i\in I}F(\nu)^{*}\mathfrak P_{i}$ defined as in example \ref{proper pullback}.

 As $(\hat f,V)$ comes from an obstruction model, $D\dbar(f)$ is injective and complementary to $V(f)$ for all $f$ in $\hat f$ and therefore all $f$ in $\hat h$, so Theorem \ref{regularity theorem} applies to $(\hat h,V)$ and implies that there is some neighborhood of $0$ in the space of simple perturbations  parametrized by $\hat h$ so that for such any  $\mathfrak P$  in this neighborhood, there is a unique small $\C\infty1$ section $\tilde \nu$ of $\hat h^{*}T_{vert}\hat{\ex B}$ vanishing at the relevant marked points so that $(\dbar -\mathfrak P)(\tilde \nu)\in V$. The fact that $(\hat f,V)$ is part of an obstruction model for $\mathcal O$ implies the following uniqueness property for $\tilde \nu$ if $\mathfrak P$ is small enough: Given any curve $h$ in $\hat h$ and section $\psi$  in $X^{\infty,\underline 1}(h)$ so that $F(\psi)$ is in $\mathcal O$, then $(\dbar -\mathfrak P)(\psi)\in V$ if and only if $\psi$ is the restriction to $h$ of $\tilde \nu$.  
 
 Denote by $\tilde X^{+}$ the pullback of $X^{+}$ over the map $\ex F(\hat h)\longrightarrow \ex F(\hat f)$, and denote by $\tilde S^{+}$ the pullback of $S^{+}$. Of course there is a natural $G_{0}$ action on $\tilde S^{+}$ and $\tilde X^{+}$ so that the following commutative diagram is $G_{0}$ equivariant.
 \[ \begin{array}{ccc}\tilde S^{+}\subset \tilde X^{+}&\longrightarrow& S^{+}\subset  X^{+}
 \\\downarrow &&\downarrow
 \\ \ex F(\hat h)&\longrightarrow& \ex F(\hat f)\end{array}\]

 This $\tilde S^{+}$ comes with two maps into $S^{+}$: one  the restriction of the map $\tilde X^{+}\longrightarrow X^{+}$, and one  the restriction of the map $\tilde X^{+}\longrightarrow \ex F(\hat h)=S^{+}$. Denote by $S^{+}_{\Delta}$ the subset of $\tilde S^{+}$ on which these two maps agree. Because these two above maps agree when composed with the relevant maps to $\ex F(\hat f)$, $\tilde S^{+}$ can be regarded as the fiber product of $ S^{+}$ with itself over $\ex F(\hat f)$ and $S^{+}_{\Delta}$ is the diagonal in this fiber product $\tilde S^{+}$. Therefore, $S^{+}_{\Delta}$ is $\C\infty1$ and the  map  $S^{+}_{\Delta}\longrightarrow S^{+}$ is an isomorphism.  A section $\tilde \nu$ of $\hat h^{*}T_{vert}\hat{\ex B}$ defines a section $\tilde \nu^{+}$ of the vector bundle $\tilde X^{+}$ so that if $\tilde \nu$ is the pullback over the map $\ex C(\hat h)\longrightarrow \ex F(\hat f)$ of some section $\nu$ of $\hat f^{*}T_{vert}\hat {\ex B}$, then $\tilde \nu^{+}$ is the pullback of $ \nu^{+}$. We can define $\tilde \nu^{+}\rvert_{\mathcal O_{s^{+}}}$ similarly to the definition of $\nu^{+}\rvert_{\mathcal O_{s^{+}}}$.
  
  As $\nu^{+}\rvert_{\mathcal O_{s^{+}}}$ is transverse to $S^{+}$ for $\nu$ small enough and $\nu^{+}\rvert_{\mathcal O_{s^{+}}}\cap S^{+}$ gives a $n$-fold cover of $\ex F(\hat f)$, $\tilde \nu^{+}\rvert_{\mathcal O_{s^{+}}}$ is transverse to $S^{+}_{\Delta}$ for $\tilde \nu$ small enough, and $\tilde \nu^{+}\rvert_{\mathcal O_{s^{+}}}\cap S^{+}_{\Delta}$ also defines a   $n$-fold cover of $\ex F(\hat f)$ with regularity $\C\infty1$. To see this, suppose that $\tilde \nu^{+}$ is the pullback of some $ \nu^{+}$. Then $\tilde \nu^{+}$ is transverse to $\tilde S^{+}$, and $\tilde \nu^{+}\cap \tilde S^{+}$ is an $n$ fold cover of $\ex F(\hat h)$ which is a pullback of $\nu^{+}\cap S^{+}$ over the map $\ex F(\hat h)\longrightarrow \ex F(\hat f)$. These $n$ sections of $\tilde S^{+}\longrightarrow \ex F(\hat h)=S^{+}$ are constant on fibers of the map $\ex F(\hat h)=S^{+}\longrightarrow \ex F(f)$, and are therefore transverse to the diagonal section $S^{+}_{\Delta}$, and when intersected with $S^{+}_{\Delta}$ give an $n$-fold section of $S^{+}_{\Delta}\longrightarrow \ex F(\hat f)$. This transversality and the fact that $\tilde \nu^{+}\rvert_{\mathcal O_{s^{+}}}\cap S^{+}_{\Delta}$ defines a   $n$-fold cover of $\ex F(\hat f)$ with regularity $\C\infty1$ is stable under perturbations of $\tilde \nu^{+}$, so it remains true for small $\tilde \nu$ which aren't the pullback of some $\nu$.

   We may consider this multiple cover of $\ex F(\hat f)$ as being a multi section  $\ex F'$ of $S^{+}=\ex F(\hat h)\longrightarrow \ex F(\hat f)$, which lifts to a $G_{0}$-equivariant  multi section of $\ex C(\hat h)\longrightarrow\ex C(\hat f)$. Restricting our $G_{0}$ invariant $\tilde \nu$ to this multi section gives locally $n$ sections $\nu_{l}$ of $\hat f^{*}T_{vert}\hat{\ex B}$ with regularity $\C\infty1$. We may similarly pullback the sections $(\tilde\nu,\dbar\nu-\mathfrak P)$ to give locally $n$ sections $(\nu_{l},\dbar'\nu_{l})$. Then
  \begin{equation}\label{nu l}(\nu,\dbar'\nu)=\sum_{l=1}^{n}\frac 1{n}t^{(\nu_{l},\dbar'\nu_{l})}\end{equation}
   is the $G_{0}$ invariant  weighed branched solution which is our `solution mod $V$'. We shall now show that this weighted branched section has the required properties if $\{\mathfrak P_{i}\}$ is small enough. Note first that close by simple perturbations $\{\mathfrak P_{i}\}$ give close by solutions $(\nu_{l},\dbar'\nu_{l})$. Also note that if we have a $\C\infty1$ family of simple perturbations $\{\mathfrak P_{i,t}\}$, Theorem \ref{regularity theorem} implies that the corresponding family of solutions $\tilde\nu_{t}$ to $(\dbar-\mathfrak P_{t})\tilde\nu_{t}\in V$ is a $\C\infty1$ family, so the corresponding weighted branched sections $(\nu_{t},\dbar'\nu_{t})$ form a $\C\infty1$ family.  
   
   If we choose $\{\mathfrak P_{i}\}$ small enough, then $F(\nu_{l})$ will not intersect $\mathcal O_{i}$ for any $i\notin I$. Therefore, the multi perturbation under study is given by 
   \begin{equation}\label{Pjl}\prod_{i}\nu_{l}^{*}\mathfrak P_{i}=\sum_{j=1}^{n}\frac 1nt^{\mathfrak P_{j,l}}\end{equation}
  where $\mathfrak P_{j,l}$ is constructed as follows: $\nu_{l}^{+}\rvert_{\mathcal O_{s^{+}}}\cap S^{+}$ is a $n$-fold cover of the open subset of $\ex F(\hat f)$ where $\nu_{l}$ is defined. By working locally, this $n$-fold cover can be thought of as $n$ local sections of $S^{+}=\ex F(\hat h)\longrightarrow \ex F(\hat f)$, which lift to $n$ local sections of $\ex C(\hat h)\longrightarrow \ex F(\hat f)$. The restriction of $\mathfrak P$ to these $n$ local sections gives the $n$ sections $\mathfrak P_{j,l}$ of $\Y{F(\nu_{l})}$ in the formula (\ref{Pjl}) above. As one of these sections of $\ex F(\hat h)\longrightarrow\ex F(\hat f)$ coincides with the multi section $\ex F'$ mentioned in the  paragraph preceding equation (\ref{nu l}) obtained using the solution $\tilde \nu$ to the equation $(\dbar-\mathfrak P)\tilde \nu\in V$, one of the sections $\mathfrak P_{l,l}$ of $\Y{F(\nu_{l})}$ has the property that $\dbar F(\nu_{l})-\mathfrak P_{l,l}=\dbar'\nu_{l}$. 

  Suppose that $f'$ is in the region of $\hat f$ where these $\nu_{l}$ in formula (\ref{nu l}) are defined and $\psi$ in $X^{\infty,\underline 1}(f')$ is small enough  so that $f=F(\psi)$ is in $\mathcal O$. If the simple perturbations $\mathfrak P_{i}$ are chosen small enough, the fact that $(\hat f,V)$ is an obstruction model will imply that if $ \prod_{i} f^{*}\mathfrak P_{i}=wt^{\mathfrak Q}+\dotsc$ where $w>0$ and $(\dbar f-\mathfrak Q)\in V$, then $\psi$ must be small - choose $\{\mathfrak P_{i}\}$ small enough that such $f$ must have $\psi^{+}\rvert_{\mathcal O_{s^{+}}}$ intersecting $S^{+}$ transversely $n$ times and $f$ is not in $\mathcal O_{j}$ for all $j\notin I$. Then $ \prod_{i} f^{*}\mathfrak P_{i}=\sum_{l=1}^{n}\frac 1nt^{\mathfrak Q_{l}}$ where the $n$ sections $\mathfrak Q_{l}$ of $\Y{f}$ are obtained as follows: The $n$  points of $\psi^{+}\rvert_{\mathcal O_{s^{+}}}\cap S^{+}$ correspond to $n$ maps of $\ex C(f)$ into $\ex C(\hat h)$ - the $n$ sections $\mathfrak Q_{l}$ are given by pulling back the simple perturbation $\mathfrak P$ over these maps. Then $\dbar f-\mathfrak Q_{l}\in V(f)$  if and only if $\psi$ is equal to the pullback under the relevant map of the solution $\tilde \nu$ to $(\dbar -\mathfrak P)\tilde \nu\in V$. Therefore, if $\{\mathfrak P_{i}\}$ is small enough, the number of $\mathfrak Q_{l}$ so that $\dbar f-\mathfrak Q_{l}\in V$ is equal to the number of $\nu_{l}$ from formula (\ref{nu l}) so that $\psi=\nu_{l}(f')$. 
  
  Similarly, if $\nu'$ is locally a section of $\hat f^{*}T_{vert}\hat{\ex B}$ vanishing on the relevant marked points so that $F(\nu')\in\mathcal O$ and $\prod_{i} F(\nu')^{*}\mathfrak P_{i}=wt^{\mathfrak Q}+\dotsc$ where $w>0$ and $(\dbar F(\nu')-\mathfrak Q)\in V$, then so long as $\{\mathfrak P_{i}\}$ is small enough, $\nu'^{+}\rvert_{\mathcal O_{s^{+}}}\cap S^{+}$ is locally a $n$-fold cover of $\ex F(\hat f)$  corresponding to $n$ sections of $\ex F(\hat h)\longrightarrow \ex F(\hat f)$ which lift locally to $n$ sections of $\ex C(\hat h)\longrightarrow \ex C(\hat f)$. Then $\mathfrak Q$ must locally correspond to the pullback of $\mathfrak P$ under one of these local maps $\ex C(\hat f)\longrightarrow \ex C(\hat h)$, and $\nu'$ must locally be the pullback of the solution $\tilde \nu$ to $(\dbar -\mathfrak P)\tilde \nu\in V$. It follows that $(\nu',\dbar\nu'-\mathfrak Q)$ must coincide locally with one of these $(\nu_{l},\dbar'\nu_{l})$ from formula (\ref{nu l}), and the weighted branched section locally equal to $\sum_{l=1}^{n}\frac 1nt^{(\nu_{l},\dbar'\nu_{l})}$ is the unique weighted branched section with the  required properties. 

Suppose that $\hat g$ is a family of curves in the subset of $\mathcal O$ projecting to the region where our $\nu_{l}$ are defined so that $\prod_{i}\mathfrak P_{i}\hat g=wt^{\dbar \hat g}+\dotsc$ and $w>(k-1)/n$. Then  using that $\hat f/G_{0}$ is a core family, after locally choosing one of the $\abs G_{0}$ maps from $\ex C(\hat g)$ to $\ex C(\hat f)$ we may pull back $X^{+}$  and $S^{+}$ to be bundles over $\ex F(\hat g)$. The corresponding bundles $X^{+}(\hat g)$ and $S^{+}(\hat g)$ may also be constructed in the same way as the original bundles using the induced trivialization on $\hat g$ from $\hat f$.  The section $\psi$ vanishing at the correct marked points so that $\hat g=F(\psi)$ corresponds to a section $\psi^{+}$ of this pulled back $X^{+}$ which is transverse to the pulled back $S^{+}$ and intersects this pulled back $S^{+}$ in an $n$-fold cover of $\ex F(\hat g)$. This $n$-fold cover of $\ex F(\hat g)$ comes with a map to $S^{+}$, corresponding to a map to $\ex F(\hat h)$ which lifts to a fiberwise holomorphic map of a $n$-fold cover of $\ex C(\hat g)$ to $\ex C(\hat h)$ so that $\prod_{i}\hat g^{*}\mathfrak P_{i}$ is determined by pulling back $\mathfrak P$ over this map, then giving the simple perturbation from each branch of the cover a weight $1/n$ and summing the result. As $w>(k-1)/n$,  locally at least $k$ of these simple perturbations must be $\dbar g$, and $\psi$ must be the pullback under each of the corresponding maps of the solution $\tilde \nu$ to $(\dbar-\mathfrak P)\tilde\nu\in V$, and its image must be contained in the subset where $\dbar'\tilde \nu=0$. It follows that around each curve in  $\hat g$, there is a map of a neighborhood into at least $k$ of the $ F(\nu_{l})$  with image contained in the subset where $\dbar'\nu_{l}=0$, and with the map $\ex C(\hat g)\longrightarrow \ex C( F(\nu_{l}))$ corresponding to our local choice of lift of the map $\ex C(\hat g)\longrightarrow \ex C(\hat f)/G_{0}$ coming from the fact that $\hat f/G_{0}$ is a core family.  Without a choice resolving this $G_{0}$-fold ambiguity, this corresponds to there being at least $k$ maps of $\hat g$ into $(\coprod_{l}F(\nu_{l}))/G_{0}$.

  \stop

\

\subsection{Construction of virtual moduli space}\label{virtual class}

\

\

\label{final theorem proof}
We now begin a construction of a virtual moduli space for the moduli stack of holomorphic curves which is an  oriented $\C\infty 1$ weighed branched substack of $\Msw$.  This will include the  proof of Theorem \ref{final theorem} stated on page \pageref{final theorem}. 

\

Any stable holomorphic curve with at least one smooth component in a basic family of targets $\hat {\ex B}\longrightarrow \ex G$ is contained in some $\C\infty 1$ obstruction model $(\hat f/G,V)$, and any obstruction model covers an  open neighborhood in the moduli stack of holomorphic curves. If $\ex G$ is compact and Gromov Compactness holds for $\hat{\ex B}\longrightarrow \ex G$ in the sense of Definition \ref{gcompactness}, then the substack of holomorphic curves in $\Msw_{g,[\gamma],\beta}$  may be covered by a finite number of extendible obstruction models.

The rough idea of how the virtual class is constructed is that the $\dbar$ equation is perturbed in some neighborhood of
the holomorphic curves in $\Msw$ to achieve `transversality', and a $\C\infty 1$ solution set. `Transversality' is easy to achieve locally with a simple perturbation parametrized by an obstruction model. For such a simple perturbation to be defined independent of coordinate choices, it must be viewed as a multi-perturbation in the sense of example \ref{proper pullback} on page \pageref{proper pullback}. One problem is that for a simple perturbation to give a $\C\infty1$ multi-perturbation restricted to a particular family, that family must meet the domain of definition of the simple perturbation properly in the sense of Definition \ref{proper meeting} on page \pageref{proper meeting}. 

  Restrict our obstruction models to satisfy the requirements of Theorem \ref{Multi solution} as follows. Each of the obstruction models we start off with has an extension $(\hat f'_{i}/G_{i},V_{i})$ on $\mathcal O'_{i}\subset\Msw$. This substack $\mathcal O_{i}'$ can be viewed as corresponding to a neighborhood of $0$ in $X^{\infty,\underline 1}(\hat f_{i})$. We may assume that this neighborhood is convex, and denote by $c\mathcal O'_{i}$ the open substack corresponding to the above neighborhood multiplied by $c$ (i.e. sections $\nu$ so that $\frac 1c\nu$ is in the above neighborhood). The fact that $(\hat f'_{i}/G_{i},V_{i})$ is an obstruction model for $\mathcal O'_{i}$ implies that  any holomorphic curves that are in $\mathcal O'_{i}$ are actually contained in the the family $\hat f'_{i}$, so all holomorphic curves in $\mathcal O'_{i}$ are contained inside $\frac 12\mathcal O'_{i}$. We may assume that  $(\hat f'_{i}/G_{i},V)$ itself is an extendible obstruction model, and that the closure $\overline{(\frac 34\mathcal O'_{i}-\frac 12\mathcal O'_{i})}$  of $(\frac 34\mathcal O'_{i}-\frac 12\mathcal O'_{i})$  contains no holomorphic curves.
  Define an open neighborhood of the part of the stack of holomorphic curves under study by 
 \[\mathcal O:=\bigcup_{i}\frac 12\mathcal O'_{i}-\bigcup_{i}\overline{(\frac 34\mathcal O'_{i}-\frac 12\mathcal O'_{i})}\]
 
 The substack $\mathcal O$ meets $\frac 12 \mathcal O'_{i}$ with the core family $\hat f'_{i}/G_{i}$ properly in the sense of Definition \ref{proper meeting}. We may restrict our original obstruction model family  to a $G_{i}$ invariant sub family $\hat f_{i}$ which is  compactly contained inside $\hat f_{i}'\rvert_{\mathcal O}$ so that $\hat f_{i}$  still contains the same set of holomorphic curves as our original obstruction model.  Use the notation $\mathcal O_{i}$ to refer to the restriction of  $\frac 12 \mathcal O'_{i}$ to the subset with core given by this new family $\hat f_{i}/G_{i}$. $\mathcal O$ meets all these new $\mathcal O_{i}$ properly, so item \ref{ft1} from Theorem \ref{final theorem} holds, and any compactly supported $\C\infty1$ simple perturbation parametrized by $\hat f_{i}$ defines a $\C\infty1$ multi-perturbation on $\mathcal O$ as in Example \ref{proper pullback}.

 Theorem \ref{Multi solution} holds for this collection of obstruction models when we use $(\hat f_{i}'\rvert_{\mathcal O}/G_{i},V_{i})$ on the corresponding restriction of $\frac 12\mathcal O_{i}'$ for the extensions of our obstruction models $(\hat f_{i}/G_{i},V_{i})$ on $\mathcal O_{i}$. It follows that item \ref{ft2} from Theorem \ref{final theorem} holds. 

In particular, for a collection of compactly supported $\C\infty1$ simple perturbations $\mathfrak P_{i}$ parametrized by $\hat f_{i}$,  let $\theta$ denote the multi perturbation on $\mathcal O$ so that 
\[\theta(\hat f):=\prod_{i}\hat f^{*}\mathfrak P_{i}\]
where $\hat f^{*}\mathfrak P_{i}$ is as in example \ref{proper pullback}.
Then for some convex $\C\infty1$ neighborhood $U$ of $0$ in the space of collections of compactly supported $\C\infty1$ simple perturbations $\{\mathfrak P_{i}\}$, for each of our obstruction models $(\hat f_{i}/G,V_{i})$, there exists a unique $\C\infty1$ weighted branched section $(\nu,\dbar'\nu)$ of $f_{i}^{*}T_{vert}\hat{\ex B}\oplus V_{i}$  so that 
  locally on $\ex F(\hat f_{i})$, 
\begin{equation}\label{nu_{k}}(\nu,\dbar'\nu)=\sum_{k=1}^{n}\frac 1n t^{(\nu_{k},\dbar'\nu_{k})}\end{equation}
where $\nu_{k}$ vanishes on marked points, and  $F(\nu_{k})$ is a family of curves in $\mathcal O_{i}\cap\mathcal O$ so that
\begin{equation}\label{theta nu}\theta(F(\nu_{k}))=\sum_{j=1}^{n}\frac 1nt^{\mathfrak P_{k,j}}\end{equation}
where \[\dbar F(\nu_{k})-\mathfrak P_{k,k}=\dbar'\nu_{k}\]
 Moreover, given any curve $f\in \mathcal O_{i}\cap\mathcal O$, if $\theta(f)=\sum_{l} w_{l}t^{\mathfrak Q_{l}}$, then the sum of the weights $w_{l}$ so that $\dbar f-\mathfrak Q_{l}$ is in $V_{i}$ is equal to $\frac 1n$ times the number of the above $\nu_{k}$ so that $f$ is contained in the family $F(\nu_{k})$. 

Say that $\theta$ is transverse to $\dbar$ on a sub family $C\subset \hat f_{i}$ if each of the sections $\dbar'\nu_{k}$ of $V_{i}$ are transverse to the zero section of $V_{i}$ on $C$. If $\theta'$ indicates the multi perturbation corresponding to a collection $\{\mathfrak P_{i}'\}$ of simple perturbations close in $\C\infty1$ to $\{\mathfrak P_{i}\}$, then Theorem \ref{Multi solution} implies that the sections $\nu_{k}'$ corresponding to $\nu_{k}$ will be $\C\infty1$ close to $\nu_{k}$, and the corresponding sections $\dbar'\nu_{k}'$ of $V_{i}$ will also be $\C\infty1$ close to the original sections. It follows that the subset of $U$ consisting of collections of perturbations $\{\mathfrak P_{i}\}$ so that $\dbar$ is transverse to $\theta$ on any particular  compact sub family $C\subset \hat f_{i}$ is open in the $\C\infty1$ topology. If we choose $\mathfrak P'_{i}-\mathfrak P_{i}$ to consist of a simple perturbation which take values in $V_{i}$ and leave the other simple perturbations unchanged, then the sections $\nu'_{k}$ from equation (\ref{nu_{k}}) will be equal to the sections $\nu_{k}$, and $\dbar\nu'_{k}$ will be modified by $(\mathfrak P_{i}'-\mathfrak P_{i})$ acted on by some element of the group $G_{i}$. It follows that the subset of $U$ so that $\theta$ is transverse to $\dbar$ on a given compact subset $C$ is dense and open in $\C\infty1$, so item \ref{ft3} from Theorem \ref{final theorem} holds. Similarly, to prove item \ref{ft 9} later on, note that given a map of a compact exploded manifold $\ex G'$ into $\ex G$, the set of perturbations for which the map to $\ex G$ from  intersection of $\dbar'\nu_{k}$ with the zero section restricted to $C$  is transverse to $\ex G'\longrightarrow \ex G$ is open and dense.
  
Say that $\theta$ is fixed point free on the sub family $C\subset \hat f_{i}$ if none of the curves in $F(\nu_{k})$ over $C$ have smooth parts with a non trivial automorphism group.  If $C$ is  compact, the set of such curves  within some $F(\nu_{k})$ over $C$  is  compact. If $\theta'$ is the multi perturbation corresponding to a close collection of simple perturbations, then Theorem \ref{Multi solution} implies that  the set of corresponding curves in $F(\nu_{k}')$ over $C$ are close  to the original set, so if $\theta$ is fixed point free on $C$, $\theta'$ is fixed point free on $C$ if the new simple perturbations are chosen close enough in $\C\infty1$. If $C$ is  compact, it is covered by a finite number of  compact subsets on which the sections $\nu_{k}$ from equation (\ref{nu_{k}}) are defined. Theorem \ref{Multi solution} implies that for any close by modification $\nu_{k}'$ of $\nu_{k}$, there exists a small modification of $\mathfrak P_{i}$ to $\mathfrak P_{i}'$ so that $\nu_{k}'$ is the solution corresponding to the modified multi perturbation $\theta'$. If the relative dimension of $\hat{\ex B}\longrightarrow \ex G$ is greater that $0$, $\nu_{k}'$ may be chosen so that $F(\nu_{k}')$ contains no curves who's smooth parts have non trivial automorphism group. We may proceed with a finite number of modifications to make each $\nu_{k}'$ fixed point free so that each modification is small enough that it doesn't change the fact that the sections previously concentrated on are fixed point free. Therefore, item \ref{ft4} from Theorem \ref{final theorem} holds and the subset of our space of perturbations $U$ so that $\theta$ is fixed  point free is open and dense. (This is of course not the case when the relative dimension of $\hat{\ex B}\longrightarrow \ex G$ is zero - in other words, $\hat {\ex B}$ is just a family of points; then all our perturbations are trivial, and the nature of the moduli space of holomorphic curves is easily deduced from Theorem \ref{moduli}.)
 Similarly, if the relative dimension of $\hat {\ex B}\longrightarrow\ex G$ is not zero, then given any two perturbations $\theta$, $\theta'$ in the set under consideration which are fixed point free on $C$, a generic family $\C\infty1$ family  of perturbations $\theta_{t}$ in the set under consideration joining $\theta$ to $\theta'$ is fixed point free on $C$.   

As $\ex G$ is compact, and Gromov compactness holds for $\hat{\ex B}\longrightarrow \ex G$,  we may cover  the moduli stack of holomorphic curves in $\Msw_{g,[\gamma],\beta}(\hat {\ex B})$ by $G_{i}$ invariant  open sub families $\hat f^{\circ}_{i}\subset \hat f_{i}$ with closures $C_{i}\subset \hat f_{i}$ which are  compact. Let $\mathcal O^{\circ}_{i}$ denote the subset of $\mathcal O_{i}$ with core $\hat f_{i}^{\circ}/G_{i}$. Let $\mathcal O^{\circ}$  denote the  union of $\mathcal O_{i}^{\circ}$. This is some open substack which contains the  holomorphic curves in $\Msw_{g,[\gamma],\beta}(\hat{\ex B})$. Theorem \ref{Multi solution} implies that if the collection of simple perturbations $\{\mathfrak P_{i}\}$ is small enough in $\C\infty1$, then if any curve $f$ in $\mathcal O$ over $\ex G'$ satisfies $\theta(f)=wt^{\dbar f}+\dotsc$ where $w>0$, then $f$ is in $\mathcal O^{\circ}$. So if we choose our open set of perturbations $U$  small enough, item \ref{ft5} from Theorem \ref{final theorem}  holds.

Say that $\theta$ is fixed point free and transverse to $\dbar$ if $\theta$ is fixed point free and transverse to $\dbar$ on all of the above subfamilies $C_{i}$. The preceding argument implies that $\theta$ is fixed point free and transverse to $\dbar$ in this sense for an open dense subset of perturbations in $U$. 

Suppose that $\theta$ is fixed point free and transverse  to $\dbar$. Then define a weighted branched substack $\Mod_{g,[\gamma],\beta}$ of $\mathcal O^{\circ}$ as follows: Given a curve $f$ in the region where equation \ref{nu_{k}} holds and section $\psi$ in $X^{\infty,\underline 1}(f)$ so that $F(\psi)\in\mathcal O^{\circ}$, on some open neighborhood of of $F(\psi)$ in $\mathcal O^{\circ}$, 

\[\Mod_{g,[\gamma],\beta}:= \sum_{k=1}^{n}\frac 1nt^{\hat g_{k}}\]

where  if $\nu_{k}(f)=\psi$ and $\dbar'\nu_{k}(f)=0$, then $\hat g_{k}$ is the family which is the subset of  $F(\nu_{k})$ given by the intersection of $\dbar'\nu_{k}$ with the zero section, restricted to some neighborhood of $F(\psi)$, and if $\nu_{k}(f)\neq\psi$ or $\dbar'\nu_{k}\neq 0$, then $\hat g_{k}$ is the empty substack. Note that as $\theta$ is fixed point free, and $\ex f/G$ is a core family,  restricting to small enough neighborhoods in the family $\hat f_{i}$, the family $F(\nu_{k})$ and therefore the family $\hat g_{k}$ is a substack of the moduli stack of curves. In other words, given any $\C\infty1$ family of curves $\hat h$ consisting of curves in $\hat g_{k}$, there exists a unique $\C\infty1$ map $\hat h\longrightarrow \hat g$.

The weighted branched substack $\Mod_{g,[\gamma],\beta}$ has the following three properties which make it well defined:
\begin{enumerate}
\item\label{modp1} If around any curve $f$ in $\mathcal O^{\circ}$, $\Mod_{g,[\gamma],\beta}$ is equal to $\sum_{k}w_{k}t^{\hat g_{k}}$ and $\theta(f)=\sum_{j}w'_{j}t^{\mathfrak Q_{j}}$, then the sum of the weights $w_{k}$ so that $f$ is in $\hat g_{k}$ is equal to the sum of the weights $w'_{j}$ so that $\dbar f=\mathfrak Q_{j}$.
\item \label{modp2}If locally $\Mod_{g,[\gamma],\beta}$ is equal to $wt^{\hat g}+\dotsc$, then $\theta(\hat g)=wt^{\dbar \hat g}+\dotsc$.
\item \label{modp3}If a family $\hat g$ in $\mathcal O^{\circ}$ containing a curve $f$ satisfies 
\[\theta(\hat g)=w t^{\dbar \hat g}\] 
then if $\mathcal M_{g,[\gamma],\beta}$ is equal to $\sum_{k}w_{k}t^{\hat g_{k}}$ on a neighborhood of $f$, then the sum of the $w_{k}$ so that there is a map from some neighborhood of $f$ in $\hat g$ to $\hat g_{k}$ is at least $w$. 
\end{enumerate}
Each of these properties follow from Theorem \ref{Multi solution}.
Below we shall show that these three properties define $\Mod_{g,[\gamma],\beta}$ uniquely as a complete weighted branched substack of $\mathcal O^{\circ}$ with total weight $1$. Suppose that we have two complete weighted branches substacks on a neighborhood of $f$ satisfying the above three properties. By choosing the neighborhood of $f$ sufficiently small, property \ref{modp1} implies that we may write the two weighted branched substacks in the form of 
\[w_{0}t^{\emptyset}+\sum_{k}w_{k}t^{\hat g_{k}}\text{ and }w_{0}t^{\emptyset}+
\sum_{l}w'_{l}t^{\hat g'_{l}}\]
so that each of the families $\hat g_{k}$ and $\hat g'_{l}$ are connected and contain $f$, and so that there is no map of any nonempty open subset  of one of the $\hat g_{k}$ into another $g_{k'}$ or map of any nonempty open subset of   $\hat g'_{l} $ into another $\hat g'_{l'} $. On the other hand, property \ref{modp2} implies that $\theta(\hat g'_{l})=w'_{l}t^{\dbar g'_{l}}+\dotsc$, so property \ref{modp3} implies that around every curve in $\hat g'_{l}$, there is a neighborhood with a map into some number of $\hat g_{k}$ so that the sum of the corresponding weights $w_{k}$ is at least $w'_{l}$. Suppose that there is a map of an open subset of $\hat g'_{l}$ into $\hat g_{k}$. Then the same properties imply that this may be composed with a map from some open subset of $\hat g_{k}$ into some $\hat g'_{l'}$ to get a map of an open subset of $\hat g'_{l}$ into $\hat g'_{l'}$, which implies that $l'=l$. Similarly, our open subset of $\hat g_{l}$ must have a map only to $\hat g_{k}$, and as $\hat g_{k}$ represents a substack, this map is unique. It follows that there are unique maps $\hat g_{k}\longrightarrow \hat g'_{l}$  and $\hat g'_{l}\longrightarrow \hat g_{k}$. As $\hat g'_{l}$ and $\hat g_{k}$ represent  substacks these two maps must be inverses of one another, and $\hat g_{k}=\hat g'_{l}$. Property \ref{modp3} then implies that $w_{k}=w'_{l}$.
 Similarly, all the other families and weights must be equal, so the two weighted branched substacks are equal.

We've seen that $\Mod_{g,[\gamma],\beta}$ is a complete weighted branched substack of $\mathcal O^{\circ}$ of some fixed dimension.  $\Mod_{g,[\gamma],\beta}$ also  has a well defined orientation relative to $\ex G$. This orientation is determined as follows:

The core family $\hat f_{i}/G_{i}$ comes with a collection of sections corresponding to marked points which when taken together give a section $s:\ex F(\hat f_{i})\longrightarrow \ex F(\hat f_{i}^{+l})$ so that $ev^{+l}(\hat f_{i}):\ex F(\hat f_{i}^{+l})\longrightarrow \M\times \sfp{\hat {\ex B}}{\ex G}l$ is an equidimensional embedding in a neighborhood of $s$. The canonical orientation of $\hat {\ex B}$ relative to $\ex G$ given by the almost complex structure, and the orientation of $\M$ given by the complex structure give an orientation to $\M\times\sfp{\hat {\ex B}}{\ex G}l$ relative to $\ex G$. Give $\hat f_{i}^{+l}$ the orientation relative to $\ex G$ so that this map $ev^{+l}(\hat f_{i})$ is oriented in a neighborhood of the image of $s$, and give $\hat f_{i}$ the corresponding orientation relative to $\ex G$ so that the complex fibers of $\ex F(\hat f_{i}^{+k})\longrightarrow \ex F(\hat f_{i}^{k-1})$ are positively oriented.  

Recall that the vector bundle $V_{i}$ over $\ex F(\hat f_{i})$ is oriented  relative to $\ex F(\hat f_{i})$ as follows: restricted to a curve $f$ in $\hat f_{i}$, we may identify $V_{i}(f)$ with the cokernel of the  injective  operator $D\dbar(f):X^{\infty,\underline 1}(f)\longrightarrow Y^{\infty,\underline 1}(f)$, and orient this using a homotopy of $D\dbar(f)$ to a complex map as in Remark \ref{orientation} on page \pageref{orientation}.

The orientation of $\hat f_{i}$ relative to $\ex G$ and the orientation of $V_{i}$ relative to $\hat f_{i}$ give an orientation to the $\hat g_{k}$ relative to $\ex G$ by considering $\hat g_{k}$ as the intersection of the section $\dbar'\nu_{k}$ of $V_{i}$ with the zero section. (The order of intersection does not matter as $V_{i}$ is always an even dimensional vector bundle because  the index of $D\dbar$ restricted to $X^{\infty,\underline 1}(f)$ is even as noted in item \ref{fr2} of Theorem \ref{f replacement}.) We must see why this construction gives a well defined orientation for $\Mod_{g,[\gamma],\beta}$ relative to $\ex G$ - in other words why we will get the same orientation using a different obstruction model. As a first step, we may replace the family $\hat f_{i}$ with the family $F(\nu_{k})$ which actually contains $\hat g_{k}$, and do our calculation of orientations at a curve $f$ in $\hat g_{k}$. This will not change the orientations constructed as above, as item \ref{rt1} of Theorem \ref{regularity theorem} imply that $V_{i}$ will remain complementary to the image of $D\dbar$. We may add a collection of $l'$ extra marked points and extend $F(\nu_{k})$  to a family $\hat h$ with extra parameters corresponding to the image of these extra marked points  so that $ev^{+(l+l')}(\hat h)$ is an equidimensional embedding in a neighborhood of the  section $s':\ex F(\hat h)\longrightarrow \ex F(\hat h^{+(l+l')})$ corresponding to all of these marked points.  Denote by $X'(f)$ the complex subspace of $X^{\infty,\underline 1}(f)$ consisting of sections which vanish at the extra marked points. The tangent space to the extra parameter space at the curve $f$ can also be identified with $X^{\infty,\underline 1}(f)/X'(f)$. The orientation of $\hat h$ relative to $\ex G$ given by $ev^{+(l+l')}(\hat h)$ agrees with the orientation from $\hat f_{i}$ and the orientation from the almost complex structure on the extra parameter space. 
Again use the method of Remark \ref{orientation} to orient of the cokernel of $D\dbar(f):X'(f)\longrightarrow Y^{\infty,\underline 1}(f)$. As noted in Remark \ref{orientation}, the orientation of this cokernel is compatible with the short exact sequence 
\[X^{\infty,\underline 1}/X'\longrightarrow Y^{\infty,\underline 1}(f)/D\dbar(f)(X')\longrightarrow Y^{\infty,\underline 1}(f)/D\dbar(f)(X^{\infty,\underline 1}(f))=V_{i}(f)\]
 Therefore, the orientation on $M_{g,[\gamma],\beta}$ we obtain does not depend on the choice of marked points in our obstruction model. Theorem \ref{regularity theorem} implies that  all obstruction models containing $f$ with the same set of marked points are homotopic in some neighborhood of $f$ as all other choices can be changed continuously. The orientation of $\hat f_{i}$ and the orientation of $V_{i}(f)$ given above do not change under homotopy, and the multisection of $V_{i}$ used to define $\Mod_{g,[\gamma],\beta}$ must change continuously under homotopy, but remain transverse to the zero section (and always intersect the zero section at $\Mod_{g,[\gamma],\beta}$). Therefore the orientation we obtain on    on $\Mod_{g,[\gamma],\beta}$ is well defined. Therefore item \ref{ft6} from Theorem \ref{final theorem} holds.

 It is clear from the construction of $\Mod_{g,[\gamma],\beta}$ that that its support in $\Msw_{g,[\gamma],\beta}$ is  compact, as it is a finite union of  compact subsets, so item \ref{ft7} of Theorem \ref{final theorem} is true.

We must now verify item \ref{ft8} of Theorem \ref{final theorem} which implies that that $\Mod_{g,[\gamma],\beta}$ gives a well defined cobordism class of finite dimensional weighted branched substacks oriented relative to $\ex G$. In particular, given any construction of a virtual moduli space $\mathcal M'_{g,[\gamma],\beta}$ defined using another small enough multi-perturbation $\theta'$ which is fixed point free and  transverse to $\dbar$, defined on some other open neighborhood of the holomorphic curves in $\Msw_{g,[\gamma],\beta}(\hat{\ex B})$ using different choices of obstruction models, we must construct $\Mod_{g,[\gamma],\beta}(\hat{\ex B}\times S^{1})$ so that its restriction to two different points $p_{1}$ and $p_{2}$ in $S^{1}$ give $\mathcal M_{g,[\gamma],\beta}$ and $\mathcal M'_{g,[\gamma],\beta}$ respectively. 
 
 Choose disjoint open intervals $I_{i}\subset S^{1}$  and compactly contained open intervals $I_{i}'\subset I_{i}$ containing $p_{i}$. We may take the product of any obstruction model on $\Msw(\hat{\ex B})$ with $I_{i}$ to obtain an obstruction model on $\Msw(\hat{\ex B}\times I_{i})\subset\Msw(\hat {\ex B}\times S^{1})$. Multiply the obstruction models used to define $\Mod_{g,[\gamma],\beta}$ and $\Mod_{g,[\gamma],\beta}$ by $I_{1}$ and $I_{2}$ respectively, and choose a finite set of obstruction models to cover the rest of the holomorphic curves in $\Msw_{g,[\gamma],\beta}(\hat{\ex B}\times S^{1})$ which project to subsets of $S^{1}$ not intersecting $I'_{i}$. Any modification of these obstruction models required for item \ref{ft1} may be chosen not to affect our product obstruction bundles on $I_{i}'$. 
 
For each of our original simple perturbations $\mathfrak P_{i}$ parametrized by $\hat f$, we may choose a compactly supported simple perturbation $\mathfrak P_{i,t}$ parametrized by $\hat f \times I_{1}$ so that $\mathfrak P_{i}=\mathfrak P_{i,t}$ for $t$ in a neighborhood of $p_{1}$, and so that the transversality and fixed point free conditions hold on the interval $I'_{1}$, and so that $\mathfrak P_{i,t}$ outside of $I'_{1}$ is small enough for the rest of the construction of $\Mod_{g,[\gamma],\beta}(\hat{\ex B})$ to proceed. We may carry out a similar procedure for the simple perturbations used to define $\Mod'_{g,[\gamma],\beta}$. Then $\Mod_{g,[\gamma],\beta}(\hat{\ex B})$  constructed using these perturbations restricts to be $\Mod_{g,[\gamma],\beta}(\hat{\ex B})$ over $p_{1}$ and $\Mod'_{g,[\gamma],\beta}(\hat{\ex B})$ over $p_{2}$.

 To prove item \ref{ft 9} from Theorem \ref{final theorem}, consider a $\C\infty1$ map of a compact exploded manifold $\ex G'$ to $\ex G$. As noted above, for an open dense subset of the space of simple perturbations, the map $\Mod_{g,[\gamma],\beta}\longrightarrow \ex G$ is transverse to $\ex G'\longrightarrow \ex G$. Given any obstruction model $(\hat f/G, V)$ for $\mathcal O_{i}\subset \Msw(\hat {\ex B})$, the fiber product of $\hat f\longrightarrow \ex G$ with $\ex G'\longrightarrow \ex G$ gives a family of curves $\hat f'$ in $\hat{\ex B}'$, the fiber product of $\hat{\ex B}\longrightarrow \ex G$ with $\ex G'\longrightarrow \ex G$. This is the inverse image of $\hat f$ under a natural map \[\Msw(\hat{\ex B}')\longrightarrow \Msw(\hat {\ex B})\] given by composing families of curves in $\hat{\ex B}'$ with the map $\hat{\ex B}'\longrightarrow \hat{\ex B}$. The action of $G$ on $\hat f$ gives an action of $G$ on $\hat f'$, and $V$ and the trivialization associated with the obstruction model pull back to $\hat{\ex B}'$ similarly to give an obstruction model $(\hat f'/G,V')$ on $\mathcal O_{i}'\subset\Msw(\hat{\ex B}')$ which is the inverse image of $\mathcal O_{i}'$. Note that as $\mathcal O$ meets $\mathcal O_{i}$ properly, the inverse image $\mathcal O'$ of $\mathcal O$ in $\Msw(\hat{\ex B}')$ also meets $\mathcal O_{i}'$ properly.  We can also pull back any compactly supported simple perturbation parametrized by $\hat f$ to a compactly supported simple perturbation parametrized by $\hat f'$. The multiperturbation defined on $\mathcal O'$ by the pullback of the simple perturbations used to define $\Mod_{g,[\gamma],\beta}(\hat{\ex B})$ is the pullback of the multiperturbation used to define $\Mod_{g,[\gamma],\beta}$, so the virtual moduli space it defines is the pullback of $\Mod_{g,[\gamma],\beta}$, which can also be regarded as the fiber product of $\ex G'\longrightarrow \ex G$ with $\Mod_{g,[\gamma],\beta}\longrightarrow \ex G$. If the simple perturbations used to define $\Mod_{g,[\gamma],\beta}$ are small enough, their pullback  will also be small enough to define a virtual moduli space in $\mathcal O'$, and  their pullback  will automatically satisfy the transversality and fixed point free requirements. Therefore, we may use the pullback of these simple perturbations to define the virtual moduli space within $\mathcal O'$. 
 
 If there is some set $\{\gamma_{i}\}$ of tropical curves in $\totb{\hat{\ex B}'}$ and maps $\beta_{i}:H^{2}(\hat{\ex B}')\longrightarrow \mathbb R$ so that the inverse image of $\Msw_{g,[\gamma],\beta}(\hat{\ex B})$ is $\coprod_{i}\Msw_{g,[\gamma_{i}],\beta_{i}}(\hat{\ex B}')$, then $\mathcal O'$ will be a neighborhood of the substack of holomorphic curves in $\coprod_{i}\Msw_{g,[\gamma_{i}],\beta_{i}}$. (A particular case of interest when more than one $\gamma_{i}$ and $\beta_{i}$ is required is when $\ex G'$ is equal to two points but $\ex G$ is connected.) Then each $\Mod_{g,[\gamma_{i}],\beta_{i}}$ can be constructed using the pulled back obstruction models and simple perturbations, and $\coprod_{i}\Mod_{g,[\gamma_{i}],\beta_{i}}$ is equal to the fiber product of $\Mod_{g,[\gamma],\beta}\longrightarrow \ex G$ with $\ex G'\longrightarrow \ex G$.

\

In the case of a single target $\ex B$,
$\Mod_{g,[\gamma],\beta}$ is 
a  finite dimensional $\C\infty1$ oriented weighted branched substack of $\Msw(\ex B)$.
This should be thought of as giving a virtual class for a component of  the moduli space of holomorphic curves in $\ex B$, which is a cobordism class of $\C\infty 1$ finite dimensional oriented weighted branched substacks of the moduli stack of $\C\infty1$ curves in $\ex B$. The above discussion implies that this virtual class behaves well in a family of targets $\hat {\ex B}\longrightarrow\ex G$, so enumerative invariants of holomorphic curves such as Gromov Witten invariants behave well in connected families of targets in the exploded category.

\section{Representing Gromov Witten invariants using differential forms}\label{numerical GW}

In this section, we define numerical Gromov Witten invariants by integrating differential forms over the virtual moduli space $\Mod_{g,[\gamma],\beta}$ constructed in section \ref{virtual class}. 

If $\ex X$ is an exploded manifold or orbifold, a $\C\infty1$ map from the moduli stack of $\omega$-positive  $\C\infty1$ curves $\Msw$ to $\ex X$  is a map from $\Msw$ to $\ex X$ considered as a stack. In particular, given any $\C\infty1$ family $\hat f$ of curves, it is a $\C\infty1$ map $\ex F(\hat f)\longrightarrow \ex X$ so that given any map of $\C\infty1$ families $\hat f\longrightarrow \hat g$, the following is a commutative diagram.
\[\begin{array}{ccc}\ex F(\hat f)&\longrightarrow &\ex X
\\ \downarrow & & \downarrow\id
\\ \ex F(\hat g)&\longrightarrow & \ex X
\end{array}\]

For example, $ev^{0}:\Msw\longrightarrow \M$ is a $\C\infty1$ map. Given any $\C\infty1$ obstruction model, $(\hat f/G,V)$, the projection from an open neighborhood of $\hat f$ in $\Msw$ to $\ex F(\hat f)/G$ is also a $\C\infty1$ map. Below, we shall construct a $\C\infty1$ evaluation map $\EV$ on $\Msw$. Numerical invariants can be extracted from the virtual moduli space $\Mod_{g,[\gamma],\beta}$ constructed in section \ref{virtual class} by integrating over $\Mod_{g,[\gamma],\beta}$ differential forms pulled back from $\C\infty1$ maps from $\Msw_{g,[\gamma],\beta}$ to finite dimensional exploded manifolds or orbifolds.

\subsection{The evaluation map EV}  \label{construction of EV}

\

\

Given a punctured holomorphic curve $f$ in a manifold $M$, evaluation at a puncture of $f$ gives a point in $M$. If $f$ is instead a holomorphic curve in an exploded manifold $\ex  B$, evaluation at a puncture of $f$ may not give a point in $\ex B$. To remedy this, we shall define the `ends' of $\ex B$ as follows:

 Given an exploded manifold $\ex B$, let $\End (\ex B)$ denote the moduli space of maps $h:\et 1{(0,\infty)}\longrightarrow \ex B$ with the equivalence relation 
 \[ h_{1}\equiv h_{2}\text{ if }  h_{1}(\tilde z)= h_{2}(c\tilde z)\text{ or } h_{2}(\tilde z)= h_{1}(c\tilde z) \]
for some constant $c\in \mathbb C^{*}\e{\mathbb R}$. The above defines $\End (\ex B)$ as a set. We can give $\End(\ex B)$ the structure of an exploded manifold so that any smooth or $\C\infty1$ map of a $\et 1{(0,\infty)}$ bundle over $\ex A$ into $\ex B$ is associated to a smooth or $\C\infty1$ map $\ex A\longrightarrow \End(\ex B)$. 
 
In particular, let $v\in \mathbb Z^{m}$ be a nonzero integer vector so that there is an infinite ray contained in the polytope $ P\subset\mathbb R^{m}$ in the direction of $v$. To $v$ we may associate a connected component $\End_{v}(\et mP)$ of $\End(\et mP)$ as follows: By a coordinate change, we may assume that $v=(0,\dotsc,0,k)$. Let $P_{v}$ be the image of $P$ under the projection $\mathbb R^{m}\longrightarrow\mathbb R^{m-1}$ which forgets the last coordinate. Then our connected component, $\End_{v}(\et mP)$ is equal to $\et {m-1}{P_{v}}$. A map $ h(\tilde z)=(c_{1},\dotsc,c_{m-1},c_{m}\tilde z^{k})$ then corresponds to the point $(c_{1},\dotsc,c_{m-1})\in \et {m-1}{P_{v}}$. Given any smooth or $\C\infty1$ map of a $\et 1{(0,\infty)}$ bundle over  $\ex A$ to $\et mP$ which is of the above type restricted to fibers, we can project to the first $k-1$ coordinates to get a smooth or $\C\infty1$ map $\ex A\longrightarrow \et {m-1}{P_{v}}$. 

All connected components of $\End(\et mP)$ are in the above form apart from the component $\End_{0}(\et mP)$ which is  equal to   $\et mP$. A map from a $\et 1{(0,\infty)}$ bundle over $\ex A$ into $\et mP$ which is trivial on the fibers is associated to the same map from $\ex A$ to $\et mP$.

 Then 
\[\End(\et mP)=\coprod_{v}\End_{v}(\et mP)\]
where the disjoint union is over all  integer vectors $v\in\mathbb Z^{m}$ so that there is an infinite ray in the direction of $v$ contained in $P$. 

 We may similarly define $\End (\mathbb R^{n}
\times \et mP)$ to be $\mathbb R^{n}\times \End(\et mP)$, and given any open subset $U$ of $\mathbb R^{n}\times \et mP$, $\End(U)$ can be identified with the open subset of $\End(\mathbb R^{n}\times \et mP)$ corresponding to maps to the open subset $U$. This construction is clearly functorial, given any smooth map $f:U\longrightarrow V$ between open subsets of $\mathbb R^{n}\times \et mP$, there exists a smooth map 

\[\End(f):\End(U)\longrightarrow \End(V)\]
sending the point corresponding to  $ h$ to the point corresponding to $ h\circ f$. The functoriality of the construction of the exploded structure on open subsets of coordinate charts implies that we may give $\End(\ex B)$ a well defined exploded structure by giving the subset of $\End(\ex B)$  corresponding to a coordinate chart $\mathbb R^{n}\times \et mP$ the structure of $\End(\mathbb R^{n}\times \et mP)$.

Given a $\C\infty1$ curve $f$ in $\ex B$ with $n$ labeled punctures,  restricting $f$ to the copy of $\et 1{(0,\infty)}$ around a puncture gives a class of maps $ h:\et 1{(0,\infty)}\longrightarrow \ex B$ which specifies a point in $\End(\ex B)$. Together, this gives a point $\EV(f)\in (\End(\ex B))^{n}$. Recall that $\Msw_{g,[\gamma],\beta}$ is the moduli space of $\C\infty1$ curves with genus $g$, homology class $\beta$ and 
tropical part isotopic to a particular tropical curve $\gamma$ in $\totb{\ex B}$. If $\gamma$ has $n$ infinite ends, we may define $\End_{\gamma}(\ex B)$ to be the connected component of $(\End(\ex B))^{n}$ which contains $\EV(f)$ for any curve $f$ in $\Msw_{g,[\gamma],\beta}$.
 
Applying this construction to each curve in a $\C\infty1$ family of curves $\hat f$ in $\Msw_{g,[\gamma],\beta}$    gives $\C\infty1$ map 
\[\EV:\ex F(\hat f)\longrightarrow \End_{\gamma}(\ex B)\]
which defines a $\C\infty1$ map
\[\EV:\Msw_{g,[\gamma],\beta}(\ex B)\longrightarrow \End_{\gamma}(\ex B)\]

\subsection{Integration of forms over virtual class}

\

\

We can define numerical Gromov Witten invariants by pulling back differential forms on $\M_{g,[\gamma]}\times\End_{\gamma}(\ex B)$ using  $ev^{0}\times\EV$ and integrating over the virtual moduli space $\Mod_{g,[\gamma],\beta}$ defined in section \ref{virtual class}.

We shall use the following class of differential forms discussed in \cite{dre}.

\begin{defn}[$\Omega^{*}(\ex B)$] Let  $\Omega^{k}(\ex B)$ be the vector space of $\C\infty1$ differential $k$ forms $\theta$ on $\ex B$ so that for all integral vectors $v$, the differential form $\theta$ vanishes on $v$, and for all maps $f:\et 1{(0,\infty)}\longrightarrow \ex B$, the differential form $\theta$ vanishes on all vectors in the image of $df$. 

Denote by $\Omega^{k}_{c}(\ex B)\subset \Omega^{k}(\ex B)$ the subspace of forms with complete support. (A form has complete support if  the set where it is non zero is contained inside a complete subset of $\ex B$ - in other words, a compact subset with tropical part consisting only of complete polytopes.)

 Denote the homology of $(\Omega^{*}(\ex B),d)$ by $H^{*}(\ex B)$, and the homology of $(\Omega^{*}_{c}(\ex B),d)$ by $H^{*}_{c}(\ex B)$.
\end{defn}

We use $\Omega^{*}(\ex B)$ instead of all $\C\infty1$ differential forms on $\ex B$ in order to use a version of Stokes' theorem proved in \cite{dre}.

\begin{defn}[Refined forms] 
  A refined form $\theta\in\ro^{*}(\ex B)$ is choice $\theta_{p}\in\bigwedge T^{*}_{p}(\ex B)$ for all $p\in\ex B$ so that given any point $p\in \ex B$, there exists an open neighborhood $U$ of $p$ and a complete, surjective, equidimensional submersion
 \[r:U'\longrightarrow U\]
  
 so that there is a form $\theta'\in\Omega^{*}(U')$ which is the pullback of $\theta$ in the sense  that if $v$ is any vector on $U'$ so that $dr(v)$ is a vector based at $p$, then 
  \[\theta'(v)=\theta_{p}(dr(v))\] 
 
 A refined form  $\theta\in\ro^{*}(\ex B)$ is completely supported if there exists some complete subset $V$ of an exploded manifold $\ex C$ with a map $\ex C\longrightarrow \ex B$ so that $\theta_{p}=0$ for all $p$ outside the image of $V$. Use the notation $\ro^{*}_{c}$ for completely supported refined forms. 
 
 Denote the homology of $(\ro^{*}(\ex B),d)$ by $\rh^{*}(\ex B)$ and  $(\ro^{*}_{c}(\ex B),d)$ by $\rh^{*}_{c}(\ex B)$.
 
 \end{defn}

The Poincare dual to a map $\ex C\longrightarrow \ex B$ as defined in \cite{dre} is correctly viewed as a refined differential form.

\begin{defn}[Differential forms generated by functions]\label{rof} A differential form is generated by functions if it is locally equal to a form constructed from $\C\infty1$ functions using the operations of exterior differentiation and wedge products. Use the notation $\rof^{*}(\ex B)\subset\ro^{*}(\ex B)$ for the set of refined forms on $\ex B$ which are locally equal to some differential form which is generated by functions on a refinement. Similarly, let $\rof_{c}(\ex B)=\ro^{*}_{c}(\ex B)\cap\rof^{*}(\ex B)$.   
\end{defn}

Differential forms generated by functions will be important in the gluing formula for Gromov Witten invariants from Theorem \ref{gluing formula}. Examples of differential forms generated by functions are the Poincare dual to a point, the Chern class defined using the Chern Weil construction, and any smooth differential form on a smooth manifold.

\

Our virtual moduli space $\Mod_{g,[\gamma],\beta}$ is an oriented weighted branched substack of $\Msw_{g,[\gamma],\beta}$. 
Locally, restricted to an open substack $\mathcal O\subset\Msw_{g,[\gamma],\beta}$ it is equal to 
\begin{equation}\label{mod form}\Mod_{g,[\gamma],\beta}=\sum_{i}w_{i}t^{\hat f_{i}}\end{equation}
where $\hat f_{i}$ are finite dimensional $\C\infty1$ families which are proper substacks of $\mathcal O$, $w_{i}$ are positive numbers, and $t$ is just a dummy variable.

\begin{defn}
A partition of unity on $\Mod_{g,[\gamma],\beta}$ subordinate to an open cover $\{\mathcal O\}$ is a countable set of nonnegative $\C\infty1$ functions  $\rho_{i}:\Mod_{g,[\gamma],\beta}\longrightarrow \mathbb R$ locally given by the restriction to $\Mod_{g,[\gamma],\beta}$ of $\C\infty1$ functions on open subsets of $\Msw$ so that $\rho_{i}$ is compactly supported inside some $\mathcal O\cap \Mod_{g,[\gamma],\beta}$, and so that $\sum_{i}\rho_{i}=1$.
\end{defn}
Assuming Gromov compactness holds for $\hat {\ex B}$, construct a partition of unity on $\Mod_{g,[\gamma],\beta}(\hat{\ex B})$ subordinate to an open cover $\{\mathcal O\}$ in which $\Mod_{g,[\gamma],\beta}$ has the form of equation \ref{mod form} as follows. For any curve $f\in \Mod_{g,[\gamma],\beta}\cap \mathcal O$, the construction of $\Mod_{g,[\gamma],\beta}$ using obstruction models implies that on some open substack $\mathcal O'\subset\mathcal O$ we may  choose some $\C\infty1$ positive function $\rho':\mathcal O'\longrightarrow \mathbb R$ so that $\rho'(f)>0$, and $\rho'$ restricted to $\mathcal O'\cap \Mod_{g,[\gamma],\beta}$ is compactly supported. (In particular, if $(\hat f/G,V$) is the obstruction model used, we may choose $\rho'$ to to come from a compactly supported $G$-invariant $\C\infty1$ function on $\ex F(\hat f)$).  We may choose some countable collection  $\{ \rho'_{i}\}$ of these $ \rho'$ so that the supports of $  \rho_{i}'$ form a  locally finite cover of $\Mod_{g,[\gamma],\beta}$. Extend each $\rho_{i}'$ to be zero wherever it is not yet defined on $\Mod_{g,[\gamma],\beta}$. Then define
\[\rho_{j}=\rho'_{j}\lrb{\sum_{i} \rho'_{i}}^{-1}\] 
Each of these functions $\rho_{i}$ is locally given by the restriction of a $\C\infty1$ function from an open substack of the stack of $\C\infty1$ curves. As $\Mod_{g,[\gamma],\beta}$ is a $\C\infty1$ substack of $\Msw_{g,[\gamma],\beta}$, $\rho_{i}$ is a $\C\infty1$ function on $\Mod_{g,[\gamma],\beta}$, and together $\{\rho_{i}\}$ form a partition of unity subordinate to $\{\mathcal O\}$. In the case that $\hat {\ex B}$ is compact, this is a finite partition of unity.

\

Let $\Phi:\Msw_{g,[\gamma],\beta}\longrightarrow \ex X$ be a $\C\infty1$ map (such as $ev^{0}\times \EV$). Given any form $\theta$  in $\Omega^{*}(\ex X)$ or $\ro^{*}(\ex X)$, we can define 

\begin{equation}\label{intdef1}\int_{\Mod_{g,[\gamma],\beta}}\Phi^{*}\theta:=\sum_{i}\int_{\Mod_{g,[\gamma],\beta}}\rho_{i}\Phi^{*}\theta\end{equation}

And define 
\begin{equation}\label{intdef2}\int_{\Mod_{g,[\gamma],\beta}}\rho_{i}\Phi^{*}\theta:=\sum_{j} w_{j}\int_{\ex F(\hat f_{j})}\rho_{i}\Phi^{*}\theta\end{equation}
where restricted to the support of $\rho_{i}$,
\[\Mod_{g,[\gamma],\beta}=\sum_{j}w_{j}t^{\hat f_{j}}\]
The integrals on the right hand side of equation \ref{intdef2} are of forms in $\Omega^{*}_{c}(\ex F(\hat f_{j}))$ or $\ro^{*}_{c}(\ex F(\hat f_{j}))$, so \cite{dre} implies that these integrals are well defined.
Clearly, if $\sum_{k}w'_{k}t^{\hat f'_{k}}=\sum_{j}w_{j}t^{\hat f_{k}}$, then 
\[\sum_{k} w'_{k}\int_{\ex F(\hat f'_{k})}\rho_{i}\Phi^{*}\theta=
\sum_{j} w_{j}\int_{\ex F(\hat f_{j})}\rho_{i}\Phi^{*}\theta\]
so the definition of $\int_{\Mod_{g,[\gamma],\beta}}\rho_{i}\Phi^{*}\theta$ is defined independent of the particular local representation of $\Mod_{g,[\gamma],\beta}$ chosen. As usual, the linearity of the integral implies independence of this definition of choice of partition of unity: in particular, if $\{\rho'_{j}\}$ is another partition of unity appropriate for defining the integral, then 
\[\sum_{i}\int_{\Mod_{g,[\gamma],\beta}}\rho_{i}\Phi^{*}\theta=\sum_{i,j}\int_{\Mod_{g,[\gamma],\beta}}\rho'_{j}\rho_{i}\Phi^{*}\theta=\sum_{j}\int_{\Mod_{g,[\gamma],\beta}}\rho'_{j}\Phi^{*}\theta
\] 
 Therefore, $\int_{\Mod_{g,[\gamma],\beta}}\Phi^{*}\theta$ is well defined for any $\theta$ in $\Omega^{*}(\ex X)$ or $\ro^{*}(\ex X)$. 
 
 \
 
 Note that $\int_{\Mod_{g,[\gamma],\beta}}\Phi^{*}\theta$ depends only on the cohomology class of $\theta$ in $\rh^{*}(\ex X)$. In particular, $\int_{\Mod_{g,[\gamma],\beta}}\Phi^{*}d\alpha$ vanishes, because linearity implies that it is equal to $\sum_{i}\int_{\Mod_{g,[\gamma],\beta}}\Phi^{*}d(\rho_{i}\alpha)$. Each of these integrals  is some sum of integrals of compactly supported exact forms in coordinate charts, which vanish because of the version of Stokes' theorem proved in \cite{dre}.
 
 \
 
 More importantly, Theorem \ref{well defined} below states that if $\theta$ is closed, $\int_{\Mod_{g,[\gamma],\beta}}\Phi^{*}\theta$ is independent of the choices involved in defining $\Mod_{g,[\gamma],\beta}$.
  
If $\theta$ is closed and generated by functions, the contribution to the integral of $\theta$ from individual tropical curves is also well defined. Given a tropical curve $\gamma_{0}$ mapping to $\totb{\ex B}$, let $\Mod_{g,[\gamma],\beta}\rvert_{\gamma_{0}}$ be the restriction of $\Mod_{g,[\gamma],\beta}$ to the subset consisting of curves with tropical part $\gamma_{0}$.
It follows from the definition of integration given in \cite{dre} that the integral $\int_{\Mod_{g,[\gamma],\beta}}\Phi^{*}\theta$ breaks up into a finite sum of integrals  over $\Mod_{g,[\gamma],\beta}\rvert_{\gamma_{0}}$ for different $\gamma_{0}$. The following theorem proves that if $\theta$ is closed and generated by functions, then 
$\int_{\Mod_{g,[\gamma],\beta}\rvert_{\gamma_{0}}}\Phi^{*}\theta$ is well defined independent of the choices involved in defining $\Mod_{g,[\gamma],\beta}$. The gluing formula of Theorem \ref{gluing formula} gives a way of calculating some of these integrals for particular $\gamma_{0}$ using relative invarants.

 \begin{thm}\label{well defined}
 If $(\ex B,J)$ is a basic exploded manifold for which Gromov compactness holds,  $\Phi:\Msw_{g,[\gamma],\beta}(\ex B)\longrightarrow \ex X$ is any $\C\infty1$ map, and $\theta\in\ro^{*}(\ex X)$ is closed, then
 \[\int_{\Mod_{g,[\gamma],\beta}}\Phi^{*}\theta\] 
 is independent of the choices made in the definition of $\Mod_{g,[\gamma],\beta}$.

If  $\theta\in\rof^{*}_{c}(\ex X)$ is generated by functions, and $\gamma$ is any tropical curve mapping to  $\totb{\ex B}$, the integral
\[\int_{\Mod_{g,[\gamma],\beta}\rvert_{\gamma}}\Phi^{*}\theta\]
is well defined independent of the choices made in defining $\Mod_{g,[\gamma],\beta}$.

 \end{thm}
 
 \pf

This follows directly from Stokes' theorem and item \ref{ft8} of Theorem \ref{final theorem}. In particular, if $\Mod_{g,[\gamma],\beta}'$ is the result of other choices, then $\Mod_{g,[\gamma],\beta}$ and $\Mod'_{g,[\gamma],\beta}$ are cobordant in the following sense:  there exists a construction of $\Mod_{g,[\gamma],\beta}(\ex B\times S^{1})$ in a trivial family of targets $\ex B\times S^{1}\longrightarrow S^{1}$ which restricts to our two different moduli spaces at two different points $p_{1}$ and $p_{2}$  of $S^{1}$.

The projection $\ex B\times S^{1}\longrightarrow \ex B$ forgetting $S^{1}$ gives a map $\pi:\Msw(\ex B\times S^{1})\longrightarrow \Msw(\ex B)$. The form $(\Phi\circ\pi)^{*}\theta$ is a closed form on each coordinate chart. Denote the restriction of $ \Mod_{g,[\gamma],\beta}$ to curves over an interval joining  $p_{1}$ and $p_{2}$ by $\tilde \Mod$. We may orient this interval so that the boundary of $\tilde\Mod$ is $\Mod'_{g,[\gamma],\beta}-\Mod_{g,[\gamma],\beta}$. As argued above, we may construct a finite partition of unity $\{\rho_{j}\}$ over $\tilde \Mod$ so that each $\rho_{j}$ is compactly supported within a subset of $\tilde \Mod$ which is equal to $\sum w_{i}t^{\hat h_{i}}$. Then the version of Stokes's theorem proved in \cite{dre} implies that 
\[\int_{\tilde \Mod}d(\rho_{j}(\Phi\circ\pi)^{*}\theta)=\int_{\partial\tilde \Mod}\rho_{j}(\Phi\circ\pi)^{*}\theta=\int_{\Mod'_{g,[\gamma],\beta}}\rho_{j}\Phi^{*}\theta-\int_{\Mod_{g,[\gamma],\beta}}\rho_{j}\Phi^{*}\theta\]

Therefore, 
\[\int_{\Mod'_{g,[\gamma],\beta}}\Phi^{*}\theta-\int_{\Mod_{g,[\gamma],\beta}}\Phi^{*}\theta=\int_{\tilde \Mod}d(\Phi\circ\pi)^{*}\theta=0\]
as required.

Now consider the case where $\theta$ is generated by functions. In this case, we can apply Stokes' theorem to the integral over $\tilde \Mod\rvert_{\gamma}$ of $d\rho_{j}(\Phi\circ\pi)^{*}\theta$. The integral over $\tilde \Mod$ is a weighted sum of integrals over refinements $U'$ of coordinate charts where $\rho_{j}(\Phi\circ\pi)^{*}\theta$ is compactly supported; the integral over $\tilde\Mod\rvert_{\gamma}$ replaces each of these integrals with the integral restricted to  a strata of $U'$ with tropical part equal to a point. This in turn is equal to the integral of a form $d\alpha$ over the tropical completion of this strata (see definition \ref{tropical completion} on page \pageref{tropical completion}), where on the subset where our new and old coordinate charts agree, $\alpha$ is equal to $\rho_{j}(\Phi\circ\pi)^{*}\theta$. The fact that $\theta$ is generated by functions is required for $\alpha$ to be in $\Omega^{*}$ of this new coordinate chart. This new form $\alpha$ is completely supported, so we may apply the version of  Stokes' theorem from \cite{dre} to get
\[\int_{\tilde \Mod\rvert_{\gamma}}d(\rho_{j}(\Phi\circ\pi)^{*}\theta)=\int_{\partial\tilde \Mod\rvert\gamma}\rho_{j}(\Phi\circ\pi)^{*}\theta=\int_{\Mod'_{g,[\gamma],\beta}\rvert_{\gamma}}\rho_{j}\Phi^{*}\theta-\int_{\Mod_{g,[\gamma],\beta}\rvert_{\gamma}}\rho_{j}\Phi^{*}\theta\]

Therefore, if $\theta\in\rof^{*}(\ex X)$ is closed, then
\[\int_{\Mod_{g,[\gamma],\beta}\rvert_{\gamma}}\Phi^{*}\theta\]
is independent of the choices made in the definition of $\Mod_{g,[\gamma],\beta}$.

\stop

  \begin{thm}\label{pd class}
 If $(\ex B,J)$ is a basic almost complex exploded manifold with taming form $\omega$ so that  Gromov compactness holds, $\ex X$ is an oriented exploded manifold or orbifold  and $\Phi:\Msw_{g,[\gamma],\beta}(\ex B)\longrightarrow \ex X$ is any $\C\infty1$ map, then there exists a closed form $\eta_{g,[\gamma],\beta}\in\ro^{*}_{c}(\ex X)$ Poincare dual to the map $\Phi:\Mod_{g,[\gamma],\beta}\longrightarrow \ex X$ in the sense that 
 \[\int_{\Mod_{g,[\gamma],\beta}}\Phi^{*}\theta=\int_{\ex X}\theta\wedge\eta_{g,[\gamma],\beta}\]
  
for all closed $\theta\in\ro^{*}(\ex X)$. 

The class of $\eta_{g,[\gamma],\beta}$ in $\rh^{*}_{c}(\ex X)$ defined by the procedure in the proof below is independent of choices made in defining it and defining $\Mod_{g,[\gamma],\beta}$.

 \end{thm}

\pf

Extend $\Phi$ to a $\C\infty1$ submersion $\psi:\Mod_{g,[\gamma],\beta}\times\mathbb R^{n}\longrightarrow\ex X$ so that $\Phi(p)=\Phi(p,0)$. Choose a compactly supported closed form $\eta_{0}$  on $\mathbb R^{n}$ with integral $1$, then let 
\[\eta_{g,[\gamma],\beta}=\psi_{!}\eta_{0}\]
More explicitly, choose a finite partition of unity $\{\rho_{i}\}$ so that restricted to the support of $\rho_{i}$, $\Mod_{g,[\gamma],\beta}$ is equal to $\sum_{j} w_{i,j}t^{\hat f_{i,j}}$, and $\psi$ is equal to a $H$-equivariant map from $H\times \coprod \ex F(\hat f_{i,j})\times \mathbb R^{n}$ into some coordinate chart $(U,H)$ on $\ex X$.
 
 \[\eta_{g,n,E}:=\frac 1{\abs H}\sum_{h\in H}\sum_{i,j}w_{i,j}(\psi_{h,i,j})_{!}\rho_{i}\eta_{0}\]
 
 where the map $(\psi_{h,i,j})_{!}$ indicates integration along the fiber of $\psi$ restricted to $(h,\ex F(\hat f_{i,j})\times \mathbb R^{n})$. (Integration along the fiber is discussed in \cite{dre}. The resulting form $\eta_{g,n,E}$ is  a closed form in $\ro^{*}_{c}(\ex X)$.) Using the defining property of integration along the fiber discussed in \cite{dre} gives
 
 \[\begin{split}\int_{X}\theta\wedge\eta_{g,[\gamma],\beta} &=\frac 1{\abs H}\sum_{h\in H}\sum_{i,j}w_{i,j}\int_{X}\theta\wedge(\psi_{h,i,j})_{!}\rho_{i}\eta_{0}
\\ &=\sum_{i,j}w_{i,j}\int_{\ex F(\hat f_{i,j})\times\mathbb R^{n}}(\frac 1{\abs H}\sum_{h\in H}\psi_{h,i,j}^{*}\theta)\wedge \rho_{i}\eta_{0}
 \\&=\int_{\Mod_{g,[\gamma],\beta}\times\mathbb R^{n}}(\psi^{*}\theta)\wedge\eta_{0}\end{split}\]
Using Stokes theorem (with details expanded as in the proof of Theorem \ref{well defined}) allows us to deform the map $\psi(p,x)$ to $\psi(p,0)=\Phi(p)$ without affecting the above integral, so
 \[\int_{\Mod_{g,[\gamma],\beta}\times\mathbb R^{n}}(\psi^{*}\theta)\wedge\eta_{0}=\int_{\Mod_{g,[\gamma],\beta}}\Phi^{*}\theta\]

We must now verify that  $[\eta_{g,[\gamma],\beta}]\in \rh_{c}^*(\ex X)$ is independent of all choices. First, note that as $\psi_{!}d=d\psi_{!}$, the cohomology class $[\eta_{g,[\gamma],\beta}]$ does not depend on the choice of $\eta_{0}$.

Theorem \ref{final theorem} part \ref{ft8} implies that given a different construction of the virtual moduli space $\Mod_{g,[\gamma],\beta}'$ and extension of $\phi$ to a submersion $\psi':\Mod_{g,[\gamma],\beta}'\times\mathbb R^{n'}\longrightarrow \ex X$ resulting in a form $\eta'_{g,[\gamma],\beta}$, there is a construction of $\Mod_{g,[\gamma],\beta}(\ex B\times S^{1})$ which restricts to $\Mod_{g,[\gamma],\beta}$ and $\Mod'_{g,[\gamma],\beta}$ at two different points $p_{1}$ and $p_{2}$ of $S^{1}$. The map $\Phi$  extends to a map \[\hat \Phi:\Mod_{g,[\gamma],\beta}(\ex B\times S^{1})\longrightarrow \ex X\times S^{1}\] We may  also extend $\hat \Phi$ to a submersion \[\hat \psi:\Mod_{g,[\gamma],\beta}(\ex B\times S^{1})\times \mathbb R^{N}\longrightarrow \ex X\times S^{1}\] so that the following conditions hold:
\begin{itemize}
\item $\hat \psi^{-1}(\ex X\times \{p_{1}\})=\Mod_{g,[\gamma],\beta}\times \mathbb R^{N}$, and there exists a linear projection from $\mathbb R^{N}$ to $\mathbb R^{n}$ so that restricted to the inverse image of $\ex X\times\{p_{1}\}$, $\hat\psi$ factors as this projection followed by $\psi$.
\item $\hat \psi^{-1}(\ex X\times \{p_{2}\})=\Mod'_{g,[\gamma],\beta}\times \mathbb R^{N}$, and there exists a linear projection from $\mathbb R^{N}$ to $\mathbb R^{n'}$ so that restricted to the inverse image of $\ex X\times\{p_{2}\}$, $\hat\psi$ factors as this projection followed by $\psi'$.
\end{itemize}

Let $\hat \eta$ be a compactly supported, closed differential form on $\mathbb R^{N}$ with integral $1$. Then $\hat \psi_{!}\hat \eta$ is a compactly supported closed differential form in $\ro^{*}(\ex X\times S^{1})$ which restricts to have the same class as $\eta_{g,[\gamma],\beta}$ on $X\times\{p_{1}\}$ and the same class as $\eta'_{g,[\gamma],\beta}$ on $X\times \{p_{2}\}$. It follows that $\eta_{g,[\gamma],\beta}$ and $\eta'_{g,[\gamma],\beta}$ represent the same cohomology class in $\rh^{*}_{c}(\ex X)$.

\stop

The next theorem establishes that Gromov Witten invariants do not change in families:

\begin{thm}\label{family invariance}
Suppose that $\pi_{\ex G}:(\hat{\ex B},J)\longrightarrow \ex G$ is a compact basic family of almost complex exploded manifolds in which Gromov compactness holds.  Suppose further that there is a commutative diagram of $\C\infty1$ maps

\[\begin{array}{ccc}\Msw_{g,[\gamma],\beta}(\hat{\ex B})&\xrightarrow{\hat \Phi}&\hat{\ex X} 
\\  \downarrow &&\downarrow \pi
\\ \ex G&\xrightarrow{\id}&\ex G
\end{array}\]

Then there exists a closed form \[\hat\eta_{g,[\gamma],\beta}\in\ro^{*}(\hat{\ex X})\] so that for any $p\in\ex G$ if a collection $\{(\gamma_{i}, \beta_{i})\}$ of tropical curves and homology classes in $\pi_{\ex G}^{-1}(p)$  satisfies the condition that the restriction of $\Msw_{g,[\gamma],\beta}(\hat {\ex B})$ to  $\Msw(\pi_{\ex G}^{-1}(p))$ is equal to $\coprod_{i}\Msw_{g,[\gamma_{i}],\beta_{i}}(\pi_{G}^{-1}(p))$, then the restriction of $\hat \eta_{g,[\gamma],\beta}$ to $\pi^{-1}(p)$ is  a form in the same cohomology class as $\sum_{i}\eta_{g,[\gamma_{i}],\beta_{i}}$.  

These $\eta_{g,[\gamma_{i}],\beta_{i}}$ should be  constructed as in  Theorem \ref{pd class} where $\pi^{-1}(p)$ plays the role of $\ex X$, $(\ex B,J)$ is $\pi_{\ex G}^{-1}(p)$, and $\Phi$ is the restriction of $\hat\Phi$ to this fiber.

\end{thm}
 

Construct $\eta_{g,[\gamma],\beta}$ as a Poincare dual to $\hat\Phi$ as in Theorem \ref{pd class}.
 In particular,
 we can extend $\hat\Phi:\Mod_{g,[\gamma],\beta}\longrightarrow \hat{\ex X}$ to a submersion 
 \[\hat\psi:\Mod_{g,[\gamma],\beta}\times \mathbb R^{n}\longrightarrow \hat{\ex X}\] 
 so that given any regular value $p$ of $\Mod_{g,[\gamma],\beta}\longrightarrow\ex G$, $\hat\psi$   restricted to the fiber over $p$ is a submersion into $\pi^{-1}(p)$, and so that $\hat\psi$ restricted $\Mod_{g,[\gamma],\beta}\times 0$ is equal to $\hat \Phi$. Choose a compactly supported form $\eta_{0}$ on $\mathbb R^{n}$ with integral $1$, then let
 \[\hat\eta_{g,[\gamma],\beta}:=\hat\psi_{!}\eta_{0}\]
 where integration along the fiber is defined using a finite partition of unity $\{\rho_{i}\}$; so if $\hat{\Mod}_{g,[\gamma],\beta}$ is equal to $\sum_{i,j}w_{i,j}t^{\hat f_{i,j}}$ on the support of $\rho_{i}$, then 
 \[\hat\psi_{!}\eta_{0}:=\sum_{i,j}w_{i,j}\lrb{\hat\psi\rvert_{\ex F(\hat f_{i,j})}}_{!}(\rho_{i}\eta_{0})\] 
 The integration along the fiber on the right hand side of the above equation is as defined in \cite{dre}. (One complication is that $\ex F(\hat f_{i,j})$ is only oriented relative to $\ex G$, but the fibers of $\hat \psi$ are all contained in the fibers of the projection to $\ex G$, so an orientation of $\ex G$ is not necessary.) The resulting form is in $\ro_{c}^*(\hat {\ex X})$.

It is proved in \cite{dre} that integration along the fiber is compatible with fiber products. Suppose that $p$ is a regular value of $ \Mod_{g,[\gamma],\beta}\longrightarrow \ex G$. Theorem \ref{final theorem} item \ref{ft 9} tells us that so long as the perturbations used to define $\Mod_{g,[\gamma],\beta}$ are small enough, then $\coprod_{i}\Mod_{g,[\gamma_{i}],\beta_{i}}(\pi_{\ex G}^{-1}(p))$ may be constructed to be equal to the restriction of $\Mod_{g,[\gamma],\beta}$ to the inverse image of $p$. Let $\hat f'_{i,j}$ indicate the restriction of $\hat f_{i,j}$  to the fiber over $p$, let $\psi$ indicate the restriction of $\hat\psi$ to the fiber over $p$ and let $\ex X$ be the fiber over $p$ of the map $\pi_{\ex G}:\hat{\ex X}\longrightarrow \ex G$. Then the following diagram is commutative

\[\begin{array}{ccc}\ro^{*}\lrb{\ex F(\hat f'_{i,j})\times\mathbb R^{n}}&\xrightarrow{\psi_{!}}&\ro^{*}(\ex X)
\\ \uparrow &&\uparrow
\\ \ro^{*}\lrb{\ex F(\hat f_{i,j})\times\mathbb R^{n}}&\xrightarrow{\hat \psi_{!}} &\ro^{*}( \hat{\ex X})
\end{array}\]

It follows that we may construct $\eta_{g,[\gamma_{i}],\beta_{i}}$ as in Theorem \ref{pd class} so that  the restriction of our $\hat \eta_{g,[\gamma],\beta}$ to the fiber over $p$ is $\sum_{i}\eta_{g,[\gamma_{i}],\beta_{i}}$. 

As in the proof of Theorem \ref{pd class}, the cohomology class of $\eta_{g,[\gamma],\beta}$ in $\rh^{*}_{c}(\hat{\ex X})$ does not depend on the construction of $\Mod_{g,[\gamma],\beta}$. Part \ref{ft 9} of Theorem \ref{final theorem} tells us that given any $p\in \ex G$, we may construct $ \Mod_{g,[\gamma],\beta}'$ so that $p$ is a regular value, so the corresponding construction for $ \Mod_{g,[\gamma],\beta}'$ would yield a form $\hat\eta'_{g,[\gamma],\beta}$ which restricts to the correct class in $\rh^{*}_{c}(\ex X)$. 

\stop

\begin{example}[Gromov Witten invariants of a compact symplectic manifold]
\end{example}

Suppose that $(M,\omega)$ is a compact symplectic manifold. We may choose a smooth almost complex structure $J$ on $M$ compatible with $\omega$. It is proved in \cite{cem} that $M,J,\omega$ satisfies Gromov compactness in our sense. Saying that a curve in $M$ has tropical part in  the isotopy class $[\gamma]$ is equivalent to labeling the punctures of that curve, so we may use the notation $\Mod_{g,n,\beta}$ for $\Mod_{g,[\gamma],\beta}$ where $\gamma$ has $n$ infinite ends. In this case, $\End_{\gamma}M$ is equal to $ M^{n}$. Then $\totl{ev^{0}}\times\EV$ gives a 
$\C\infty1$ map 
\[\Mod_{g,n,\beta}\longrightarrow \bar{\mathcal M}_{g,n}\times M^{n}\]
where $\bar{\mathcal M}_{g,n}$ is Deligne Mumford space. The Gromov Witten invariants of $M$ are given by  integrating closed  forms on $\bar{\mathcal M}_{g,n}\times M^{n}$ over $\Mod_{g,n,\beta}$, which is the same as integrating forms on $\bar{\mathcal M}_{g,n}\times M^{n}$ against a closed form $\eta_{g,n,\beta}$ on $\bar{\mathcal M}_{g,n}\times M^{n}$. The cohomology class of $\eta_{g,n,\beta}$ for all $(g,n,\beta)$ encapsulates all Gromov Witten invariants of $(M,\omega)$.  

Theorem \ref{family invariance} implies that these Gromov Witten invariants are independent of our choice of $J$. In particular, let $\eta'_{g,n,\beta}$ be the form constructed using another choice  $J'$ of compatible almost complex structure. Then there is a smooth family $J_{t}$  of compatible almost complex structures  on $M$ for $t\in\mathbb R/2\pi$ so that $J=J_{0}$ and $J'=J(1)$. It is noted in Appendix \ref{compactness} that Gromov compactness holds for such a family. Theorem \ref{family invariance} tells us that there is a closed differential form 
 on $\bar{\mathcal M}_{g,n}\times M^{n}\times[0,1]$ which restricts at one boundary to have the same cohomology class as $\eta_{g,n,\beta}$, and which restricts to the other boundary to have the same cohomology class as $\eta'_{g,n,\beta}$. It follows that $\eta_{g,n,\beta}$ and $\eta'_{g,n,\beta}$ represent the same cohomology class. 

\

Now suppose that we have a connected basic family of almost complex exploded manifolds $(\hat {\ex B},J)\longrightarrow \ex G$ with a family of taming forms so that Gromov compactness holds. Suppose further our manifold $(M,J)$ is one fiber of our family. 

Theorem \ref{family invariance} implies that Gromov Witten invariants of $(M,J)$ may be calculated from the Gromov Witten invariants of other fibers of $\hat{\ex B}$. In particular, given any two points in a coordinate chart of $\ex G$, there exists a map of a complete exploded manifold $\ex G'$ to $\ex G$  with image containing our two points, and contained in a subset of that coordinate chart with bounded tropical part. The fiber product $\hat {\ex B'}$ of $\hat{\ex B}\longrightarrow \ex G$ with $\ex G'\longrightarrow \ex G$ is then a complete family for which Gromov compactness holds. As $\hat{\ex B}'$
 is pulled back from a subset of a single coordinate chart on $\ex G$ with bounded tropical part, the proof of invariance of cohomology in families from \cite{dre} implies that we may canonically identify the cohomology of each fiber, and  that the map from $H^{*}(\hat{\ex B}')$ to $H^{*}$ of any fiber is surjective and has a fixed kernel. Therefore, the homology of  fibers includes canonically into the homology of $\hat{\ex B}'$, so we may use the same notation $\beta$ for a homology class in any fiber or its pushforward to the cohomology of $\hat{\ex B}'$. Similarly, as $\hat {\ex B}'$ is pulled back from $\hat{\ex B}$, a connected family which contains one fiber with tropical part equal to a point, the isotopy classes of tropical curves contained in fibers of $\hat {\ex B'}$ will always just be determined by the number of ends, so we may continue using the notation $\Mod_{g,n,\beta}$ in place of $\Mod_{g,[\gamma],\beta}$.

 Then Theorem \ref{family invariance} tells us that there exists a closed form \[\hat \eta_{g,n,\beta}\in \ro^{*}_{c}(\bar{\mathcal M}_{g,n}\times \sfp {\hat{\ex B}'}{\ex G'}n)\] which restricts to each fiber to have the same cohomology class as $\eta_{g,n,\beta}$.

It is proved in \cite{dre} that $H^{*}$ cohomology does not change in connected families, and that the identification is locally given by identifying the restriction of forms in $H^{*}$ of the total space. It follows that with this natural identification, the restriction of $\hat \eta_{g,n,\beta}$ to any fiber represents the same element of the dual of $H^{*}$ of the fiber.

 Therefore, Gromov Witten invariants of $M$ can be calculated in any fiber of $\hat{\ex B}$. (Note that this identification of cohomology of the fiber and therefore of Gromov Witten invariants is only locally canonical - if $\totl{\ex G}$ is not simply connected, then different paths in $\totl{\ex G}$ may correspond to different identifications.) In the next section, we shall derive some gluing formulas which make the calculation of Gromov Witten invariants in exploded manifolds with nontrivial tropical part easier.

\section{Gluing relative invariants}\label{gluing relative invariants}

\

\subsection{Tropical completion}

\

\

In order to state gluing theorems for Gromov Witten invariants, we shall need the notion of tropical completion.

\begin{defn}[Tropical completion in a coordinate chart] 
The tropical completion of a strata $S$ of a polytope $P$ in $\mathbb R^{m}$ is a polytope $\check P_{S}\subset \mathbb R^{m}$ which is the union of all rays in $\mathbb R^{m}$ which begin in $S$ and intersect $P$ in more than one point.

Given a coordinate chart $U=\mathbb R^{n}\times\et mP$ and a strata $S$ of $P$, if $U_{S}$ indicates the strata of $U$ corresponding to $S$, then the tropical completion of $U_{S}\subset U$ is defined to be
\[\check U_{S}:=\mathbb R^{n}\times \et m{\check P_{S}}\]

\end{defn}

Tropical completion in coordinate charts is functorial:

In particular, given a map $\totb f:P\longrightarrow Q$ sending a strata $S\subset P$ into a strata $S'\subset Q$, then there is a unique map $\totb{\check f}: \check P_{s}\longrightarrow \check Q_{s} $ which restricted to $P\subset \check P_{S}$ is equal to $\totb f$.

Similarly,  given a smooth or $\C\infty1$ map of coordinate charts
 \[f:U\longrightarrow U'\] sending  a strata $U_{S}\subset U$ into  $U'_{S'}\subset U'$ there is a unique map 
\[\check f:\check U_{S}\longrightarrow \check U'_{S'}\]
so that $\check f$ restricted to $U_{S}$ is equal to $f$. Of course, the tropical part of $\check f$ is equal to the map $\totb {\check f}$ above. This map $\check f$ is smooth or $\C\infty1$ if $f$ is, and is an isomorphism onto an open subset of $\check U'_{S'}$ if $f$ is an isomorphism onto an open subset of $U'$. There is therefore a functorial construction of the tropical completion $\check {\ex B}_{S}$ of any strata $\ex B_{S}$ of an exploded manifold $\ex B$ as follows:

\begin{defn}[tropical completion]\label{tropical completion}The tropical completion of a strata $\ex B_{S}\subset \ex B$ of an exploded manifold is an exploded manifold $\check{\ex B}_{S}$ constructed as follows: Let $\{(U,U_{S})\}$ be the set of all possible coordinate charts $U$ of $\ex B$ which send a strata $U_{S}$ to $\ex B_{S}$, and let $\{\phi\}$ indicate the set of all possible inclusions $(U,U_{S})\longrightarrow (U',U'_{S'})$ of these coordinate charts. Then $\check {\ex B}_{S}$ is an exploded manifold with $\{\check U_{S}\}$ as the set of all possible coordinate charts on $\check {\ex B}_{S}$, and $\{\check \phi\}$ as the set of all possible inclusions of these coordinate charts on $\ex B$.

\end{defn}

Note that $\check {\ex B}_{S}$ always contains a copy of $\ex B_{S}$ as a dense subset. If $\ex B$ is basic, then $\check {\ex B}_{S}$ also contains a copy of the closure $\bar{\ex B}_{S}$ of $\ex B_{S}\subset\ex B$.
The construction of tropical completions is functorial in the sense that given a map
\[f:\ex A\longrightarrow \ex B\]
which sends a strata $\ex A_{S}$ to $\ex B_{T}$, there is a unique map 
\[\check f:\check{\ex A}_{S}\longrightarrow \check{\ex B}_{T}\]
which restricts to be $f$ on $\ex A_{T}\subset \check{\ex A}_{T}$.

Given any tensor $\theta$ on $\ex B$ such as an almost complex structure, metric, or differential form, there is a unique tensor $\check \theta$ on $\check{\ex B}_{S}$ which restricts to be $\theta$ on ${\ex B_{S}}$. Therefore, if $\ex B$ is an almost complex exploded manifold, we may talk about holomorphic curves in $\check {\ex B}_{S}$.

\subsection{$\gamma$-decoration}

\

\

A tropical curve $\gamma$ in $\totb{\ex B}$ is a continuous map with domain equal to a complete graph so that each edge of the graph has an integral affine structure and $\gamma$ restricted to each edge is an integral affine map.
Let $\gamma$ be a tropical curve and let $f:\ex C\longrightarrow \ex B$ be a $\C\infty1$ curve with tropical part equal to $\gamma$.  Label the strata of $\ex C$ corresponding to a vertex $v$ of $\gamma$  by $\ex C_{v}$, and the strata corresponding to an edge $e$ of $\gamma$ by $\ex C_{e}$. Indicate by ${\ex B}_{v}$ or $\ex B_{e}$ the  strata of $\ex B$ which  $f$ sends  $\ex C_{v}$ or $\ex C_{e}$ to. Taking tropical completions gives maps
\[\check f_{v}:\check{\ex C}_{v}\longrightarrow \check{\ex B}_{v}\]
\[\check f_{e}:\check {\ex C}_{e}\longrightarrow \check{\ex B}_{e}\]
Use the notation $\gamma_{v}$ for the tropical part of $\check f_{v}$ and $\gamma_{e}$ for the tropical part of $\check f_{e}$.
Note that $f$ is holomorphic if and only if $\check f_{v}$ is holomorphic for all vertices $v$ of $\gamma$. 


If the edge $e$ is adjacent to $v$, then $\check{\ex C}_{v}$ contains $\ex C_{e}$ and if $v$ is at both ends of $e$, $\check{\ex C}_{v}$ contains two copies of $\ex C_{e}$. In general, if $e$ is an internal edge of $\gamma$, there are two copies of $\ex C_{e}$ in the union of all $\check{\ex C}_{v}$. To distinguish these two copies, choose an orientation for each internal edge of $\gamma$. Let $e_{0}$ denote the strata of $\check{\ex C}_{v}$ which corresponds to $v$ being an incomming end of the edge $\gamma$, and let $e_{1}$ be the strata of $\check{\ex C}_{v}$ which corresponds to $v$ being an outgoing end of the edge $e$.

 The map  $\check f_{e}$ can be obtained from  $\check f_{v}$ by taking the tropical completion of $\check f_{v}$ using the strata $e_{1}$ or $e_{0}$. As $\check {\ex C}_{e}$ is equal to $\ex T$, the moduli space $\mathcal M_{e}$ of possibilities for $\check f_{e}$ is finite dimensional and equal to the quotient of $\check{\ex B}_{e}$ by a $\ex T$ action. Even though this $\ex T$ action is not be free when $\check f_{e}$ is not injective, we shall simply treat $\mathcal M_{e}$ as an exploded manifold instead of as a stack. (There will be some adjustments to fiber products over $\mathcal M_{e}$ that we shall need to do below because of this simplification.) Use the notation $\mathcal M_{\gamma_{e}}$ for the restriction of $\mathcal M$ to curves with tropical part $\gamma_{e}$.

\begin{defn} Denote by $\Msw_{\gamma}(\ex B)$ the moduli stack of $\C\infty1$ curves in $\ex B$ with an isomorphism between their tropical part and $\gamma$, similarly, denote the substack of $\Msw_{\gamma}(\ex B)$ consisting of curves with genus $g$ and homology class $\beta$ by $\Msw_{g,\gamma,\beta}$. 

Use the notation $\Msw_{\gamma_{v}}$ for $\Msw_{\gamma_{v}}(\check{\ex B}_{v})$, and $\Msw_{\gamma}$ for $\Msw_{\gamma}(\ex B)$. 
\end{defn}

Note that in the case that $\gamma$ has automorphisms,  a map of families in $\Msw_{\gamma}$ must be compatible with the identification of tropical parts with $\gamma$.  If $\gamma$ has automorphism group of size $k$, then $\Msw_{\gamma}$ should be thought of as a $k$-fold cover of the  substack $\Msw(\ex B)\rvert_{\gamma}$ of $\Msw(\ex B)$ consisting of curves with tropical parts which are isomorphic to $\gamma$.
 
 \
 
 \begin{defn}[$\gamma$-decoration]
 Consider the domain of  $\gamma$ as a graph with an affine structure on the edges, and consider the tropical part of the domain of a family of curves in $\Msw(\ex B)$ as a union of affine polytopes glued along faces. (This forgets the integral part of the integral affine structure on the tropical part of the domain.) Define a $\gamma$-decorated tropical curve to be a tropical curve in $\ex B$ with a continuous affine map of its domain to the domain of $\gamma$ which is a homeomorphism restricted to the inverse image of the interior of all edges of $\gamma$ and which is an integral affine isomorphism restricted to exterior edges of $\gamma$.   
 
 Define a $\gamma$-decorated curve to be a curve with a $\gamma$ decorated tropical part.  Consider the $\gamma$-decorated moduli space of  $\C\infty1$ families of curves $\hat f\in \Msw(\ex B)$ with an affine map of the tropical part of the domain of $\hat f$ to the domain of $\gamma$ which makes each individual curve into a $\gamma$-decorated curve.   $\Msw_{\gamma}$ is a substack of this $\gamma$-decorated stack of curves. Define  $\bMs_\gamma$ to be the closure of $\Msw_{\gamma}$ in this $\gamma$-decorated
stack of curves. Similarly, define $\bMs_{g,\gamma,\beta}$ to be the closure of $\Msw_{g,\gamma,\beta}$ in this $\gamma$-decorated stack of curves.
\end{defn}

$\bMs_\gamma$ can be thought of as a kind of multiple cover of the substack of $\Msw$ consisting of the closure of all the strata containing curves with tropical part $\gamma$. If $\gamma$ has no internal edges, then $\bMs_{g,\gamma,\beta}=\Msw_{g,[\gamma],\beta}$. In particular, $\bMs_{g,\gamma_{v},\beta}=\Msw_{g,[\gamma_{v}],\beta}$.

Forgetting the $\gamma$-decoration gives a map 

\[\bMs_\gamma\longrightarrow \Msw\]

  Given any curve in $\Msw_{\gamma}$, taking tropical completions for all vertices $v$ in $\gamma$ give curves in $\Msw_{\gamma_{v}}$ which in turn give curves in $\mathcal M_{\gamma_{e}}$ which agree with each other. When $\ex B$ is basic, this tropical completion map can be extended in an obvious way to a map 
  \[\bMs_\gamma\longrightarrow \bMs_{\gamma_{v}}\]

  In particular for a family $\hat f\in \bMs_\gamma$,  consider the closure of the collection of strata in the domain of $\hat f$ that are sent to $v$. Our map $\hat f$ restricted to this subset lands in the closure $\bar{\ex B}_{v}$ of $\ex B_{v}$ inside $\ex B$, which includes in $\check {\ex B}_{v}$ because $\ex B$ is basic. We can then canonically extend this restricted map to a $\C\infty1$ family of curves $\hat f_{v}$ in $\check{\ex B}_{v}$ by extending all the edges of the tropical part of curves in $\hat f$ with only one end to have infinite length. To distinguish this map from the tropical completion map defined earlier, call this process $\gamma$-decorated tropical completion. 
  
  We can apply $\gamma$-decorated tropical completion to appropriate $\gamma$-decorated tropical curves analogously: Suppose that $\gamma'$ is a $\gamma$-decorated tropical curve so that all strata of $\gamma'$ which are sent to $v\in\gamma$ have their image inside $\totb{\bar{\ex B}_{v}}$. (Note that while $\bar{\ex B}_{v}$ is the closure of the strata $\ex B_{v}\subset \ex B$, $\totb{\bar{\ex B}_{v}}$ consists of all strata of $\totb{\ex B}$ which contain $\totb{\ex B_{v}}$ in their closures using the topology on $\totb{\ex B}$.) Then $\gamma'$ restricted to the inverse image of $\totb {\bar{\ex B}_{v}}$ intersected with all strata of $\gamma'$ which are sent to $v$ or an edge adjacent to $v$ is a tropical curve in $\totb{\bar{\ex B}_{v}}\subset\totb{\check{\ex B}_{v}}$ which has some edges which have finite length, but are only attached to a vertex at one end. Increasing the length of these edges to be infinite gives $\gamma'_{v}$, the $\gamma$-decorated tropical completion of $\gamma'$.

   Similarly, if $\ex B$ is basic, then given any edge $e$ of $\gamma$, there is a map 
  \[\bMs_\gamma\longrightarrow \mathcal M_{e}\]
  extending the map 
  \[\Msw_{\gamma}\longrightarrow \mathcal M_{\gamma_{e}}\]
   Both of these maps are given by restricting curves in $\hat f$ to the inverse image of the edge $e$, which gives a family  of maps all of the form $\et 1{(0,l)}\longrightarrow \bar {\ex B}_{e}$. These can be extended canonically to maps $\ex T\longrightarrow \check{\ex B}_{e}$ which are in $\mathcal M_{e}$, giving a family $\hat f_{e} $ in $\mathcal M_{e}$. 
   
 For each inclusion of $v$ as an end of the edge $e$, we may associate an edge of $\gamma_{v}$ with $e$, and there is a map 
 \begin{equation}\label{bnm}\bMs_{\gamma_{v}}\longrightarrow{\mathcal M_{e}}\end{equation}
defined in exactly the same way as   the map $\bMs_\gamma\longrightarrow \mathcal M_{e}$. Of course, the result of concatenating the maps $\bMs_\gamma\longrightarrow \bMs_{\gamma_{v}}\longrightarrow \mathcal M_{e}$ is the map $\bMs_\gamma\longrightarrow \mathcal M_{e}$ defined above.

\subsection{Gluing theorems}

\

\

Choose an orientation on each of the internal edges of a tropical curve $\gamma$. Putting together all the above maps from (\ref{bnm}) for the ingoing ends of internal edges gives
   \[\EV_{0}:\prod_{v}\bMs_{\gamma_{v}}\longrightarrow \prod_{e}\mathcal M_{e}\]
   and using the outgoing ends of internal edges gives
   \[\EV_{1}:\prod_{v}\bMs_{\gamma_{v}}\longrightarrow \prod_{e}\mathcal M_{e}\]

 A curve in $\prod_{v}\bMs_{\gamma_{v}}$ is a curve  $\check f_{v}$ in  $\bMs_{\gamma_{v}}$ for all $v$. In order for these curves to glue together to a curve $ f$ in $\bMs_{\gamma}$, it is necessary that $\EV_{0}$ and $\EV_{1}$ agree on these curves, and it is also necessary that the tropical parts $\totb {\check f_{v}}$ glue together to form a ($\gamma$ decorated) tropical curve $\totb f$ in $\totb{\ex B}$ so that applying $\gamma$-decorated tropical completion to $\totb f$ gives $\totb{\check f_{v}}$ for all $v$. Two things might go wrong  with this tropical gluing if $\EV_{0}$ and $\EV_{1}$ agree: The vertices of $\totb{\check f_{v}}$ might not all be inside the image of $\bar{\ex B}_{v}\subset \check {\ex B}$, or two vertices that must be glued together might require an edge of negative length in order to join them.
   
%

\

\

 \begin{lemma} \label{fp bundle} Suppose that 
 \begin{itemize}
 \item  $\hat f$ is a $\C\infty1$ family of curves in $\prod_{v}\bMs_{\gamma_{v}}$ which represents a substack of $\prod_{v}\bMs_{\gamma_{v}}$.
 \item  $EV_{0}$ and $EV_{1}$ are identical restricted to $\hat f$. 
 \item The tropical part of every curve in $\hat f$ may be obtained by applying $\gamma$-decorated tropical completion to some $\gamma$-decorated tropical curve.
 \end{itemize}
 Then the inverse image of $\hat f$ in $\bMs_{\gamma}$ is a substack of $\bMs_{\gamma}$ represented by a  $\C\infty1$ family $\hat f'$, and the map
   
\[\ex F(\hat f')\longrightarrow \ex F(\hat f)\]  
is a proper submersion with fibers equal to the product of $\et 1{(0,\infty)}$ for each internal edge of $\gamma$ sent to a point, and a set with $m_{e}$ elements for each internal edge $e$ of $\gamma$ with nonzero multiplicity $m_{e}$.

  \end{lemma}

\pf

In what follows, we will construct $\hat f'$. For each vertex $v$ of $\gamma$, we have a family $\hat f_{v}$ of curves in $\bMs_{\gamma_{v}}$ parametrized by $\ex F(\hat f)$. We must glue all these families together over the edges corresponding to internal edges of $\gamma$. For each internal edge $e$ of $\gamma$, let $\ex C_{e_{0}}$ and $\ex C_{e_{1}}$ indicate the subsets of the appropriate $\ex C(\hat f_{v})$ corresponding to the edges $e_{0}$ and $e_{1}$. (Each of these $\ex C_{e_{i}}$ is a $\et 1{(0,\infty)}$ bundle over $\ex F(\hat f)$.) $\hat f$ applied to $\ex C_{e_{i}}$ gives a commutative diagram 
\[\begin{array}{ccc}\ex C_{e_{i}}&\longrightarrow &\check {\ex B}_{e}
\\\downarrow &&\downarrow
\\ \ex F(\hat f)&\longrightarrow &\mathcal M_{e}\end{array}\]

 where the map $\ex F(\hat f)\longrightarrow \mathcal M_{e}$ is given by $\EV_{i}$ followed by projection to $\mathcal M_{e}$.  In the case that the edge $e$ is not sent to a point, consider $\check {\ex B}_{e}$ as a $\ex T$ bundle over $\mathcal M_{e}$, and pull back this bundle over this map $\ex F(\hat f)\longrightarrow \mathcal M_{e}$ to obtain a $\ex T$ bundle $\check{\ex B}'_{e}$ over $\ex F(\hat f)$. Our commutative diagram then factors into a bundle map $\hat f_{e_{i}}$ over $\ex F(\hat f)$ and a pullback diagram:
 \[\begin{array}{ccccc}\ex C_{e_{i}}&\xrightarrow{\hat f_{e_{i}}}&\check{\ex B}'_{e}&\longrightarrow &\check {\ex B}_{e}
\\\downarrow &&\downarrow&&\downarrow
\\ \ex F(\hat f)&\xrightarrow{\id}&\ex F(\hat f)&\longrightarrow  &\mathcal M_{e}\end{array}\]

 This bundle map looks in local coordinates like \[\hat f_{e_{i}}(p,\tilde z_{i})= \lrb{p,f_{i}(p)\tilde z_{i}^{(-1)^{i}m_{e}}}\] where $m_{e}$ is the multiplicity of the edge $e$. In the case that our edge $e$ is sent to a point, $\mathcal M_{e}=\check{\ex B}_{e}$, $\check {\ex B}_{e}'=\ex F(\hat f)$ and the equivalent of $\hat f_{e_{i}}$ is just the projection to $\ex F(\hat f)$.

Consider the fiber product of $\ex C_{e_{0}}$ and $\ex C_{e_{1}}$ over the maps $\hat f_{e_{i}}$. (Recall that we are assuming  that $\EV_{0}$ and $\EV_{1}$ coincide on $\ex F(\hat f)$, so $\check{\ex B}'_{e}$ does not depend on $i$.)
Given any family of curves $\hat g$ inside the inverse image of $\hat f$ in $\bMs_{\gamma}$, the fact that $\hat f$ represents a substack implies that there is a unique map from the image of $\hat g$ in $\prod_{v}\bMs_{\gamma_{v}}$ to $\hat f$, which gives a canonical $\C\infty1$ map 
$\ex F(\hat g)\longrightarrow \ex F(\hat f)$. For each internal edge $e$ of $\gamma$, this map lifts naturally to a map $\ex C(\hat g)_{e}\longrightarrow \ex C_{e_{0}}$ and a map $\ex C(\hat g)_{e}\longrightarrow \ex C_{e_{1}}$. As $\hat f_{e_{i}}$ composed with each of these maps is equal to $\hat g$, this gives a natural map $\ex C(\hat g)_{e}\longrightarrow \ex C_{e_{0}}\fp{\hat f_{e_{0}}}{\hat f_{e_{1}}}\ex C_{e_{1}}$. 

Take the quotient of $\ex C_{e_{0}}\fp{\hat f_{e_{0}}}{\hat f_{e_{1}}}\ex C_{e_{1}}$ by the equivalence relation $(p,\tilde z_{0},\tilde z_{1})\cong(p,c\tilde z_{0},c^{-1}\tilde z_{1})$. The result of this quotient is an exploded manifold $\ex F_{e}$ with a proper submersion $\ex F_{e}\longrightarrow \ex F(\hat f)$ with fibers equal to $\et 1{(0,\infty)}$ if the edge $e$ is sent to a point, and fibers equal to a set with $m_{e}$ elements if $e$ is an edge with multiplicity $m_{e}$. (This uses our assumptions that $\EV_{0}$ and $\EV_{1}$ coincide and that the tropical part of every curve in $\hat f$ is in the image of  $\gamma$-decorated tropical completion.) Our map $\ex C(\hat g)_{e}\longrightarrow \ex C_{e_{0}}\fp{\hat f_{e_{0}}}{\hat f_{e_{1}}}\ex C_{e_{1}}$ induces a natural map $\ex F(\hat g)\longrightarrow \ex F_{e}$ which lifts our map $\ex F(\hat g)\longrightarrow \ex F(\hat f)$. Define $\ex F(\hat f')$ to be the fiber product of all $\ex F_{e}$ over $\ex F(\hat f)$. So far, we have that $\ex F(\hat f')\longrightarrow \ex F(\hat f)$ is a proper submersion with fibers as specified in our lemma, and given any family $\hat g$ in the inverse image of $\hat f$ in $\bMs_{\gamma}$, we have constructed a canonical $\C\infty1$ map $\ex F(\hat g)\longrightarrow\ex F(\hat f')$.

Let us now construct $\ex C(\hat f')$ and the map $\hat f'$ itself. 
Let $\gamma_{v}'$ indicate  subset of $\gamma$ obtained by removing the middle $1/3$ of each edge of $\gamma$ and taking the connected component containing $v$. This $\gamma_{v}'$ can be considered as a subset of both $\gamma$ and $\gamma_{v}$. The restriction $\hat f'_{\gamma_{v}'}$ of $\hat f'$ to the inverse image of $\gamma_{v}'$ shall be equal to the restriction to the inverse image of  $\gamma_{v}'$ of the pullback of $\hat f_{v}$ over the map $\ex F(\hat f')\longrightarrow \ex F(\hat f)$.

We must now describe the restriction $\hat f'_{e}$ of $\hat f'$ to each internal edge $e$, and specify how $\hat f'_{e}$ is attached to $\hat f'_{\gamma_{v}'}$. The domain of $\hat f'_{e}$ shall be the pullback of  $\ex C_{e_{0}}\fp{\hat f_{e_{0}}}{\hat f_{e_{1}}}\ex C_{e_{1}}\longrightarrow \ex F_{e}$ over the map $\ex F(\hat f')\longrightarrow \ex F_{e}$, so the following is a pullback diagram:
\[\begin{array}{ccc}\ex C(\hat f'_{e})&\longrightarrow & \ex C_{e_{0}}\fp{\hat f_{e_{0}}}{\hat f_{e_{1}}}\ex C_{e_{1}}
\\\downarrow &&\downarrow
\\ \ex F(\hat f')&\longrightarrow& \ex F_{e}
\end{array}\]
 The pullback of $\hat f$ over the map $\ex C_{e_{0}}\fp{\hat f_{e_{0}}}{\hat f_{e_{1}}}\ex C_{e_{1}}\longrightarrow \ex C_{e_{0}}$ is equal to the pullback of $\hat f$ over the map $\ex C_{e_{0}}\fp{\hat f_{e_{0}}}{\hat f_{e_{1}}}\ex C_{e_{1}}\longrightarrow \ex C_{e_{1}}$, and can be considered to be a map to $\ex B$ instead of $\check{\ex B}_{e}$ because of our assumption that the tropical part of $\hat f$ is in the image of $\gamma$ decorated tropical completion. Let $\hat f'_{e}:\ex C(\hat f'_{e})\longrightarrow \ex B$ be the pullback of this map over $\ex F(\hat f')\longrightarrow \ex F_{e}$. 

 We can attach $\hat f_{e}'$ to $\hat f'_{\gamma_{v}}$ at the end of $e$ corresponding to $e_{i}$ by considering $\hat f_{e}'$ as pulled back from the map $\ex C_{e_{0}}\fp{\hat f_{e_{0}}}{\hat f_{e_{1}}}\ex C_{e_{1}}\longrightarrow \ex C_{e_{i}}\subset \ex C(\hat f_{v})$, so the definition of $\hat f'$ restricted to the corresponding $1/3$ of the edge of $e$ is the same, wheather $\hat f'_{e}$ or $\hat f'_{\gamma'_{v}}$ is used. 

We have now constructed a family $\hat f'$ in the inverse image of $\hat f$  in $\bMs_{\gamma}$. Reconsider our family $\hat g$ of curves in the inverse image of $\hat f$ in $\bMs_{\gamma}$. We must show that there exists a unique map $\hat g\longrightarrow \hat f'$. As $\hat f$ is a substack,  we have a unique map from the $\gamma$-decorated tropical completion of  $\hat g$ to $\hat f$. Any map $\hat g\longrightarrow \hat f'$ must be compatible with the corresponding unique map of the $\gamma$-decorated tropical completion of $\hat f'$ into $\hat f$. The maps $\hat g_{v}\longrightarrow \hat f_{v}$ give for each internal edge $e$ of $\gamma$ a map $\ex C(\hat g)_{e}\longrightarrow \ex C_{e_{0}}$ and $\ex C(\hat g)_{e}\longrightarrow\ex C_{e_{1}}$ which together define a map $\ex C(\hat g)_{e}\longrightarrow\ex C_{e_{0}}\fp{\hat f_{e_{0}}}{\hat f_{e_{1}}}\ex C_{e_{1}}$ so that $\hat g$ restricted to $\ex C(\hat g)_{e}$ is the pullback of $\hat f$ under the above map followed by  projection to $\ex C_{e_{i}}$. Any map $\hat g\longrightarrow \hat f$ must give a commutative diagram

\[\begin{array}{ccc}\ex C(\hat g)_{e}&\longrightarrow &\ex C(\hat f)_{e}
\\ \downarrow&& \downarrow 
\\ \ex C_{e_{0}}\fp{\hat f_{e_{0}}}{\hat f_{e_{1}}}\ex C_{e_{1}}&\xrightarrow{\id}&\ex C_{e_{0}}\fp{\hat f_{e_{0}}}{\hat f_{e_{1}}}\ex C_{e_{1}}\end{array}\]

Therefore the fiber product of the corresponding maps $\ex F(\hat g)\longrightarrow \ex F_{e}$ over $\ex F(\hat f)$ specify the only possible map $\ex F(\hat g)\longrightarrow \ex F(\hat f')$ which can come from a map $\hat g\longrightarrow\hat f$. The maps $\ex C(\hat g)_{e}\longrightarrow \ex F(\hat g)\longrightarrow \ex F(\hat f')$ and $\ex C(\hat g)_{e}\longrightarrow\ex C_{e_{0}}\fp{\hat f_{e_{0}}}{\hat f_{e_{1}}}\ex C_{e_{1}}$ then specify a map $\ex C(\hat g)_{e}\longrightarrow \ex C(\hat f')_{e}$, so that $\hat g$ restricted to $\ex C(\hat g)_{e}$ is equal to the pullback of $\hat f'$. Similarly, the maps $\ex C(\hat g)_{\gamma'_{v}}\longrightarrow \ex F(\hat g)\longrightarrow \ex F(\hat f')$  and $\ex C(\hat g)_{\gamma'_{v}}\longrightarrow \ex C(\hat  f_{v})$ specify a map $\ex C(\hat g)_{\gamma'_{v}}\longrightarrow \ex C(\hat f')$ so that the pullback of $\hat f'$ is $\hat g$. As with the definition of $\hat f'$, the construction of this map on an edge $e$ is identical to the construction coming from the vertices at either side, so we have a unique $\C\infty1$ map $\hat g\longrightarrow \hat f'$.

\stop

   \begin{thm} \label{gluing cobordism} Let $\gamma$ be a tropical curve in $\ex B$ with genus $g_{\gamma}$, and suppose that Gromov compactness holds for $\ex B$ and $\check{\ex B}_{v}$ for all vertices $v$ of $\gamma$. Choose a genus $g$ and an energy $E$. Then  we may construct the virtual moduli spaces of holomorphic curves
    $\Mod_{g_{v},[\gamma_{v}],E_{v}}(\check {\ex B}_{v})$ for all vertices $v$ of $\gamma$ and $g_{v}\leq g-g_{\gamma}$ and $E_{v}\leq E$ so that the following holds:
   \begin{itemize} \item The maps $\EV_{0}$ and $\EV_{1}$ are transverse applied to  $\prod_{v}\Mod_{g_{v},[\gamma_{v}],E_{v}}$ whenever $\sum_v g_{v}=g-g_{\gamma}$ and $\sum_v E_{v}=E$.
   \item The pullback of the virtual moduli space of holomorphic curves in $\Msw(\ex B)$ to $\bMs_{g,\gamma,E}$  is cobordant to the pullback to $\bMs_{g,\gamma,E}$ of the virtual moduli space of curves in $\prod_{v}\bMs_{\gamma_{v}}$.
   \end{itemize}
If holomorphic curves in $\bMs_{g,\gamma,E}$ have the same number of automorphisms as their image in $\Msw_{g,[\gamma],E}$, the $\Mod_{g_{v},[\gamma_{v}],E_{v}}(\check{\ex B}_{v})$ are constructed using the zero perturbation and  $\EV_{0}$ and $\EV_{1}$ are transverse, then the virtual moduli space of holomorphic curves in $\ex B$ may be constructed so that the pullbacks of the two different virtual moduli spaces to $\bMs_{g,\gamma,E}$ are equal.
   \end{thm}

The proof of this theorem is contained in appendix \ref{gluing cobordism proof}. Note that the inverse image in $\bMs_{\gamma}$ of  any curve in $\prod_{v}\Msw_{g_{v},[\gamma_{v}],E_{v}}$ has genus $\sum_v{g_{v}}+g_{\gamma}$ and energy $\sum_{v} E_{v}$ so $EV_{0}$ and $EV_{1}$ are transverse on the part of the virtual moduli space which pulls back to $\bMs_{g,\gamma,E}$.

\

Note that when $\ex B$ is a symplectic manifold, Theorem \ref{gluing cobordism} applied to tropical curves with a single internal edge implies the splitting and genus reduction axioms of Kontsevich and Manin  stated in \cite{KM}. These splitting and genus reduction axioms require using Gromov Witten invariants defined using $(\totl {ev^{0}},EV):\Mod_{g,n,\beta}\longrightarrow \overline{\mathcal M}_{g,n}\times \ex B^{n}$.  We shall first prove a  gluing theorem for Gromov Witten invariants defined using only $EV$, then prove a generalization of the splitting and genus reduction axioms in Theorem \ref{splitting}.

\begin{thm}\label{gluing formula} Let $\ex B$ be a basic exploded manifold for which Gromov compactness holds. Let  $\gamma$ be a tropical curve in  $\totb{\ex B}$ with genus $g_{\gamma}$, and let
$\theta\in \rof^{*}_{c}(\End_{[\gamma]} \ex B)$ be a closed differential form. Suppose further that for all vertices $v$ of $\gamma$, Gromov compactness holds for $\check {\ex B}_{v}$. Then
\[\sum_{g,E}\hbar^{g}\lambda ^{E}\int_{\Mod_{g,[\gamma],E}\rvert_{\gamma}}\EV^{*}\theta=k_{\gamma}\hbar^{g_{\gamma}}\int_{\prod_{e}\mathcal M_{\gamma_{e}}}\theta \bigwedge_{v\in \gamma}\eta_{v} \]

where 
\begin{itemize}
\item The equality is in the ring of formal series $\sum_{g,E} c_{g,E}\hbar^{g}\lambda^{E}$ where $g\in \mathbb N$, $E\in [0,\infty)$, so that given any bounded subset of $\mathbb N\times [0,\infty)$, there are only a finite number of nonzero coefficients $c_{g,E}\in \mathbb R$ with $(g,E)$ in this bounded subset.
\item $k_{\gamma}$ is $0$ if $\gamma$ has any internal edges which are sent to points, and otherwise 
$k_{\gamma}$ is the product of the multiplicity of the internal edges of $\gamma$ divided by the number of automorphisms of $\gamma$ as a tropical curve with labeled ends.
\item $\prod_{e}\mathcal M_{\gamma_{e}}$ indicates the product over each edge $e$ of $\gamma$, of the moduli space of curves with tropical part equal to $e$ (considered as a manifold, not an orbifold when $e$ has multiplicity greater than $1$).
\item The $\theta$ on the right hand side indicates the pull back of $\theta$ over the map $\prod_{e}\mathcal M_{\gamma_{e}}\longrightarrow \End_{[\gamma]} \ex B$ which is independent of $\mathcal M_{\gamma_{e}}$ for each internal edge $e$ of $\gamma$, and which is the product of the inclusions $\mathcal M_{\gamma_{e}}\subset \End{\ex B}$ for each of the external edges $e$  of $\gamma$.
\item For each vertex $v$ of $\gamma$, $\eta_{v}$ is a formal series with coefficients which are  differential forms on $\prod_{e} \mathcal M_{\gamma_{e}}$ constructed as follows:

From $\bMs_{\gamma_{v}}$, there is an evaluation map \[\EV_{v}:\bMs_{\gamma_{v}}\longrightarrow \End_{[\gamma_{v}]}{\check{\ex B}_{v}}= \prod_{e'\subset\gamma_{v}}\mathcal M_{e'}\]
 Let $\eta_{g,[\gamma_{v}],E}$ be the Poincare dual in the sense of Theorem \ref{pd class} to $\EV_{v}$ applied  $\Mod_{g,[\gamma_{v}],E}$. This is a refined differential form on the product of $\mathcal M_{e'}$ for each edge $e'$ of $\gamma_{v}$. Then $\eta_{v}$ is the pullback of $\sum_{g,E}\hbar^{g}\lambda^{E}\eta_{g,\gamma_{v},E}$ over the product of the maps $\prod_{e}\mathcal M_{\gamma_{e}}\longrightarrow \mathcal M_{e'}$ which is the projection onto $\mathcal M_{\gamma_{e}}$ for the edge $e$ of $\gamma$ corresponding to $e'$ followed by the inclusion $\mathcal M_{\gamma_{e}}\subset\mathcal M_{e'}$.

\end{itemize}

\end{thm}

\pf Let us calculate $\int_{ \Mod_{g,[\gamma],E}\rvert_{\gamma}}\EV^{*}\theta$ using Theorem \ref{gluing cobordism}. 
 Note that the pullback of $ \Mod_{g,[\gamma],E}\rvert_{\gamma}$ to $\bMs_{\gamma}$ is an oriented $\abs{\Aut \gamma}$-fold multiple cover of $ \Mod_{g,[\gamma],E}\rvert_{\gamma}$, where $\Aut \gamma$ indicates the group of automorphisms of $\gamma$. (This arises from the difference between specifying that a curve has tropical part which is equal to $\gamma$, and specifying a particular isomorphism of the tropical part of a curve with $\gamma$.)

We may therefore exchange $\int_{ \Mod_{g,[\gamma],E}\rvert_{\gamma}}\EV^{*}\theta$ with the integral of $\abs{\Aut \gamma}^{-1}\EV^{*}\theta$ over the pullback of $ \Mod_{g,[\gamma],E}\rvert_{\gamma}$ to $\bMs_{\gamma}$. This space is equal to the restriction of the pullback of $ \Mod_{g,[\gamma],E}$ to $\bMs_{\gamma}$ to curves with tropical part gamma. As $\theta$ is generated by functions, we may apply the proof of Theorem \ref{well defined} to the cobordism from Theorem \ref{gluing cobordism}, and exchange this integral for an integral over the restriction to curves with tropical part $\gamma$ of the pullback to $\bMs_{g,\gamma,E}$ of the virtual moduli space of curves in $\prod_{v}\bMs_{\gamma_{v}}$. If $\gamma$ has some internal edges which are sent to a point, then this integral will be $0$ because $\EV^{*}\theta$ can not depend on the coordinate corresponding to the length of that edge. In the case that $\gamma$ has no internal edges sent to a point, Lemma \ref{fp bundle} implies that this space is an $m$-fold cover of the fiber product of the appropriate virtual moduli spaces in $\prod_{v}\Msw_{\gamma_{v}}$ over the maps $\EV_{1}$ and $\EV_{0}$, where $m$ is the product of the multiplicities of the internal edges of $\gamma$. This fiber product has an obvious evaluation map to the product of $\mathcal M_{\gamma_{e}}$ for all external edges $e$ of $\gamma$ which pulls back to $\bMs_{\gamma}$ be the evaluation map $\EV$ referred to above; call this new evaluation map $\EV$ as well.

In summary, so far we have that 

\[\int_{ \Mod_{g,[\gamma],E}\rvert_{\gamma}}\EV^{*}\theta=k_{\gamma}\int_{X}\EV^{*}\theta\]
where $k_{\gamma}=m\abs{\Aut\gamma}^{-1}$ and $X$ is the fiber product over  $\EV_{0}$ and $\EV_{1}$ of the virtual moduli space of holomorphic curves in $\prod_{v}\Msw_{\gamma_{v}}$ which have total genus $g-g_{\gamma}$ and total energy $E$. This fiber product $X$ can also  be constructed as follows:

For each vertex $v$ of $\gamma$, there is an evaluation map 
\[\EV_{v}:\bMs_{\gamma_{v}}\longrightarrow \prod_{e'\subset \gamma_{v}} \mathcal M_{e'}\]
The map we are calling $\EV$ is equal to $\prod_{v}\EV_{v}$ followed by projection to the product of those $\mathcal M_{e'}$ which correspond to external edges  of $\gamma$. The map $\EV_{1}$ is $\prod_{v}\EV_{v}$ followed by projection to the product of those $\mathcal M_{e'}$ which correspond to incoming ends of internal edges of $\gamma$, and $\EV_{0}$ is equal to $\prod_{v}\EV_{v}$ followed by projection to the product of those $\mathcal M_{e'}$ which correspond to outgoing ends of internal edges of $\gamma$. Let $\iota$ indicate the inclusion
\[\iota:\prod_{e\subset \gamma}\mathcal M_{e}\longrightarrow \prod_{v,e'\subset \gamma_{v}}\mathcal M_{e'}\]
which is the product of the identification of $\mathcal M_{e}$ with the appropriate $\mathcal M_{e'}$ for each external edge $e$ of $\gamma$, and the inclusion of $\mathcal M_{e}$ as the diagonal in $\mathcal M_{e_{1}}\times \mathcal M_{e_{0}}$ for each internal edge ${e}$ of $\gamma$. Then $X$ is also equal to the fiber product of the map $\iota$ with the map $\prod_{v}\EV_{v}$ applied to the virtual moduli space of curves in $\prod_{v}\Msw_{\gamma_{v}}$ with total genus $g-g_{\gamma}$ and total energy $E$.

Note that  $\theta$ is defined on the product of $\mathcal M_{\gamma_{e}}$, for each external edge $e$ of $\gamma$ but not on the product of $\mathcal M_{e}$. 
The image of $\EV$ applied to $ \Mod_{g,[\gamma],E}\rvert_{\gamma}$ is contained inside the subset of $\End_{[\gamma]} \ex B$ which is the product of $\mathcal M_{\gamma_{e}}$ for all external edges $e$ of $\gamma$. This subset is the inverse image of a point $p$ in the tropical part of $\End_{[\gamma]} \ex B$. As $\theta\in \rof^{*}_{c}(\End_{[\gamma]} \ex B)$, there is  a neighborhood  $U$ of $p$ and and a refinement $U'$ of $U$ so that $\theta$ is equal to a differential form  on $U'$ which may be constructed from $\C\infty1$ functions using the operations of exterior differentiation and wedge products. The tropical completion of the strata of $U'$ containing $p$ is a refinement of the product of $\mathcal M_{e}$ for all external edges $e$ of $\gamma$. Let $\check\theta$ indicate the tropical completion of $\theta$ on this refinement of the product of $\mathcal M_{e}$. As $\theta$ was generated by functions, and $\mathcal M_{e}$ is complete,  $\check \theta$ is a completely supported differential form on this refinement of $\mathcal M_{e}$. Of course, as $\theta$ and $\check \theta$ are equal restricted to the product of $\mathcal M_{\gamma_{e}}$ for all external edges $e$ of $\gamma$, we may use $\check \theta$ instead of $\theta$ to prove our formula.

Let $\eta_{g_{v},[\gamma_{v}],E_{v}}$ be the Poincare dual in the sense of Theorem \ref{pd class} to $\EV_{v}$ applied to the virtual moduli space $\Mod_{g_{v},[\gamma_{v}],E_{v}}(\check{\ex B}_{v})$, so with a slight abuse of notation, $\bigwedge_v \eta_{g_{v},[\gamma_{v}],E_{v}}$ is the Poincare dual to the virtual moduli space $\prod_{v}\Mod_{g_{v},\gamma_{v},E_{v}}$.  To complete our proof, we must show that
\begin{equation}\label{fppd}\int_{X}\EV^{*}\check \theta=\sum_{\{g_{v},E_{v}\}\vert \sum_{v} g_{v}+g_{\gamma}=g, \ \sum_{v} E_{v}=E} \int_{\prod_{e\subset \gamma}\mathcal M_{\gamma_{e}}}\check\theta\wedge\iota^{*}\bigwedge_{v}\eta_{g_{v},[\gamma_{v}],E_{v}}\end{equation}
This follows if restricted to the component $X'$ of $X$ so that the corresponding curves in $\bMs_{\gamma_{v}}$ have genus $g_{v}$ and energy  $E_{v}$, we have
\begin{equation}\label{gf2}\int_{X'}\EV^{*}\check\theta= \int_{\prod_{e\subset \gamma}\mathcal M_{\gamma_{e}}}\check\theta\wedge\iota^{*}\bigwedge_{v}\eta_{g_{v},\gamma_{v},E_{v}}\end{equation}

To see this, recall from the proof of Theorem \ref{pd class} that  $\bigwedge_{v}\eta_{g_{v},[\gamma_{v}],E_{v}}$ is constructed as follows: Extend $\prod_{v}\EV_{v}$ to a submersion $\psi$ from $\prod_{v}\Mod_{g_{v},[\gamma_{v}],E_{v}}\times \mathbb R^{m}$. Then choose a compactly supported form $\eta_{0}$ on $\mathbb R^{m}$ with integral $1$, and then $\bigwedge_{v}\eta_{g_{v},[\gamma_{v}],E_{v}}=\psi_{!}\eta_{0}$ (the result of integrating along the fibers of $\psi$ the pullback of $\eta_{0}$ to  $\prod_{v}\Mod_{g_{v},[\gamma_{v}],E_{v}}\times \mathbb R^{m}$, interpreted as in the proof of Theorem \ref{pd class}).

At this stage, we may use the fact proved in \cite{dre} that integration along the fiber behaves well under fiber products. In particular, consider the fiber product diagram:

\[\begin{array}{ccc}
(\prod_{e}\mathcal M_{e})\fp\iota\psi(\prod_{v}\Mod_{g_{v},[\gamma_{v}],E_{v}}\times \mathbb R^{m})&\xrightarrow {\psi'}&\prod_{e}\mathcal M_{e}
\\\downarrow \iota'&&\downarrow \iota
\\ \prod_{v}\Mod_{g_{v},[\gamma_{v}],E_{v}}\times \mathbb R^{m} &\xrightarrow{\psi}&\prod_{e'}\mathcal M_{e'}

\end{array}\]

Then $\psi'_{!}\lrb{(\iota^{'})^{*}\eta_{0}}=\iota^{*}(\psi_{!}\eta_{0})=\iota^{*}\bigwedge_{v}\eta_{g_{v},[\gamma_{v}],E_{v}}$. So,

\begin{equation}\label{gf3}\int_{\prod_{e\subset \gamma}\mathcal M_{\gamma_{e}}}\check\theta\wedge\iota^{*}\bigwedge_{v}\eta_{g_{v},[\gamma_{v}],E_{v}}=\int_{(\psi')^{-1}(\prod_{e}\mathcal M_{\gamma_{e}})} (\psi')^{*}\check\theta\wedge\iota'^{*}\eta_{0}\end{equation}

Consider the tropical part of our virtual moduli space $\prod_{v}\Mod_{g_{v},[\gamma_{v}],E_{v}}$. The tropical part of each coordinate chart on this moduli space may be identified with a complete cone so that $0$ corresponds to the curves with tropical part equal to $\gamma_{v}$ in $\bMs_{\gamma_{v}}$. The image of this point $0$ under $\totb{\psi}$ is of course the point corresponding to $\prod_{e'}\mathcal M_{\gamma _{e'}}\subset \prod_{e}\mathcal M_{e'}$, which is equal to the image under $\totb{\iota}$ of the point $p$ corresponding to $\prod_{e}\mathcal M_{\gamma_{e}}$. Our form $\check \theta$ is a closed form in $\Omega^{*}_{c}$ of some refinement of $\prod_{e}\mathcal M_{e}$ with tropical part given by cones centered on $p$. It follows that $(\psi')^{*}\check \theta$ is a form in $\Omega^{*}_{c}$ of some refinement of $(\prod_{e}\mathcal M_{e})\fp\iota\psi(\prod_{v}\Mod_{g_{v},[\gamma_{v}],E_{v}}\times \mathbb R^{m})$ for which the tropical part is some cone around the point corresponding to curves in $\bMs_{\gamma_{v}}$ with tropical part equal to $\gamma_{v}$, and $(\iota')^{*}\eta_{0}$ is also a form in $\Omega^{*}_{c}$ of the same refinement. Therefore their integral over the entire space is equal to their integral over the subset with tropical part $0$,
\begin{equation}\label{gf4}\int_{(\psi')^{-1}(\prod_{e}\mathcal M_{\gamma_{e}})} (\psi')^{*}\check\theta\wedge\iota'^{*}\eta_{0}=
\int_{(\prod_{e}\mathcal M_{e})\fp\iota\psi(\prod_{v}\Mod_{g_{v},[\gamma_{v}],E_{v}}\times \mathbb R^{m})} (\psi')^{*}\check\theta\wedge\iota'^{*}\eta_{0}\end{equation}
Similarly, 
\begin{equation}\label{gf5}\int_{X'}\EV^{*}\check\theta=\int_{\prod_{e}\mathcal M_{e}\fp\iota{\prod_{v}\EV_{v}}\prod_{v}\Mod_{g_{v},[\gamma_{v}],E_{v}}}\EV^{*}\check\theta\end{equation}

 By considering the family of maps $\psi_{t}:=\psi(\cdot,t\cdot)$ for $t\in[0,1]$, we may deform $\psi$ through a family of maps to the map $\psi_{0}$  which is projection to $\prod_{v}\Mod_{g_{v},[\gamma_{v}],E_{v}}$ followed by $\prod_{v}\EV_{0}$. Each of these maps remains transverse to $\iota$, so we get a family of fiber product diagrams
 
\[\begin{array}{ccc}
(\prod_{e}\mathcal M_{e})\fp\iota{\psi_{t}}(\prod_{v}\Mod_{g_{v},[\gamma_{v}],E_{v}}\times \mathbb R^{m})&\xrightarrow {\psi_{t}'}&\prod_{e}\mathcal M_{e}
\\\downarrow \iota_{t}'&&\downarrow \iota
\\ \prod_{v}\Mod_{g_{v},[\gamma_{v}],E_{v}}\times \mathbb R^{m} &\xrightarrow{\psi_{t}}&\prod_{e'}\mathcal M_{e'}

\end{array}\]

Using Stokes' theorem with details expanded as in Theorem \ref{well defined} then gives that

\begin{equation}\label{gf6}\begin{split}\int_{(\prod_{e}\mathcal M_{e})\fp\iota\psi(\prod_{v}\Mod_{g_{v},[\gamma_{v}],E_{v}}\times \mathbb R^{m})} &(\psi')^{*}\check\theta\wedge\iota'^{*}\eta_{0}
\\ &=\int_{(\prod_{e}\mathcal M_{e})\fp\iota{\psi_{0}}(\prod_{v}\Mod_{g_{v},[\gamma_{v}],E_{v}}\times \mathbb R^{m})} (\psi_{0}')^{*}\check\theta\wedge\iota_{0}'^{*}\eta_{0}
\end{split}\end{equation}

Note that $\psi_{0}$ is equal to projection of $\prod_{v}\Mod_{g_{v},[\gamma_{v}],E_{v}}\times \mathbb R^{m}$ to $\prod_{v}\Mod_{g_{v},[\gamma_{v}],E_{v}}$ composed with $\prod_{v}\EV_{v}$, therefore associativity for the orientation of fiber products (discussed in \cite{dre}) gives that
 \[\lrb{\prod_{e}\mathcal M_{e}}\fp\iota{\psi_{0}}(\prod_{v}\Mod_{g_{v},[\gamma_{v}],E_{v}}\times \mathbb R^{m})=\lrb{\lrb{\prod_{e}\mathcal M_{e}}\fp\iota{\prod_{v}\EV_{v}}\prod_{v}\Mod_{g_{v},[\gamma_{v}],E_{v}}}\times \mathbb R^{m}\]
 Furthermore, $\iota'^{*}_{0}\eta_{0}$ is simply equal the pullback of our form $\eta_{0}$ via the obvious projection to $\mathbb R^{m}$. Therefore, 
 
\begin{equation}\label{gf7} \int_{(\prod_{e}\mathcal M_{e})\fp\iota{\psi_{0}}(\prod_{v}\Mod_{g_{v},[\gamma_{v}],E_{v}}\times \mathbb R^{m})} (\psi_{0}')^{*}\check\theta\wedge\iota_{0}'^{*}\eta_{0}=\int_{\lrb{\prod_{e}\mathcal M_{e}}\fp\iota{\prod_{v}\EV_{v}}\prod_{v}\Mod_{g_{v},[\gamma_{v}],E_{v}}}\EV^{*}\check\theta
\end{equation}

Subsituting equation (\ref{gf5}) into (\ref{gf7}), (\ref{gf6}), (\ref{gf4}), then (\ref{gf3}) gives the required equation (\ref{gf2}) which as noted earlier, completes the proof.

\stop

The following is a generalization of Kontevich and Mannin's splitting and genus reduction axioms stated in \cite{KM}.

\begin{thm}[Splitting and genus reduction]\label{splitting}Let $\ex B$ be a basic exploded manifold for which Gromov compactness holds with tropical part equal to a cone, and let $\eta_{g,[\gamma],\beta}$ be the Poincare dual to the map 
\[(\totl{ev_{0}},EV):\Mod_{g,[\gamma],\beta}(\ex B)\longrightarrow \overline{\mathcal M}_{g,[\gamma]}\times \End_{[\gamma]}{\ex B} \]
 
Choose an internal edge $e$  of $\gamma$ which is sent to a point in $\totb{\ex B}$
 and let $\gamma_{i}$ for $i\in I$ indicate the set of tropical curves obtained by cutting $\gamma$ at this internal edge and extending the new ends to be infinite. ($I$ has two elements if $e$ separates $\gamma$, and otherwise has one element.) Choose some $g_{i}$ for all $i\in I$ so that $g_{i}>1$ if $\gamma_{i}$ has less than $3$ external edges, and 
 \[\sum g_{i}=g+\abs I -2\] 
 
 Let \[\phi:\prod_{i\in I}\overline{\mathcal M}_{g_{i},[\gamma_{i}]}\times \End_{[\gamma]}\ex B\longrightarrow \overline{\mathcal M}_{g,[\gamma]}\times \End_{[\gamma]}\ex B\] be the product of the identity on $\End_{[\gamma]}\ex B$ with the  map of Deligne Mumford spaces obtained by identifying the two marked points in curves in  $\overline{\mathcal M}_{g,[\gamma_{i}]}$ which correspond to the  ends of the $\gamma_{i}$'s obtained by cutting $e$.
 
 Note that $\prod_{i}\End_{[\gamma_{i}]}\ex B=\End_{\gamma}\ex B\times \ex B^{2}$. Let $\Delta$ be the Poincare dual to the diagonal in $\ex B^{2}$, pulled back to $\overline{\mathcal M}_{g_{i},[\gamma_{i}]}\times \End_{[\gamma]}\ex B\times \ex B^{2}$.
 \[\text{Let }\pi:\prod_{i\in I}\overline{\mathcal M}_{g_{i},[\gamma_{i}]}\times \End_{[\gamma]}\ex B \times\ex B^{2}\longrightarrow \prod_{i\in I}\overline{\mathcal M}_{g_{i},[\gamma_{i}]}\times \End_{[\gamma]}\ex B \]
be the obvious projection.

Then if $\theta$ is any closed differential form in $\ro^{*}(\prod_{i\in I}\overline{\mathcal M}_{g_{i},[\gamma_{i}]}\times \End_{[\gamma]}\ex B)$,
\[\sum_{\{\beta_{i},i\in I\}}\int \pi^{*}(\theta)\wedge\Delta\wedge\prod_{i\in I}\eta_{g_{i},[\gamma_{i}],\beta_{i}}=\int\theta\wedge\phi^{*}\eta_{g,[\gamma],\beta}\]
 where the sum is over choices of $\beta_{i}$ for all $i\in I$ so that $\sum_{i}\beta_{i}=\beta$. 
 \end{thm}

\pf

As $\totb{\ex B}$ is a cone, we may assume without losing generality that $e$ is the only internal edge of  $\gamma$, and that all vertices are contained in the smallest strata of $\totb{\ex B}$. Then $\bMs_{g,\gamma,\beta}$ is an $(\Aut \gamma)$-fold cover of its image in  $\Msw_{g,[\gamma],\beta}$. We are only interested in the component of $\bMs_{g,\gamma,\beta}$ which has genus $g_{i}$ at the $i$th vertex of $\gamma$, so we should modify this a little. In particular, let $ \gamma'$ indicate $\gamma$ decorated with genus $g_{i}$ at the $i$th vertex, and let $\bMs_{g,\gamma',\beta}$ indicate the component of $\bMs_{g,\gamma,\beta}$ that has genus $g_{i}$ at the $i$th vertex of $\gamma$. Then $\bMs_{g,\gamma',\beta}$ is an $(\Aut \gamma')$-fold cover of its image in $\Msw_{g,[\gamma],\beta}$.
 $\abs{\Aut\gamma'}=2$ if our internal edge $e$ is a loop, or if $e$ is the only edge of $\gamma$ and $g_{1}=g_{2}$. Otherwise $\gamma'$ has no automorphisms. As the $\gamma_{i}$ have no internal edges or automorphisms, $\bMs_{g_{i},\gamma_{i},\beta_{i}}=\Msw_{g_{i},[\gamma_{i}],\beta_{i}}$.  Note that  the map 
\[\prod_{i\in I}\overline{\mathcal M}_{g_{i},[\gamma_{i}]}\longrightarrow \overline{\mathcal M}_{g,[\gamma]}\]
is also an $(\Aut \gamma')$-fold cover of its image which is a boundary strata of $\overline{\mathcal M}_{g,[\gamma]}$.

In the notation of Theorem \ref{gluing cobordism}, we may construct the virtual moduli spaces $\Mod _{g_{i},[\gamma_{i}],E_{i}}$ so that $EV_{0}$ and $EV_{1}$ are transverse maps $\prod_{i}\Mod_{g_{i},[\gamma_{i}],E_{i}}\longrightarrow \ex B$ when $\sum E_{i}= \beta(\omega)$. As we are gluing along an edge sent to a point, it makes sense to talk about summing homology classes $\beta_{i}$, and the inverse image in $\bMs_{\gamma}$ of  $\prod_{i}\Mod_{g_{i},[\gamma_{i}],\beta_{i}}$  is in $\bMs_{\sum g_{i}+2-\abs I,\gamma',\sum\beta_{i}}$. The above transversality therefore implies that $EV_{0}$ and $EV_{1}$ are transverse maps  $\prod_{i}\Mod_{g_{i},[\gamma_{i}],\beta_{i}}\longrightarrow \ex B$ when $\sum \beta_{i}= \beta$.

 Let $X$ be the disjoint union of  the  fiber product of $EV_{0}$ with $EV_{1}$ for each choice of $\{\beta_{i}\}$ summing to $\beta$, and $\psi:X\longrightarrow \prod_{i\in I}\overline{\mathcal M}_{g_{i},[\gamma_{i}]}\times \End_{[\gamma]}\ex B$ be the corresponding evaluation map.

Then a similar argument to the proof of equation (\ref{fppd}) in the proof of theorem \ref{gluing formula} gives
\[\sum_{\beta_{i}}\int_{\prod_{i\in I}\overline{\mathcal M}_{g_{i},[\gamma_{i}]}\times \End_{[\gamma_{i}]}\ex B} \pi^{*}(\theta)\wedge\Delta\wedge\prod_{i\in I}\eta_{g_{i},[\gamma_{i}],\beta_{i}}=\int_{X}\psi^{*}\theta\]

Lemma \ref{fp bundle} implies that the inverse image of $\coprod_{\beta_{i}}\prod_{i}\Mod_{g_{i},[\gamma_{i}],\beta_{i}}$ in $\bMs_{\gamma}$ is equal to a $\et 1{(0,\infty)}$ bundle $\hat X$ over $X$. In fact,  the subset of $\M_{g,[\gamma]}$ corresponding to $\prod_{i}\M_{g_{i},[\gamma_{i}]}$   is equal to an $\et 1{(0,\infty)}$ bundle  over $\prod_{i}\M_{g_{i},[\gamma_{i}]}$, and the inverse image of $\prod_{i}\Mod_{g_{i},[\gamma_{i}],\beta_{i}}$ in $\bMs_{\gamma}$ is the pullback of this bundle over an  evaluation map. A similar statement holds for any family in the image of $\bMs_{\gamma}$ in $\prod_{i}\Msw_{g_{i},[\gamma_{i}]}$.

Construct a closed differential form $\alpha$ in $\ro^{2}_{c}$ of this $\et 1{(0,\infty)}$ bundle over $\prod_{i}\M_{g_{i},[\gamma_{i}]}$ with integral $1$ on these $\et 1{(0,\infty)}$ fibers, and extend $\alpha$ to be a form on  $\M_{g,[\gamma]}$ by setting it equal to $0$ everywhere else. Denote by $\hat \psi$ the evaluation map $\bMs_{\gamma'}\longrightarrow \prod_{i\in I}\overline{\mathcal M}_{g_{i},[\gamma_{i}]}\times \End_{[\gamma]}\ex B $

\[\int_{X}\psi^{*}\theta=\int_{\hat X}\hat \psi^{*}\theta\wedge (ev^{0})^{*}\alpha\]

Let $\bar{\Mod}_{g,\gamma',\beta}$ indicate the inverse image of $\Mod_{g,[\gamma],\beta}$ in $\bMs_{g,\gamma',\beta}$. Theorem \ref{gluing cobordism} implies that $\hat X$ is cobordant to $\bar{\Mod}_{g,\gamma',\beta}$, so 

\[\int_{\hat X}\hat \psi^{*}\theta\wedge (ev^{0})^{*}\alpha=\int_{\bar{\Mod}_{g,\gamma',\beta}}\hat\psi^{*}\theta\wedge(ev^{0})^{*}\alpha \]

(Note that this cobordism is compact but not quite complete, as it is an $\et 1{(0,\infty)}$ bundle over a complete cobordism. As $(ev^{0})^{*}\alpha$ has complete support on this cobordism, this is not a problem for applying the version of  Stoke's theorem from \cite{dre}.) Let $\theta'$ be any differential form on $\M_{g,[\gamma],\beta}\times \End_{[\gamma]}\ex B$ which   pulls back to $\theta$ under the map \[\phi:\prod_{i\in I}\overline{\mathcal M}_{g_{i},[\gamma_{i}]}\times \End_{[\gamma]}\ex B\longrightarrow \overline{\mathcal M}_{g,[\gamma]}\times \End_{[\gamma]}\ex B\]
and which is closed in a neighborhood of the image of $\phi$. Note that $\bar{\Mod}_{g,\gamma',\beta}$ is an $(\Aut \gamma')$-fold cover of its image in $\Mod_{g,[\gamma],\beta}$, which contains the support of $(ev^{0})^{*}\alpha$. Then
\[\int_{\bar{\Mod}_{g,\gamma,\beta}}\hat\psi^{*}\theta\wedge(ev^{0})^{*}\alpha=\abs{\Aut \gamma'}\int_{\Mod_{g,[\gamma],\beta}}(\totl{ev^{0}},EV)^{*}\theta'\wedge(ev^{0})^{*}\alpha\]
Let $\alpha'$ be a Poincare dual of the image of $\phi$, supported in a connected neighborhood of the image of $\phi$ where $\theta'$ is closed.  Then the pullback of $\alpha'$ to $\M_{g,\gamma}\times \End_{[\gamma]}\ex B$ is in the same cohomology class as the pull back of $\alpha$ to $\M_{g,\gamma}\times \End_{[\gamma]}\ex B$. Then

\[\begin{split}\abs{\Aut \gamma'}\int_{\Mod_{g,[\gamma],\beta}}(\totl{ev^{0}},EV)^{*}\theta'\wedge(ev^{0})^{*}\alpha&=\abs{\Aut \gamma'}\int_{\Mod_{g,[\gamma],\beta}}(\totl{ev^{0}},EV)^{*}(\theta'\wedge\alpha')
\\ &= \abs{\Aut \gamma'}\int \theta'\wedge\alpha'\wedge\eta_{g,[\gamma],\beta}
\\ &=\int \theta \wedge\phi^{*}\eta_{g,[\gamma],\beta}\end{split}\]
where the last step uses that $\phi$ is an $(\Aut \gamma')$-fold cover of its image.

\stop

\section{Example computations}

\subsection{Curves in $\ex T^{n}$}\label{curves in Tn}

\

\

The easiest nontrivial curves to analyze are zero genus curves in $\ex T^{n}$. It is easily verified that the only stable zero genus curves in $\ex T^{n}$ with less than $3$ punctures are the maps $\ex T\longrightarrow \ex T^{n}$. It follows that for the case of zero genus curves,  we need only analyze curves in $\ex T^{n}$ with tropical parts which are equal to zero genus graphs with no univalent or bivalent vertices.

 As noted in \cite{iec}, the tropical part of any curve in $\ex T^{n}$ obeys the `balancing' or `conservation of momentum'  condition familiar to tropical geometers.  In particular, let $\gamma$ be the tropical part  of a curve in $\ex T^{n}$. This is a map of a graph into $\mathbb R^{n}$ considered as a tropical part of $\ex T^{n}$. The edges of this graph have an integral affine structure, and $\gamma$ restricted to each edge is an integral affine map.  If we choose any orientation on the domain of $\gamma$, we may define a momentum for each edge to be the image under the derivative of $\gamma$ of the unit vector pointing in the positive direction of our edge. This momentum is a vector in $\mathbb Z^{n}$. The conservation of momentum condition states that at any vertex of $\gamma$, the sum of the momentum of the incoming edges is equal to the sum of the momentum of the outgoing edges. Call any infinite edge of $\gamma$  an end of $\gamma$, and make the convention that if we do not specify an  orientation for the ends of $\gamma$, they are oriented to be outgoing. (Of course, the conservation of momentum condition implies that the sum of the momentum of the ends of $\gamma$ is $0$.)

 Let us now examine the moduli space of zero genus holomorphic curves in $\ex T^{n}$ that have $3$ ends with momentum $a,\ b$ and $-a-b$. For such a curve to be stable, its tropical part must be equal to a graph with a single vertex and $3$ ends. Therefore, the domain of our curve is uniquely isomorphic to the explosion of $\mathbb CP^{1}$ relative to the three points $0$,$1$, $\infty$ so that $0$, $1$ and $\infty$ correspond to our first, second and third ends. Restricted to $\mathbb CP^{1}-\{0,1,\infty\}$, our holomorphic curves are curves of the form:
 \[\lrb{c_{1}z^{a_{1}}(z-1)^{b_{_{1}}},\dotsc,c_{n}z^{a_{n}}(z-1)^{b_{n}}}\] 
so our moduli space is parametrized by $(c_{1},\dotsc,c_{n})\in \ex T^{n}$. Let us now trace through the steps of the construction of the virtual moduli space to see that this explicit moduli space is in fact our virtual moduli space.

In this case, we may cover our entire moduli space by a single core family. Let $\hat f$ indicate the above family of curves parametrized by $\ex T^{n}$. In this case, the group $G$ of automorphisms from the definition of core families is the trivial group. Choose a point $p$ in $\ex CP^{1}-\{0,1,\infty\}$, then let $s$ indicate the section of $\ex C(\hat f)\longrightarrow \ex F(\hat f)$ corresponding to the point $p$. The criteria \ref{crit1}, \ref{crit2} and  \ref{crit3} from Theorem \ref{core criteria}  
are easily seen to be satisfied by $(\hat f/G,s)$. We must check criterion \ref{crit4}, so we need that there exists a neighborhood of the section $s:\ex F(\hat f)\longrightarrow \ex C(\hat f)$ on which  
\[ev^{+1}(\hat f):\ex C(\hat f)\longrightarrow \M_{0,4}\times \ex T^{n}\] 
is an equidimensional embedding, and we need to check that the tropical part of $ev^{+1}\circ s$ is a compete map, and restricted to any polytope in $\ex F(f)$ is an isomorphism onto a strata of the image
 in $\M_{0,4}\times\ex T^{n}$
 under $ev^{+1}$ of some open neighborhood of $f$ in $\Msw$. In this case, $\ex F(\hat f)$ is equal to $\ex T^{n}$, $\ex C(\hat f)$ is equal to $\M_{0,4}\times \ex T^{n}$, and $ev^{+1}(\hat f)$ restricted to $(\mathbb  CP-\{0,1,\infty\})\times \ex T^{n}\subset \M_{0,4}\times \ex T^{n}$ is given by the formula
\[ev^{+1}(\hat f)(z,c_{1},\dotsc,c_{n})=\lrb{z,c_{1}z^{a_{1}}(z-1)^{b_{_{1}}},\dotsc,c_{n}z^{a_{n}}(z-1)^{b_{n}}}\]
This is an isomorphism onto an open subset of $\M_{0,4}\times \ex T^{n}$,  so criterion \ref{crit4} is satisfied.  To complete the description of our core family, we need only define a map
\[F:\hat f^{*}T\ex T^{n} \longrightarrow \ex T^{n}
\]
satisfying criterion \ref{crit5}. Standard coordinates on the tangent bundle of $\ex T^{n}$ given by the real and imaginary parts of $\tilde z_{i}\frac{\partial}{\partial \tilde z_{i}}$ identify $T\ex T^{n}$ with a trivial $\mathbb C^{n}$ bundle. Then we may define $F:\mathbb C^{n}\times \ex C(\hat f)$ by 
\[F(v_{1},\dotsc,v_{n},z):=(e^{v_{1}}\hat f(z)_{1},\dotsc, e^{v_{n}}\hat f(z)_{n})\]
This map $F$ satisfies the conditions of criterion \ref{crit5}, so Theorem \ref{core criteria}  implies that $(\hat f,s,F)$ is a core family which covers our entire moduli space of holomorphic curves.

Next, we shall verify that $\hat f$ is trivially an obstruction model when we use the zero dimensional vector bundle in place of $V$. To do this, we need to  verify that for all curves $f\in\hat f$, $D\dbar(f):X^{\infty,\underline 1} ( f)\longrightarrow Y^{\infty,\underline 1}( f)$ is a bijection. In this case $D\dbar(f)$ is the usual $\dbar$ operator. Standard complex analysis implies that any holomorphic $\C\infty1$ map $\ex C(f)\longrightarrow \mathbb C^{n}$ must be constant, as restricted to the smooth part of $\ex C(f)$, it is a bounded holomorphic map from $\mathbb CP^{1}\setminus \{0,1,\infty \}$. As $X^{\infty,\underline 1}(f)$ consists of $\C\infty1$ maps $\ex C(f)\longrightarrow \mathbb C^{n}$ which vanish at a particular marked point, it follows that $D\dbar(f)$ is injective restricted to $X^{\infty,\underline 1}(f)$. A dense subset of $Y^{\infty,\underline 1}( f)$ is the space of smooth maps from $\mathbb CP^{1}(f)$ into $\mathbb C^{n}$ which vanish in a neighborhood of $0$, $1$, and $\infty$. Cauchy's integral formula then gives that $D\dbar(f)$ is surjective onto this dense subset, and therefore  surjective onto $Y^{\infty,\underline 1}$ as Theorem \ref{f replacement} tells us that $D\dbar(f)$ has a closed image. Therefore, for every curve $f$ in $\hat f$, $D\dbar(f):X^{\infty,\underline 1}(f)\longrightarrow Y^{\infty,\underline 1}(f)$ is a bijection, so $\hat f$ is an obstruction model which covers our entire moduli space. As in this case $D\dbar$ is a complex map, the orientation of this trivial obstruction bundle is the positive one, so our virtual moduli space has the same orientation as $\hat f$, which has the orientation given by the complex structure on $\ex F(\hat f)$. 

We may follow the construction of the virtual moduli space from section \ref{virtual class} using our obstruction model $\hat f$, and using trivial perturbations to get that $\hat f$ parametrizes the virtual moduli space, which in this case is the actual moduli space of holomorphic curves.

We now have a parametrization of our moduli space $\mathcal M$ of holomorphic curves by $\ex T^{n}$. To compute the map $\EV$, we need coordinates on the moduli spaces $\mathcal M_{a}$, $\mathcal M_{b}$
 and $\mathcal M_{-a-b}$ of possibilities for the ends of our curves. The curves in $\mathcal M_{a}$ are maps $\ex T\longrightarrow \ex T^{n}$ in the form
 \[\tilde z\mapsto(c_{1}\tilde z^{a_{1}},\dotsc,c_{n}\tilde z^{a_{n}})\]
 where $(c_{1},\dotsc,c_{n})\in \ex T^{n}$
  In the special case that $a=0$, $\mathcal M_{a}$ is equal to $\ex T^{n}$, otherwise, $\mathcal M_{a}$ is isomorphic to $\ex T^{n-1}$ (where as in section \ref{gluing relative invariants} and \ref{construction of EV}, we ignore the orbifold structure on $\mathcal M_{a}$ when $a$ is a nontrivial multiple of another integral vector). We may identify exploded functions $\mathcal M_{a}\longrightarrow \ex T$ to be exploded functions on  $\ex T^{n}$ which are constant on each curve in $\mathcal M_{a}$. For example, if we change coordinates so that $a=(\abs a,0,\dotsc,0)$, then $\tilde z_{2},\dotsc,\tilde z_{n}$ give coordinates on $\mathcal  M_{a}$. This defines a map $\pi_{a}:\ex T^{n}\longrightarrow \mathcal M_{a}$. Straightforward computation then gives that 
  $\EV:\mathcal M\longrightarrow \mathcal M_{a}\times \mathcal M_{b}\times\mathcal M_{-a-b}$ is given  by 
  \[\EV(c_{1},\dotsc,c_{n})=\lrb{\pi_{a}((-1)^{b_{1}}c_{1},\dotsc,(-1)^{b_{n}}c_{n}),\pi_{b}(c_{1},\dotsc,c_{n}),\pi_{-a-b}(c_{1},\dotsc,c_{n})}\]
 
Note that in this case, our moduli space and $\EV$ may be read off from the corresponding tropical problem. Our moduli space is parametrized by $\ex T^{n}$, which corresponds to the image of a chosen extra  marked point on our domain. The corresponding moduli space of tropical curves is parametrized by $\totb{\ex T^{n}}$, which corresponds to the position of the vertex of our tropical curve. Similarly,  $\mathcal M_{a}$ is either $\ex T^{n-1}$ or $\ex T^{n}$, and the tropical part of $\mathcal M_{a}$  is equal to the corresponding tropical moduli space of lines in the direction of $a$. Up to multiplication by a constant, $\EV$ is just given by the product of the three obvious projection maps $\EV:\ex T^{n}\longrightarrow \mathcal M_{a}\times\mathcal M_{b}\times \mathcal M_{-a-b}$, and the tropical part  $\totb{\EV}$ is equal to to corresponding obvious tropical projections.

Now Lemma \ref{fp bundle} and Theorem \ref{gluing cobordism} allow us to compute the part of the virtual moduli space of curves in $\ex T^{n}$ consisting of curves with tropical part equal to a trivalent graph $\gamma$ with genus equal to the genus of the corresponding curve. We shall see that the nature of this moduli space can be read off from the corresponding tropical moduli space.

Choose an oriented trivalent graph $\Gamma$ with labeled free ends and assign a vector  $\alpha_{e}\in\mathbb Z^{n}$ to each edge $e$ of $\Gamma$ so that the sum of $\alpha_{e}$ for edges $e$ leaving a vertex is equal to the sum of $\alpha_{e'}$ for edges $e'$ entering a vertex. An automorphism of $\Gamma$ is a isomorphism of the graph $\Gamma$ to itself which fixes the free ends of $\Gamma$, so that if an edge $e$ is sent to $e'$, then $\alpha_{e}=\alpha_{e'}$ the map is orientation preserving, and $\alpha_{e}=-\alpha_{e'}$ if the map is orientation reversing.  Let $\Aut \Gamma$ indicate the group of automorphisms of $\Gamma$. Say that a tropical curve $\gamma$ in $\totb{\ex T^{n}}$ has the combinatorial type of $\Gamma$ if there is an isomorphism of $\Gamma$ with the domain of our tropical curve preserving the labeling of free ends so that on each edge $e$, the unit vector in the positive direction as determined by the orientation of $\Gamma$ is sent to $\alpha_{e}$.

Now consider the moduli space of tropical curves $\gamma$ with the combinatorial type of $\Gamma$. To specify such a tropical curve, we may specify the length $l_{e}$ of each internal edge $e$ of $\gamma$, and specify the position $x_{v}$ of each vertex of $\gamma$. Together, $(l_{e},x_{v})$ give coordinates on some space $(0,\infty)^{\#e}\times\mathbb R^{n(\#v)}$ where $\#e$ is the number of internal edges of $\Gamma$ and $\#v$ is the number of vertices of $\Gamma$. The group $\Aut \Gamma$ acts in an obvious way on this space. Let $v(e_{0})$ be the vertex of $\Gamma$ which is attached to the initial end of $e$, and $v(e_{1})$ be the vertex of $\Gamma$ attached to the final end of $e$. Then the moduli space of tropical curves with the combinatorial type of $\Gamma$ is equal to the quotient by $\Aut\Gamma$ of the subset of $(0,\infty)^{\#e}\times\mathbb R^{n(\#v)}$ satisfying the equations:
\begin{equation}\label{tropeqn}x_{v(e_{0})}+l_{e}\alpha_{e}-x_{v(e_{1})}=0\end{equation}
so our tropical moduli space is the inverse image of $0$ under some integral linear map
\begin{equation}\label{defa}A:(0,\infty)^{\#e}\times\mathbb R^{n(\#v)}\longrightarrow \mathbb R^{n(\#e)} \end{equation}
This map $A$ will be important for describing the corresponding moduli space of holomorphic curves. The corresponding virtual moduli space of holomorphic curves will be empty if $A$ is not transverse to $0$, and otherwise the tropical part of the corresponding moduli space will be $\abs{\mathbb Z^{n(\#e)}/A(\mathbb Z^{\#e+n(\#v)})}$ disjoint copies of $A^{-1}(0)/\Aut \Gamma$.

 Note that in the case that the genus of $\Gamma$ is $0$, then $\Aut\Gamma$ is the trivial group,  $A$ is always transverse to zero, and $A(\mathbb Z^{\#e+n(\#v)})=\mathbb Z^{n(\#e)}$, so in this case  we shall see that the tropical part of our moduli space of curves is equal to $A^{-1}(0)$.
 
 Consider the problem of finding the moduli space of holomorphic curves in $\ex T^{n}$ with genus equal to the genus of $\Gamma$, and with tropical part with the combinatorial type of $\Gamma$. Suppose $A^{-1}(0)$ is nonempty, so such a tropical curve $\gamma$ exists. Then the pullback of our (virtual) moduli space of holomorphic curves to $\bMs_{\gamma}$ is an $\abs{\Aut\Gamma}$-fold cover of its image in $\Msw(\ex B)$. In light of Theorem \ref{gluing cobordism}, we shall first study the fiber product of our moduli space of curves in  $\prod_{v}\bMs_{\gamma_{v}}$ over the maps $\EV_{0}$ and $\EV_{1}$. 
 
 As the total genus of our curve is equal to the genus of $\gamma$, we are interested in zero genus curves in $\bMs_{\gamma_{v}}$. All our curves in $\bMs_{\gamma_{v}}$ therefore have domain equal to the explosion of $\mathbb CP^1$ relative to $3$ points. We may choose a point on the domain, and parametrize our moduli space by its image $\tilde z_{v}$ in $\ex T^{n}$. So $\tilde z_{v}$ for all vertices  $v$ in $\gamma$ give coordinates on our moduli space in $\prod_{v}\bMs_{\gamma_{v}}$ which we have now identified with $\ex T^{n(\#v)}$. If $\pi_{e}$ denotes the standard projection  $\mathbb T^{n}\longrightarrow\mathcal M_{e}$, then our  map $\EV_{i}$ is equal to the product of  $\pi_{e}(\sigma_{e_{i}}\tilde z_{v(e_{i})})$ where $\sigma_{e_{i}}\in \{1,-1\}^{n}$.    
The result of taking the fiber product of our moduli space over $\EV_{0}$ and $\EV_{1}$ is therefore equivalent to the equations
\[\pi_{e}(\sigma_{e_{0}}\tilde z_{v(e_{0})})=\pi_{e}(\sigma_{e_{1}}\tilde z_{v(e_{1})})\text{ for all edges }e \text{ of }\gamma\]
 These equations are equivalent to the requirement that there exists some $\tilde c_{e}\in\ex T$ so that \begin{equation}\label{nontropeqn}\tilde z_{v(e_{0})}\tilde c_{e}^{\alpha_{e}}\tilde z^{-1}_{v(e_{0})}=\sigma_{e_{1}}\sigma_{e_{0}}\end{equation} In this case, the condition that the tropical part of the corresponding curve is in the image of $\gamma$-decorated tropical completion is equivalent to the condition that $\tilde c_{e}\in \et{1}{(0,\infty)}$.
 
  Note at this stage that the equations \ref{tropeqn} and  our tropical map $A$ from \ref{defa} can be obtained from the tropical part of these equations \ref{nontropeqn} by equating the tropical part of $\tilde z_{v}$ with $x_{v}$ and the tropical part of $\tilde c_{e}$ with $l_{e}$. If $A$ is surjective,  $\EV_{0}$ and $\EV_{1}$ are transverse, so if $\Gamma$ has no automorphisms, Theorem \ref{gluing cobordism} tells us that we may construct our virtual moduli space to coincide with the actual moduli space of holomorphic curves in the region we are studying. (In the case that $\Gamma$ has nontrivial automorphisms, our moduli space of actual holomorphic curves will still be transversely cut out, but some holomorphic curves will have automorphisms, so some multiperturbation is required to get rid of these automorphisms). On the other hand, if A is not surjective, the fact that our moduli space is some power of $\ex T$ implies that  any perturbation to make the maps $\EV_{0}$ and $\EV_{0}$ transverse will result in $\EV_{0}$ and $\EV_{1}$ not intersecting at all. Therefore in this case, the virtual moduli space restricted to curves with the correct genus and tropical part is empty (or at least cobordant to the empty set; in fact, dimension considerations imply that this part of the virtual moduli space must always be empty even though this part of the actual moduli space of holomorphic curves may not be empty.) 
 
The coordinates $\tilde z_{v}$ and $\tilde c_{e}$ are coordinates for a core family in $\bMs_{\gamma}$ containing all our holomorphic curves, where $\tilde c_{e}$ are the gluing coordinates, and $\tilde z_{v}$ correspond to the image of the obvious marked points. Therefore, our moduli space of holomorphic curves in $\bMs_{\gamma}$ is given by the subset of $(\et 1{(0,\infty)})^{\#e}\times\ex T^{n(\#v)}$ where 
\[\tilde z_{v(e_{0})}\tilde c^{\alpha_{e}}\tilde z^{-1}_{v(e_{0})}=\sigma_{e_{1}}\sigma_{e_{0}}\text{ for all edges } e\text{ in }\Gamma\]
 In other words, our moduli space is the inverse image of some point under a map $(\et 1{(0,\infty)})^{\#e}\times\ex T^{n(\#v)}\longrightarrow \ex T^{n(\#e)}$ with tropical part equal to $A$. Therefore, the nature of this part of the moduli space of holomorphic curves can be read off simply from considering the corresponding moduli space of tropical curves. 
 
 In the next section, we prove that the virtual moduli space of curves in a refinement of $\ex B$ is the corresponding refinement of the virtual moduli space of curves in $\ex B$. Therefore, the above computations apply to any refinement of $\ex T^{n}$ such as the explosion of a toric manifold relative to its toric boundary divisors.

\subsection{Refinements}

\

\

  Recall  from \cite{iec} that given any refinement $\ex B'\longrightarrow \ex B$ and any map $\ex A\longrightarrow \ex B$, taking the fiber product of $\ex B'$ with $\ex A$ over $\ex B$ gives a refinement $\ex A'\longrightarrow \ex A$ with a corresponding map $\ex A'\longrightarrow \ex B'$. This functorial construction may be applied to any $\C\infty1$ curve $ f$  in $\ex B$ to obtain a corresponding curve $ f'$  in $\ex B'$. 
  
  For a $\C\infty1$ family of curves $\hat f$ in $\ex B$, the construction is more complicated, because the fiber product of $\hat f$ with $\ex B'\longrightarrow \ex C$ will produce a map $\ex C'\longrightarrow \ex B'$ which may not be the total space of any family of curves, as there may not be a refinement $\ex  F'$ of $\ex F(\hat f)$ so that $\ex C'\longrightarrow \ex F'$ is a family (the derivative applied to integral vectors in $T\ex C'$ may not be surjective onto the integral vectors in $T\ex F'$). 
  
  \begin{lemma}
  Given a $\C\infty1$ family $\hat f$ of curves in $\ex B$ and a refinement $\ex B'\longrightarrow \ex B$, there exists a unique $\C\infty1$ family $\mathcal R(\hat f)$ of curves in $\ex B'$ fitting into the commutative diagram:
  \[\begin{array}{ccccc}\ex F(\mathcal R(\hat f))&\longleftarrow&\ex C(\mathcal R(\hat f))&\longrightarrow &\ex B'
  \\ \downarrow && \downarrow &&\downarrow
  \\ \ex F(\hat f)&\longleftarrow&\ex C(\hat f)&\longrightarrow &\ex B\end{array}\]
$\mathcal R(\hat f)$ satisfies the universal property that given any $\C\infty1$ family $\hat g$ of curves in $\ex B'$ and a commutative diagram 
  \[\begin{array}{ccccc}\ex F(\hat g)&\longleftarrow&\ex C(\hat g)&\longrightarrow &\ex B'
  \\ \downarrow && \downarrow &&\downarrow
  \\ \ex F(\hat f)&\longleftarrow&\ex C(\hat f)&\longrightarrow &\ex B\end{array}\]
there exists a unique factorization of the above maps into
  \[\begin{array}{ccccc}\ex F(\hat g)&\longleftarrow&\ex C(\hat g)&\longrightarrow &\ex B'
  \\ \downarrow && \downarrow &&\downarrow
  \\
  \ex F(\mathcal R(\hat f))&\longleftarrow&\ex C(\mathcal R(\hat f))&\longrightarrow &\ex B'
  \\ \downarrow && \downarrow &&\downarrow
  \\ \ex F(\hat f)&\longleftarrow&\ex C(\hat f)&\longrightarrow &\ex B\end{array}\]

  \end{lemma}
  
  \pf
  We may use the universal property of $\mathcal R(\hat f)$ to reduce to the case that $\ex F(\hat f)$ is covered by a single coordinate chart $ U$. Let $\ex C'$ indicate the fiber product of $\ex C(\hat f)$ with $\ex B'$ over $\ex B$.
  
  Consider lifts of our coordinate chart $U$ on $\ex F(\hat f)$ to a coordinate charts $\tilde U_{i}$ on $\ex C(\hat f)$. The refinements of $\tilde U_{i}$ from taking the fiber product of $\hat f$ with $\ex B'\longrightarrow \ex B$ correspond to  subdivisions of the polytopes $\totb{\tilde U_{i}}$. The projection of all the polytopes in these subdivision to $\totb U$ corresponds to a subdivision of $\totb U$, and hence a refinement $U'$ of $U$. 
 Let $\ex F(\mathcal R(\hat f))$ be $U'$, $\ex C(\mathcal R(\hat f))$ be the fiber product of $\ex C'$ with $U'$ over $U$, and let $\mathcal R(\hat f)$ be the composition of the maps $\ex C(\mathcal R(\hat f))\longrightarrow \ex C'\longrightarrow \ex B'$. The map $\ex C(\mathcal R(\hat f))\longrightarrow \ex F(\mathcal R(\hat f))$ is now a family of curves,  because the derivative is now surjective on integral vectors. So we have constructed a $\C\infty1$ family of curves $\mathcal R(\hat f)$ which comes with a commutative diagram
  
   \[\begin{array}{ccccc}\ex F(\mathcal R(\hat f))&\longleftarrow&\ex C(\mathcal R(\hat f))&\longrightarrow &\ex B'
  \\ \downarrow && \downarrow &&\downarrow
  \\ \ex F(\hat f)&\longleftarrow&\ex C(\hat f)&\longrightarrow &\ex B\end{array}\]
 
 Now suppose that we have any $\C\infty1$ family  of curves $\hat g$ which fits into the commutative diagram 
  
  \[\begin{array}{ccccc}\ex F(\hat g)&\longleftarrow&\ex C(\hat g)&\longrightarrow &\ex B'
  \\ \downarrow && \downarrow &&\downarrow
  \\ \ex F(\hat f)&\longleftarrow&\ex C(\hat f)&\longrightarrow &\ex B\end{array}\]
  
  Again, we may use the universal property we are trying to prove to reduce to the case that $\ex F(\hat g)$ is covered by a single coordinate chart. As $\ex C'$ is a fiber product  of $\ex B'$ with $\ex C(\hat f)$, we get a map $\ex C(\hat g)\longrightarrow \ex C'$ and the following commutative diagram.
  \[\begin{array}{ccc}{\ex C(\hat g)}&\longrightarrow & {\ex C'}
  \\ \downarrow &&\downarrow
  \\ {\ex F(\hat g)}&\longrightarrow&{\ex F(\hat f)}\end{array}\]

  The inverse image of every polytope in the tropical part of $\ex C'$ is a polytope in the tropical part of $\ex C(\hat g)$. As $\ex C(\hat g)\longrightarrow \ex F(\hat g)$ is a family, the projection of each of these polytopes to $\ex F(\hat g)$ is a polytope in the tropical part of $\ex F(\hat g)$. Therefore, the inverse image in the tropical part of $\ex F(\hat g)$ of the projection to the tropical part of $\ex F(\hat g)$ of any polytope in the tropical part of $\ex C'$ is a polytope in the tropical part of $\ex F(\hat g)$. It follows that the map $\ex F(\hat g)\longrightarrow \ex F(\hat f)$ factors through $\ex F(\mathcal R(\hat f))$, so we get a commutative diagram   
  \[\begin{array}{ccc}{\ex C(\hat g)}&\longrightarrow & {\ex C'}
  \\ \downarrow &&\downarrow
  \\ {\ex F(\mathcal R(\hat f))}&\longrightarrow&{\ex F(\hat f)}\end{array}\]
Then using the fact that $\ex C(\mathcal R(\hat f))$ is defined as the fiber product of $\ex C'$ with $\ex F(\mathcal R(\hat f))$ over $\ex F(\hat f)$, we get the required commutative diagram

  \[\begin{array}{ccccc}\ex F(\hat g)&\longleftarrow&\ex C(\hat g)&\longrightarrow &\ex B'
  \\ \downarrow && \downarrow &&\downarrow
  \\
  \ex F(\mathcal R(\hat f))&\longleftarrow&\ex C(\mathcal R(\hat f))&\longrightarrow &\ex B'
  \\ \downarrow && \downarrow &&\downarrow
  \\ \ex F(\hat f)&\longleftarrow&\ex C(\hat f)&\longrightarrow &\ex B\end{array}\]

The uniqueness of this diagram is automatic as the maps $\ex F(\mathcal R(\hat f))\longrightarrow \ex F(\hat f)$ and $\ex C(\mathcal R(\hat f))\longrightarrow \ex C(\hat f)$ are bijective.
    
  \stop

 \begin{thm} \label{refinement theorem} Given a refinement $\ex B'\longrightarrow \ex B$, the virtual moduli space $\mathcal M$ of holomorphic curves in $\ex B$ and  the virtual moduli space of holomorphic curves in $\ex B'$ may be constructed so that the virtual moduli space of holomorphic curves in $\ex B'$ is the refinement $\mathcal R(\mathcal M)$ of $\mathcal M$. 
 
\end{thm}

\pf 
First, note that given a stable holomorphic curve $f'$ in $\ex B'$, composing with $\ex B'\longrightarrow\ex B$ gives a holomorphic curve in $\ex B$ which may not be stable, because it may be a refinement of another holomorphic curve. Let $f$ indicate the underlying stable holomorphic curve in $\ex B$. We have the following commutative diagram
\[\begin{array}{ccc} \ex C(f')&\xrightarrow{f'} &\ex B'
\\ \downarrow &&\downarrow
\\ \ex C(f)&\xrightarrow{f}&\ex B\end{array}\]
There is a corresponding unique map of $f'\longrightarrow \mathcal R(f)$ which is a refinement map because the composition $\ex C(f')\longrightarrow \ex C(\mathcal R(f))\longrightarrow \ex C(f)$ is a refinement. As $f'$ is stable it follows that $f'=\mathcal R(f)$. Therefore the moduli stack of stable holomorphic curves in $\ex B'$ is $\mathcal R$ applied to the moduli stack of stable holomorphic curves in $\ex B$. 
 
Now suppose that $(\hat f/G,V)$ is an obstruction model on $\Msw(\ex B)$. We shall examine what happens when the functor $\mathcal R$ is applied to such an obstruction model for $\Msw(\ex B)$.  

 Implicit in $(\hat f/G,V)$, we have a core family $(\hat f/G,F,\{s_{i}\})$ for an open substack $\mathcal O$ of $\Msw(\ex B)$. 
Note that the group $G$ still acts on $\mathcal R(\hat f)$. The sections $\{s_{i}\}$ of $\ex C(\hat f)\longrightarrow\ex F(\hat f)$ all consist of marked points in the smooth part of curves in $\hat f$, so they correspond to sections of $\ex C(\mathcal R(\hat f))\longrightarrow \ex F(\mathcal R(\hat f))$ which we shall again call $\{s_{i}\}$.
  The functor $\mathcal R$ applied to $\mathcal O$ is an open substack of $\Msw(\ex B)$.
Note that $\mathcal R(\hat f)^{*}T\ex B'$ is the pullback of $\hat f^{*}T\ex B$ under the map $\ex C(\mathcal R(\hat f))\longrightarrow \ex C(\hat f)$, so there is a unique $G$ invariant map $\mathcal R(F)$ which fits into the commutative diagram

\[\begin{array}{ccc}\mathcal R(\hat f)^{*}T\ex B'&\xrightarrow{\mathcal R(F)}&\ex B'
\\ \downarrow && \downarrow
\\ \hat f^{*}T\ex B&\xrightarrow{F}&\ex B
\end{array}\]

 $(\mathcal R(\hat f)/G,\mathcal R(F),\{s_{i}\})$ may not be an obstruction model for $\mathcal R(\mathcal O)$, however, it satisfies criterion \ref{core family main} from the definition of obstruction models on page \pageref{core family main}.  Given a family $\hat g$ of curves in $\mathcal R(\mathcal O)$, composition with $\ex B'\longrightarrow \ex B$ gives a family of curves in $\mathcal O$, which criterion \ref{core family main} for $(\hat f/G,F,\{s_{i}\})$ tells us comes with a unique fiberwise holomorphic map
 
  \[\begin{array}{ccc}(\ex C(\hat g),j)&\xrightarrow{\Phi_{\hat  g}}&({\ex C}(\hat f),j)/G
  \\ \downarrow & & \downarrow
  \\ \ex F(\hat g')&\longrightarrow &\ex F(\hat f)/G\end{array}\]
  and unique  $\C\infty 1$ section 
  \[\psi_{\hat g}:\ex C(\hat g)\longrightarrow\Phi_{\hat g}^{*}\lrb{\hat f^{*}T_{vert}\hat{\ex B}}\]
   which vanishes on the the pullback of the extra marked points so that 
   $F\circ\psi_{\hat g} $ is the composition of $g$ with the refinement map $\ex B'\longrightarrow \ex B$.
The universal property of $\mathcal R(\hat f)$ implies that there is a unique  lift of $\Phi_{\hat g}$ to  $\mathcal R\Phi_{\hat g}:\ex C(\hat g)\longrightarrow \ex C(\mathcal R(\hat f))$. We may also regard $\psi_{\hat g}$ as a section of $\mathcal R\Phi_{\hat g}^{*}\lrb{\hat f^{*}T_{vert}\hat{\ex B}}$. Then $\mathcal R(F)\circ\psi_{\hat g}$ is equal to 
$\hat g$. 

We shall now construct a family $\hat f'$ which should be thought of as adding some extra marked point sections to $\mathcal R(\hat f)$, and adding extra coordinates to allow the value of $\hat f'$ at these marked points to vary appropriately. In particular, we must add extra marked points where the edges of curves in  $\hat  f$ must be refined. 
If $\hat f$ was chosen small enough, we may  choose extra sections  $s'_{j}:\hat f'\longrightarrow \ex C(\mathcal R(\hat f))$ corresponding to extra marked points on the smooth part of curves in the domain of $\mathcal R(\hat f)$ disjoint from each other and $\{s_{i}\}$ so that the action of $G$ permutes these sections,  and each smooth component of each curve in $\mathcal R(\hat f)$ contains at least one of the  marked points from $\{s_{i},s'_{j}\}$.

Let $n'$ be the total number of sections $\{s_{i},s'_{j}\}$, and $n$ the number of sections in $\{s_{i}\}$. Denote by $s':\ex C(\mathcal R(\hat f))\longrightarrow \ex F(\mathcal R(\hat f^{+n'}))$ the section defined by taking all the above sections at once. We can now verify that the tropical part of $ev^{+n'}\circ s'$ is a complete map which gives an isomorphism
from any strata of $\totb{\ex F(\mathcal R(\hat f))}$ to a strata in the image of $\mathcal O$ in the tropical part of $\M^{+n'}\times (\ex B')^{n'}$ under $\totb{ev^{+n'}}$. The corresponding condition, (criterion \ref{lccb}) for $\hat f$ tells us that the tropical part of $ev^{+n}\circ s$ is a complete map so $\totb{ev^{+n'}\circ s'}$ is also complete. Criterion \ref{lccb} also tells us that  each polytope in $\totb{\ex F(\hat f)}$ may be described by taking the product of the polytopes from $\totb{\ex B}$ which contain the image of the marked points points in $\{s_{i}\}$ and a copy of $(0,\infty)$ for each internal edge, and then subjecting this polytope to the condition that there exists an appropriate tropical curve in $\totb{\ex B}$ with the corresponding data.  A polytope in $\ex F(\mathcal R(\hat f))$ is given by taking the subset of a polytope in $\ex F(\hat f)$ corresponding to  curves in $\mathcal R(\hat f)$ with a given combinatorial type.  As $\{s_{i},s_{j}'\}$ includes a marked point in each smooth component of curves in $\mathcal R(\hat f)$, the combinatorial type of curves in $\mathcal R(\hat f)$ is determined by which strata of $\totb{\ex B'}$ each  extra marked point is sent to, and which smooth component contains each of these extra marked points. Therefore, the combinatorial type of curves in $\mathcal R(\hat f)$ is determined by which strata $\totb{ev^{+n'}\circ s'}$ lands in. In other words, each polytope in $\totb{\ex F(\mathcal R(\hat f))}$ may be described by taking the product of the polytopes from $\totb{\ex B'}$ which contain the image of the marked points points in $\{s_{i},s'_{j}\}$ and a copy of $(0,\infty)$ for each internal edge, and then subjecting this polytope to the condition that there exists an appropriate tropical curve in $\totb{\ex B'}$ with the corresponding data. Therefore $\mathcal R(\hat f)$ with the section $s'$ satisifies criterion \ref{lccb}.

Given any family $\hat g\in\mathcal R(\mathcal O)$,  we may pull back $s'$ using $\mathcal R\Phi_{\hat g}^{+(n'-1)}:\ex C(\hat g^{+n'-1})\longrightarrow \ex C(\hat f^{+n'-1})$ to obtain  $\abs G$ sections $\ex F(\hat g)\longrightarrow \ex F(\hat g^{+n'})$. Composing these sections with $ev^{+n'}(\hat g)$ gives a subset of $\M^{+n'}\times (\ex B')^{n'}$. Define $\ex F(\hat f')$ to be the  subset of $\M^{+n'}\times (\ex B')^{n'}$ which is the union of the image of all such $\hat g$ from $\mathcal R(\mathcal O)$. Recall that  $G$ acts on the $n$ marked point sections by permutation. There is a corresponding action of $G$ on $\ex F(\hat f')$ given by relabeling marked points and permuting the corresponding coordinates on $\M^{+n'}\times(\ex B')^{n'}$.


 Note that criterion \ref{local core crit} implies that $\ex F(\hat f)$ may be regarded as a subset of $\M^{+n}\times \ex B^{+n}$.  Forgetting the last $n'-n$ marked points from curves in $\M^{+n'}$, and using the refinement map $\ex B'\longrightarrow \ex B$ on the first $n$ components of $(\ex B')^{n'}$ gives a map $\M^{+n'}\times (\ex B')^{n'}\longrightarrow \ex M^{+n}\times \ex B^{n}$. By construction, $\ex F(\hat f')$ lies in the inverse image of $\ex F(\hat f)$ under this map, and each fiber of the map $\ex F(\hat f')\longrightarrow \ex F(\hat f)$ is some open subset of a refinement of $(\ex B')^{n'-n}$  (times a constant in the other coordinates of $\M^{+n'}\times (\ex B')^{n'}$).  The map $\ex F(\hat f')\longrightarrow \ex F(\hat f)$ is $G$ equivariant.

 Define $\ex C(s'^{*}\hat f'^{+n'})\longrightarrow \ex F(\hat f')$ to be the restriction of  $\M^{+(n'+1)}\times \ex B^{n'}\longrightarrow \M^{+n'}\times \ex B^{n'}$ to $\ex F(\hat f')$. Note that the map $\ex F(\hat f')\longrightarrow \ex F(\hat f)$ lifts to a fiberwise holomorphic map 
 \[\begin{array}{ccc}\ex C(s'^{*}\hat f'^{+n'})&\longrightarrow &\ex C(\hat f)
 \\\downarrow&&\downarrow
 \\ \ex F(\hat f')&\longrightarrow &\ex F(\hat f)
 \end{array}\]
 We may pull $\hat f$ back over this map to obtain a family of curves $\ex C(s'^{*}\hat f'^{+n'})\longrightarrow \ex B$. As every curve in this family has extra marked points where edges meet places where $\ex B$ is refined, this family lifts without modification to a family $\ex C(s'^{*}\hat f'^{+n})\longrightarrow \ex B'$, and is therefore the pullback of $\mathcal R(\hat f)$ under a map 
 \[\begin{array}{ccc}\ex C(s'^{*}\hat f'^{+n'})&\longrightarrow &\ex C(\mathcal R(\hat f))
 \\\downarrow&&\downarrow
 \\ \ex F(\hat f')&\longrightarrow &\ex F(\mathcal R(\hat f))
 \end{array}\]

  Define $\ex C(\hat f')\longrightarrow \ex F(\hat f')$  by removing the extra $n'$ marked points from $\ex C(s'^{*}\hat f'^{+n'})$ (without modifying any components that become unstable when these marked points are removed.) Of course, $\ex C(\hat f')$ comes with $n'$ extra marked points which we may remember the location of.  The above diagram then factors into

\[\begin{array}{ccccc}\ex C(s'^{*}\hat f'^{+n'})&\longrightarrow&\ex C(\hat f')&\longrightarrow &\ex C(\mathcal R(\hat f))
 \\\downarrow&&\downarrow&&\downarrow
 \\ \ex F(\hat f')&\longrightarrow&\ex F(\hat f')&\longrightarrow &\ex F(\mathcal R(\hat f))
 \end{array}\]
 where the right hand square is a pullback diagram of abstract families of curves. There is a lift of the $G$ action on $\ex F(\hat f')$ to a $G$ action on $\ex C(\hat f')$ which permutes marked point sections so that the right hand square consists of $G$-equivariant maps. Denote by $\hat f'_{0}:\ex C(\hat f')\longrightarrow \ex B'$ the pullback of $\mathcal R(\hat f)$. We shall modify $\hat f'_{0}$ near the extra marked points below, after pulling back  the map $F:\hat f^{*}T\ex B\longrightarrow \ex B$ to $(\hat f'_{0})^{*}T\ex B'$.

  Define a $G$ invariant map $\hat f':\ex C(\hat f')\longrightarrow \ex B'$ as follows: Choose a neighborhood of each extra marked point section in $\ex C(\mathcal R(\hat f))$ so that no two neighborhoods intersect. Pull these neighborhoods back to $\ex C(\hat f')$ and  modify the pullback of $\mathcal R(f)$ in the neighborhood surrounding  each extra marked point  so that the   projection of $\ex C(\hat f'^{+n})\longrightarrow \ex F(\hat f')\subset \M^{+n'}\times (\ex B')^{n'}$ onto the $i$th copy of $\ex B'$ is equal to evaluation at the $i$th extra marked point.  This means that that if $s'$ indicates the section $\ex F(\hat f')\longrightarrow \ex F(\hat f^{+n'})$ coming from these extra marked points, then $ev^{+n}(\hat f')\circ s'$ is equal to the inclusion of $\ex F(\hat f)$ into $\mathcal M^{+n'}\times (\ex B')^{n'}$.  
 More specifically, we may choose $
 \hat f'$ as follows: Regard $\mathcal R (F)$ as giving a map $(\hat f'_{0})^{*}T\ex B'\longrightarrow \ex B'$. Choose $\hat f'$  so that it is in $\mathcal R(\mathcal O)$ and so that $\hat f'=\mathcal R(F)\circ \psi$ for some $G$ equivariant  section $\psi$ of $(\hat f_{0}')^{*}T\ex B'$ which vanishes outside the neighborhoods of our extra marked points, and which vanishes inside the neighborhood around a marked point if it vanishes at that marked point. In particular, $\hat f'$ contains $\mathcal R(\hat f)$ as a subfamily corresponding to those curves on which $\psi$ is identically $0$.
     
     There is a unique  $G$ equivariant map $F':\hat f'^{*}T\ex B'\longrightarrow \ex B'$ satisfying criterion \ref{F crit}  given by a fiberwise affine map from $\hat f'^{*}T\ex B'$ to $(\hat f'_{0})^{*}T\ex B'$ followed by $\mathcal R(F)$, in particular, 
     \[F'(\phi):=\mathcal R(F)(D\mathcal R(F)(\psi)^{-1}(\phi)+\psi)\]
     
   Now $(\hat f',F',\{s_{i},s_{j}'\}) $ is a core family for $\mathcal R(\mathcal O)$. 
     
  We can now check that $\hat f'$ satisfies criterion \ref{local core crit} from the definition of a core family: 
As $\hat f$ is a core family, $ev^{+n}(\hat f)$ restricted to some neighborhood of the image of $s$ is an equidimensional embedding, so $ev^{+n}(\mathcal R(\hat f))$ is also an equidimensional embedding restricted to the same neighborhood of $s$, and $ev^{+n'}(\mathcal R(\hat f))$ composed with the projection $\M^{+n'}\times (\ex B')^{n'}\longrightarrow \M^{+n'}\times (\ex B')^{n}$ is also an equidimensional embedding when restricted to the lift of this neighborhood. On the other hand, $ev^{+n'}(\hat f')\circ s'$ is an embedding with derivative that is surjective onto the fibers of the projection $\M^{+n'}\times (\ex B')^{n'}\longrightarrow \M^{+n'}\times (\ex B')^{n}$. As the domain and range have the same dimension, it follows that in some neighborhood of the image of $s'$, $ev^{+n'}(\hat f')$ is an equidimensional embedding, so $\hat f'$ with the section $s'$ satisfies criterion \ref{lcca}. $\hat f'$ with $s'$ also satisfies criterion \ref{lccb} because $\mathcal R(\hat f)$ with $s'$ does, and $\hat f'$ and $\mathcal R(\hat f)$ have the same tropical part.

$(\hat f',F',\{s_{i},s'_{j}\})$ also satisfies criterion \ref{core family main} for being a core family. In particular, let $\hat g$ be a family in $\mathcal R(\mathcal O)$, then as discussed already, there exists a unique fiberwise holomorphic map

  \[\begin{array}{ccc}(\ex C(\hat g),j)&\xrightarrow{\mathcal R\Phi_{\hat  g}}&{\ex C}(\mathcal R(\hat f))/G
  \\ \downarrow & & \downarrow
  \\ \ex F(\hat g')&\longrightarrow &\ex F(\mathcal R(\hat f))/G\end{array}\]
  and unique $\C\infty 1$ section $\psi_{\hat g}$ which vanishes on the pullback of the extra marked points 
  \[\psi_{\hat g}:\ex C(\hat g)\longrightarrow\mathcal R\Phi_{\hat g}^{*}\lrb{\hat f^{*}T_{vert}\hat{\ex B}}\]
   so that   $\mathcal R(F)\circ\psi_{\hat g}=\hat g$. 
We may lift $\mathcal R\Phi_{\hat g}$ to a fiberwise holomorphic map 

\[\begin{array}{ccc}(\ex C(\hat g),j)&\xrightarrow{\mathcal R\Phi_{\hat  g}}&{\ex C}(\mathcal R(\hat f))/G
  \\ \downarrow & & \downarrow
  \\ \ex F(\hat g')&\longrightarrow &\ex F(\mathcal R(\hat f))/G\end{array}\]
The other criteria for $(\hat f',F',\{s_{i},s'_{j}\})$ to be a core family are easily seen to be satisfied, so $(\hat f',F',\{s_{i},s'_{j}\})$ is a core family for $\mathcal R(\mathcal O)$.

To describe an obstruction model with core family $(\hat f',F',\{s_{i},s'_{j}\})$, we also need a trivialization in the sense of definition \ref{trivialization def}. For this we may pull back the trivialization from $\hat f$ to $\hat f'_{0}$, then use the induced trivialization on $\hat f'$. The induced trivialization is described in \cite{reg}. The  construction agrees with our construction of $F'$ from $F$. We may use this induced trivialization to pullback $V$ from $(\hat f,V)$  to get a $G$-invariant pre obstruction model $(\hat f',V)$. Restricted to some neighborhood of $\mathcal R(\hat f)\subset \hat f'$,  $D\dbar$ is injective and has image intersecting $V$ on only at $0$, so on this neighborhood we may extend $V$ to a $G$ invariant pre obstruction model $(\hat f', V')$ so that $D\dbar$ is injective and has image complementary to $V'$.  Theorem \ref{regularity theorem} implies that we may modify $(\hat f'/G,V')$ to be a core family on some neighborhood of $\mathcal R(\hat f)$ in $\Msw$.  To avoid further notational complications, simply call this modified core family $(\hat f'/G,V')$. Note that this modified core family $\hat f'$ contains $\mathcal R(\hat f)$ because $\dbar\mathcal R(\hat f)$ is a section of $V\subset V'$. By choosing $\mathcal O$ small enough to begin with, we may assume that  $(\hat f'/G,V')$ is a core family for $\mathcal R(\mathcal O)$. 

Therefore, we may choose a locally finite cover of the moduli space of stable holomorphic curves in $\ex B$ with core families $(\hat f/G,V)$ and a corresponding cover of the moduli space of stable holomorphic curves in $\ex B'$ with the corresponding  core families $(\hat f'/G,V')$, and use these to construct the virtual moduli space as in section \ref{virtual class}. After choosing the compact subsets $C$ of the core families  $\hat f$ in which simple perturbations will be supported,  choose an open neighborhood $\mathcal U$ of the moduli space of holomorphic  curves with the property that $\mathcal R(\mathcal U)\cap \mathcal R(\mathcal O)$ projects to a subset of $\hat f'$ which when intersected with the inverse image of $C$ under the map $\ex F(\hat f')\longrightarrow \ex F(\hat f)$, is contained in a compact subset.  By using small simple perturbations supported inside $C\subset\ex F(\hat f)$, we may arrange that $\mathcal M$ is in $\mathcal U$. Then we may use the pullback of the same simple perturbations under the maps $\ex C(\hat f')\longrightarrow \ex C(\hat f)$ and $\ex B'\longrightarrow \ex B$, then cut them off so that they are compactly supported in $\ex F(\hat f')$, but still describe the same multi perturbation on $\mathcal R(\mathcal U)$. If small enough perturbations are used, then all solutions must be in $\mathcal R(\mathcal U)$, so the solution set of the corresponding multiperturbation will be $\mathcal R(\mathcal M)$. 

It remains to check that these perturbations give the required transversality in $\mathcal R(\mathcal U)$, and that the orientation on $\mathcal R(\mathcal M)$ agrees with the orientation on the corresponding virtual moduli space of curves in $\ex B'$. Consider a curve $f'$ in $\hat f'$ which projects to a curve $f$ in $\hat f$ under the map $\ex C(\hat f')\longrightarrow \ex C(\hat f)$. So $\ex C( f')$ is some refinement of $\ex C(f)$, and using our trivialization we may regard $ f'^{*}T\ex B'$ as the pullback of $ f^{*}T\ex B$ and $Y( f')$ as the pullback of $Y(f)$. On the edges of $\ex C(f)$,  $f^{*}T\ex B$ and $Y(f)$  are naturally trivial bundles, and a section of $f'^{*}T\ex B'$ or $Y(f')$ is the pullback of a section of $f^{*}T\ex B$ or $Y(f)$ if and only if that section is constant on the extra smooth components of $\ex C(f')$ where edges have been refined.  As sections have to vanish on edges in order to be in $Y^{\infty,\underline 1}(f)$, it follows that $Y^{\infty,\underline 1}(f')$ is equal to the product of a space of sections for each smooth component, so it is equal to the product of   $Y^{\infty,\underline 1}(f)$ with a space of sections for each extra smooth component. Each extra smooth component is a twice punctured sphere, and if $f'=\mathcal R(f)$, $D\dbar$ on sections of $f'^{*}T\ex B'$ restricted to this sphere is just the usual $\dbar$ equation on maps to $\mathbb C^{k}$, so has kernel equal to the constant sections and is surjective. For $f'$ close enough to $\mathcal R(f)$, it follows that $D\dbar$ restricted to sections of $f'^{*}T\ex B'$ which vanish on the marked points corresponding to $\{s_{i}\}$ ( but non necessarily $\{s_{j}'\}$)  is injective and has image complementary to $V$ in $Y^{\infty,\underline 1}(f')$. It follows that $\dbar \hat f'$ as a section of $V'$ is transverse to $V\subset V'$ (and has intersection equal to $\mathcal R(\hat f)$).
If the simple perturbations used to define the moduli space are chosen small enough, it follows that $\dbar$ of the corresponding solutions given by Theorem \ref{Multi solution} are also transverse to $V\subset V'$. As their intersection with $V$ is equal to lift of the corresponding solutions in the core family $(\hat f/G,V)$, it follows that the transversality  required for the construction of the virtual moduli space in section \ref{virtual class} will then hold. 

As argued in section \ref{virtual class}, adding extra marked points to a core family and adding extra parameters corresponding to the image of those extra marked points does not affect the orientation on the virtual moduli space. As $\hat f'$ is constructed from $\mathcal R(\hat f)$ simply by adding these extra parameters, it follows that the orientation on $\mathcal R(\mathcal M)$ as a virtual moduli space is the same as the orientation coming from $\mathcal M$. 

\stop

\begin{example}\label{toric like divisors}\end{example} Suppose that $M^{2n}$ is a connected compact symplectic manifold with codimension $2$ embedded symplectic submanifolds $N_{i}$ which intersect each other symplectically orthogonally. Suppose further that some $N_{j}$ has a neighborhood in $M$ so that restricted to this neighborhood, $M$ with the submanfolds $N_{i}$ is symplectomorphic to an open subset of a toric symplectic manifold with its toric boundary divisors.
 
 Then the explosion of $(M,\{N_{i}\})$ discussed in \cite{iec} gives an exploded manifold $\ex M$ with tropical part $\totb{\ex M}$ equal to a cone. The open neighborhood of $N_{j}$ discussed above corresponds to an open subset of $\ex M$ which is isomorphic to a refinement of $\ex T^{n}$, so for computing the contribution of a tropical curve to Gromov Witten invariants using Theorem \ref{gluing formula} or Theorem \ref{gluing cobordism}, vertices contained in the ray corresponding to $N_{j}$ or any adjacent positive dimensional strata  can be regarded as being inside a refinement of $\ex T^{n}$. 
 
 Note that considering the ray corresponding to $ N_{j}$ together with all adjacent positive dimensional strata as being a subset of the tropical part of a refinement of $\ex T^{n}$ allows us to put a natural affine structure on these strata which extends over faces. If more of our symplectic submanifolds $N_{i}$ obey the same condition of having a neighborhood equal to an open subset of a toric symplectic manifold with its toric boundary strata, then we can consider part of the tropical part of $\ex M$ corresponding to these $N_{i}$ and all adjacent positive dimensional strata as a subdivision of a cone with an integral affine structure.  In light of Theorem \ref{refinement theorem}, the subdivision of this cone is not important for the computation of Gromov Witten invariants, but the integral affine structure is important.

\begin{example} There does not exist a compact 4-dimensional symplectic manifold $M$ which contains $3$ embedded symplectic spheres  which intersect each other once symplectically orthogonally, and which have self intersection numbers 1,1 and 2.\end{example} To see why this is not possible, suppose that such a manifold existed, and let $\ex M$ be the explosion of $M$ relative to these three symplectic submanifolds discussed in \cite{iec}. Example \ref{toric like divisors} implies that we may regard the tropical part of $\ex M$ as a subdivision of a two dimensional integral affine cone which in this case has monodromy around $0$. 
In particular, removing the strata corresponding to the sphere with self intersection 2, the tropical part or $\ex M$ should be regarded as $\mathbb R^{2}$ minus a `cut' along the ray generated by $(1,0)$, subdivided by the rays generated by $(0,1)$ and $(-1,-1)$. This integral affine structure can be continued over the cut by identifying it and adjacent strata with the union of the cone in $\mathbb R^{2}$ generated by $(1,0)$ and $(0,1)$ with the cone generated by $(1,0)$ and $(-2,-1)$. In particular, with this integral affine structure, a straight line entering the bottom of the cut in direction $(0,1)$ will exit in direction $(1,1)$.

Consider the moduli space of holomorphic curves in $\ex M$ with genus $0$ and tropical part with $3$ punctures having momentum $(0,1)$, $(0,-1)$ and $0$ respectively. Note that tropical curves in $\ex M$ which are the image of holomorphic curves obey a conservation of momentum condition. The only such tropical curves obeying our conditions are therefore refinements of vertical lines contained in the left hand side of $\mathbb R^{2}$ identified with the tropical part of $\ex M$. Section \ref{curves in Tn} together with Theorem \ref{refinement theorem} imply that the integral over the corresponding virtual moduli space  of pullback under the evaluation map at the third puncture of the Poincare dual to a point with tropical part contained in the left hand side of $\mathbb R^{2}$ is $1$, but the corresponding integral for the Poincare dual to  a point with tropical part contained in the right hand side of $\mathbb R^{2}$ must be $0$. If such a manifold $M$ existed, Theorem \ref{pd class} would imply that two integrals are be equal, so no such symplectic manifold $M$ exists.

 \subsection{The case of a nice cokernel}
 
 \
 
 \
 
 Sometimes, (a connected component of) the moduli space of holomorphic curves will be an exploded  manifold (or orbifold) with dimension $k$ greater than expected, and come with a natural $k$-dimensional obstruction bundle. In these cases, integrating the pullback of closed forms over the virtual moduli space is equivalent to integrating the wedge product of the Euler class of this bundle with the pullback of these forms  over the moduli space of holomorphic curves. (The case when $k=0$ is the case of `transversality'.) 
 
 Up to this point, we have defined the linearization of the $\dbar$ operator in the context of a family with a trivialization. We now give a more natural definition of the tangent space to the image of $\dbar$ at a holomorphic curve $f$ projected to $Y^{\infty,\underline 1}(f)$ which is a closed subspace $\tdbar f\subset Y^{\infty,\underline 1}(f)$ with finite codimension. The annihilator of $\tdbar f$ can be regarded as the cokernel of the linearized $\dbar$ operator at $f$. 
 
 \begin{defn}(Tangent space to $\dbar$) Let $\hat f$ be a family of curves in $\Msw$ containing a holomorphic curve $f$, let $\nabla$ be a $\C\infty1$ connection on $Y(\hat f)$, and let $v$ be a vector field on $\ex C(\hat f)$ which projects to a vector field on $\ex F(\hat f)$. Then $\nabla_{v}\dbar \hat f\rvert_{\ex C(f)}$ is in  $Y^{\infty,\underline 1}(f)$. 
Let 
\[\tdbar f:=\{\nabla_{v}\dbar \hat f\rvert_{\ex C(f)}\}\subset Y^{\infty,\underline 1}(f)\] 
 for all choices of $\hat f$, $v$ and $\nabla$ satisfying the conditions above.
 \end{defn}  
 
 \begin{lemma}\label{tdbar}Given any stable holomorphic curve $f$ in a basic exploded manifold $\ex B$, $\tdbar f$ is a closed linear subspace of $Y^{\infty,\underline 1}(f)$ with finite codimension. 
 
 If $\hat g$ is a family of stable holomorphic cures in $\ex B$ so that $\tdbar f$ has fixed codimension for all curves $f$ in $\hat g$, then there is a  $\C\infty1$ vector bundle $E$ over $\ex F(\hat g)$ with fiber over $f$ given by  
 \[E(f):=Y^{\infty,\underline 1}(f)/\tdbar f\]
 so that sections in $Y^{\infty,\underline 1}(\hat g)$ project to $\C\infty1$ sections of $E$.
   \end{lemma}

\pf Let $f$ be a stable holomorphic curve. If $\ex C( f)=\ex T$, then $Y^{\infty,\underline 1}(f)$ is zero dimensional, so this lemma holds trivially. In all other cases, Theorem \ref{construct obstruction model} tells us that there is an obstruction model $(\hat f/G,V)$ containing $f$. Lemma \ref{f replacement} tells us that the image of $D\dbar(f)$ is a closed linear subspace of $Y^{\infty,\underline 1}(f)$ and has finite codimension. Of course, $\tdbar f$ contains this image of $D\dbar(f)$. To obtain the entire $\tdbar f$, we may restrict to families $\hat f'$ parametrized by $\mathbb R$ and contained in the open substack for which $(\hat f/G,V)$ is an obstruction model. Assume that our holomorphic curve $f$ is the curve over $0$.  Using the fact that $\hat f/G$ is a core family in the sense of Definition \ref{core family}, we get a map
  \[\begin{array}{ccc}(\ex C(\hat f'),j)&\xrightarrow{\Phi_{\hat f'}}&({\ex C}(\hat f),j)/G
  \\ \downarrow & & \downarrow
  \\ \ex F(\hat f')&\longrightarrow &\ex F(\hat f)/G\end{array}\]
  and unique $\C\infty 1$ section 
  \[\psi_{\hat f'}:\ex C(\hat f')\longrightarrow\Phi_{\hat f'}^{*}\lrb{\hat f^{*}T{\ex B}}\]
    which vanishes on the pullback of marked points, so that 
   
\[\hat f'=F\circ\psi_{\hat f'}\]

We may resolve the $G$-ambiguity of the above map $\Phi_{\hat f'}$ and choose some lift to a map to $\ex C(\hat f)$. We may also lift the resulting map to a map 
 \[\Phi:\ex C(\hat f')\longrightarrow \ex C(\hat f\times \mathbb R)\] which is the identity on the $\mathbb R$ factor when remembering that $\ex F(\hat f')=\mathbb R$. Then there exists a section $\psi_{t}$ of $(\hat f\times \mathbb R)^{*}T\ex B$ so that 
\[\psi_{\hat f'}=\Phi^{*}\psi_{t}\]
Regarding $t$ as the coordinate on $\mathbb R$, we may think of  $\psi_{t}$ as defining a family of sections of $X^{\infty,\underline 1}$. We may assume that  $\psi_{0}$ is the zero section.  $\Phi$ defines an inclusion of $\hat f'$ into $F(\psi_{t})$.
Similarly, given any vectorfield $v$ on $\ex C(\hat f')$ which projects to a vectorfield on $\ex F(\hat f')$, there exists a vectorfield $v'$ on $\ex C(\hat f\times \mathbb R)$ projecting to a vectorfield on $\ex C(\hat f'\times \mathbb R)$ so that 
\[\nabla_{v}\dbar\hat f'\rvert_{\ex C(f)}=\nabla_{v'}\dbar \psi_{t}\rvert_{\ex C(f)}\] 
In the above, we may use any connection $\nabla$, because $f$ is holomorphic and the above expression is actually independent of choice of $\nabla$. We may write $v'$ as $s\frac \partial{\partial t}+v''(t)$, where  $s$ is a function of $\ex F(\hat f)$, and  for each $t$, $v''(t)$ may be regarded as a vectorfield on $\ex C(\hat f)$ which projects to a vectorfield on $\ex F(\hat f)$. 
Then
\[\begin{split}\nabla_{v'}\dbar \psi_{t}\rvert_{\ex C(f)}&=\nabla_{s\frac\partial{\partial t}}\dbar\psi_{t}\rvert_{\ex C(\hat f)}
+\nabla_{v''(t)}\dbar \psi_{t}\rvert_{\ex C(\hat f)}
\\&= D\dbar(f)\lrb{s(f)\frac \partial{\partial t}\psi_{t}(f)}+\nabla_{v''(0)}\dbar\hat f\rvert_{\ex C(\hat f)}\end{split}\] 
As we may construct our family $\hat f'$ so that $s(f)\frac \partial{\partial t}\psi_{t}(f)$ and $v''(0)\rvert_{\ex C(f)}$ are whatever we like, it follows that $\tdbar f$ is equal to the linear span of the image of $D\dbar(f)$ and the linear subspace of $V(f)\subset Y^{\infty, \underline 1}(f)$ which is the image of the derivative of the section $\dbar:\ex F(\hat f)\longrightarrow \ex V$. Therefore, $\tdbar f$ is indeed a closed linear subspace of $Y^{\infty,\underline 1}(f)$ which has finite codimension.

Now suppose that $\hat g$ is a family of stable holomorphic curves so that for all $f$ in $\hat g$, $\tdbar f$ has codimension $k$. To prove that the vector spaces $Y^{\infty,\underline 1}(f)/\tdbar f$ for all $f$ in $\hat g$ have a natural structure of a $k$ dimensional $\C\infty1$ vector bundle, it suffices to work locally, so we may assume that $\hat g$ is contained in our obstruction model $(\hat f/G,V)$, and that there is a map of curves $\hat g\longrightarrow \hat f$. As $D\dbar(f)$ is complementary to $V$, our assumption on the codimension of $\tdbar f$ and the above characterization of $\tdbar f$ imply that the image of the derivative of the section $\dbar:\ex F(\hat f)\longrightarrow V$ has codimension $k$ at all curves $f$ in $\hat g\subset \hat f$. As the image of the derivative of this section $\dbar$ has constant codimension, its image $V'\subset V$ is a $\C\infty1$ sub vector bundle of $V$ restricted to the image of $\hat g$. We may pull back $V$ and $V'$ to $\C\infty1$ vector bundles $V(\hat g)$ and $V'(\hat g)$ over $\ex F(\hat g)$. 

As $V(f)$ is complementary to the image of $D\dbar(f)$, the inclusion $V(f)\longrightarrow Y^{\infty,\underline 1}(f)$ induces an isomorphism 
\[V(f)/V'(f)\longrightarrow Y^{\infty,\underline 1}(f)/\tdbar f\]
Our desired $\C\infty1$ vector bundle $E$ is equal to $V(\hat g)/V'(\hat g)$.
Note that $(\hat g,V(\hat g))$ is a pre obstruction model so that $D\dbar(f)$ is injective and complementary to $V(f)$ for all $f$ in $\hat g$. Applying the linear version of Theorem \ref{regularity theorem} part \ref{rt2}, (proved separately in \cite{reg}) gives that for any section $\theta$ in $Y^{\infty,\underline 1}(\hat g)$, there exists a unique section $\nu$  in $X^{\infty,\underline 1}(\hat g)$ and $\C\infty1$ section $v$ of $V(\hat g)$ so that 
\[D\dbar(\hat g)(\nu)=\theta-v\]
Then the section of $E$ corresponding to $\theta$ is equal to the section of $E$ corresponding to $v$, which gives a $\C\infty1$ section of $V(\hat g)/V'(\hat g)$. Therefore, the $\C\infty1$ structure on $E$ given by $V(\hat g)/V'(\hat g)$ has the desired property.

 \stop

 \begin{thm}\label{obstruction bundle} Suppose that a connected component $\mathcal M$ of the moduli space of curves in a basic exploded manifold $\ex B$ for which Gromov compactness holds satisfies the following:
 
 \begin{itemize}
 \item For all curves $f$ in $\mathcal M$, $\tdbar f\subset Y^{\infty,\underline 1}(f)$ has codimension $k$.
 \item Either $\mathcal M$ is an orientable exploded orbifold with dimension $k$ greater than the expected dimension of the virtual moduli space, or $k=0$.
  
 \end{itemize}
 
 Then $\mathcal M$ is a complete orientable exploded manifold or orbifold, and there is a $\C\infty1$ vector bundle $E$ over $\mathcal M$ with the fibers
 \[E(f):=Y^{\infty,\underline 1}(f)/\tdbar f\]
 and the natural $\C\infty1$ structure from Lemma \ref{tdbar}. Given an orientation of $\mathcal M$, there is a natural orientation for $E$ so that the following holds:
 
  Given any $\C\infty1$ map $\psi:\Msw\longrightarrow X$ and closed differential form $\alpha\in\ro^{*}(X)$, the integral of $\psi^{*}\alpha$ over the component of the virtual moduli space corresponding to $\mathcal M$ is equal to 
  \[\int_{\mathcal M}\psi^{*}\alpha\wedge e(E)\]
  where $e(E)$ is the Euler class of the vector bundle $E$.
 
 \end{thm}
 
 \pf
 
 Note that we may construct the virtual moduli space of curves using one set of obstruction models covering $\mathcal M$ and another disjoint set of obstruction models which cover all other holomorphic curves. The construction of the virtual moduli space with these two sets of obstruction models is then completely independent, so it makes sense to talk about the component of the virtual moduli space corresponding to $\mathcal M$ as that component which is contained in the image of the union of the obstruction models used to cover $\mathcal M$.
 
 Now consider one of our obstruction models $(\hat f/G,V)$ used to cover $\mathcal M$. We shall assume that we have chosen these obstruction models small enough that the only holomorphic curves in some extension of $\hat f$ are in $\mathcal M$. The holomorphic curves in $\hat f$ therefore form a $G$-fold cover $U$ of  an open subset of $\mathcal M$.
 
  If $k=0$, then the characterization of $\tdbar f$ given in the proof of Lemma \ref{tdbar} implies that section $\dbar:\ex F(\hat f)\longrightarrow V$ is transverse to the zero section, and it follows that $\mathcal M$ must be a complete exploded orbifold of the dimension expected for the virtual moduli space. We may then give $\mathcal M$ the orientation from the oriented intersection of $\dbar$ with the zero section using the orientation on $\ex F(\hat f)$ and the relative orientation of $V$. As noted in section \ref{virtual class}, this gives a well defined global orientation on $\mathcal M$. 
 
  If on the other hand $\mathcal M$ is of dimension $k$ larger than the expected dimension, the set of holomorphic curves  in $\ex F(\hat f)$ is a closed submanifold of  $U\subset \ex F(\hat f)$ with dimension $k$ larger than the dimension of $\ex F(\hat f)$ minus the codimension of $V\longrightarrow \ex F(\hat f)$. The characterization of $\tdbar f$ from the proof of Lemma \ref{tdbar} implies that we may choose a   $k$ dimensional $\C\infty1$ sub bundle $E_{\hat f}\subset V$ so that for all holomorphic  $f$, $E_{\hat f}(f)$ is complementary to $\tdbar f$. ($E_{\hat f}$ is naturally isomorphic to the pullback of our obstruction bundle $E$  to $\hat f$.) It follows that in some neighborhood $U\subset \ex F(\hat f)$, the intersection of the section $\dbar$ with $E_{\hat f}$ is transverse. As this intersection is also the same dimension as $U$, it must be equal to U ( at least, it must be equal to $U$ when  restricted to some neighborhood of $U$.) Therefore, the $U$ is actually a complete sub exploded manifold of $\ex F(\hat f)$.
 
  As Gromov compactness holds for $\ex B$, we already know that $\mathcal M$ is compact, so it follows that $\mathcal M$ is complete. Given an orientation of $\mathcal M$, we may orient $E_{\hat f}$ using the exact sequence
  \[0\longrightarrow T_{f}\mathcal M\longrightarrow T_{f}\ex F(\hat f)\xrightarrow{d\dbar} V(f)\longrightarrow E_{\hat f}(f)\longrightarrow 0\]
 
  So restricted to a small enough open neighborhood of the holomorphic curves in $\hat f$, the inverse image of $\mathcal M$ in $\ex F(\hat f)$ is equal to the oriented sub manifold $U:=\dbar^{-1}(E_{\hat f})$. That this gives a global orientation for $E$ may be proved in the same way as the proof that the orientation for the virtual moduli space is well defined in section \ref{virtual class}.
 
Choose compactly supported simple perturbations parametrized by $\hat f$ so that the resulting multiperturbation $\mathfrak P$ restricted to the family of curves $\mathcal M$ gives a multisection of $Y(\mathcal M)$ which  corresponds to a multisection $s$ of $E$ which is transverse to the zero section.

 Now consider perturbing the $\dbar$ equation using $\epsilon$ times the above multiperturbation $\mathfrak P$. Apply Theorem \ref{Multi solution} and let $(\nu_{\epsilon},\dbar'\nu_{\epsilon})$ be the corresponding family of multisections which are the solutions mod $V$ parametrized by $\hat f$. Theorem \ref{Multi solution} tells us that for $\epsilon$ small enough, this is a $\C\infty1$ family of multisections. 
 
 We have that locally around a curve $f$ in $\hat f/G$,  
 \[(\nu_{\epsilon},\dbar'\nu_{\epsilon})=\sum_{i=1^{n}}\frac 1nt^{(\nu_{i,\epsilon},\dbar'\nu_{i,\epsilon})}\]
 where 
 \[\mathfrak P(\nu_{i,\epsilon})=\sum_{j=1}^{n}\frac 1nt^{\mathfrak P_{i,j,\epsilon}}\]
 and $\dbar'\nu_{i,\epsilon}$ is a section of $V$ satisfying
 \[\dbar'\nu_{i,\epsilon}:=\dbar \nu_{i,\epsilon}-\epsilon\mathfrak P_{i,i,\epsilon}\]
Restrict $\hat f$ to a small enough compactly contained $G$-invariant open subset so that the above holds, and simply relabel this subset $\hat f$ to avoid notational complications. Each $\mathfrak P_{i,j,\epsilon}$  corresponds to a choice of map of $F(\nu_{i,\epsilon})$ to an obstruction model, which removes the ambiguity from the automorphisms of that obstruction model.  The proof of Theorem \ref{Multi solution} involves locally extending the choices of such maps for $\hat f$ to choices for $F(\nu)$ so that the index $j$ of $\mathfrak P_{i,j,\epsilon}$ corresponds to a particular local choice for $\hat f$. As $\nu_{i,\epsilon}=0$, it follows that 
\[\mathfrak P(\hat f)=\sum_{i}\mathfrak P_{i,i,0}\]

 As $\dbar'\nu_{i,0}=\dbar\hat f$, and $\dbar\hat f$ is transverse to $E_{f}$,   $\dbar'\nu_{i,\epsilon}$  remains transverse to $E_{\hat f}\subset V$ for $\epsilon$ small enough, and we may parametrize the intersection with $E_{\hat f}$ with a weighted branched map from $U\subset \ex F(\hat f)$ (and this intersection with $E_{\hat f}$ will converge to $U$ in $\C\infty1$ as $\epsilon\to 0$.)

  Recall that on $\mathcal M$, there is a map $\pi_{E}$ from  $Y^{\infty,\underline 1}(\hat M)$ to $\C\infty1$ sections of $E$ corresponding to the identification $E(f):=Y^{\infty,\underline 1}(f)/\tdbar f$.
  Choose some $\C\infty1$ bundle map $\pi_{E_{\hat f}}:V\longrightarrow E_{\hat f}$ which is a projection that which on curves $f$ in $U$ has kernel equal to the image of the derivative of $d\dbar$, so on these curves $\pi_{E_{\hat f}}$ is equal to the restriction of $\pi_{E}$ to $V(f)$.
  
   Note that  $\pi_{E_{\hat f}}\dbar \hat f$ vanishes to first order on $U$. Similarly, note that on $U$, $\pi_{E_{\hat f}}(f)$ applied to the image of $D\dbar(f)$ is $0$, therefore restricted to $\epsilon=0$ and $U$, 
   \[\frac\partial{\partial \epsilon}\pi_{E_{f}}\dbar'\nu_{i,\epsilon}=-\pi_{E}\mathfrak P_{i,i,0}\]
   
   To summarize, we have that $\dbar'\nu_{i,\epsilon}$ is transverse to $E_{f}\subset V$ for $\epsilon$ small enough, and restricted to the intersection of $\dbar'\nu_{i,\epsilon}$ with $E_{f}$, $\dbar'\nu_{i,0}=0$ and the partial derivative with respect to $\epsilon$ of $\dbar'\nu_{i,\epsilon}$ is $-\pi_{E}\mathfrak P_{i,i,0}:=-s_{i}$ where the pullback to $U$ of our  multisection $s$ of $E$ is given by  
   \[\sum_{i=1}^{n}\frac 1nt^{s_{i}}\]
   
    For $\epsilon$ small enough, we may therefore approximate the integral of the pullback of any differential form over  $\dbar'\nu_{i,\epsilon}=0$ with the integral over the intersection of $s_{i}$ with the zero section.

For $\epsilon$ small enough,  $\dbar'\nu_{i,\epsilon}$ is transverse to the zero section, so we may change our simple perturbations slightly for each $\epsilon$  so that the resulting family of solutions to $\dbar'\nu_{i,\epsilon}=0$ is fixed point free, so the resulting solutions fit together to form the virtual moduli space.   Given any closed form  $\alpha\in\ro^{*}X$, the integral of  $\psi^{*}\alpha$ over (our component of) the virtual moduli space is independent of $\epsilon$, and locally given by the sum of the integrals over the intersection of $\dbar\nu_{i,\epsilon}$ with $0$ divided by $n\abs G$. Similarly, the integral of $\psi^{*}\alpha$ over the weighted branched sub exploded orbifold of $\mathcal M$  defined by the intersection of 
 our multisection $s$ of $E$  with the zero section is locally given by the sum of integrals over $s_{i}\cap 0$ divided by $n\abs G$.  The above arguments imply that as we may make the above two integrals locally as close as we like by choosing $\epsilon$ small, therefore the fact that this integral is independent of $\epsilon$ if $\alpha$ is closed implies that these two integrals are equal. Therefore, the integral of $\psi^{*}\alpha$ over our component of the virtual moduli space is equal to

\[\int_{\mathcal M}\psi^{*}\alpha\wedge e(E)\]
  where $e(E)$ is the Euler class of the vector bundle $E$.
\stop

 \begin{example}[Curves mapping to a point]
 \end{example} 
 
 Consider the component of the moduli space of holomorphic curves in $\ex B$ consisting of curves with genus $g$ and $n$ punctures which map to a point in $\ex B$. This component of the moduli space of holomorphic curves is equal to $\M_{g,n}\times \ex B$, which has dimension $g$ times the dimension of $\ex B$ greater than expected, and has a nice obstruction bundle which we shall now describe. 
 
  Given a particular curve $f$ with domain $\ex C$ in $\M_{g,n}$ and a point $p$ in $\ex B$, let $E^{*}(f)$ be the space of holomorphic sections of $T^{*}\ex C\otimes_{\mathbb C} T^{*}_{p}\ex B$ which vanish on any external edge of $\ex C$. Given any $\theta\in E^{*}(f)$ and $\C\infty1$ section $\alpha$ of $Y(f)$ which vanishes on edges of $\ex C$, we may regard $\alpha\wedge \theta$ as a two form on $\ex C$ which vanishes on all edges of $\ex C$. As $\alpha\wedge \theta$ vanishes on edges of $\ex C$, 
  \[\alpha\mapsto \int_{\ex C}\alpha \wedge \theta\]
  gives a linear functional on the space of $\C\infty1$ sections of $Y(f)$ which vanish on edges of $\ex C$, so we may regard $E^{*}(f)$ as a linear subspace of the dual of $Y^{\infty,\underline 1}(f)$. We shall now check that $E^{*}(f)$ is a subspace of the cokernel of $D\dbar$.

   Given any $\C\infty1$ map $\nu:\ex C\longrightarrow T_{p}\ex B$, we may regard $\nu\theta$ as a $\C\infty1$ one form on $\ex C$ which vanishes on external edges of $\ex C$. Then $\dbar\nu\wedge \theta$ is equal to $d(\nu\theta)$. As $\nu\theta$ may not vanish on integral vectors, we can not apply the version of Stokes' theorem proved in \cite{dre} directly, however $\int_{\ex C}\dbar \nu\wedge \theta$ does vanish. This is because the usual Stokes' theorem applied to each smooth component of $\ex C$ gives that the integral of $d (\nu\theta)$ is equal to the sum of the limits of the integral of $\nu\theta$ over suitably oriented loops around punctures as those loops are sent into the edges of $\ex C$. As $\nu\theta$ vanishes on external edges of $\ex C$, the contribution to the integral from external edges disappears. On the other hand, the contribution from each end of an internal edge cancels out, so the sum of the integral of $d(\nu\theta)$ over all smooth components of $\ex C$ is $0$. 
   
    Note that $\theta$ on the smooth part of $\ex C$ is a holomorphic one form with values in $T_{p}\ex B$ that has simple poles with opposite residues at each side of a node corresponding to an internal edge of $\ex C$, and which is bounded and hence smooth  at punctures corresponding to external edges of $\ex C$. The dimension of $E^{*}(f)$ is equal to $g(\dim \ex B)$. The only holomorphic maps from $\ex C$ to $T_{p}\ex B$ are the constant maps, so the kernel of the $\dbar$ operator has dimension $\dim \ex B $. The index of the $\dbar$ operator acting on the space of $\C\infty1$ maps from $\ex C$ to $T_{p}\ex B$ is $(1-g)(\dim \ex B)$, so $E^{*}(f)$ is the cokernel of this $\dbar$ operator. Note that in this case $\tdbar f$ is equal to the image of this $\dbar$ operator, so we may apply Theorem \ref{obstruction bundle}, and the relevant obstruction bundle $E$ is has fibers $E(f)$ dual to $E^{*}(f)$.

  If we give $\M_{g,n}\times \ex B$ the orientation from its almost complex structure, the orientation of $E^{*}(f)$ given in the proof of Theorem \ref{obstruction bundle} is the orientation from its complex structure, as in this case, the linearization of the map $\dbar$ is complex.
  Integrating the pullback of a closed differential form over the virtual moduli space of curves mapping to points in $\ex B$ with genus $g$ and  $n$ marked points is  therefore equivalent to integrating that form against the Euler class of $E$ on $\M_{g,n}\times \ex B$. 
%

 Note that the pullback of this bundle $E$ over the map which forgets one marked point gives the equivalent bundle on $\M_{g,n}\times \ex B$, so all cases follow from $\M_{0,3}$, $\M_{1,1}$ and $\M_{g,0}$ where $g\geq 2$. For dimension reasons, the Euler class of $E$ over $\M_{g,0}\times \ex B$ will be $0$ when $g\geq 2$ and $\dim \ex B>6$, so these curves will not contribute to Gromov Witten invariants.

\appendix

\section{Construction and properties of $\hat f^{+n}$}\label{construct f^{+n}}

In this section we fill in the details of Definition \ref{f^{+n}} from page \pageref{f^{+n}}, and construct the family of curves $\hat f^{+n}$ with $n$ extra marked points from a given family of curves $\hat f$.  As the definition is inductive, with $\hat f^{+n}=(\hat f^{+n-1})^{+1}$, we shall describe $\hat f^{+1}$.
This is some family of cuves
\[\begin{array}{ccc}\hat{\ex C}^{+1}&\xrightarrow{\hat f^{+1}} &\sfp {\hat {\ex B}}{\ex G}2
\\ \downarrow & & \downarrow
\\ \hat{\ex C}&\xrightarrow {\hat f}&\hat{\ex B}
\end{array}\]
that fits into the following diagram
\[\begin{array}{ccc}
\hat{\ex C}^{+1}&&
\\ \downarrow & \hat f^{+1}\searrow &
\\ \hat{\ex C}\fp{\pi_{\ex F}}{\pi_{\ex F}}\hat{\ex C}  &\longrightarrow & \hat{\ex B}\fp{\pi_{\ex G}}{\pi_{\ex G}}{\hat {\ex B}}
\\ \downarrow & & \downarrow
\\ \hat{\ex C}&\xrightarrow {\hat f}&\hat{\ex B}
\\ \downarrow\pi_{\ex F} & & \downarrow\pi_{\ex G}
\\ \ex F & \longrightarrow & \ex G
\end{array}\]

The total space of the domain, $\hat{\ex C}^{+1}$ is constructed by `exploding' the diagonal of $\sfp{\hat {\ex C}}{\ex F}2$ as follows: 

Consider the diagonal map $\Delta:\hat {\ex C}\longrightarrow \sfp{\hat {\ex C}}{\ex F}2$. The image of the tropical part of this map $\totb\Delta$ defines a subdivision of the tropical part of   $\sfp{\hat {\ex C}}{\ex F}2$, which determines a unique refinement $\hat {\ex C}'\longrightarrow \sfp{\hat{\ex C}}{\ex F}2$. Note that the diagonal map to  this refinement $\hat{\ex C}'$ is still defined, 
\[\begin{array}{ccc} & & \hat{\ex C}'
\\ &\nearrow &\downarrow
\\\hat{\ex C}&\xrightarrow{\Delta}&\sfp{\hat{\ex C}}{\ex F}2
\end{array}\]
and a  neighborhood of the image of the diagonal in $\ex C'$ is equal to a neighborhood of $0$ in a $\mathbb C$ bundle over $\hat {\ex C}$.

 Now `explode' the image of the diagonal in $\hat{\ex C}'$ to make $\hat {\ex C}^{+1}\longrightarrow \hat {\ex C}'$ as follows: 
We may choose coordinate charts on $\hat{\ex C}'$ so that any coordinate chart intersecting the image of the diagonal is equal to some subset of  $\mathbb C\times U$ where $U$ is a coordinate chart on $\hat{\ex C}$, the projection to $\hat {\ex C}$ is the obvious projection to $U$, the complex structure on the fibers of this projection is equal to the standard complex structure on $\mathbb C$, and  the image of the diagonal is $0\times U$. Replace these charts with the corresponding subsets of $\et 11\times U$, and leave coordinate charts that do not intersect the image of the diagonal unchanged. Any transition map between coordinate charts of the above type is of the form $(z,u)\mapsto (g(z,u)z,\phi(u))$ where $g(z,u)$ is $\mathbb C^{*}$ valued. In the corresponding `exploded' charts, the corresponding transition map is given by $(\tilde z,u)\mapsto(g(\totl{\tilde z},u)\tilde z,\phi(u))$. The transition maps between other charts can remain unchanged. This defines $\hat{\ex C}^{+1}$. The map $\hat{\ex C}^{+1}\longrightarrow \hat {\ex C}'$ is given in the above coordinate charts by $(\tilde z,u)\mapsto (\totl{\tilde z},u)$. Composing this with the refinement map $\hat{\ex C}'\longrightarrow \sfp{\hat {\ex C}}{\ex F}2$ then gives a degree one fiberwise holomorphic map 
\[\begin{array}{ccc}\hat{\ex C}^{+1}&\longrightarrow &\sfp{\hat{\ex C}}{\ex F}2
\\ \downarrow &&\downarrow
\\ \hat{\ex C}&\xrightarrow{\id}&\hat{\ex C}\end{array}\] 
The map $\hat f^{+1}:\hat{\ex C}^{+1}\longrightarrow\sfp {\hat {\ex B}}{\ex G}2$ is given by the above constructed map $\hat {\ex C}^{+1}\longrightarrow \sfp{\hat{\ex C}}{\ex F}2$ composed with the map 
\[\sfp{\hat{\ex C}}{\ex F}2\longrightarrow \sfp {\hat {\ex B}}{\ex G}2\] which is $\hat f$ in each component. All the above maps are smooth or $\C\infty 1$ if $\hat f$ is.

\

The above construction is functorial. Given a map of families $\hat f\longrightarrow \hat g$, there is an induced map $\hat f^{+1}\longrightarrow \hat g^{+1}$. To see this, consider the naturally induced map 
\[\begin{array}{ccc}
\sfp{\ex C(\hat f)}{\ex F(\hat f)}2&\xrightarrow{\phi\times \phi} &\sfp{\ex C(\hat g)}{\ex F(\hat g)}2
\\ \downarrow & & \downarrow
\\ \ex C(\hat f) &\xrightarrow{\phi} &\ex C(\hat g)
\\ \downarrow & & \downarrow
\\ \ex F(\hat f)& \longrightarrow &\ex F(\hat g)
\end{array}\]

As $\phi\times \phi$ sends the diagonal to the diagonal, this map lifts to the refinement referred to in the above construction. As $\phi\times\phi$ is holomophic on fibers and sends the diagonal to the diagonal, in the special coordinates on the refinement used in the above construction of the form $(z,u)$ and $(w,v)$, the map $\phi\times\phi$ is of the form 
   \[\phi\times\phi(z,u)=(h(z,u)z,\phi(u))\]
   where $h(z,u)$ is $\mathbb C^{*}$ valued. Then the map $\phi^{+1}:\ex C(\hat f^{+1})\longrightarrow \ex C(\hat g^{+1})$ is given in the corresponding exploded coordinates by
   \[\phi^{+1}(\tilde z,u)=(h(\totl{\tilde z},u)\tilde z,\phi(u))\]
We then get a map 
\[\begin{array}{ccc}
\ex C(\hat f^{+1})&\xrightarrow{\phi^{+1}}&\ex C(\hat g^{+1})
\\ \downarrow & & \downarrow
\\ 
\sfp{\ex C(\hat f)}{\ex F(\hat f)}2&\xrightarrow{\phi\times \phi} &\sfp{\ex C(\hat g)}{\ex F(\hat g)}2
\\ \downarrow & & \downarrow
\\ \ex C(\hat f) &\xrightarrow{\phi} &\ex C(\hat g)
\\ \downarrow & & \downarrow
\\ \ex F(\hat f)& \longrightarrow &\ex F(\hat g)
\end{array}\]
The map $\phi^{+1}$ is clearly compatible with the maps $\hat f^{+1}$ and $\hat g^{+1}$, so $\phi^{+1}$ is a map of families $\hat f^{+1}\longrightarrow \hat g^{+1}$. It follows that the construction of $\hat f^{+n}$ is functorial for all $n$.

\

We can apply a similar construction to moduli stacks of curves. Let $\gamma$ be a tropical curve in $\totb{\ex B}$ and $\gamma'$ be obtained from $\gamma$ by adding an infinite edge which maps to a point in $\ex B$. Then given any family $\hat f$ in $\Msw_{g,[\gamma],\beta}$, $\hat f^{+1}$ is a family in $\Msw_{g,[\gamma'],\beta}$. Conversely, if $2g$ plus the number of edges of $\gamma$ is at least $3$, or if $\beta\neq0$, then given any family $\hat f'$ in $\Msw_{g,[\gamma'],\beta}$, forgetting the extra marked point and removing unstable components gives a family $\hat f$ in $\Msw_{g,[\gamma],\beta}$ so that there is a map $\hat f'\longrightarrow \hat f^{+1}$.  Therefore, it makes sense to refer to $\Msw_{g,[\gamma'],\beta}$ as $(\Msw_{g,[\gamma],\beta})^{+1}$. Of course the inverse image of the stack of holomorphic curves under the map $(\Msw_{g,[\gamma],\beta})^{+1}\longrightarrow \Msw_{g,[\gamma],\beta}$ is the stack of holomorphic curves in $(\Msw_{g,[\gamma],\beta})^{+1}$. 

Similarly, construct $\mathcal M^{+1}_{g,[\gamma],\beta}$ as follows. If on some open subset $O$ of $\Msw_{g,[\gamma],\beta}$, $\mathcal M_{g,[\gamma],\beta}$ is equal to $\sum w_{i}t^{\hat f_{i}}$, then on $O^{+1}$, $\mathcal M^{+1}_{g,[\gamma],\beta}$ is equal to $\sum w_{i}t^{\hat f_{i}^{+1}}$.
Let $\gamma$ be a tropical curve in $\totb{\ex B}$ and $\gamma'$ be obtained from $\gamma$ by adding an infinite edge that maps to a point in $\totb{\ex B}$. Suppose that either the homology class $\beta\neq 0$ or $2g$ plus the number of external edges of $\gamma$ is at least $3$. Then $\Mod_{g,[\gamma'],\beta}$ is cobordant to $\Mod_{g,[\gamma],\beta}^{+1}$. This fact implies that our Gromov Witten invariants satisfy the `fundamental class' and `divisor' axioms of \cite{KM}. Its proof involves a slight modification of Theorem \ref{Multi solution} to allow simple perturbations from $\Msw_{g,[\gamma],\beta}$ to  be pulled back and used in defining a virtual moduli space in $\Msw_{g,[\gamma'],\beta}$.

\section{Proof of Theorem \ref{gluing cobordism}}\label{gluing cobordism proof}

This section is dedicated to the proof of Theorem \ref{gluing cobordism}. In particular, we must prove that  given a tropical curve $\gamma$ in $\ex B$ with genus $g_{\gamma}$, and an energy $E$ and genus  $g$, the virtual moduli space of holomorphic curves in $\check {\ex B}_{v}$ for all vertices $v$ of $\gamma$ may be constructed so that the maps $\EV_{0}$ and $\EV_{1}$ are transverse applied to  $\prod_{v}\Mod_{g_{v},[\gamma_{v}],E_{v}}$ whenever $\sum_{v}g_{v}+g_{\gamma}=g$ and $\sum_{v}E_{v}=E$, and we must prove that the pullback of the virtual moduli space of holomorphic curves in $\Msw(\ex B)$ to $\bMs_{g,\gamma,E}$  is cobordant to the pullback to $\bMs_{g,\gamma,E}$ of the virtual moduli space of curves in $\prod_{v}\bMs_{\gamma_{v}}$.

We shall also show that if the virtual moduli space $\prod_{v}\Mod_{g_{v},[\gamma_{v}],E_{v}}$   is constructed using the zero perturbation and  $\EV_{0}$ and $\EV_{1}$ are transverse whenever $g_{\gamma}+\sum g_{v}=g$ and $\sum E_{v}=E$, then the virtual moduli space of holomorphic curves in $\ex B$ may be constructed so that the pullbacks of the two different virtual moduli spaces to $\bMs_{g,\gamma,E}$ are equal.

We have our two maps 
\[\Msw\longleftarrow\bMs_{\gamma}\longrightarrow \prod_{v}\bMs_{\gamma_{v}}\]
described in section \ref{gluing relative invariants}. A curve in $\bMs_{\gamma}$ is holomorphic if and only if it has a holomorphic image in $\prod_{v}\bMs_{\gamma_{v}}$ under the above map, as having a holomorphic image in $\prod_{v}\bMs_{\gamma_{v}}$ means that it is holomorphic on each strata. It is also obvious that a curve in $\bMs_{\gamma} $ is holomorphic if and only if it has a holomorphic image in $\Msw$, as the $\gamma$-decoration has nothing to do with being holomorphic. Therefore the subset of holomorphic curves in $\bMs_{\gamma}$ is equal to the inverse image  of the set of holomorphic curves in $\prod_{v}\bMs_{\gamma_{v}}$ and also equal to the inverse image of the set of holomorphic curves in $\Msw$.

We must now deal with the following problem: the multiperturbations we used to define the virtual moduli space in $\Msw$ will pull back to multiperturbations on $\bMs_{\gamma}$ that look different from the pullback of the multiperturbations used to define the virtual moduli space in $\bMs_{\gamma_v}$. 
We are forced to consider   multiperturbations  slightly more general than those defined using simple perturbations parametrized by core families in  $\Msw(\ex B)$. In particular, we must construct the virtual moduli space of holomorphic curves in $\bMs_{\gamma}$ using perturbations general enough to include multiperturbations pulled back from both $\prod_{v}\bMs_{\gamma_v}$ and $\Msw(\ex B)$.

As in the case of $\Msw(\ex B)$, we have core families and obstruction models covering the moduli space of holomorphic curves in $\bMs_{\gamma}$. If we were to follow the construction of section \ref{virtual class}, we would use multiperturbations obtained from simple perturbations which are parametrized by obstruction models in a compactly supported way. If we are to allow the pullback of mulitperturbations from $\bMs_{\gamma_v}$, we must give up this `compact support' property. Now we begin the construction of virtual moduli spaces.

\begin{itemize}
\item for each holomorphic curve in $\Msw(\ex B)$ choose a core family $\hat f/G$ containing it. (We will specify extra conditions later which will amount to $\hat f/G$ being small enough and having enough marked points.) 
 Construct these core families  $\hat f/G$ on $\ex B$ using the method in the proof of Proposition \ref{smooth model family}.

 After a choice which fixes the $G$-fold ambiguity, the element of $\ex F(\hat f)$ corresponding to a close by curve  $g$ is given by
\begin{enumerate}
 \item The complex structure of  $\ex C(g)$ restricted to each smooth component of $g$.
\item The (transverse) intersection of $g$ with some codimension $2$ submanifolds of $\ex B$.  (This transverse intersection always occurs in the smooth components of $\ex C(g)$)
\item The image under $g$ of some extra marked points on smooth components of $\ex C(g)$. The position of these extra marked points on a given smooth component of $\ex C(g)$ is determined by the complex structure of that smooth component and the transverse intersections of that smooth component with the above submanifolds.
\item \label{gp} A gluing parameter in $\et 1{(0,\infty)}$ corresponding to each internal edge of $\ex C(g)$.
\end{enumerate}

\item The inverse image of $\hat f$ in $\bMs_{\gamma}$ consists of some number of families. Each of these families can be considered as a sub family $\hat f_{l}$ of $\hat f$ with a $\gamma$-decoration in $\abs{G/G_{l}}$ ways, where $G_{l}\leq G$ is the subgroup of $G$ which preserves the $\gamma$-decoration of $\hat f_{l}$. Call $\hat f_{l}$ a lift of $\hat f$ to $\bMs_{\gamma}$. We may consider $\hat f_{l}/G_{l}$ as a core family on $\bMs_{\gamma}$. When our core family is constructed around a curve $f$ not in the image of $\bMs_{\gamma}$, reduce the size of $\hat f$ so that it does not intersect the image of $\bMs_{\gamma}$. When our core family $\hat f$  is constructed around a holomorphic curve $f$ in the image of $\bMs_{\gamma}$, by choosing our core family small enough, we may arrange that all lifts $\hat f_{l}$ contain $f$.
  
\item  Now given a vertex $v$ of $\gamma$ and a lift $\hat f_{l}$ of $\hat f$,  we can construct a core family  $\hat f_{l,v}/G_{l,v}$ around the the image of $ f_{l}$ under the map $\bMs_{\gamma}\longrightarrow \bMs_{\gamma_v}$.   Here, we may use the same choices as the choices for the strata of  $\ex C(f)$ which are associated with $v$ when we lift $f$ to $\bMs_{\gamma}$. 
  This means that there is a map from $\ex F(\hat f_{l})$ to $\ex F(\hat f_{l,v})$ given by simply restricting to the relevant coordinates. 

Now, so long as  our original core family was constructed with enough marked points on the components corresponding to $v$, we can make $\hat f_{l,v}/G_{l,v}$ into an obstruction model $(\hat f_{l,v}/G_{l,v},V_{l,v})$. Here, we may need to shrink the size of $\hat f_{l,v}$ and therefore shrink the size of $\hat f$ in order to keep the map $\ex F(\hat f_{l})\longrightarrow \ex F(\hat f_{l,v})$ well defined.

\item So long as  $\hat f$ was chosen small enough,  $\hat f_{l}/G_{l}$ can be made into an obstruction model $(\hat f_{l}/G_{l},V_{l})$ where $V_{l}=\oplus_{v}V_{l,v}\oplus_{e}V_{e}$ so that 
\begin{itemize} 
\item The pre-obstruction bundle $(\hat f_{l}, V_{l,v})$ above is the pullback to $\hat f_{l}$ of the pre-obstruction bundle $(\hat f_{l,v},V_{l,v})$. 

More explicitly, applying our map $\bMs_{\gamma}\longrightarrow \bMs_{\gamma_v}$ to $\hat f_{l}$ gives a family $(\hat f_{l})_{v}$ close to $\hat f_{l,v}$ in $\bMs_{\gamma_v}$, which comes with a map $\ex C((\hat f_{l})_{v})\longrightarrow \ex C(\hat f_{l,v})$ compatible with the map $\ex F(\hat f_{l})\longrightarrow \ex F(\hat f_{l,v})$ given by restricting to the relevant coordinates. Using the trivialization associated with the pre obstruction bundle $(\hat f_{l,v},V_{l,v})$, we may pull back $V_{l,v}$ to give a pre obstruction bundle $((\hat f_{l})_{v},V_{l,v})$. Each element $(f,x)$ of this $V_{l,v}$ is a curve $f$ in $(\hat f_{l})_{v}$ along with a section in $Y(f)$ (which of course, vanishes on the edges of $f$). By shortening some infinite edges, we may consider $f$ as part of a curve in $\hat f_{l}$, and then extend our section in $Y(f)$ to be zero everywhere else on this curve in $\hat f_{l}$. We may similarly extend our pre obstruction model $((\hat f_{l})_{v},V_{l,v})$ to become a pre obstruction model $(\hat f_{l},V_{l,v})$.

 \item The pre-obstruction bundle $(\hat f_{l},V_{e})$ is the pullback of a pre-obstruction bundle $(\hat f_{l,v},V_{e})$ where $v$ is the vertex attached the incoming end of $e$ (where we have chosen an orientation on the internal edges of $\gamma$)
\begin{itemize}
\item $V_{e}\oplus V_{l,v}$ is complementary to the image of $D\dbar$ restricted to sections in $X^{\infty,\underline 1}(\hat f_{l,v})$ which vanish on the edge $e_{1}$ corresponding to our incoming end of $e$.
\item $V_{e}$  is $G_{l,v}$ invariant, and  contained in the image of $D\dbar$ restricted to sections in $X^{\infty,\underline 1}(\hat f_{l,v})$ which are compactly supported inside some small neighborhood of  $e_{1}$.
\end{itemize}
\end{itemize}

Note that for any curve $f_{l}$ in $\hat f_{l}$, $(\oplus_{v}V_{l,v}\oplus_{e} V_{e})(f_{l})$ is the same dimension as the cokernel of $D\dbar(\hat f_{l})$, and that if we make $\hat f_{l}$ small enough, $D\dbar(\hat f_{l})$ will be transverse to $\oplus_{v}V_{l,v}\oplus_{e}V_{e}$. Note also that $\oplus_{v}V_{l,v}\oplus_{e}V_{e}$ is $G_{l}$ invariant, so we may apply Theorem \ref{regularity theorem} to modify $(\hat f_{l}/G_{l},V)$ into an obstruction model (which we shall still refer to as $(\hat f_{l}/G_{l},V)$).

\item Now choose a finite cover of the set of holomorphic curves in $\bMs_{g,\gamma,E}$  by our obstruction models $(\hat f_{l}/G_{l},V_{l})$ on $\bMs_{\gamma}$.
\item As in the construction of the virtual moduli space in section \ref{virtual class}, by shrinking the open sets for which $\hat f_{l}/G_{l}$ and $\hat f_{l,v}/G_{l,v}$ are core families, we may arrange that there is an open neighborhood $\mathcal O$ of the set of holomorphic curves in $\bMs_{\gamma}$ satisfying our energy and genus bound so that
 \begin{itemize}\item  $\mathcal O$ meets the set for which  $\hat f_{l}/G_{l}$ is a core family properly in the sense of definition \ref{proper meeting},
\item and the image of $\mathcal O$ in $\bMs_{\gamma_v}$ meets the sets for which $\hat f_{l,v}/G_{l,v}$ is a core family properly.
\end{itemize}

\item We need the above proper meeting properties in order to apply the following slight modification of Theorem \ref{Multi solution} to our context.

\begin{thm} \label{mms}
Given
\begin{itemize}
\item a finite collection of  core families $\hat f'_{i}/G_{i}$ for the substacks $\mathcal O'_{i}$ of $\bMs_{\gamma}$ or $\bMs_{\gamma_{v_{i}}}$
\item  an open substack $\mathcal O$ of $\bMs_{\gamma}$ which  meets $\mathcal O'_{i}$ properly for all $\mathcal O'_{i}\subset\bMs_{\gamma}$ (definition \ref{proper meeting}), and which has an image in $\bMs_{\gamma_{v_{i}}}$ which meets  $\mathcal O'_{i}$ properly if $\mathcal O'_{i}\subset\bMs_{\gamma_{v_{i}}}$
\item an obstruction model $(\hat f/G,V)$ for the substack $\mathcal O$
 \item compactly contained $G_{i}$ invariant sub families $\hat f_{i}\subset \hat f'_{i}$
 \end{itemize}
 then given any collection of  $\C\infty1$ simple perturbations $\mathfrak P_{i}$ parametrized by  $\hat f'_{i}$ which are compactly supported in $\hat f_{i}$ and are small enough in $\C\infty1$, and a sufficiently small simple perturbation $\mathfrak P_{0}$ parametrized by $\hat f$ there exists a solution mod $V$ on $\hat f$  which is a $G$-invariant $\C\infty1$ weighted branched section $(\nu,\dbar'\nu)$ of $\hat f^{*}T{\ex B}\oplus V$  with weight $1$ (see example \ref{wbv} below definition \ref{multi}) so that the following holds:

Locally on $\ex F(\hat f)$, 
\[(\nu,\dbar'\nu)=\sum_{l=1}^{n}\frac 1n t^{(\nu_{l},\dbar'\nu_{l})}\]
where $\dbar'\nu_{l}$ is a section of $V$, and  $\nu_{l}$ is in $X^{\infty,\underline 1}(\hat f)$ so that $F(\nu_{l})$ is in $\mathcal O$ and the multiperturbation defined by our perturbations $\mathfrak P_{i}$ is equal to 

\[\prod_{i}F(\nu_{l})^{*}\mathfrak P_{i}=\sum_{j=1}^{n} \frac 1nt^{\mathfrak P_{j,l}}\]
 
 so that 
 \[\dbar'\nu_{l}=\dbar F(\nu_{l})-\mathfrak P_{l,l}\]
 
  The weighted branched section $(\nu,\dbar'\nu)$ is the unique weighted branched section of $\hat f^{*}T_{vert}\hat {\ex B}\oplus  V$ with weight $1$ satisfying the following two conditions:
 \begin{enumerate}
 \item
   Given any curve $f$ in $\hat f$ and  any section $\psi$ in $X^{\infty,\underline 1}(f)$ so that $F(\psi)\in \mathcal O$, if $\prod_{i} F(\psi)^{*}\mathfrak P_{i}=\sum w_{k}t^{\mathfrak Q_{k}}$, and near $f$,  $(\nu,\dbar'\nu)=\sum w'_{l}t^{(\nu_{l},\dbar'\nu_{l})}$ then the sum of the weights $w_{k}$ so that $\dbar F(\psi)-\mathfrak Q_{k}$ is in $V$ is equal to  the sum of the weights $w'_{l}$  so that $\nu_{l}(f)=\psi$.
 \item For any  locally defined section $\psi$ in $X^{\infty,\underline 1}(\hat f)$, if the multi perturbation 
 $\prod_{i} F(\psi)^{*}\mathfrak P_{i}=wt^{\mathfrak Q}+\dotsc$, and $\dbar F(\psi)-\mathfrak Q$ is a section of $V$, then locally, $(\nu,\dbar'\nu)=wt^{(\psi,\dbar F(\psi)-\mathfrak Q)}+\dotsc$.
 \end{enumerate}

This weighted branched section determines the solutions to the perturbed $\dbar$ equation on $\mathcal O$ in the following sense:  
Given any family $\hat g$ in $\mathcal O$, if $\prod_{i}\hat g^{*}\mathfrak P_{i}=wt^{\dbar\hat g}+\dotsc$,  then around each curve in $\hat g$ which projects to the region where the above $\nu_{l}$ are defined, there is a connected open neighborhood in $\hat g$ with at least $nw$  different  maps to $\lrb{\coprod_{l} F(\nu_{l})}/G$.

\
 
 If $\{\mathfrak P_{i}'\}$ is another collection of simple perturbations satisfying the same assumptions as $\{\mathfrak P_{i}\}$ then the sections $(\nu_{l}',\dbar'\nu_{l}')$ corresponding to $(\nu_{l},\dbar'\nu_{l})$, with the correct choice of indexing can be forced to be  as close to $(\nu_{l},\dbar'\nu_{l})$ as we like in $\C\infty1$ by choosing $\{\mathfrak P_{i}'\}$ close to $\{\mathfrak P_{i}\}$  in $\C\infty1$. If $\{\mathfrak P_{i,t}\}$ is a $\C\infty1$ family of simple perturbations satisfying the same assumptions as $\{\mathfrak P_{i}\}$, then the corresponding family of solutions mod $V$, $(\nu_{t},\dbar'\nu_{t})$ form a $\C\infty1$ family of weighted branched sections.

\end{thm}

\pf The proof of this theorem is analogous to the proof of Theorem \ref{Multi solution} - we just need to deal differently with the $\hat f'_{i}$ which are core families on $\bMs_{\gamma_{v_{i}}}$ instead of $\bMs_{\gamma}$.

 As in the proof of Theorem \ref{Multi solution}, we may reduce our proof to the case that for all $i$, some neighborhood of $\hat f$ in $\mathcal O$ is contained in $\mathcal O'_{i}$, or has image in $\Msw_{\gamma_{v_{i}}}$ contained in $\mathcal O'_{i}$.

Use $\mathcal O_{i}$ to denote the restriction of $\mathcal O_{i}'$ to the subset with core $\hat f_{i}/G_{i}$.

As in the proof of Theorem \ref{Multi solution}, we will extend $\hat f$ to a  family $\hat h$ which can be regarded as parametrizing the simple perturbations $\mathfrak P_{i}$ for all $i\in I$ and use the resulting unique solution $\tilde \nu$ to the corresponding perturbed $\dbar$ equation over $\hat h$ to construct the weighted branched section of $\hat f^{*}T_{vert}\hat{\ex B}$ which is our `solution' with the required properties.

 Use the notation 
 \[s^{i}:\ex F(\hat f'_{i})\longrightarrow  \ex F(\hat f_{i}'^{+n_{i}})\]
for the map coming from the extra marked points on the core family $\hat f'_{i}$.

When we are dealing with curves with a $\gamma$-decoration, we shall use a slightly different evaluation map. Let $\M_{\gamma}$ indicate the moduli stack of abstract curves with a $\gamma$ decoration, and $\M_{\gamma}^{+n}$ indicate the moduli stack of curves with a $\gamma$ decoration and $n$ extra punctures. So there is a natural evaluation map 
\[ev^{+n_{i}}(\hat f_{i}'):\ex F(\hat f_{i}'^{+n_{i}})\longrightarrow \M_{\gamma}^{+n_{i}}\times {\ex B}^{n_{i}}\]
in the case that $\hat f'_{i}$ is in $\bMs_{\gamma}$, and a corresponding evaluation map
\[ ev^{+n_{i}}(\hat f_{i}'):\ex F(\hat f_{i}'^{+n_{i}})\longrightarrow \M \times (\check{\ex B}_{v_{i}})^{n_{i}}\]
in the case that $\hat f_{i}$ is in $\bMs_{\gamma_{v_{i}}}$.

Each of these evaluation maps
 have the property that they are  equidimensional embeddings in a neighborhood of the section $s^{i}$. There exists an open neighborhood $\mathcal O_{s^{i}}$ of the family of curves $(s^{i})^{*}\hat f_{i}'^{+n}$ so that  given any curve $f$ in $\mathcal O'_{i}$, if $ f^{+n_{i}}\rvert_{\mathcal O_{s^{i}}}$ indicates the restriction of the family $f^{+n_{i}}$ to $\mathcal O_{s^{i}}$, then $ev^{+n_{i}}(f)(\ex F(f^{+n_{i}}\rvert_{\mathcal O_{s^{i}}}))$ intersects $ev^{+n_{i}}(\hat f')(s^{i}(\ex F(\hat f'_{i})))$ transversely exactly $\abs {G_{i}}$ times, corresponding to the $\abs {G_{i}}$ maps from $\ex C(f)$ into $\ex C(\hat f'_{i})$. 
 
 In the case that $\hat f_{i}$ is in $\bMs_{\gamma_{v_{i}}}$, we wish to pull back $\mathcal O_{s^{i}}$ to $(\bMs_{\gamma})^{+n_{i}}$, (the result of adding $n_{i}$ extra punctures to the moduli stack of curves in $\bMs_{\gamma}$.) The extra punctures given by $s^{i}$ are all contained in the smooth components of $\hat f_{i}$, so we may ensure that the same is true for the extra punctures on curves in $\mathcal O_{s^{i}}$. We may extend our map $\bMs_{\gamma}\longrightarrow \bMs_{\gamma_{v_{i}}}$ in an obvious way to a map to $(\bMs_{\gamma_{v_{i}}})^{+n_{i}}$ from the subset of $(\bMs_{\gamma})^{+n_{i}}$ where all the extra marked points are contained in strata sent to the vertex $v_{i}$. To reduce the amount of extra notation required, refer to the inverse image of the above $\mathcal O_{s^{i}}$ under this map as $\mathcal O_{s^{i}}$.

Consider the family $\hat f^{+(n-1)}:\ex F(\hat f^{+n})\longrightarrow {\ex B}^{n}$. Use the notation 
$X^{+(n-1)}$ to denote the vector bundle over $\ex F(\hat f^{+n})$ which is the pullback under $\hat f^{+n-1}$ of $T\ex B^{n}$.

Any section $\nu$ of $\hat f^{*}T{\ex B}$ corresponds in an obvious way to a section $\nu^{+(n-1)}$ of $X^{+(n-1)}$, and the map $F:\hat f^{*}T{\ex B}\longrightarrow {\ex B} $ corresponds to a $\C\infty 1$ map  
\[F^{+(n-1)}:X^{+(n-1)}\longrightarrow  \ex B^{n}\]
so that
\[F^{+(n-1)}(\nu^{+(n-1)})=\lrb{F(\nu)}^{+(n-1)}\]
Use the notation $\nu^{+(n_{i}-1)}\rvert_{\mathcal O_{s^{i}}}$ to denote the restriction of $\nu^{+(n_{i}-1)}$ to the subset $\ex F(F(\nu)^{+n_{i}}\rvert_{\mathcal O_{s^{i}}})\subset \ex F(\hat f^{+n_{i}})$ 

Define a map \[EV^{+n}:X^{+(n-1)}\longrightarrow  \M_{\gamma}^{+n}\times{\ex B}^{n}\] so 
that $EV^{+n}$ is equal to the natural map coming from the complex structure of curves in $\ex C(\hat f^{+n})\longrightarrow \ex F(\hat f^{+n})$ on the first component, and $F^{+(n-1)}$ on the second component. So
\[EV^{+n}( \nu^{+(n-1)}(\cdot))=ev^{+n}(F(\nu))(\cdot)\] 
If we restrict to the subset of $X^{+(n-1)}$ where the extra marked points are in components sent to $v_{i}$, then there is an analogous map $EV^{+n}_{v_{i}}$ to $\M_{\gamma_{v_{i}}}^{+n}\times(\check{\ex B}_{v_{i}})^{n}$.

 Use the notation $\nu(g)^{+(n_{i}-1)}\rvert_{\mathcal O_{s^{i}}}$ for the restriction of  $\nu^{+(n_{i}-1)}\rvert_{\mathcal O_{s^{i}}}$ to the inverse image of a curve $g\in \hat f$. In the case that $\hat f_{i}$ is in $\bMs_{\gamma}$, for any section $\nu$ small enough in $C^{1,\delta}$,
  the map $EV^{+n}$ restricted to $\nu(g)^{+(n_{i}-1)}\rvert_{\mathcal O_{s^{i}}}$  intersects $ev^{+n_{i}}(\hat f'_{i})(s^{i}(\ex F(\hat f'_{i})))$ transversely in exactly $\abs {G_{i}}$ points.  Denote by $S_{i}$ the subset of $X^{+(n_{i}-1)}$ which is the pullback of the image of the section $s^{i}$:
\[S_{i}:=(EV^{+n_{i}})^{-1}\lrb{ev^{+n_{i}}(\hat f'_{i})(s^{i}(\ex F(\hat f'_{i})))}\subset X^{+(n_{i}-1)}\]

In the case that $\hat f_{i}$ is in $\bMs_{\gamma_{v_{i}}}$, for any section $\nu$ small enough in $C^{1,\delta}$,
  the map $EV^{+n}_{v_{i}}$ restricted to $\nu(g)^{+(n_{i}-1)}\rvert_{\mathcal O_{s^{i}}}$  intersects $ev^{+n_{i}}(\hat f'_{i})(s^{i}(\ex F(\hat f'_{i})))$ transversely in exactly $\abs {G_{i}}$ points.  Denote by $S_{i}$ the subset of $X^{+(n_{i}-1)}$ which is the pullback of the image of the section $s^{i}$:
\[S_{i}:=(EV^{+n_{i}}_{v_{i}})^{-1}\lrb{ev^{+n_{i}}(\hat f'_{i})(s^{i}(\ex F(\hat f'_{i})))}\subset X^{+(n_{i}-1)}\]

Close to the zero section in $X^{+(n_{i}-1)}$, $S_{i}$ has regularity $\C\infty1$, and for sections $\nu$ small enough in $C^{1,\delta}$, $S_{i}$ is transverse to $\nu^{+(n_{i}-1)}\rvert_{\mathcal O_{s^{i}}}$ and  $\nu^{+(n_{i}-1)}\rvert_{\mathcal O_{s^{i}}}\cap S_{i}$ is a  $\abs {G_{i}}$-fold  multisection $\ex F(\hat f)\longrightarrow S_{i}\subset X^{+(n_{i}-1)}$.

Use the notation $\tilde X^{+(n-1)}$ for the pullback along the map $\ex C(\hat f^{+n})\longrightarrow \ex F(\hat f^{+n})$ of the vector bundle $X^{+(n-1)}$, $\tilde S_{i}$ for the inverse image of $S_{i}$ in $\tilde X^{+(n_{i}-1)}$, and 
$\tilde \nu^{+(n_{i}-1)}_{\rvert_{\mathcal O_{s^{i}}}}$ for the pullback of $\nu^{+(n_{i}-1)}\rvert_{\mathcal O_{s^{i}}}\subset X^{+(n_{i}-1)}$ to a section of $\tilde X^{+(n_{i}-1)}$. Considering $\tilde S_{i}\longrightarrow S_{i}$ as the pullback of $\ex C(\hat f^{+n} )\longrightarrow \ex F(\hat f^{+n})$ gives $\tilde S_{i}\longrightarrow S_{i}$ the structure of a family of curves. Forgetting the extra $n$ marked points gives a family of curves $\check S_{i}\longrightarrow S_{i}$ which is the pullback of $\ex C(\hat f)\longrightarrow \ex F(\hat f)$.

 As discussed in the proof of Theorem \ref{Multi solution}, in the case that $\hat f_{i}$ is in $\bMs_{\gamma}$, there is a map 
 \[\begin{array}{ccc}\check S_{i}&\xrightarrow{\Phi_{i}} &\ex C(\hat f_{i}')
 \\\downarrow &&\downarrow
 \\ S_{i}&\longrightarrow &\ex F(\hat f_{i}')
 \end{array}\] 
 
  so  $\Phi_{i}$ determines the maps $\Phi_{F(\nu)}$ from the definition of the core family $\hat f_{i}'/G_{i}$ in the following sense: For $\nu$ small enough, the intersection of $ \nu^{+(n_{i}-1)}\rvert_{\mathcal O_{s^{i}}}$
 with $ S_{i}$ is transverse, and is a $\abs {G_{i}}$-fold cover  of $\ex F(\hat f)$, which lifts to a $\abs {G_{i}}$-fold cover of $\ex C(\hat f)$ which is a subset of $\check S$. Then $\Phi_{i}$ gives a map of our $\abs {G_{i}}$-fold cover of $\ex C(\hat f)$ into $\ex C(\hat f'_{i})$, which corresponds to the map $\Phi_{F(\nu)}:\ex C(\hat f)\longrightarrow \ex C(\hat f'_{i})/G_{i}$.
 
 In the case that $\hat f_{i}$ is in $\bMs_{\gamma_{v_{i}}}$, denote by $(\check S_{i})_{v_{i}}$ the restriction of $\check{S_{i}}$ to the subset which is sent to $v_{i}$ using the $\gamma$ decoration. Then there is an analogous map
 
 \[\begin{array}{ccc}(\check S_{i})_{v_{i}}&\xrightarrow{\Phi_{i}} &\ex C(\hat f_{i}')
 \\\downarrow &&\downarrow
 \\ S_{i}&\longrightarrow &\ex F(\hat f_{i}')
 \end{array}\] 
 
 so that the following holds: For $\nu$ small enough, the intersection of $ \nu^{+(n_{i}-1)}\rvert_{\mathcal O_{s^{i}}}$
 with $ S_{i}$ is transverse, and is a $\abs {G_{i}}$-fold cover  of $\ex F(\hat f)$, which lifts to a $\abs {G_{i}}$-fold cover of $\ex C(\hat f)$ which is a subset of $\check S$. Then $\Phi_{i}$ gives a map of a subset of our $\abs {G_{i}}$-fold cover of $\ex C(\hat f)$ into $\ex C(\hat f'_{i})$ which can also be constructed as follows: Applying the map $\bMs_{\gamma}\longrightarrow \bMs_{\gamma_{v_{i}}}$ to the family $F(\nu)$ gives a family $F(\nu)_{v_{i}}$ in $\bMs_{\gamma_{v_{i}}}$. Then the map $\Phi_{F(\nu)_{v_{i}}}$ can be regarded as giving a map of a $G_{i}$ fold cover of $\ex C(F(\nu)_{v_{i}})$ to $\ex C(\hat f_{i}')$. Pull this $G_{i}$-fold cover and map back over the inclusion into $\ex C(F(\nu)_{v_{i}})$ of the strata of $\ex C(\hat f)$ which are sent to $v_{i}$.   The resulting map is $\Phi_{i}$.

%
%
%

Denote by $X^{+}$ the fiber product of $X^{+(n_{i}-1)}$ over $\ex F(\hat f)$ for all $i $,
and  denote by $\ex F^{+}$ the fiber product of $\ex F(\hat f^{+n_{i}})$ over $\ex F(\hat f)$ for all $i $;
    so $X^{+}$ is a vector bundle over $\ex F^{+}$.
        A $\C\infty1$ section $\nu$ of $\hat f^{*}T_{vert}\hat {\ex B}$ corresponds in the obvious way to a $\C\infty1$ section $\nu^{+}$ of $X^{+}$ which is equal to $\nu^{+(n_{i}-1)}$ on each $X^{+(n_{i}-1)}$ factor. Similarly, denote by $\nu^{+}\rvert_{\mathcal O_{s^{+}}}$ the open subset of $\nu^{+}$ inside $\nu^{+(n_{i}-1)}\rvert_{\mathcal O_{s^{i}}}$ on each $X^{+(n_{i}-1)}$ factor.
          Denote by $ S^{+}\subset X^{+}$ the  subset corresponding to all $S_{i}\subset X^{+(n_{i}-1)}$  restricted to a neighborhood of the zero section small enough that $S^{+}$ is $\C\infty1$.
We can choose $S^{+}$ so that pulling back $(\hat f,V)$ over the map $S^{+}\longrightarrow \ex F(\hat f)$ gives an allowable pre obstruction model $(\hat h,V)$. Note that $\ex C(\hat h)$ is some open subset of the fiber product of $\check S_{i}$ over $\ex C(\hat f)$ for all $i $, so the maps $\Phi_{i}$ induce  maps
\[\begin{array}{ccc}\ex C(\hat h)&\xrightarrow{\Phi_{i}}& \ex C(\hat f'_{i})
\\ \downarrow & & \downarrow
\\\ex F(\hat h)&\longrightarrow &\ex F(\hat f'_{i})\end{array}\]
(When $\hat f_{i}$ is in $\bMs_{\gamma_{v_{i}}}$, note that the above map $\Phi_{i}$ is only defined on the subset of $\ex C(\hat h)$ consisting of strata which are sent to $v_{i}$ using the $\gamma$ decoration.)

Pulling a simple perturbation $\mathfrak P_{i}$ parametrized by $\hat f_{i}'$ back over the map $\Phi_{i} $ gives a simple perturbation $\Phi_{i}^{*}\mathfrak P_{i}$ parametrized by $\hat h$. In the case that $\hat f_{i}$ is in $\bMs_{\gamma_{v}}$, this simple perturbation is only defined on the strata of $\ex C(\hat h)$ sent to $v_{i}$, but we may extend the simple perturbation to be $0$ everywhere else. Now the rest of the proof is identical to the corresponding part of the proof of Theorem \ref{Multi solution}.

 Use the notation
\[\mathfrak P:=\sum_{i}\Phi_{i}^{*}\mathfrak P_{i}\]
   If $\nu$ is any small enough section of $\hat f^{+}T_{vert}\hat{\ex B}$, then the multi perturbation $\prod_{i }F(\nu)^{*}\mathfrak P_{i}$  can be constructed as follows: If $\nu$ is small enough, then $\nu^{+}\rvert_{\mathcal O_{s^{+}}}$ is transverse to $S^{+}$, and the intersection of $\nu^{+}\rvert_{\mathcal O_{s^{+}}}$ with $S^{+}$ is a $n$-fold cover of $\ex F(\hat f)$ in $\ex F(\hat h)$ which lifts to a multiple cover  of $\ex C(\hat f)$ inside $\ex C(\hat h)$ (where $n=\prod_{i }\abs{ G_{i}}$). Together these give the domain for a family of curves $F(\nu)'$ which is a  
 $n$-fold multiple cover of $F(\nu)$.  Restricting $\mathfrak P$ to $F(\nu)'$ then gives a section of $\Y{F(\nu)'}$, which corresponds to a $n$-fold multi section of $\Y{F(\nu)}$. Locally, giving each of these $n$ sections a weight $1/n$ gives a weighted branched section of $\Y{F(\nu)}$ with total weight $1$ which is equal to the multi perturbation   $\prod_{i }F(\nu)^{*}\mathfrak P_{i}$.

 As $(\hat f,V)$ is an obstruction model, Theorem \ref{regularity theorem} applies to $(\hat h,V)$ and implies that there is some neighborhood of $0$ in the space of simple perturbations  parametrized by $\hat h$ so that for such any  $\mathfrak P$  in this neighborhood, there is a unique small  section $\tilde \nu\in X^{\infty,\underline 1}(\hat h)$ so that $(\dbar -\mathfrak P)\tilde \nu=\in V$. The fact that $(\hat f,V)$ is part of an obstruction model for $\mathcal O$ implies the following uniqueness property for $\tilde \nu$ if $\mathfrak P$ is small enough: Given any curve $h$ in $\hat h$ and section $\psi$ in $X^{\infty,\underline 1}$  so that $F(\psi)$ is in $\mathcal O$, then $(\dbar -\mathfrak P)\psi\in V$ if and only if $\psi$ is the restriction to $h$ of $\tilde \nu$.  
 
 Denote by $\tilde X^{+}$ the pullback of $X^{+}$ over the map $\ex F(\hat h)\longrightarrow \ex F(\hat f)$, and denote by $\tilde S^{+}$ the pullback of $S^{+}$. 
\[ \begin{array}{ccc}\tilde S^{+}\subset \tilde X^{+}&\longrightarrow& S^{+}\subset  X^{+}
 \\\downarrow &&\downarrow
 \\ \ex F(\hat h)&\longrightarrow& \ex F(\hat f)\end{array}\]
 
 This $\tilde S^{+}$ comes with two maps into $S^{+}$: one  the restriction of the map $\tilde X^{+}\longrightarrow X^{+}$, and one  the restriction of the map $\tilde X^{+}\longrightarrow \ex F(\hat h)=S^{+}$. The significance of these two maps is as follows: $S^{+}$ is used to parametrize the different possible maps into $\ex C(\hat f_{i})$. The restriction of the map $\tilde X^{+}\longrightarrow X^{+}$ determines the maps into $\ex C(\hat f_{i})$ that the families of curves corresponding to sections of $X^{+}$ have. The restriction of the map $\tilde X^{+}\longrightarrow \ex F(\hat h)=S^{+}$ specifies the maps to $\ex C(\hat f_{i})$ used to define our simple perturbation parametrized by $\hat h$.
 
  Denote by $S^{+}_{\Delta}$ the subset of $\tilde S^{+}$ on which the above two maps agree. Because these two above maps agree when composed with the relevant maps to $\ex F(\hat f)$, $\tilde S^{+}$ can be regarded as the fiber product of $ S^{+}$ with itself over $\ex F(\hat f)$ and $S^{+}_{\Delta}$ is the diagonal in this fiber product $\tilde S^{+}$. Therefore, $S^{+}_{\Delta}$ is $\C\infty1$ and the  map  $S^{+}_{\Delta}\longrightarrow S^{+}$ is an isomorphism.  A section $\tilde \nu$ of $\hat h^{*}T{\ex B}$ defines a section $\tilde \nu^{+}$ of the vector bundle $\tilde X^{+}$ so that if $\tilde \nu$ is the pullback over the map $\ex C(\hat h)\longrightarrow \ex F(\hat f)$ of some section $\nu$ of $\hat f^{*}T{\ex B}$, then $\tilde \nu^{+}$ is the pullback of $ \nu^{+}$. We can define $\tilde \nu^{+}\rvert_{\mathcal O_{s^{+}}}$ similarly to the definition of $\nu^{+}\rvert_{\mathcal O_{s^{+}}}$.
  
  As $\nu^{+}\rvert_{\mathcal O_{s^{+}}}$ is transverse to $S^{+}$ for $\nu$ small enough, and $\nu^{+}\rvert_{\mathcal O_{s^{+}}}\cap S^{+}$ gives a $n$-fold cover of $\ex F(\hat f)$, $\tilde \nu^{+}\rvert_{\mathcal O_{s^{+}}}$ is transverse to $S^{+}_{\Delta}$ for $\tilde \nu$ small enough, and $\tilde \nu^{+}\rvert_{\mathcal O_{s^{+}}}\cap S^{+}_{\Delta}$ also defines a   $n$-fold cover of $\ex F(\hat f)$ with regularity $\C\infty1$. To see this, suppose that $\tilde \nu^{+}$ is the pullback of some $ \nu^{+}$. Then $\tilde \nu^{+}$ is transverse to $\tilde S^{+}$, and $\tilde \nu^{+}\cap \tilde S^{+}$ is an $n$ fold cover of $\ex F(\hat h)$ which is a pullback of $\nu^{+}\cap S^{+}$ over the map $\ex F(\hat h)\longrightarrow \ex F(\hat f)$. These $n$ sections of $\tilde S^{+}\longrightarrow \ex F(\hat h)=S^{+}$ are constant on fibers of the map $\ex F(\hat h)=S^{+}\longrightarrow \ex F(f)$, and are therefore transverse to the diagonal section $S^{+}_{\Delta}$, and when intersected with $S^{+}_{\Delta}$ give an $n$-fold section of $S^{+}_{\Delta}\longrightarrow \ex F(\hat f)$. This transversality and the fact that $\tilde \nu^{+}\rvert_{\mathcal O_{s^{+}}}\cap S^{+}_{\Delta}$ defines a   $n$-fold cover of $\ex F(\hat f)$ with regularity $\C\infty1$ is stable under perturbations of $\tilde \nu^{+}$, so it remains true for small $\tilde \nu$ which aren't the pullback of some $\nu$.

   We may consider this multiple cover of $\ex F(\hat f)$ as being a multi section  $\ex F'$ of $S^{+}=\ex F(\hat h)\longrightarrow \ex F(\hat f)$, which lifts to a multi section of $\ex C(\hat h)\longrightarrow\ex C(\hat f)$. Restricting $\tilde \nu$ to this multi section gives locally $n$ sections $\nu_{l}$ of $\hat f^{*}T{\ex B}$ with regularity $\C\infty1$. We may similarly pullback the sections $(\tilde\nu,\dbar\tilde \nu-\mathfrak P)$ to give locally $n$ sections $(\nu_{l},\dbar'\nu_{l})$. Then
  \begin{equation}\label{2nu l}(\nu,\dbar'\nu)=\sum_{l=1}^{n}\frac 1{n}t^{(\nu_{l},\dbar'\nu_{l})}\end{equation}
   is the  weighed branched solution which is our `solution mod $V$'. We shall now show that this weighted branched section has the required properties if $\{\mathfrak P_{i}\}$ is small enough. Note first that close by simple perturbations $\{\mathfrak P_{i}\}$ give close by solutions $(\nu_{l},\dbar'\nu_{l})$. Also note that if we have a $\C\infty1$ family of simple perturbations $\{\mathfrak P_{i,t}\}$, Theorem \ref{regularity theorem} implies that the corresponding family of solutions $\tilde\nu_{t}$ to $(\dbar-\mathfrak P_{t})\tilde\nu_{t}\in V$ is a $\C\infty1$ family, so the corresponding weighted branched sections $(\nu_{t},\dbar'\nu_{t})$ form a $\C\infty1$ family.  
   
If $\{\mathfrak P_{i}\}$ are chosen small enough, the multi perturbation under study is given by 
   \begin{equation}\label{2Pjl}\prod_{i}\nu_{l}^{*}\mathfrak P_{i}=\sum_{j=1}^{n}\frac 1nt^{\mathfrak P_{j,l}}\end{equation}
  where $\mathfrak P_{j,l}$ is constructed as follows: $\nu_{l}^{+}\rvert_{\mathcal O_{s^{+}}}\cap S^{+}$ is a $n$-fold cover of the open subset of $\ex F(\hat f)$ where $\nu_{l}$ is defined. By working locally, this $n$-fold cover can be thought of as $n$ local sections of $S^{+}=\ex F(\hat h)\longrightarrow \ex F(\hat f)$, which lift to $n$ local sections of $\ex C(\hat h)\longrightarrow \ex F(\hat f)$. The restriction of $\mathfrak P$ to these $n$ local sections gives the $n$ sections $\mathfrak P_{j,l}$ of $\Y{F(\nu_{l})}$ in the formula (\ref{2Pjl}) above. As one of these sections of $\ex F(\hat h)\longrightarrow\ex F(\hat f)$ coincides with the multi section $\ex F'$ mentioned in the  paragraph preceding equation (\ref{2nu l}) obtained using the solution $\tilde \nu$ to the equation $(\dbar-\mathfrak P)\tilde \nu\in V$, one of the sections $\mathfrak P_{l,l}$ of $\Y{F(\nu_{l})}$ has the property that $\dbar F(\nu_{l})-\mathfrak P_{l,l}=\dbar'\nu_{l}$. 
  
  Suppose that $f$ is  some curve in $\hat f$  where these $\nu_{l}$ in formula (\ref{2nu l}) are defined, and $\psi\in X^{\infty,\underline 1}$ is small enough that $F(\psi)\in\mathcal O$. If the simple perturbations $\mathfrak P_{i}$ are chosen small enough, the fact that $(\hat f,V)$ is an obstruction model will imply that if $ \prod_{i} F(\psi)^{*}\mathfrak P_{i}=wt^{\mathfrak Q}+\dotsc$ where $w>0$ and $(\dbar F(\psi)-\mathfrak Q)\in V$, then $\psi$ must be small - choose $\{\mathfrak P_{i}\}$ small enough that such $f$ must have $\psi^{+}\rvert_{\mathcal O_{s^{+}}}$ intersecting $S^{+}$ transversely $n$ times and our reduction to the case that the families corresponding to sections of $X(\hat f)$ are in $\mathcal O'_{i}$ is valid. Then $ \prod_{i} F(\psi)^{*}\mathfrak P_{i}=\sum_{l=1}^{n}\frac 1nt^{\mathfrak Q_{l}}$ where the $n$ sections $\mathfrak Q_{l}$ of $\Y{F(\psi)}$ are obtained as follows: The $n$  points of $\psi^{+}\rvert_{\mathcal O_{s^{+}}}\cap S^{+}$ correspond to $n$ maps of $\ex C(F(\psi))$ into $\ex C(\hat h)$ - the $n$ sections $\mathfrak Q_{l}$ are given by pulling back the simple perturbation $\mathfrak P$ over these maps. Then $(\dbar F(\psi)-\mathfrak Q_{l})\in V$  if and only if $\psi$ is equal to the pullback under the relevant map of the solution $\tilde \nu$ to $(\dbar -\mathfrak P)\tilde \nu\in V$. Therefore, if $\{\mathfrak P_{i}\}$ is small enough, the number of $\mathfrak Q_{l}$ so that $(\dbar F(\psi)-\mathfrak Q_{l})\in V$ is equal to the number of $\nu_{l}$ from formula (\ref{2nu l}) so that $\psi =\nu_{l}(f)$. 
  
  Similarly, if $\nu'$ is locally a section of $\hat f^{*}T{\ex B}$ vanishing on the relevant marked points so that $F(\nu')\in\mathcal O$ and $\prod_{i} F(\nu')^{*}\mathfrak P_{i}=wt^{\mathfrak Q}+\dotsc$ where $w>0$ and $(\dbar F(\nu')-\mathfrak Q)\in V$, then so long as $\{\mathfrak P_{i}\}$ is small enough, $\nu'^{+}\rvert_{\mathcal O_{s^{+}}}\cap S^{+}$ is locally a $n$-fold cover of $\ex F(\hat f)$  corresponding to $n$ sections of $\ex F(\hat h)\longrightarrow \ex F(\hat f)$ which lift locally to $n$ sections of $\ex C(\hat h)\longrightarrow \ex C(\hat f)$. Then $\mathfrak Q$ must locally correspond to the pullback of $\mathfrak P$ under one of these local maps $\ex C(\hat f)\longrightarrow \ex C(\hat h)$, and $\nu'$ must locally be the pullback of the solution $\tilde \nu$ to $(\dbar -\mathfrak P)\tilde \nu\in V$. It follows that $(\nu',\dbar\nu'-\mathfrak Q)$ must coincide locally with one of these $(\nu_{l},\dbar'\nu_{l})$ from formula (\ref{2nu l}), and the weighted branched section locally equal to $\sum_{l=1}^{n}\frac 1nt^{(\nu_{l},\dbar'\nu_{l})}$ is the unique weighted branched section with the  required properties.  The fact that our weighted branched section is the unique one satisfying these properties implies that it is $G$-invariant. (This may also be checked directly by noting the the construction is $G$-invariant.)
  
  The relationship of our solution to the solutions to the perturbed $\dbar$ equation is as follows: Consider a family $\hat g$ in $\mathcal O$ so that $\prod_{i}\mathfrak P_{i}\hat g=wt^{\dbar \hat g}+\dotsc$ and $w>(k-1)/n$.  By working locally on $\hat g$, we may assume without losing generality that $\hat g$ projects to the region in $\hat f$ where equation \ref{2nu l} holds, and that we may resolve the $G$-fold ambiguity of the map $\ex C(\hat g)\longrightarrow \ex C(\hat f)/G$ to a map $\ex C(\hat g)\longrightarrow \ex C(\hat f)$.
   We may then pull back $X^{+}$  and $S^{+}$ to be bundles over $\ex F(\hat g)$. The corresponding bundles $X^{+}(\hat g)$ and $S^{+}(\hat g)$ may also be constructed in the same way as the original bundles using the induced trivialization on $\hat g$ from $\hat f$.  The section $\psi$ vanishing at the correct marked points so that $\hat g=F(\psi)$ corresponds to a section $\psi^{+}$ of this pulled back $X^{+}$ which is transverse to the pulled back $S^{+}$ and intersects this pulled back $S^{+}$ in an $n$-fold cover of $\ex F(\hat g)$. This $n$-fold cover of $\ex F(\hat g)$ comes with a map to $S^{+}$, corresponding to a map to $\ex F(\hat h)$ which lifts to a fiberwise holomorphic map of a $n$-fold cover of $\ex C(\hat g)$ to $\ex C(\hat h)$ so that $\prod_{i}\hat g^{*}\mathfrak P_{i}$ is determined by pulling back $\mathfrak P$ over this map, then giving the simple perturbation from each branch of the cover a weight $1/n$ and summing the result. As $w>(k-1)/n$,  locally at least $k$ of these simple perturbations must be $\dbar\hat g$, and $\psi$ must be the pullback under each of the corresponding maps of the solution $\tilde \nu$ to $(\dbar-\mathfrak P)\tilde\nu\in V$, and its image must be contained in the subset where $\dbar'\tilde \nu=0$. It follows that around each curve in  $\hat g$, there is a map of a neighborhood into at least $k$ of the $ F(\nu_{l})$  with image contained in the subset where $\dbar'\nu_{l}=0$, and with the map $\ex C(\hat g)\longrightarrow \ex C( F(\nu_{l}))$ corresponding to our local choice of lift of the map $\ex C(\hat g)\longrightarrow \ex C(\hat f)/G$ coming from the fact that $\hat f/G$ is a core family.  Without a choice resolving this $G$-fold ambiguity, this corresponds to there being at least $k$ maps of $\hat g$ into $(\coprod_{l}F(\nu_{l}))/G$.

  \stop

\item By pulling back generic simple perturbations parametrized by $\hat f_{l,v}$ for all $v$, transversality can be achieved in any compact subset of an obstruction model $(\hat f_{l}/G_{l},V_{l})$ as follows:

\begin{lemma}\label{transversality conditions} Suppose that we have compactly supported simple perturbations parametrized by our $\hat f_{l,v}$ which are small enough that Theorems \ref{Multi solution} and \ref{mms} may be applied in the following two contexts.  Over $\ex F(\hat f_{l,v})$, Theorem \ref{Multi solution} gives a solution to the corresponding perturbed $\dbar$ equation mod $V_{l,v}$ which is locally equal to a weighted sum of sections  in $X^{\infty,\underline 1}(\hat f_{l,v})$ each of which has a corresponding section $\dbar'_{v}$ of $V_{l,v}$. Similarly, over $\ex F(\hat f_{l})$, Theorem \ref{mms} gives a solution to the corresponding perturbed $\dbar$ equation mod $V_{l}$ which is locally equal to a weighted sum of sections in $X^{\infty,\underline 1}(\hat f_{l})$, each of which has a corresponding section $\dbar'$ of $V_{l}$.

Let $C$ be a subfamily of  $\hat f_{l}$, and let $C'$ be the subfamily of $\prod_{v}\hat f_{l,v}$ which is the intersection of $\prod_{v}\hat f_{l,v}$ with the  image of $C$ in $\prod_{v}\bMs_{\gamma_{v}}$.
Then if our simple perturbations are small enough, and on $C'$,
\begin{enumerate}\item\label{tc1} $\dbar'_{v}$ is transverse to the zero section in $V_{l,v}$,
\item\label{tc2} the solutions to the equations $\dbar'_{v}=0$ are transverse when mapped to $\prod_{e}\mathcal M_{e}$,
\item\label{tc3} and none of the solutions to the equations $\dbar'_{v}=0$ have any automorphisms,
\end{enumerate}

then restricted to $C$, $\dbar'$ is transverse to the zero section in $V_{l}$.
\end{lemma}

\pf  Recall that $V_{l}=\oplus_{v}V_{l,v}\oplus_{e} V_{e}$, where $V_{e}$ is considered as pulled back from a pre-obstruction model $(\hat f_{l,v},V_{e})$ where $v$ is the vertex attached to one  of the ends of $e$. Use the notation $S_{v}$ to indicate the set of edges $e$ for which $V_{e}$ is considered as coming from a pre-obstruction model $(\hat f_{l,v},V_{e})$.  Instead of a single simple perturbation parametrized by $\hat f_{l,v}$, we may consider a family of simple perturbations parameterized by the total space of the bundle $\oplus_{e\in S_{v}}V_{e}$ where at $(f,x)$ our new simple perturbation is equal to the old simple perturbation plus  $x$.  Close to the zero section in $\oplus_{e\in S_{v}}V_{e}$, we may then apply Theorem \ref{Multi solution} with this new family of simple perturbations, to obtain the local moduli space $M'_{v}$ of solutions to $\dbar'\in V_{l,v}\oplus_{e\in S_{v}}V_{e}$. This moduli space $M'_{v}$ is locally equal to a weighted sum of families  together with sections $\dbar'$ of $V_{l,v}\oplus_{e\in S_{v}}V_{e}$. It follows from our application of Theorem \ref{Multi solution} that these sections $\dbar'$ are transverse to $V_{l,v}\subset V_{l,v}\oplus_{e\in S_{v}}V_{e}$, therefore condition \ref{tc1} implies that these sections $\dbar'$ are transverse to the zero section.

Condition \ref{tc2} implies that near $\dbar'=0$, these moduli spaces $M'_{v}$ will be transverse when mapped by $\EV_{1}$ and $\EV_{0}$ to $\prod_{e}\mathcal M_{e}$. Lemma \ref{fp bundle} implies that close to the subset where $\dbar'=0$, the result of applying Theorem \ref{mms} to $(\hat f_{l},V_{l})$ is a bundle of gluing choices over the result of taking the fiber product of our $M_{v}$ over $\prod_{e}\mathcal M_{e}$ and restricting to the subset with tropical parts that can be obtained by $\gamma$-decorated tropical completion. The sections $\dbar'$ of $V_{l}$ are equal to the pullback of the fiber product of the individual sections $\dbar'$ of $V_{l,v}\oplus_{e\in S_{v}}V_{e}$, so the fact that these individual sections $\dbar'$ are transverse to the zero section combined with the transversality \ref{tc2} for the solutions to  $\dbar'=0$ when taking our fiber product imply that $\dbar'$ is transverse to the zero section in $V_{l}$.

\stop

\begin{lemma}\label{EV transversality} If $\ex B$ is not zero dimensional, then given any compact subfamily $C\subset \hat f_{l}$, the conditions \ref{tc1}, \ref{tc2}, and \ref{tc3} of Lemma \ref{transversality conditions} are satisfied for a generic choice of simple perturbations parametrized by $\hat f_{l,v}$ for all $v$.
\end{lemma}

\pf

Clearly, the conditions \ref{tc1}, \ref{tc2} and \ref{tc3} are open conditions, and Theorem \ref{final theorem} states that the set of perturbations satisfying conditions \ref{tc1} and \ref{tc3} is dense, so we need only show that perturbations satisfying  condition \ref{tc2} are dense. This is easy: start off with some $\C\infty1$ perturbation satisfying condition \ref{tc1} and \ref{tc3}, apply Theorem \ref{Multi solution} and let $\sum_{i=1}^{n}w_{i}t^{\prod_{v}\nu_{i,v}}$ indicate the product for all $v$ of the $\C\infty1$ multi sections which solve the equation $\dbar'_{v}\in V_{l,v}$. Assume that for each fixed $i<k$, the the sections $\nu_{i,v}$ correspond to families satisfying the transversality condition \ref{tc2}. We may choose a $\C\infty1$ section $\nu_{v,k}'$ in  $X(\hat f_{l,v})$ which  is as close as we like to $\nu_{v,k}$ in $\C\infty1$, so that these sections $\nu_{v,k}'$ restricted to the sub families of $\hat f_{l,v}$ where $\dbar'\hat f_{l,v}=0$ correspond to families which are transverse when mapped to $\prod_{e}\mathcal M_{e}$ using $\EV_{0}$ and $\EV_{1}$. Then if $\nu'_{v,k}$ is sufficiently close to $\nu_{v,k}$, modifying our original perturbation by $\dbar\nu'_{v,k}-\dbar\nu_{v,k}$ will give modified solution sections $\nu_{v,i}'$ which satisfy the conditions \ref{tc1} and \ref{tc3}, and which also satisfy transversality condition $\ref{tc2}$ for $i\leq k$. Continuing this argument for larger $k$, it follows that the set of perturbations satisfying condition \ref{tc2} is also dense, so a generic small perturbation parametrized by $\hat f_{l,v}$ for all $v$ will satisfy conditions \ref{tc1}, \ref{tc2} and \ref{tc3}.

\stop 

 \item Choose some finite set of extra obstruction models which together with our $(\hat f_{l,v}/G_{l,v},V_{l,v})$ cover the set of holomorphic curves in $\prod_{v}\bMs_{g_{v},\gamma_{v},E_{v}}$ where $\sum_v{g_{v}}=g-g_{\gamma}$ and $\sum_{v}E_{v}=E$, choosing our new obstruction models to avoid some neighborhood of the image of the set of holomorphic curves in $\bMs_{g,\gamma,E}$. We may now construct the virtual moduli space of holomorphic curves in $\prod_{v}\bMs_{g_{v},\gamma_{v},E_{v}}$ using these obstruction models and the method of section \ref{virtual class}. Lemma \ref{EV transversality} implies that we may construct our virtual moduli space so that  the maps $\EV_{1}$ and $\EV_{0}$ will be transverse, and the pullback of this component of the virtual moduli space to $\bMs_{\gamma}$ will be $\C\infty1$. Note that   pulling the multiperturbation used to define our virtual moduli space back to $\bMs_{g,\gamma,E}$ will only involve the perturbations parametrized by our $\hat f_{l,v}$ in some neighborhood of the holomorphic curves, so for the purposes of describing the pullback, we may ignore the other perturbations, and we may apply Lemma \ref{transversality conditions}.
 
\item Cover the moduli space of holomorphic curves in $\Msw_{g,[\gamma],E}(\ex B)$ by a  finite collection of obstruction models so that any obstruction model which intersects the image of $\bMs_{\gamma}$ uses one of our core families $\hat f/G$ discussed earlier. Then construct the virtual moduli space of holomorphic curves $\Mod_{g,[\gamma],E}(\ex B)$ using the method of section \ref{virtual class}. 

Recall that there are $\abs G/\abs {G_{l}}$ different lifts of $\hat f$ which give the same core family as $\hat f_{l}/\abs {G_{l}}$. Let $\hat f_{L}$ indicate the disjoint union of these lifts. The group $G$ acts on $\hat f_{L}$ so that the map $\hat f_{L}\longrightarrow \hat f$ is $G$-equivariant. We can consider $\hat f_{L}/G$ as a core family for $\bMs_{\gamma}$. Around this  core family $\hat f_{L}/G$, the pullback of the multiperturbation defined by a simple perturbation parametrized by $\hat f$ is the same as the multiperturbation defined using the simple perturbation parametrized by $\hat f_{L}$ which is just the pullback of the original perturbation using the map $\hat f_{L}\longrightarrow \hat f$. Let $(\hat f_{L},V_{L})$ indicate the pre obstruction model which is the disjoint union of the corresponding $(\hat f,V_{l})$. Then $(\hat f_{L}/G,V_{L})$ is an obstruction model on $\bMs_{\gamma}$. Using generic perturbations in our construction of the virtual moduli space of curves in $\ex B$, we may achieve transversality for the corresponding perturbed $\dbar$ equation on any compact subset of these obstruction models $(\hat f_{L}/G,V_{L})$.

\item We may now construct a cobordism in $\bMs_{g,\gamma,E}$ between the pullback of our two different moduli spaces by perturbing the $\dbar$ equation with a family of multiperturbations parameterized by $s\in [0,1]$. (We shall verify that the orientation of this cobordism is compatible with the orientations on the pullback of our different moduli spaces afterwards. ) In particular, multiply our generic multiperturbation pulled back from $\prod_{v}\bMs_{\gamma_{v}}$ by a cutoff function $\rho(s):[0,1]\longrightarrow[0,1]$ which is $1$ when $t=0$ and $0$ when $t=1$. (In other words, if the original pulled back multiperturbation on a family was  $\sum_{i} w_{i}t^{\mathfrak P_{i}}$, then use $\sum_{i} w_{i}t^{\rho(s)\mathfrak P_{i}}$.) Multiply (in the sense of definition \ref{multi}) the resulting multiperturbation with the pullback of a family of multiperturbations from $\Msw(\ex B)$ given by a generic  family of small simple perturbations parametrized by our obstruction models on $\Msw(\ex B)$ which vanishes when $s=0$, and is the multiperturbation used to define the virtual moduli space when $s=1$.

Recall that in our construction of virtual moduli spaces, we choose compact subsets of obstruction models $C\subset \hat f/G$ who's interiors still cover the moduli space of holomorphic curves, and only consider transversality in these compact subsets of obstruction models. Choose one set of these compact subsets for our $\hat f_{l}/G_{l}$, and another set of these compact subsets for our $\hat f_{L}/G$.

Lemma \ref{transversality conditions} together with Lemma \ref{EV transversality} imply that at $s=0$ we have transversality on the relevant compact subsets of $\hat f_{l}/G_{l}$.
It then follows from Theorem \ref{mms} that for $s$ small, we have transversality for the $\dbar$ equation on our  compact subset $\hat f_{l}/G_{l}$. Then Theorem \ref{mms} implies that with our generic family of simple perturbations parametrized by  $\hat f$, we may achieve transversality for the other values of $s$ within our compact subsets of $\hat f_{L}/G_{L}$. Similarly, so long as $\ex B$ has dimension greater than zero, we may ensure that our solutions have no automorphisms. It follows that our two virtual moduli spaces  of curves in $\bMs_{g,\gamma,E}$ are cobordant. 

Of course, in the case that the conditions of Lemma \ref{transversality conditions} are satisfied by the zero perturbation, and the holomorphic curves in $\bMs_{g,\gamma,E}$ have the same number of automorphisms as their image in $\Msw_{g,[\gamma],E}$ then  Lemma \ref{transversality conditions} implies that we may construct our virtual moduli space on $\ex B$ using perturbations which pull back to $\bMs_{g,\gamma,E}$ to give the zero perturbation, so in this case, the pullback of our two virtual moduli spaces coincide. 

\item We must verify that the orientation on our cobordism is compatible with the orientations pulled back from $\prod_{v}\Mod_{g_{v},[\gamma_v],E_{v}}$ and $\Mod_{g,[\gamma],E}(\ex B)$. As usual, we shall orient our moduli space by orienting the obstruction models $(\hat f_{l}/G_{l},V_{l})$. First, note that we can orient $\hat f_{l}$ exactly the same way as we orient $\hat f$ - using the fact that it is a core family. Similarly, we may orient $\hat f_{l,v}$ for all $v$ using the fact that these are core families too. These orientations are related by equating  $\ex F(\hat f_{l})$ with the product of a coordinate in $\et 11$ for each internal edge of $\gamma$ and $\ex F(\hat f_{l,v})$ for all vertices $v$ in $\gamma$. (Each of these factors is even dimensional, so the order of this product does not matter.) 

As usual, we orient $V_{l}$ and $V_{l,v}$ by identifying them with the cokernel of the relevant linearized $\dbar$ operator and using the method given in Remark \ref{orientation} to orient this cokernel using a homotopy of $D\dbar$ to complex operator. This method of constructing orientations for our moduli space in $\bMs_{g,\gamma,E}$ will clearly give the same orientation as pulling back our orientation for the moduli space in $\Msw(\ex B)$. 

We now have an orientation of $\hat f_{l,v}$ and an orientation of $V_{l,v}$. These define an orientation on our virtual moduli space  $\Mod_{g_{v},[\gamma_v],E_{v}}$. In particular, applying Theorem \ref{Multi solution} to $(\hat f_{l,v}/G_{l,v},V_{l,v})$ gives a solution to the perturbed $\dbar$ equation mod $V_{l,v}$ which consists of a weighted sum of families parametrized by $\ex F(\hat f_{l,v})$ along with sections $\dbar'_{v}$ of $V_{l,v}$; our virtual moduli space is given by a weighted sum of the oriented intersection of these $\dbar'_{v}$ with the zero section.  To pull back this oriented virtual moduli space to $\bMs_{g,\gamma,E}$, we take the fiber product over $\prod_{e}\mathcal M_{e}$ (which has an almost complex structure which orients it), restrict to the subset with tropical part in the image of $\gamma$-decorated tropical completion then take the bundle of gluing choices over this fiber product (which has has a fiberwise complex structure with orients it). As all factors in this fiber product are even dimensional, we need not specify the order in which we take fiber products to specify the orientation. 

As $V_{l,v}$ is an even dimensional bundle, we may swap the order of the fiber product and the intersection of $\dbar'$ with  the zero section, so another way of specifying the orientation of the pullback of  $\prod_{v}\Mod_{g_{v},[\gamma_v],E_{v}}$ is as follows:

Theorem \ref{Multi solution} gives us a weighted branched sum of families parametrized by $\ex F(\hat f_{l,v})$ together with sections $\dbar'_{v}$ of $V_{l,v}$. Recall from the proof of Lemma \ref{transversality conditions} that we may also apply the version of Theorem \ref{Multi solution} for families to obtain the local moduli space of solutions to $\dbar'_{v}\in V_{l,v}\oplus_{e\in S_{v}}V_{e}$ which will be a weighted branched sum of families parametrized by the total space of the bundle $\oplus_{e\in S_{v}}V_{e}$ over $\ex F(\hat f_{l,v})$ together with sections $\dbar ''_{v}$ of $V_{l,v}\oplus_{e\in S_{v}}V_{e}$. These sections $\dbar''_{v}$ are the identity on the $V_{e}$ coordinates, and the intersection of $\dbar''_{v}$ with $V_{l,v}\subset V_{l,v}\oplus_{e\in S_{v}}V_{e}$ is equal to our moduli space of solutions mod $V_{l,v}$. Giving $V_{e}$ any orientation we like, the oriented intersection of $\dbar''_{v}$ with the zero section will be equal to the oriented intersection of $\dbar'_{v}$ with the zero section.

As noted in the proof of Lemma \ref{transversality conditions}, we can take the fiber product of these families over $\prod_{e}\mathcal M_{e}$, restrict to the subset with tropical part in the image of $\gamma$-decorated tropical completion, and then take the cover corresponding to gluing choices over this fiber product to obtain the local moduli space of solutions to $\dbar'\in V_{l}$  which is obtained by applying Theorem \ref{mms} to  $(\hat f_{l},V_{l})$. We must now verify that we obtain the correct orientation on this moduli space of solutions to $\dbar'\in V_{l}$. Recall that $V_{l,v}\oplus V_{e}$ is identified with the dual of the cokernel of $D\dbar$ restricted to sections in $X^{\infty,\underline 1}(\hat f_{l,v})$ which vanish on the edge corresponding to a chosen end of $e$, and that $V_{e}$ is in the image of $D\dbar$. As $D\dbar$ has no kernel, applying $D\dbar^{-1}$ to $V_{e}$, then restricting the resulting sections to a point on the relevant edge will give an isomorphism between $V_{e}$  and the tangent space of $\ex B$ at that point, which has an almost complex structure. We may also choose $V_{e}$ so that $D\dbar^{-1}(V_{e})$ is complex. 
It follows that the orientation on $V_{e}$ determined using a homotopy of $D\dbar$ to a complex operator (in other words the method of Remark \ref{orientation})  is equal to the orientation determined by identifying $V_{e}$ with some pullback of the tangent space to $\ex B$.  Note in particular, this orientation agrees with an orientation of the solution of $\dbar'_{v}\in V_{l,e}\oplus_{e\in S_{v}}V_{e}$ given by parametrizing solutions by $\ex F(\hat f_{l,v})$ and the restriction of our maps to the relevant  ends corresponding to $S_{v}$. This orientation then agrees with the given orientation of $\hat f_{l}$ after taking fiber products over $\prod_{e}\mathcal M_{e}$ and taking the bundle of gluing choices over the result. In other words, the orientation from taking the fiber product of the solutions of $\dbar_{v}'\in V_{l,v}\oplus_{e\in S_{v}}V_{e}$ is equal to the orientation of the solutions of $\dbar'\in V_{l}$.

 The orientation of $V_{l}$ using the method of Remark \ref{orientation} will agree with the orientation of $\oplus_{v}V_{l,v}\oplus_{e}V_{e}$ because of the following: We may extend the domain of $D\dbar$ to include sections which do not agree along edges. Then the cokernel will be represented by $\oplus_{v}V_{l,v}$, and $\oplus_{e}V_{e}$ will be in the image of $D\dbar$. Again, we may assume that we've chosen $V_{e}$ so that $D\dbar^{-1}(V_{e})$ is complex, so the homotopy of $D\dbar$ to a complex operator will give $V_{e}$ the orientation corresponding to the complex structure on $D\dbar^{-1}(V_{e})$, which is the same as the orientation we gave $V_{e}$ when orienting $V_{l,v}\oplus V_{e}$ using the method of Remark \ref{orientation}.  
 
 Therefore, the pullback of the orientation of  $\prod_{v}\Mod_{g_{v},[\gamma_v],E_{v}}$ to $\bMs_{g,\gamma,E}$ is equal to the orientation determined within $\bMs_{g,\gamma,E}$ using $(\hat f_{l},V_{l})$.

This completes the proof of Theorem \ref{gluing cobordism}.

\end{itemize}

\section{Gromov Compactness}\label{compactness}

In this section, we give some examples of exploded manifolds $\ex B$ with an almost complex structure $J$ tamed by as symplectic form $\omega$ so that Gromov compactness holds in the following sense:

\begin{defn}\label{agcompactness}
Say that Gromov compactness holds for $(\ex B,J)$ if   the substack of   holomorphic curves  in $\Msw_{g,[\gamma],\beta}(\ex B)$ is compact (in the topology on $\Msw_{g,[\gamma],\beta}$ described in \cite{cem}), and  there are only a finite number of $\beta$ with $\beta(\omega)\leq E$  so that there is a holomorphic curve in $\Msw_{g,[\gamma],\beta}$. 

Say that Gromov compactness holds for a family $\hat {\ex B}\longrightarrow\ex G$ if the following holds: given any compact exploded manifold $\hat {\ex G}'$ with a map $\hat {\ex G}'\longrightarrow \ex G$, let $\hat {\ex B}'\longrightarrow  \ex G'$ be the pullback of our original family $\hat{\ex B}\longrightarrow \ex G$. Then the substack of holomorphic curves in $\Msw_{g,[\gamma],\beta}(\hat {\ex B}')$ is compact and there are only a finite number of $\beta$ with $\beta(\omega)\leq E$ so that there is a holomorphic curve in $\Msw_{g,[\gamma],\beta}(\hat {\ex B}')$.

\end{defn}

To use the results of \cite{cem}, we must construct a strict taming of $J$ containing $\omega$. This involves some work, but is done in the following cases in \cite{cem}.

\begin{thm}\label{compactness examples} Gromov compactness holds in the following cases

\begin{enumerate}
\item \label{explosion compactness}If $\expl M$ is the explosion of a compact complex manifold with normal crossing divisors which are embedded submanifolds, and $\omega$ is a symplectic form on $M$ taming the complex structure, then Gromov compactness holds.

\item The paper \label{symplectic explosion compactness}\cite{cem} constructs an exploded manifold $\ex M$ from a compact symplectic manifold $M$ with orthogonally intersecting codimension $2$ symplectic submanifolds. This construction is similar to the construction of the explosion functor, so the smooth part of $\ex M$ is $M$, and the tropical part of $\ex M$ is the dual intersection complex which has a vertex for each connected component of $M$, a ray for each of our codimension $2$ submanifolds, and an $n$-dimensional face $[0,\infty)^{n}$ for each $n$-fold intersection. There exists an almost complex structure on $\ex M$ tamed by a symplectic form $\omega$ corresponding to the symplectic structure on $M$ so that Gromov compactness holds. 

\item\label{standard simplex compactness} Assume $\ex B$ is a basic exploded manifold so that each polytope in the tropical part of $\ex B$ is a standard simplex. Suppose that $J$ is an almost complex structure on $\ex B$ so that for each exploded function $\tilde z$ locally defined on $\ex B$, $\tilde z^{-1}\dbar\tilde z$ is a smooth function. (For example any integrable complex structure satisfies this condition.) If $\omega$ is a symplectic form on $\ex B$ which tames $J$, then Gromov Compactness holds.

If $J$ is a civilized  almost complex structure on $\ex B$ tamed by $\omega$ which does not satisfy the above condition, then it may be modified to an almost complex structure still tamed by $\omega$ which does

\end{enumerate}

\end{thm}  
  
  Note that in particular, item \ref{standard simplex compactness} of Theorem \ref{compactness examples} implies that Gromov compactness holds for any tamed almost complex structure on a compact symplectic manifold. 
  
  Gromov compactness obviously holds for the product of two exploded manifolds with almost complex structures tamed by symplectic forms for which Gromov compactness holds. As well as the above examples, a useful base case is the following:
  
  \begin{lemma}Gromov compactness holds for $\ex T$ with the standard complex structure and the zero symplectic form.
  \end{lemma} 
  
 \pf A refinement $\ex T'$ of $\ex T$ is the explosion of $\mathbb CP^{1}$ relative to $0$ and $\infty$. Item \ref{explosion compactness} of Theorem \ref{compactness examples} then tells us that Gromov compactness holds for this refinement when we use a standard toric form on $\mathbb CP^{1}$ for $\omega$. Therefore, the set of holomorphic curves in $\Msw_{g,[\gamma],E}(\ex T')$ is compact. The energy of a holomorphic curve in $\ex T'$ with tropical part isotopic to $[\gamma]$ is determined by the sum of the momentum of the ends of $\gamma$, so is determined by $\gamma$. Therefore, the set of holomorphic curves in $\Msw_{g,[\gamma]}(\ex T')$ is compact. This is just a refinement of the set of holomorphic curves in $\Msw_{g,[\gamma],0}(\ex T)$, so the set of holomorphic curves in $\Msw_{g,[\gamma],0}$ is compact, and Gromov compactness holds for $\ex T$.
 
 \stop
 
 \
   
   As a corollary, Gromov compactness holds for the product of $\ex T^{n}$ with any of the examples from Theorem \ref{compactness examples}. In particular, if $\ex B$ is of type \ref{standard simplex compactness} from Theorem \ref{compactness examples}, then the tropical completion of any strata of $\ex B$ will be in the form of $\ex T^{n}$ times something of type \ref{standard simplex compactness}, so Gromov compactness will also hold for the tropical completion of any strata of $\ex B$.

\

One way to prove Gromov compactness in a family, $\hat {\ex B}$ is to construct a family of strict tamings and use \cite{cem}. With minimal work, we can also use the above examples to construct examples of families  in which Gromov compactness holds:

 \begin{example} Gromov compactness holds for any trivial family of symplectic manifolds with a family of almost complex structures which need not be trivial.
 \end{example}
We must prove in particular that for any compact trivial family of symplectic manifolds, Gromov compactness holds. Let $F$ be the base of our compact family. We may embed $F$ into some compact symplectic manifold $N$, then choose a tamed almost complex structure on $M\times N$ so that fibers of $M\times N$ are holomorphic and so that the pullback of the corresponding fiberwise almost complex structure to $M\times F$ is our original family of almost complex structures. Then Gromov compactness for $M\times N$ implies Gromov compactness for our family.

\begin{example}Suppose that Gromov compactness holds for $(\ex B,J,\omega)$ and there is a family $\ex B\longrightarrow \ex G$ with $J$-holomorphic fibers. Let $\hat {\ex B}\longrightarrow \ex G'$  be a family obtained by taking the fiber product of some map $\ex G'\longrightarrow \ex G$ with $\ex B\longrightarrow \ex G$, let $\omega'$ be the pullback of $\omega$ under the map $\hat{\ex B}\longrightarrow \ex B$, and let $J'$ be the pullback of the fiberwise almost complex structure obtained by restricting $J$ to fibers of $\ex B\longrightarrow \ex G$. Then Gromov compactness holds for $(\hat {\ex B},J',\omega')$.
\end{example} 

In particular, the above example  can be used to show that the explosion of many normal crossing degenerations from algebraic geometry give families of exploded manifolds in which Gromov compactness holds.

\bibliographystyle{plain}
\bibliography{ref.bib}
\end{document}